\documentclass[10pt,leqno,english,a4paper]{amsart}
\usepackage{color}
\usepackage{hyperref}
\usepackage{amssymb}
\usepackage[square,numbers]{natbib}
\usepackage[all]{xy}

\theoremstyle{plain}
\newtheorem{mainthm}{Theorem}

\newtheorem{thm}{Theorem}
\numberwithin{thm}{section}

\newtheorem{thmfact}[thm]{Fact}
\swapnumbers
\newtheorem{lemm}[subsubsection]{Lemma}
\newtheorem{prop}[subsubsection]{Proposition}
\newtheorem{claim}[subsubsection]{Claim}
\newtheorem{fact}[subsubsection]{Fact}
\newtheorem{obsv}[subsubsection]{Observation}

\numberwithin{subsection}{section}
\numberwithin{subsubsection}{subsection}

\renewcommand{\labelenumi}{(\alph{enumi})}

\DeclareMathOperator{\Hom}{\mathrm{Hom}}

\DeclareMathOperator{\DGHom}{\mathit{Hom}}

\DeclareMathOperator{\HopfHom}{\mathit{HopfHom}}
\DeclareMathOperator{\HopfEnd}{HopfEnd}
\DeclareMathOperator{\PrimEnd}{PrimEnd}
\DeclareMathOperator{\Alg}{Alg} \DeclareMathOperator{\CoAlg}{CoAlg}

 \DeclareMathOperator{\Mod}{Mod}
\DeclareMathOperator{\HopfMod}{HopfMod}
\DeclareMathOperator{\Op}{Op} \DeclareMathOperator{\dg}{dg}
\DeclareMathOperator{\gr}{gr} 
 
 \DeclareMathOperator{\sk}{sk}
 \DeclareMathOperator{\skcell}{sk}
\DeclareMathOperator{\ckcell}{ck} 
\DeclareMathOperator{\itr}{i} \DeclareMathOperator{\str}{s}
\DeclareMathOperator{\ctr}{c}

\DeclareMathOperator{\dec}{dec} 
\DeclareMathOperator{\tr}{tr} \DeclareMathOperator{\ev}{ev}
\DeclareMathOperator{\Set}{Set} \DeclareMathOperator{\Gr}{Gr}
\DeclareMathOperator{\pk}{p}

\DeclareMathOperator{\cofibration}{\xymatrix@M=0pt@C=16pt{*{\hspace{4pt}}\ar@{>->}[r]&}}
\DeclareMathOperator{\wecofib}{\xymatrix@M=0pt@C=16pt{*{\hspace{4pt}}\ar@{>->}[r]^{\sim}&}}
\DeclareMathOperator{\fibration}{\xymatrix@W=0pt@M=0pt@C=16pt{\ar@{->>}[r]&}}
\DeclareMathOperator{\wefib}{\xymatrix@W=0pt@M=0pt@C=16pt{\ar@{->>}[r]^{\sim}&}}

\DeclareMathOperator{\N}{\mathbb{N}}
\DeclareMathOperator{\F}{\mathbb{F}}
\DeclareMathOperator{\Z}{\mathbb{Z}}

\DeclareMathOperator{\A}{\mathcal{A}}
\DeclareMathOperator{\B}{\mathcal{B}}
\DeclareMathOperator{\C}{\mathcal{C}}

\DeclareMathOperator{\E}{\mathcal{E}}
\DeclareMathOperator{\FOp}{\mathcal{F}}

\DeclareMathOperator{\K}{\mathcal{K}}

\renewcommand{\P}{\mathcal{P}}
\renewcommand{\O}{\mathcal{O}}
\DeclareMathOperator{\Q}{\mathcal{Q}}
\DeclareMathOperator{\R}{\mathcal{R}}
\DeclareMathOperator{\SOp}{\mathcal{S}}

\DeclareMathOperator{\HopfOp}{\mathrm{HopfOp}}
\DeclareMathOperator{\PrimOp}{\mathrm{PrimOp}}
\DeclareMathOperator{\Id}{Id} \DeclareMathOperator{\id}{id}
\DeclareMathOperator{\pr}{pr}
\DeclareMathOperator{\sgn}{\mathrm{sgn}}
\DeclareMathOperator{\shuffle}{\mathit{Sh}}

\DeclareMathOperator{\colim}{\mathrm{colim}}
\DeclareMathOperator{\coker}{\mathrm{coker}}

\DeclareMathOperator{\mylim}{\mathrm{lim}}

\DeclareMathOperator{\op}{\mathrm{op}}

\newcommand{\ptensor}{\underline{p}}

\DeclareMathOperator{\I}{\mathbb{I}}
\newcommand{\Cube}[1]{\mathbb{D}^{#1}}
\newcommand{\IntCube}[1]{\mathbb{D}'{}^{#1}}

\newcommand{\x}[1]{\mathrm{#1}}

\newcommand{\Match}[1]{\mathbb{M}#1}
\newcommand{\Latch}[1]{\mathbb{L}#1}

\newcommand{\push}[1]{} 
\newcommand{\pull}[1]{} 

\title[The universal Hopf operads of the bar construction]{The universal Hopf operads\\of the bar construction}
\date{9 January 2007}

\author{Benoit Fresse}
\address{Universit\'e de Lille 1,
UMR 8524 du CNRS, 59655 Villeneuve d'Ascq C\'edex (France)}
\email{Benoit.Fresse@math.univ-lille1.fr}
\urladdr{http://math.univ-lille1.fr/\~{ }fresse}
\thanks{Research supported in part by the ANR grant JCJC06-143080.
The author enjoyed a stay at the Institut Mittag Leffler (Sweden)
during the preparation of this memoir.}


\begin{document}

\maketitle

\begin{abstract}
The goal of this memoir is to prove that the bar complex B(A) of an
E-infinity algebra A is equipped with the structure of a Hopf
E-infinity algebra, functorially in A. We observe in addition that
such a structure is homotopically unique provided that we consider
unital operads which come equipped with a distinguished 0-ary
operation that represents the natural unit of the bar complex. Our
constructions rely on a Reedy model category for unital Hopf
operads. For our purpose we define a unital Hopf endomorphism operad
which operates functorially on the bar complex and which is
universal with this property. Then we deduce our structure results
from operadic lifting properties. To conclude this memoir we hint
how to make our constructions effective and explicit.
\end{abstract}

\newpage\tableofcontents

\newpage\begin{center}{\scshape Main theorems}\end{center}\medskip

\begin{list}{}{\setlength{\labelsep}{2em}\setlength{\labelwidth}{4em}\setlength{\itemindent}{-2em}}

\item[\ref{thm:LambdaModuleModelStructure}] The Reedy model
  category of $\Lambda_*$-modules, \pageref{thm:LambdaModuleModelStructure}

\item[\ref{thm:UnitalOperadModelStructure}] The Reedy model
  category of unital dg-operads, \pageref{thm:UnitalOperadModelStructure}

\item[\ref{thm:HopfLambdaModuleModelStructure}] The Reedy model
  category of Hopf $\Lambda_*$-modules,
  \pageref{thm:HopfLambdaModuleModelStructure}

\item[\ref{thm:HopfOperadModelStructure}] The Reedy model
  category of unital Hopf operads, \pageref{thm:HopfOperadModelStructure}

\item[\ref{thm:WHopfOperad}] The Boardman-Vogt construction
  recall, \pageref{thm:WHopfOperad}

\item[\ref{thm:BoardmanVogtOperadicCellularDecomposition}] The
  operadic cellular structure of the Boardman-Vogt construction,
  \pageref{thm:BoardmanVogtOperadicCellularDecomposition}

\item[\ref{thm:BarHopfEndomorphismOperadRecall}] The universal
  definition of the Hopf endomorphism operad of the bar construction,
  \pageref{thm:BarHopfEndomorphismOperadRecall}

\item[\ref{thm:HopfEndomorphismOperadFibration}] The Hopf
  operad of universal bar operations preserve fibrations and acyclic
  fibrations, \pageref{thm:HopfEndomorphismOperadFibration}

\item[\ref{thm:HopfOperadLifting}] The Hopf
  operad of universal bar operations for associative and commutative algebras,
  \pageref{thm:HopfOperadLifting}

\item[\ref{thm:HopfExistenceUniquenessRecall}] The existence
  and uniqueness of Hopf operad actions on the bar complex of an
  $E_\infty$-algebra, \pageref{thm:HopfExistenceUniquenessRecall}

\item[\ref{thm:HopfOperadMorphismHomotopyUniqueness}] The
  uniqueness property for the action of connected Hopf operads,
  \pageref{thm:HopfOperadMorphismHomotopyUniqueness}

\item[\ref{fact:BarOperationMorphismReduction}] Recall of
  the factorization of operad actions on the bar complex,
  \pageref{fact:BarOperationMorphismReduction}

\item[\ref{fact:OperationStructure}] The
  structure of operations on the bar complex,
  \pageref{fact:OperationStructure}

\item[\ref{fact:OperationExpansion}] The expansion of
  operations on the bar complex, \pageref{fact:OperationExpansion}

\item[\ref{thm:OperationExplicitDefinition}] The explicit equations of
  operations on the bar complex, \pageref{thm:OperationExplicitDefinition}

\item[\ref{thm:OperationRecursiveDefinition}] The explicit
  recursive construction of operations on the bar complex, \pageref{thm:OperationRecursiveDefinition}

\end{list}

\clearpage

\part*{Introduction}

This memoir is a sequel of the article~\cite{Bar1} in which we
proved that the classical bar complex $B(A)$ can be equipped with
the structure of an $E_\infty$-algebra if $A$ is an $E_\infty$
algebra.

To be precise we consider operads in the category of differential
graded modules (dg-modules for short) and an $E_\infty$-operad
refers to a dg-operad $\E$ equivalent to the operad of associative
and commutative algebras $\C$. An $E_\infty$-algebra refers to an
algebra over any given $E_\infty$-operad. Similarly one considers
$A_\infty$-operads, defined as dg-operads equivalent to the operad
of associative algebras $\A$, and $A_\infty$-algebras, defined as
algebras over an $A_\infty$-operad. For our purpose we fix a
specific model $\K$ of an $A_\infty$-operad. Namely we consider the
classical chain operad of Stasheff's associahedra. The classical
notion of an $A_\infty$-algebra, defined by a collection of higher
associative products $\mu_r: A^{\otimes r}\rightarrow A$, is
equivalent to the structure of an algebra over this $A_\infty$
operad. Recall that the bar complex is defined precisely for
algebras equipped with such operations. For any $E_\infty$-operad
$\E$ there exists a homotopically unique morphism $\K\rightarrow\E$.
Once we fix such a morphism any $\E$-algebra is provided with the
structure of a $\K$-algebra and hence has an associated bar complex
$B(A)$.

In~\cite{Bar1} we consider only the chain structure of the bar
complex $B(A)$. But, classically, one identifies $B(A)$ with the
tensor coalgebra generated by $\Sigma\bar{A}$, the suspension of the
augmentation ideal of $A$. One observes in addition that the bar
differential $\partial: B(A)\rightarrow B(A)$ is defined by a
coalgebra coderivation so that the bar complex $B(A)$ forms a
dg-coalgebra. Therefore a natural aim consists in extending the
constructions of~\cite{Bar1} in the framework of dg-coalgebras.

For this purpose we consider operads in the ground symmetric
monoidal category of dg-coalgebras, usually called Hopf operads, and
algebras over operads in this category. The algebras over a Hopf
operad $\P$ in the category of dg-coalgebras are usually referred to
as Hopf $\P$-algebras. In fact, one can incidentally forget
coalgebra structures and consider algebras over Hopf operads in the
underlying category of dg-modules and similarly for other
structures. Therefore, as a general rule, the objects defined in the
category of dg-coalgebras are referred to by the qualifier Hopf.
Otherwise we consider tacitely an underlying structure in the
category of dg-modules.

For a commutative algebra the classical shuffle product of tensors
defines a morphism of differential graded coalgebras $\smile:
B(A)\otimes B(A)\rightarrow B(A)$. Consequently, the bar complex of
a commutative algebra $B(A)$ is equipped with the structure of a
commutative Hopf algebra. The goal of this memoir is precisely to
extend this structure result to $E_\infty$-algebras.

\subsection*{Sketch of the memoir results}
Recall briefly that a Hopf operad $\P$ is defined by a collection of
dg-coalgebras $\P(r)$, acted on by the symmetric group $\Sigma_r$,
equipped with operad composition products
$\P(r)\otimes\P(n_1)\otimes\dots\otimes\P(n_r)\rightarrow\P(n_1+\dots+n_r)$
which are morphisms in the category of coalgebras. Similarly, a Hopf
$\P$-algebra consists of a dg-coalgebra $\Gamma$ equipped with
evaluation products $\P(r)\otimes\Gamma^{\otimes
r}\rightarrow\Gamma$ which are morphisms in the category of
coalgebras (we refer to~\ref{subsection:HopfOperadsIntroduction} for
detailed recalls).

Observe that the bar coalgebra $\Gamma = B(A)$ is equipped with a
natural unit morphism $\F\rightarrow B(A)$. Therefore, for the
purposes of this memoir, it is natural to consider unital Hopf
operads, equipped with a distinguished unital operation $*\in\P(0)$,
and operad actions $\P(r)\otimes B(A)^{\otimes r}\rightarrow B(A)$
for which the unital operation $*: \F\rightarrow B(A)$ agrees with
the natural unit of the bar complex $\F\rightarrow B(A)$. This unit
requirement gives the analogue of a boundary condition for the
construction of operad actions on $B(A)$. Indeed, if we restrict
ourself to non-negatively graded Hopf operads, then the unit
requirement implies the assumption of the uniqueness theorem
of~\cite{Bar1}.

As in~\cite{Bar1} we consider also connected operads which, in the
unital context, satisfy $\P(0) = \F$ and $\P(1) = \F$.

To recapitulate, we suppose given an $E_\infty$-operad $\E$ in the
category of dg-modules, equipped with a fixed operad morphism
$\K\rightarrow\P$, where $\K$ denotes Stasheff's chain operad, and
we consider the bar coalgebra $B(A)$ for $A$ an algebra over $\E$.
Then let $\Q$ denote a unital non-negatively graded Hopf
$E_\infty$-operad. Our main goal is to prove the following theorem:

\begin{mainthm}\label{thm:HopfExistenceUniqueness}\hspace*{2mm}

\begin{enumerate}
\item
The bar complex of an $\E$-algebra $B(A)$ can be equipped with the
structure of a Hopf $\Q$-algebra functorially in $A$ and so that the
unital operation $\Q(0)\rightarrow B(A)$ agrees with the natural
unit of the bar complex $\F\rightarrow B(A)$ provided that $\Q$ is a
Reedy cofibrant object in the category of unital Hopf
$E_\infty$-operads.
\item
Any such $\Q$-algebra structure where $\Q$ is connected and
non-negatively graded satisfies the requirement of the uniqueness
theorem of~\cite{Bar1}. More explicitly, if the unit condition of
claim (a) is satisfied and the operad $\Q$ is connected and
non-negatively graded, then, for a commutative algebra $A$, the
$\Q$-algebra structure of $B(A)$ reduces automatically to the
classical commutative algebra structure of $B(A)$, the one given by
the shuffle product of tensors.
\end{enumerate}
\end{mainthm}

The proof of the existence claim (a) follows the same lines of
argument as in the framework of dg-modules. Namely we introduce
first a universal unital Hopf operad, the Hopf endomorphism operad
of the bar construction $\HopfEnd_B^{\P}$, which operates
functorially on the bar complex of algebras over a given operad
$\P$. More precisely, we prove the following result:

\begin{mainthm}\label{thm:BarHopfEndomorphismOperad}
Let $\P$ denote an operad (in dg-modules) equipped with an operad
morphism $\K\rightarrow\P$, where $\K$ denotes Stasheff's chain
operad. There is a universal unital Hopf operad $\Q =
\HopfEnd_B^{\P}$ such that the bar complex of a $\P$-algebra $B(A)$
is equipped with the structure of a Hopf algebra over $\Q$,
functorially in $A\in\P\Alg$.

More precisely, the Hopf operad $\HopfEnd_B^{\P}$ operates on the
coalgebra $B(A)$ functorially in $A\in\P\Alg$ and so that the unital
operation $*: \F\rightarrow B(A)$ agrees with the unit of $B(A)$.
Furthermore, we have a one-to-one correspondence between such Hopf
operad actions and morphisms of unital Hopf operads $\rho:
\Q\rightarrow\HopfEnd_B^{\P}$.
\end{mainthm}

For the sake of completeness, we should point out that the map
$\P\mapsto\HopfEnd_B^{\P}$ defines a functor on the category of
operads under~$\K$.

As in the context of dg-modules, any coalgebra $\Gamma$ has an
associated endomorphism Hopf operad defined by
\begin{equation*}
\HopfEnd_{\Gamma}(r) = \HopfHom(\Gamma^{\otimes r},\Gamma)
\end{equation*}
where $\HopfHom(K,L)$ denotes an appropriate internal hom-object in
the category of coalgebras. By definition a Hopf $\Q$-algebra
structure is equivalent to a Hopf operad morphism
$\Q\rightarrow\HopfEnd_{\Gamma}$. One can adapt this construction
for unital coalgebras so that $\HopfEnd_{\Gamma}$ forms a unital
Hopf operad and a morphism $\rho: \Q\rightarrow\HopfEnd_{\Gamma}$
preserves the distinguished unital operations if and only if in the
equivalent Hopf $\Q$-algebra structure the unital operation $*:
\F\rightarrow\Gamma$ agrees with the unit of $\Gamma$.

The Hopf endomorphism operad of the bar construction is defined by
the coend of the bifunctors
\begin{equation*}
\HopfEnd_{B(A)}(r) = \HopfHom(B(A)^{\otimes r},B(A))
\end{equation*}
where $A$ ranges over the category of $\P$-algebras. The assertions
of theorem~\ref{thm:BarHopfEndomorphismOperad} are immediate from
this construction. Notice that this endomorphism operad is
$\Z$-graded and non-connected, like any endomorphism object.

The classical commutative algebra structure of the bar complex gives
a morphism $\nabla_c: \C\rightarrow\HopfEnd_B^{\C}$ and the
existence assertion of theorem~\ref{thm:HopfExistenceUniqueness} is
equivalent to the lifting problem:
\begin{equation*}
\xymatrix{ & & \HopfEnd_B^{\E}\ar[d] \\
\Q\ar@{-->}^{\exists ?}[urr]\ar[r] & \C\ar[r]^(0.33){\nabla_c} &
\HopfEnd_B^{\C} }.
\end{equation*}
As in~\cite{Bar1} we introduce another universal operad
$\HopfOp_B^{\P}$, the Hopf operad of natural operations of the bar
complex, that forms a suboperad of $\HopfEnd_B^{\P}$ and that agrees
with this one only if the ground field is infinite. In general this
operad is endowed with better homotopical properties than the Hopf
endomorphism operad $\HopfEnd_B^{\P}$. More specifically, we prove
the following property:

\begin{mainthm}\label{thm:HopfHomotopyInvariance}
The functor $\P\mapsto\HopfOp_B^{\P}$ maps a fibration, respectively
an acyclic fibration, of (non-unital) operads under $\K$ to a Reedy
fibration, respectively an acyclic Reedy fibration, of unital Hopf
operads.
\end{mainthm}

In this statement and in theorem~\ref{thm:HopfExistenceUniqueness}
we refer to a particular model structure on the category of Hopf
operads, the Reedy model structure, in which we have cofibrations,
fibrations and weak-equivalences endowed with the classical lifting
properties. In fact, we introduce a new model structure that differs
from the general adjoint model structures considered
in~\cite{BergerMoerdijk} and which is more appropriate in the unital
context. Usually, for a solvable operadic lifting problem
\begin{equation*}
\xymatrix{ \A\ar[d]_{i}\ar[r]^{f} & \P\ar[d]^{p} \\
\B\ar[r]^{g}\ar@{-->}[ur]^{h}
  & \Q },
\end{equation*}
the lifting $h: \B(r)\rightarrow\P(r)$ can be constructed
effectively by induction on $r\in\N$. In the unital context one can
observe that the operadic composites with a unital operation
$*\in\P(0)$ provide the underlying $\Sigma_*$-module of a unital
operad $\P$ with operations $\partial_i: \P(r)\rightarrow\P(r-1)$
that decrease the operadic arity. Consequently, in the inductive
construction, one can assume that the composites
$\B(r)\xrightarrow{h}\P(r)\xrightarrow{\partial_i}\P(r-1)$ are
specified so that the lift-component $h: \B(r)\rightarrow\P(r)$ is
indeed obtained by lifting a matching morphism
\begin{equation*}
(\mu,p): \P(r)\rightarrow\Match{\P}(r)\times_{\Match{\Q}(r)}\Q(r)
\end{equation*}
for a natural notion of matching objects associated to unital
operads. Accordingly, for our purposes it is natural to let an
operad morphism $p: \P\rightarrow\Q$ be a fibration if for all
$r\in\N$ the matching morphism $(\mu,p):
\P(r)\rightarrow\Match{\P}(r)\times_{\Match{\Q}(r)}\Q(r)$ is a
fibration in the underlying ground category (in the category of
dg-coalgebras for Hopf operads, dg-modules otherwise). We call this
class of fibrations the Reedy fibrations in order to distinguish
them from the classical fibrations. We have also a class of Reedy
cofibrations characterized by the left lifting property with respect
to acyclic Reedy fibrations as usual, where we consider the usual
weak-equivalences. We prove precisely that these class of morphisms
define a model structure on unital operads.

To be precise, we prove the axioms of a model structure only for
non-negatively graded unital Hopf operads. But we mention that an
endomorphism operad $\HopfOp_B^{\P}$ is $\Z$-graded like any
endomorphism object. In fact, we extend abusively the notion of a
fibration, respectively of an acyclic fibration, to such operads by
the left-lifting property requirement. Precisely, we let a morphism
of (possibly $\Z$-graded) unital Hopf operads be a fibration,
respectively acyclic fibration, if it satisfies the left-lifting
property with respect to acyclic cofibrations, respectively
cofibrations, of non-negatively graded unital Hopf operads.
Equivalently, we have a truncation functor $\tr_+$ which gives a
right adjoint to the category embedding of $\N$-graded unital Hopf
operads into the category of all $\Z$-graded unital Hopf operads.
Abusively, we let a morphism of unital Hopf operads be a fibration,
respectively an acyclic fibration, if its truncation defines a
fibration, respectively an acyclic fibration, in the category of
non-negatively graded unital Hopf operads.

Finally, the existence claim of
theorem~\ref{thm:HopfExistenceUniqueness} is a corollary of
theorem~\ref{thm:HopfHomotopyInvariance}. Explicitly, we observe
that the Hopf operad morphism $\nabla_c:
\C\rightarrow\HopfEnd_B^{\C}$ associated to the classical shuffle
algebra structure factorizes through $\HopfOp_B^{\C}$. Then we
consider the lifting problem
\begin{equation*}
\xymatrix{ & & \HopfOp_B^{\E}\ar[r]\ar@{->>}[d]^{\sim} & \HopfEnd_B^{\E}\ar[d] \\
\Q\ar@{-->}^{\exists}[urr]\ar[r] & \C\ar[r]^(0.33){\nabla_c} &
\HopfOp_B^{\C}\ar[r] & \HopfEnd_B^{\C} }
\end{equation*}
which has automatically a solution if $\Q$ is cofibrant as the
augmentation of an $E_\infty$-operad $\E\wefib\C$ induces an acyclic
fibration of unital Hopf operads
\begin{equation*}
\HopfOp_B^{\E}\wefib\HopfOp_B^{\C}
\end{equation*}
by theorem~\ref{thm:HopfHomotopyInvariance}.

Notice that the commutative operad forms a final object in the
category of Hopf operads. For the uniqueness claim we observe that
any morphism of unital Hopf operads $\nabla:
\Q\rightarrow\HopfOp_B^{\C}$ where $\Q$ is connected and
non-negatively graded makes the diagram
\begin{equation*}
\xymatrix{ \Q\ar[d]\ar[dr]^{\nabla} & \\ \C\ar[r]_(0.33){\nabla_c} &
\HopfOp_B^{\C} }.
\end{equation*}
commute. As a consequence, for a given $E_\infty$-operad $\E$, any
morphism of unital Hopf operads $\nabla:
\Q\rightarrow\HopfOp_B^{\E}$ such that $\Q$ is connected and
non-negatively graded defines a lifting of~$\nabla_c:
\C\rightarrow\HopfOp_B^{\C}$.

\subsection*{Toward effective constructions}
The lifting process can be made effective for certain $E_\infty$
operads. Indeed we observed in~\cite{BarEinfinityNote} that the bar
complex of an algebra over the so-called surjection operad
$\mathcal{X}$ forms a Hopf algebra over the Barratt-Eccles operad
$\E$. This structure result gives an explicit instance of a morphism
$\E\rightarrow\HopfOp_B^{\mathcal{X}}$. Recall that the surjection
operad is an instance an $E_\infty$-operad in the category of
dg-modules and the Barratt-Eccles operad is an instance of a Hopf
$E_\infty$-operad. We are particularly interested in these operads
as, on one hand, we proved in~\cite{BergerFresse} (see
also~\cite{McClureSmith}) that the surjection operad, as well as the
Barratt-Eccles operad, operates on the cochain complex of simplicial
sets $C^*(X)$, and, on the other hand, we observed in~\cite{Bar1}
that the bar complex $B(C^*(X))$ is equivalent as an $E_\infty$
algebra to $C^*(\Omega X)$, the cochain complex of the loop space
of~$X$.

Nevertheless the Barratt-Eccles operad is not cofibrant and
theorem~\ref{thm:HopfExistenceUniqueness} does not apply to this
operad. One could check that the bar complex of an algebra over the
Barratt-Eccles operad is not acted on by the Barratt-Eccles operad
itself. Therefore we consider the so-called Boardman-Vogt'
construction $\Q = W(\E)$ which provides an explicit cofibrant
replacement of the Barratt-Eccles operad and for which
theorem~\ref{thm:HopfExistenceUniqueness} applies. Indeed for the
Barratt-Eccles operad $\E$ a lifting $\nabla:
W(\E)\rightarrow\HopfOp_B^{\E}$ can be constructed by an explicit
induction process so that we can obtain effectively a Hopf
$W(\E)$-algebra structure on the bar complex of algebras over the
Barratt-Eccles operad. This morphism seems to factorize through a
simplicial decomposition $W_{\Delta} = W_{\Delta}(\E)$ (in fact the
chain operad associated to the simplicial Boardman-Vogt construction
of the simplicial Barratt-Eccles operad) of the Boardman-Vogt
dg-operad $W_{\square} = W(\E)$. To be more precise, we seem to have
closed formulas for a morphism $\nabla_\delta:
W_{\Delta}\rightarrow\HopfOp_B^{ W_{\Delta}}$ that fits a
commutative diagram
\begin{equation*}
\xymatrix{ W_{\Delta}\ar[r]^(0.35){\nabla_\delta} & \HopfOp_B^{ W_{\Delta}}\ar[d] \\
W_{\square}\ar[r]^(0.35){\nabla}\ar[u] & \HopfOp_B^{\E} }.
\end{equation*}
This result, obtained as an application of explicit formulas given
in this memoir, should be confirmed and published in a subsequent
article.

\subsection*{Further prospects}
The functorial constructions addressed in this memoir and in the
previous article have the drawback to yield only global formulas,
valid for all algebras in a category. But one would like to control
the structure of the bar construction for algebra subclasses or for
a particular subclass of operations. Our idea to address this issue
is to introduce cellular operads obtained by a limit-colimit
decomposition of an $E_\infty$-operads. To motivate this idea,
observe that a quotient $\Q$ of an operad $\P$ is associated to a
subcategory of the category of $\P$-algebras; a suboperad
$\mathcal{R}$ is associated to a larger category of algebras as the
$\mathcal{R}$-algebras support less operations. To generalize this
process, we shall consider functor operads (in the sense
of~\cite{McClureSmithCosimplicial}) on an operad in the category of
categories. In the text, we give a few indications on how our
constructions might be extended to this context.

The crux is to provide the operad of universal bar operations with
appropriate cellular structures. As such, the cell categories that
one would like to consider have to be motivated by the internal
structure of the operad of universal bar operation itself. A good
example includes the $G$-cellular operads
of~\cite{Berger,BergerEnglish} that give rise to sequences of
$E_n$-operads.

\subsection*{Memoir organization}
The first part of this memoir
(sections~\ref{section:UnitaryHopfOperads}-\ref{section:CofibrantHopfOperads})
is of a general interest. In this part, we give a comprehensive
account of Reedy model categories of unital operads. The general
theory is set in~\ref{section:UnitaryHopfOperads}. The main theorems
regarding the Reedy model category of unital operads in dg-modules
and dg-coalgebras, the model categories used in the memoir, are also
stated in this section. Section~\ref{section:CofibrantHopfOperads}
is devoted to the Boardman-Vogt construction, a construction which
returns explicit cofibrant replacements in operad categories. Though
we do not use this construction explicitly in the memoir, we give a
detailed account of it for the sake of completeness and for
subsequent references.

The goal of the second part of the memoir
(sections~\ref{section:UnitaryHopfOperads}-\ref{section:CofibrantHopfOperads}),
is to define precisely the universal Hopf operads introduced in the
memoir introduction and to study the internal structure of these
objects. The first section of this part,
section~\ref{section:CocellularCoalgebras}, is still of general
interests: we define an appropriate notion of cocellular objects in
order to obtain classes of effective fibrations in the context of
$\Z$-graded dg-coalgebras. The core of the memoir is formed by
section~\ref{section:OperadActionConstruction} in which we prove our
main theorem sketched in the introduction. Namely: we define the
internal hom object $\HopfHom(K,L)$ of coalgebras; we prove the
existence of the universal Hopf endomorphism operad
$\HopfEnd_B^{\P}$; we define the suboperad of universal bar
operations $\HopfOp_B^{\P}$ and we study the structure of these
objects; more specifically, we prove that the functor
$\P\mapsto\HopfOp_B^{\P}$ preserves fibrations and acyclic
fibrations; then we prove our main existence and uniqueness theorem.

In the last part of this memoir
(section~\ref{section:ExplicitOperadActions}), we give an elementary
interpretation of the structure results obtained
in~\ref{section:OperadActionConstruction}. More particularly, we
make explicit a recursive construction that yields an effective
$E_\infty$-structure on the bar complex. For motivations, we suggest
the reader to have a glance at the introduction of this section and
theorems~\ref{thm:OperationExplicitDefinition}-\ref{thm:OperationRecursiveDefinition}
before to start a thorough reading of the memoir.

We refer to the section introductions for detailed sketches of the
content of each part.

\clearpage

\part*{The operadic framework}

\section{On unital and Hopf operads}\label{section:UnitaryHopfOperads}

\subsection{Introduction and general conventions}
In general the objects that we consider belong to the category of
dg-modules over a ground ring $\F$\glossary{$\F$}, always assumed to
be a field, and fixed once and for all.

We adopt the conventions of~\cite{Bar1} except that in this memoir
we may find more convenient to deal with \emph{augmented unitary
algebras} instead of \emph{non-unitary algebras}. To be precise, as
we equip the bar construction $B(A)$ with a Hopf algebra structure,
we find better to assume $B_0(A) = \F$ so that the functor $A\mapsto
B(A)$ targets to a category of augmented unitary objects. On the
other hand, the complex $B(A)$ is more naturally defined for a
non-unitary algebra $A$. Therefore we shall assume on the contrary
that the source of the functor $A\mapsto B(A)$ lives in a category
of non-unitary objects.

Anyway, as stated in the memoir introduction, we shall consider
\emph{unital operads} $\P$\index{operad!unital}\index{unital!operad}
which come equipped with a distinguished \emph{unital operation}
$*\in\P(0)$\index{unital!operation} that generates $\P(0)$ so that
$\P(0) = \F
*$ (for an operad in a module category, the usual case for the
operads considered in this memoir). As explained
in~\cite[\S3.2.1]{Bar1}, this assumption $\P(0) = \F$ implies that
the ground field~$\F$ defines the initial object in the category
of~$\P$-algebras. In this context an augmented $\P$-algebra refers
simply to an algebra $A$ equipped with a fixed $\P$-algebra morphism
$\epsilon: A\rightarrow\F$. The kernel of this augmentation morphism
defines the augmentation ideal of $A$, denoted by $\bar{A}$. Recall
that $\bar{A}$ is an algebra over the \emph{reduced operad}
$\overline{\P}$\glossary{$\overline{\P}$}\index{unital
operad!reduced operad of a} such that $\overline{\P}(r) = 0$ if
$r=0$ and $\overline{\P}(r) = \P(r)$ otherwise. Furthermore, the map
$A\mapsto\bar{A}$ defines an equivalence of categories from the
category of augmented $\P$-algebras to the category of
$\overline{\P}$-algebras.

On the other hand, we consider always \emph{unitary operads}
$\P$\index{operad!unitary}\index{unitary!operad}, equipped with a
unit operation $1\in\P(1)$\index{operad!unit operation of
an}\index{unit!operation of an operad}. In the literature, the
unital and unitary terminologies are used interchangeably. The
reader should not be confused: in this memoir, we use always these
terminologies in the sense specified in this introduction.

As in~\cite{Bar1}, we may consider connected operads for which the
module $\P(1)$ is spanned by~$1\in\P(1)$\index{operad!connected}. In
general the connectedness assumption is not necessary but it can
simplify certain constructions.

Naturally, a morphism of unital operads is assumed to preserve the
unital operation $*\in\P(0)$. Accordingly, there is an initial
object in the category of unital operads $*$ defined explicitly by
$*(r) = \F$ for $r = 0,1$ and $*(r) = 0$ otherwise. The category of
augmented algebras over this initial operad is equivalent to the
category of dg-modules since any augmented algebra over $*$ has the
form $V_+ = \F\oplus V$ for a dg-module $V$. Observe that the
category of unital operads is also endowed with a terminal object
given by the commutative operad $\C$
(see~\ref{item:UnitaryOperadAdjunction}).

The aim of this section is to fix our conventions in regard to
unital and Hopf operads and to prove the fundamental properties of
the categories formed by these objects. More specifically, as
explained in the memoir introduction, we introduce a new model
structure for unital operads which gives an appropriate framework
for the constructions of Hopf operad actions on the bar complex.
Recall briefly that our model category structure is deduced from a
generalization of the classical Reedy constructions and differs from
the model structure defined in~\cite{BergerMoerdijk} like classical
Reedy model structures differ from adjoint model structures.

Here is the plan of this section. For our purposes we need more
insights into the structure of the underlying $\Sigma_*$-module of a
unital-operad and we devote the next
subsection~\ref{subsection:LambdaModulesIntroduction} to this topic.
Then we perform the construction of the Reedy model structure of
unital operads in two steps: first,
in~\ref{subsection:ModelCategoryLambdaModules} we define a model
structure at the $\Sigma_*$-module level; then
in~\ref{subsection:UnitalOperadsReedyModelStructure} we transfer
this structure to unital operads by the classical adjunction
process. In the last
subsection~\ref{subsection:HopfOperadsIntroduction} we generalize
these constructions to unital Hopf operads after recalls and
precisions on this notion.

\subsection{On $\Lambda_*$-modules and unital operads}\label{subsection:LambdaModulesIntroduction}
Classically, one observes that the operadic composites with a unital
operation $\partial_i(p) = p\circ_i *$\glossary{$\partial_i$},
defined for $i = 1,\dots,r$ if $p\in\P(r)$, provide the underlying
$\Sigma_*$-module of a unital operad $\P$ with morphisms
$\partial_i: \P(r)\rightarrow\P(r-1)$ such that
$\partial_i\partial_j =
\partial_{j-1}\partial_i$ for $i<j$. The structure of a
\emph{preoperad}, introduced in~\cite{Berger}, consists precisely of
a $\Sigma_*$-module $M$ equipped with such morphisms $\partial_i:
M(r)\rightarrow M(r-1)$. As explained in \emph{loc. cit.} one can
observe that the structure of a preoperad is equivalent to a
contravariant functor $M: \Lambda_*^{\op}\rightarrow\Mod$, where
$\Lambda_*$\glossary{$\Lambda_*$} denotes the category formed by the
finite sets $\underline{r} = \{1,\dots,r\}$ and the injective maps
$u: \{1,\dots,r\}\rightarrow\{1,\dots,s\}$. In this memoir we prefer
to adopt the terminology of a
\emph{$\Lambda_*$-module}\index{$\Lambda_*$-module} for this
structure which is more consistent with our terminology of a
\emph{$\Sigma_*$-module}\index{$\Sigma_*$-module} for a collection
of $\Sigma_r$-modules. Similarly, the set of morphisms $u:
\{1,\dots,r\}\rightarrow\{1,\dots,s\}$ in $\Lambda_*$ is denoted by
$\Lambda_{r}^{s}$\glossary{$\Lambda_{r}^{s}$} so that
$\Lambda_{r}^{r} = \Sigma_r$. In this memoir, we consider tacitely
only right $\Lambda_*$-modules, that are contravariant functors on
the category $\Lambda_*$, in contrast to left $\Lambda_*$-modules
that are covariant functors. Accordingly, the $\Lambda_*$-action on
a $\Lambda_*$-module restricts naturally to a right action of the
symmetric groups $\Sigma_*$. On the other hand, we assume by
convention that a $\Sigma_*$-module is acted on by permutations on
the left, but left and right actions are equivalent for group
actions and this convention difference does not create any actual
difficulty.

In this subsection we make explicit the relationship between
$\Sigma_*$-modules and $\Lambda_*$-modules and between non-unital
operads and unital operads. To be precise, we consider the following
categories:
\begin{enumerate}
\item
the category $\Sigma_*\Mod_0$\glossary{$\Sigma_*\Mod_0$} formed by
the \emph{non-unital
$\Sigma_*$-modules}\index{$\Sigma_*$-module!non-unital}\index{non-unital!$\Sigma_*$-module}
-- the $\Sigma_*$-modules $M$ such that $M(0) = 0$; the category
$\Sigma_*\Mod^1_0$\glossary{$\Sigma_*\Mod^1_0$} formed by the
unitary objects of
$\Sigma_*\Mod_0$\index{$\Sigma_*$-module!unitary}\index{unitary!$\Sigma_*$-module}
-- a unitary object in $\Sigma_*\Mod_0$ consists of a
$\Sigma_*$-module $M\in\Sigma_*\Mod_0$ equipped with a unit element
$1\in M(1)$; and the category
$\Sigma_*\Mod^1_0/\overline{\C}$\glossary{$\Sigma_*\Mod^1_0/\overline{\C}$}
formed by unitary objects $M\in\Sigma_*\Mod^1_0$ equipped with an
augmentation $\epsilon: M\rightarrow\overline{\C}$ over the
underlying $\Sigma_*$-module of the reduced commutative operad
$\overline{\C}$\glossary{$\epsilon$};
\item
the category
$\Lambda^{\op}_*\Mod_0$\glossary{$\Lambda^{\op}_*\Mod_0$} formed by
the \emph{non-unital
$\Lambda_*$-modules}\index{$\Lambda_*$-module!non-unital}\index{non-unital!$\Lambda_*$-module}
-- the $\Lambda_*$-modules $M$ such that $M(0) = 0$; the category
$\Lambda^{\op}_*\Mod^1_0$\glossary{$\Lambda^{\op}_*\Mod^1_0$} formed
by the unitary objects of
$\Lambda^{\op}_*\Mod_0$\index{$\Lambda_*$-module!unitary}\index{unitary!$\Lambda_*$-module}
-- a unitary object in $\Lambda^{\op}_*\Mod_0$ consists of a
$\Lambda_*$-module equipped with a unit element $1\in M(1)$; and the
category
$\Lambda^{\op}_*\Mod^1_0/\overline{\C}$\glossary{$\Lambda^{\op}_*\Mod^1_0/\overline{\C}$}
formed by unitary objects $M\in\Sigma_*\Mod^1_0$ equipped with an
augmentation $\epsilon:
M\rightarrow\overline{\C}$\glossary{$\epsilon$};
\item
the category $\Op^1_0$\glossary{$\Op^1_0$} formed by the
\emph{non-unital unitary operads
$\P$}\index{operad!non-unital}\index{non-unital!operad} -- the
unitary operads $\P$ such that $\P(0) = 0$; and the category
$\Op^1_0/\overline{\C}$\glossary{$\Op^1_0/\overline{\C}$} formed by
the operads $\P\in\Op^1_0$ equipped with an augmentation $\epsilon:
M\rightarrow\overline{\C}$ over the reduced commutative operad
$\overline{\C}$\glossary{$\epsilon$};
\item
the category $\Op^1_*$\glossary{$\Op^1_*$} formed by the
\emph{unital unitary operads $\P$} -- the unitary operads $\P$ such
that $\P(0) = \F *$ for a distinguished operation $*\in\P(0)$.
\end{enumerate}
Observe that the unit of a
$\Sigma_*$-module\index{$\Sigma_*$-module!unit element of
a}\index{unit!element of a $\Sigma_*$-module} is equivalent to a
$\Sigma_*$-module morphism $\eta: I\rightarrow M$, where $*$ denotes
the non-unital $\Sigma_*$-module such that $I(1) = \F$ and $I(r) =
0$ for $r\not=1$. Accordingly, the category of augmented unitary
$\Sigma_*$-modules $\Sigma_*\Mod^*_0/\overline{\C}$ is defined by a
comma category and similarly for the category of augmented unitary
$\Lambda_*$-modules $\Lambda^{\op}_*\Mod^*_0/\overline{\C}$.

The connections between these categories can be summarized by a
diagram of categorical adjunctions
\begin{equation*}
\xymatrix{
\Sigma_*\Mod^1_0/\overline{\C}\ar@<1mm>[r]^(0.55){\FOp}\ar@<2mm>[d]\ar@<-2mm>[d]
&
\Op^1_0/\overline{\C}\ar@<1mm>[l]\ar@<2mm>[d]\ar@<-2mm>[d] \\
\Lambda^{\op}_*\Mod^1_0/\overline{\C}\ar@<1mm>[r]^(0.55){\FOp_*}\ar[u]
& \Op^1_*\ar@<1mm>[l]\ar[u] }.
\end{equation*}
For the upper horizontal adjunction we consider the obvious
forgetful functor which maps a non-unital operad $\P$ to its
underlying $\Sigma_*$-module. The unit element is given by the unit
operation~$1\in\P(1)$. This functor has a left-adjoint defined by
the classical free operad functor $M\mapsto\FOp(M)$ -- we refer
to~\ref{FreeNonunitalOperad} for recalls and more precision. In our
context this free object is equipped with a natural operad morphism
$\epsilon: \FOp(M)\rightarrow\overline{\C}$ induced by the
$\Sigma_*$-module augmentation $\epsilon: M\rightarrow\overline{\C}$
so that $\FOp$ induces a functor $\FOp:
\Sigma_*\Mod^1_0/\overline{C}\rightarrow\Op^1_0/\overline{\C}$. For
the bottom horizontal adjunction we consider simply unital versions
of these functors -- we refer to~\ref{item:FreeUnitalOperads} for
explicit definitions. For vertical adjunctions we consider the
canonical forgetful functor from $\Lambda_*$-modules to
$\Sigma_*$-modules and the functor $\P\mapsto\overline{\P}$ which
maps a unital operad $\P$ to the associated reduced operad
$\overline{\P}$ defined in the introduction of this section. One
observes readily that these functors preserve colimits and limits
and hence have both a left and a right adjoint. At the module level
this adjunction is an instance of a general adjunction relation
\begin{equation*}
\xymatrix{
\R^{\op}_*\Mod\ar@<14pt>[r]^{\phi_!}\ar@<-14pt>[r]^{\phi_*} &
\ar@<0pt>[l]_{\phi^! = \phi^*}\SOp^{\op}_*\Mod }
\end{equation*}
associated to a morphism of small categories $\phi:
\R_*\rightarrow\SOp_*$. Recall simply that the (right)
$\SOp_*$-module associated to a (right) $\R_*$-module $M$, also
denoted by $\phi_! M = M\otimes_{\R_*}\SOp_*$, respectively $\phi_*
M = \DGHom_{\R_*}(\SOp_*,M)$, can be defined by a coend,
respectively end, formula. Namely:
\begin{multline*}
M\otimes_{\R_*}\SOp_*(y) = \int_{x\in\R_*} M(x)\otimes\SOp_{y}^{\phi(x)}, \\
\text{respectively}\quad\DGHom_{\R_*}(\SOp_*,M)(y) = \int^{x\in\R_*}
M(x)^{\SOp_{y}^{\phi(x)}},
\end{multline*}
where the notation $K\otimes V$, respectively $V^K$, refers to the
classical tensor product, respectively mapping object, of a
dg-module $V$ with a set $K$. (As explained in \cite{Bar1}, our
conventions for ends and coends are converse to the usual one: we
use superscripts for ends and subscripts for coends.)

To begin our constructions we study the connection between the
category of unital and non-unital operads.

\subsubsection{Unitary operads versus non-unital operads}\label{item:UnitaryOperadAdjunction}
In fact, if we reverse the definition of the reduced
operad~$\overline{\P}$\glossary{$\overline{\P}$}, then, for a unital
operad $\P$, we obtain the relation
\begin{equation*}
\P(r) = \begin{cases} \F & \text{if $r=0$}, \\
\overline{\P}(r) & \text{otherwise}.
\end{cases}
\end{equation*}
Furthermore, the operad composition products $\circ_i:
\P(s)\otimes\P(t)\rightarrow\P(s+t-1)$ are determined for $s,t\geq
1$ by the composition products of the reduced operads
$\overline{\P}$. As a consequence, the structure of $\P$ is
completely determined by the associated reduced operad
$\overline{\P}$ and by morphisms $\partial_i:
\overline{\P}(r)\rightarrow\overline{\P}(r-1)$ that keep track on
the composite with the unital operation $*\in\P(0)$ at the level of
the reduced operad $\overline{\P}$.

Observe in addition that a unital operad comes equipped with
augmentation morphisms $\epsilon:
\P(r)\rightarrow\F$\glossary{$\epsilon$} which map an
operation~$p\in\P(r)$ to the composite $p(*,\dots,*)\in\P(0)$ and
such that $\epsilon\partial_i = \epsilon$ for all $i$. As a
byproduct, these augmentation morphisms give rise to an operad
morphism $\epsilon: \P\rightarrow\C$, where $\C$ denotes the operad
of commutative algebras (recall that $\C(r) = \F$ for $r\in\N$), so
that any unital operad $\P$ is canonically augmented over the
commutative operad $\C$. Hence the commutative operad $\C$ defines
the terminal object in the category of unital operad as mentioned in
the introduction of this section.

Finally, our observations prove that a unital operad $\P$ is
equivalent to a non-unital operad $\overline{\P}$ equipped with
morphisms $\partial_i:
\overline{\P}(r)\rightarrow\overline{\P}(r-1)$ for $i = 1,\dots,r$
and $r\geq 2$, such that the associativity and commutativity
properties of operad composition products are satisfied for the
operations $\partial_i = -\circ_i *$, and with an augmentation
morphism $\epsilon: \overline{\P}\rightarrow\F$ that preserves all
operadic composites including the composites with the unital
operation $*$. As a byproduct, one can prove that the functor
$\P\mapsto\overline{\P}$ creates small colimits in the category of
unital operads. Explicitly, for a diagram of unital
operads~$\P_\alpha$ one observes readily that the colimit
$\colim_\alpha\overline{\P}_\alpha$ is equipped with canonical
operations $\partial_i$ and with an augmentation $\epsilon:
\colim_\alpha\overline{\P}_\alpha\rightarrow\F$ induced by the
corresponding operations of the unital operads~$\P_\alpha$. As a
consequence, this colimit is associated to a unital operad an this
operad defines necessarily the colimit of the operads $\P_\alpha$ in
the category of unital operads.

\subsubsection{The adjunction between $\Sigma_*$-modules and non-unital
operads}\label{FreeNonunitalOperad} As explained in the introduction
of this subsection, we consider the obvious forgetful functor from
the category of non-unital unitary operads $\Op^*_0$ to the category
of non-unital unitary $\Sigma_*$-modules $\Sigma_*\Mod^*_0$ which
maps an operad $\P\in\Op^*_0$ to its underlying $\Sigma_*$-module.

The left adjoint of this functor is given by a variant of the
classical free operad
$\FOp(M)$\glossary{$\FOp(M)$}\index{free!non-unital
operad}\index{operad!non-unital!free}\index{operad!free non-unital}
defined in the literature. Namely, in the unitary context, the
universal morphism $\eta: M\rightarrow\FOp(M)$ is supposed to map
the unit of~$M$ to the unit of the free operad and hence to define a
morphism in the category of unitary $\Sigma_*$-modules. Furthermore,
the adjunction relation supposes that an operad morphism $\phi_f:
\FOp(M)\rightarrow\P$ is associated to a $\Sigma_*$-module morphism
$f: M\rightarrow\P$ that preserves unit elements. Equivalently, for
the universal property, one assumes that the $\Sigma_*$-module
morphisms $f: M\rightarrow\P$ that preserve unit elements have a
unique factorization
\begin{equation*}
\xymatrix@C=4mm{ M\ar[rd]_f\ar[rr]^\phi & & \P \\ &
\FOp(M)\ar@{-->}[ru]_{\phi_f} &
\\ }
\end{equation*}
such that $\phi_f: \FOp(M)\rightarrow\P$ is an operad morphism.
Accordingly, the free operad associated to a unitary
$\Sigma_*$-module has to be defined by a quotient of the classical
free operad associated to a non-unitary $\Sigma_*$-module.
Explicitly, for a unitary $\Sigma_*$-module $M$ we make the internal
unit $1\in M(1)$ equivalent to the unit of the free operad
in~$\FOp(M)$. Up to this quotient process we refer to our
article~\cite{OperadTextbook}, from which we borrow our conventions,
for an explicit construction of the free operad that follows closely
the original constructions of~\cite{GetzlerGoerss,GinzburgKapranov}.
The definition of tree structures that occur in this construction is
also recalled in~\ref{subsection:CellMetricTrees} for the purposes
of Boardman-Vogt' $W$-construction. Recall that $\FOp(M)$ can be
defined simply by the modules spanned by formal expressions
$(\dots((x_1\circ_{i_2} x_2)\circ_{i_3}\dots x_{l-1})\circ_{i_l}
x_l$ which represent composites of generators $x_1\in
M(n_1),\ldots,x_l\in M(n_l)$.

\subsubsection{On unitary and non-unitary $\Sigma_*$-modules}
In fact, in our context, we could replace unitary $\Sigma_*$-modules
by equivalent non-unitary $\Sigma_*$-modules. Explicitly, recall
that we consider $\Sigma_*$-module which are augmented over the
reduced commutative operad $\overline{\C}$. For such
$\Sigma_*$-modules we have a canonical splitting $M(1) =
\widetilde{M}(1)\oplus\F 1$, where $\widetilde{M}(1)$ denotes the
cokernel of the morphism $\eta: \F\rightarrow M(1)$ defined by the
unit element $1\in M(1)$. Indeed, since $\overline{\C}(1) = \F$, the
morphism $\eta: \F\rightarrow M(1)$ has a canonical left-inverse
given by the augmentation $\epsilon:
M(1)\rightarrow\overline{\C}(1)$.

Consequently, if we let $\widetilde{M}$ denote the non-unitary
$\Sigma_*$-module such that $\widetilde{M}(0) = 0$,
$\widetilde{M}(1) = \coker(\eta: \F\rightarrow M(1))$ and
$\widetilde{M}(r) = M(r)$ for $r\geq 1$, then the map
$M\mapsto\widetilde{M}$ defines an equivalence from the category
$\Sigma_*\Mod^*_0/\overline{\C}$, formed by unitary
$\Sigma_*$-modules augmented over the reduced commutative operad, to
the category $\Sigma_*\Mod_0/\widetilde{\C}$, formed by non-unitary
$\Sigma_*$-modules which are augmented over the augmentation ideal
of the commutative operad. Furthermore, the free operad~$\FOp(M)$
associated to a unitary $\Sigma_*$-module~$M$ is clearly isomorphic
to the free operad~$\FOp(\widetilde{M})$ associated to the
equivalent non-unitary $\Sigma_*$-module~$\widetilde{M}$.

For our purpose it is more natural to deal with unitary objects and
we do not use this relationship. Nevertheless one can conclude
immediately from these observations that the free operad functor
$M\mapsto\FOp(M)$ satisfies the same homological properties as the
free operad $\FOp(\widetilde{M})$. For instance, the functor
$M\mapsto\FOp(M)$ maps weak-equivalences of unitary
$\Sigma_*$-modules to weak-equivalences of operads.

\subsubsection{The adjunction between $\Lambda_*$-modules and unital operads}\label{item:FreeUnitalOperads}
The functor
$M\mapsto\FOp_*(M)$\glossary{$\FOp_*(M)$}\index{free!unital
operad}\index{operad!unital!free}\index{operad!free unital}
considered in the introduction of this subsection is defined as the
left adjoint of the functor $\P\mapsto\overline{\P}$ from the
category of unital operads $\Op^1_*$ to the category of augmented
non-unital unitary $\Lambda_*$-modules
$\Lambda^{\op}_*\Mod_0/\overline{\C}$. The operad $\FOp_*(M)$
associated to a given non-unital module can also be characterized by
the usual universal property. In fact, the reduced operad associated
to $\FOp_*(M)$ can be identified with the usual free
operad~$\FOp(M)$ so that
\begin{equation*}
\FOp_*(M)(r) = \begin{cases} \F & \text{if $r=0$}, \\
\FOp(M)(r) & \text{otherwise}.
\end{cases}
\end{equation*}
Accordingly, for $s,t\geq 1$ the composition products $\circ_i:
\FOp_*(M)(s)\otimes\FOp_*(M)(t)\rightarrow\FOp_*(M)(s+t-1)$ are
given by the usual formal composition process of the free operad.
The composites with the unital operation $-\circ_i *:
\FOp_*(M)(r)\rightarrow\FOp_*(M)(r-1)$ are induced by the morphisms
$\partial_i: M(r)\rightarrow M(r-1)$ for $r\geq 2$ and by the
augmentation morphism $\epsilon: M(r)\rightarrow\F$ for $r = 1$.
Explicitly, for a generator $x\in M(r)$ we set $x\circ_1 * =
\epsilon(x)$ if $r = 1$, $x\circ_i * =
\partial_i(x)$ for $r\geq 2$ and we observe that these operations $-\circ_i *:
M(r)\rightarrow\FOp_*(M)(r-1)$ have a unique extension to
$\FOp_*(M)$ that satisfies the commutativity properties of operad
composition products.

\subsubsection{On unital $\Lambda_*$-modules}
Sometimes we can find convenient to consider \emph{unital unitary
$\Lambda_*$-modules}\index{unital!$\Lambda_*$-module}\index{$\Lambda_*$-module!unital}
-- explicitly, unital $\Lambda_*$-modules $N$ such that $N(0) = \F
*$ for a distinguished element $*\in N(0)$ such that $\partial_1(1)
=
*$. Equivalently, a unital $\Lambda_*$-module consists of a
$\Lambda_*$-module $N$ equipped with a unit morphism $\eta:
*\rightarrow N$, where $*$ denotes the underlying
$\Lambda_*$-module of the initial unital operad, such that $\eta:
\F\rightarrow N(0)$ is iso. Clearly, the underlying
$\Lambda_*$-module of a unital operad is unital.

Observe that a unital $\Lambda_*$-module $N$ is equipped with a
canonical augmentation $\epsilon:
N\rightarrow\C$\glossary{$\epsilon$} like a unital operad.
Explicitly, the map $p\mapsto p(*,\dots,*)$ defined for a unital
operad in~\ref{item:UnitaryOperadAdjunction} can be identified with
the $\Lambda_*$-module operation $\eta_0^*: N(r)\rightarrow N(0)$
associated to the initial map $\eta_0:
\emptyset\rightarrow\{1,\dots,r\}$\glossary{$\eta_0$}. If we assume
$N(0) = \F$, then this operation gives a morphism from $N$ to the
constant $\Lambda_*$-module~$\F$, that represents the underlying
$\Lambda_*$-module of the commutative operad~$\C$.

The definition of the reduced operad can also be extended to unital
$\Lambda_*$-modules. Explicitly, for a unital $\Lambda_*$-module $N$
we consider the $\Lambda_*$-module $\overline{N}$ such that
$\overline{N}(0) = 0$ and $\overline{N}(r) = N(r)$ for $r\geq 1$.
Clearly, this object forms a non-unital unitary $\Lambda_*$-module.
Since any unital $\Lambda_*$-module is augmented over $\C$, the map
$N\mapsto\overline{N}$ defines a functor $(-)^{-}:
\Lambda_*^{\op}\Mod^1_*\rightarrow\Lambda_*^{\op}\Mod^1_0/\overline{\C}$,
where $\Lambda_*^{\op}\Mod^1_*$ denotes the category of unital
unitary $\Lambda_*$-modules.

The functor $N\mapsto\overline{N}$ has clearly a left adjoint
$M\mapsto M_+$. Explicitly, the unital $\Lambda_*$-module $M_+$
associated to an augmented non-unital $\Lambda_*$-module $M$ is
defined by $M_+(0) = \F$ and $M_+(r) = M(r)$ for $r\geq 1$. The
operations $\partial_i: M_+(r)\rightarrow M_+(r)$ are induced by the
corresponding operations of~$M$ and by the augmentation $\epsilon:
M(1)\rightarrow\F$. Clearly, these functors yield inverse adjoint
equivalences of categories:
\begin{equation*}
\xymatrix{ (-)_+: \Lambda^{\op}_*\Mod^1_0\ar@<2pt>[r]^{\simeq} &
\ar@<2pt>[l]^{\simeq}\Lambda^{\op}_*\Mod^1_* :(-)^{-}}.
\end{equation*}

\subsubsection{Connected unital operads}\label{item:ConnectedUnitaryOperads}
As mentioned in the section introduction, we may consider
\emph{connected unital operads} $\P$ such that the module $\P(1)$ is
spanned by the operad unit
$1\in\P(1)$\index{operad!connected}\index{operad!unitary!connected}.
Let $\Op^*_*$\glossary{$\Op^*_*$} denote the category formed by
these objects. The category embedding $\itr^1_*:
\Op^*_*\rightarrow\Op^1_*$\glossary{$\itr^1_*$} preserves clearly
colimits and limits and hence admits both a right and a left adjoint
denoted by $\str^1_*:
\Op^1_*\rightarrow\Op^*_*$\glossary{$\str^1_*$} and $\ctr^1_*:
\Op^1_*\rightarrow\Op^*_*$\glossary{$\ctr^1_*$} respectively. We
give an explicit construction of these functors in the next
paragraph. We check in addition that the adjunction unit $\eta:
\P\rightarrow\str^1_*\itr^1_*(\P)$, respectively the adjunction
augmentation $\epsilon: \ctr^1_*\itr^1_*(\Q)\rightarrow\Q$, is an
isomorphism for all connected unital operad $\P$, respectively $\Q$.
Consequently, we have an adjunction ladder:
\begin{equation*}
\xymatrix{ \Op^*_*\ar[r]^{\itr^1_*} &
\Op^1_*\ar@<-6mm>[l]_{\ctr^1_*}\ar@<+6mm>[l]_{\str^1_*} }
\end{equation*}
such that $\ctr^1_*\itr^1_* = \Id = \str^1_*\itr^1_*$. In
forthcoming constructions non-connected operads $\Q$ can be replaced
by the associated connected object $\str^1_*(\Q)$. Therefore we
could restrict ourself to connected unital operads.

To conclude, the map $\P\mapsto\overline{\P}$ and the free operad
$M\mapsto\FOp_*(M)$ restrict to adjoint functors
\begin{equation*}
\FOp_*:
\dg\Lambda^{\op}_*\Mod^*_0/\overline{\C}\rightleftarrows\dg\Op^*_*
:(-)^{-},
\end{equation*}
where $\Lambda^{\op}_*\Mod^*_0$ denotes the category formed by the
connected non-unital unitary $\Lambda_*$-modules, the non-unital
unitary $\Lambda_*$-modules $M$ such that $M(0) = 0$ and $M(1) =
\F$.

\subsubsection{The adjunction between connected and non-connected unital
objects}\label{item:ConnectedAdjointOperads} We give an explicit
construction of the connected operad $\str^1_*(\Q)\in\Op^*_*$,
respectively $\ctr^1_*(\P)\in\Op^*_*$, associated to a unital operad
$\Q\in\Op^1_*$, respectively $\P\in\Op^1_*$. In fact, we define an
adjunction ladder for the underlying categories of unital unitary
$\Lambda_*$-modules:
\begin{equation*}
\xymatrix{ \Lambda^{\op}_*\Mod^*_*\ar[r]^{\itr^1_*} &
\Lambda^{\op}_*\Mod^1_*\ar@<-6mm>[l]_{\ctr^1_*}\ar@<+6mm>[l]_{\str^1_*}
}.\glossary{$\itr^1_*$}\glossary{$\ctr^1_*$}\glossary{$\str^1_*$}
\end{equation*}
One can check readily that the connected unital unitary
$\Lambda_*$-module $\str^1_*(\Q)\in\Mod^*_*$ associated to a unital
operad $\Q\in\Op^1_*$ forms a suboperad of~$\Q$ and similarly the
$\Lambda_*$-module $\ctr^1_*(\P)\in\Mod^*_*$ forms a quotient operad
of~$\P$. Hence, for operads, the map $\Q\mapsto\str^1_*(\Q)$,
respectively $\P\mapsto\str^1_*(\P)$, gives rise to a right,
respectively left, adjoint of the category embedding $\itr^1_*:
\Op^*_*\rightarrow\Op^1_*$.

For $r\geq 1$ and $i = 1,\dots,r$, we consider the
$\Lambda_*$-module operation $\eta_i^*: N(r)\rightarrow N(1)$
associated to the map $\eta_i:
\{1\}\rightarrow\{1,\dots,r\}$\glossary{$\eta_i$} such that
$\eta_i(1) = i$. For an operad we have equivalently $\eta_i^*(p) =
p(*,\dots,1,\dots,*)$, where unital operations $*$ are substituted
to the entries $k\not=i$ of the operation. For $r\geq 1$, let
$\str^1_*(N)(r)$ denotes the submodule of $N(r)$ defined by the
pullback diagram
\begin{equation*}
\xymatrix{ \str^1_*(N)(r)\ar@{-->}[r]\ar@{-->}[d] &
N(r)\ar[d]^{(\eta_i^*)_i} \\
\F^{\times r}\ar[r]^{1^{\times r}} & N(1)^{\times r} }
\end{equation*}
in which we consider the morphism $1: \F\rightarrow N(1)$ defined by
the unit element $1\in N(1)$. One checks readily that these
dg-modules $\str^1_*(N)(r)$ define a $\Lambda_*$-submodule of~$N$.
Furthermore, any morphism of unitary unital $\Lambda_*$-module $f:
M\rightarrow N$ such that $M$ is connected factorizes through
$\str^1_*(N)$ since we have a commutative diagram
\begin{equation*}
\xymatrix{ M(r)\ar[r]^f\ar[d]^{\eta_i^*} &
N(r)\ar[d]^{\eta_i^*} \\
M(1) = \F\ar[r]^{1} & N(1) }
\end{equation*}
for all $i = 1,\dots,r$. Therefore the functor $N\mapsto\str^1_*(N)$
satisfies the adjunction relation
$\Hom_{\Lambda^{\op}_*\Mod^1_*}(\itr^1_*(M),N)
=\Hom_{\Lambda^{\op}_*\Mod^*_*}(M,\str^1_*(N))$. According to this
construction, we have clearly $M = \str^1_*\itr^1_*(M)$ for any
connected unital unitary $\Lambda_*$-module~$M$.

The other connected $\Lambda_*$-module $\ctr^1_*(M)$ associated to a
unital unitary $\Lambda_*$-module $M$ can clearly be defined by
$\ctr^1_*(M)(r) = \F$ for $r = 0,1$ and $\ctr^1_*(M)(r) = M(r)$ for
$r\geq 2$. For $r = 2$, the operations $\partial_1,\partial_2:
\ctr^1_*(M)(2)\rightarrow\F$ are given by the augmentation morphism
$\epsilon: M(2)\rightarrow\F$. The relation $\ctr^1_*\itr^1_*(N) =
N$, that holds for a connected $\Lambda_*$-module~$N$, is also
immediate from this definition.

\subsection{The Reedy model structure for $\Lambda_*$-modules}\label{subsection:ModelCategoryLambdaModules}
The aim of this section is to prove the following theorem:

\begin{thm}\label{thm:LambdaModuleModelStructure}
The category of ($\N$ or $\Z$-graded) dg-$\Lambda_*$-modules
$\dg\Lambda^{\op}_*\Mod$ is equipped with the structure of a
cofibrantly generated model category such that a morphism $f:
M\rightarrow N$ is a weak-equivalence, respectively a cofibration,
if $f$ defines a weak-equivalence, respectively a cofibration, in
the category of dg-$\Sigma_*$-modules $\dg\Sigma_*\Mod$.
\end{thm}

In this statement we consider the category $\dg\Lambda^{\op}_*\Mod$
formed by all $\Lambda_*$-modules which do not satisfy necessarily
$M(0) = 0$. On the other hand, for our purposes we need a model
structure on the category of non-unital $\Lambda_*$-modules
$\dg\Lambda^{\op}_*\Mod_0$ and on the associated comma category of
augmented unitary $\Lambda_*$-modules
$\dg\Lambda^{\op}_*\Mod^1_0/\overline{\C}$. In fact, we use the
following claim in order to equip $\dg\Lambda^{\op}_*\Mod_0$ with
the structure of a model subcategory of~$\dg\Lambda^{\op}_*\Mod$:

\begin{prop}\label{prop:ConnectedAdjoints}
The category embedding $\itr^{\Lambda}_0:
\dg\Lambda^{\op}_*\Mod_0\hookrightarrow\dg\Lambda^{\op}_*\Mod$
admits a left and a right adjoint
\begin{equation*}
\xymatrix{ \dg\Lambda^{\op}_*\Mod_0\ar[r]^{\itr^{\Lambda}_0} &
\dg\Lambda^{\op}_*\Mod\ar@<-6mm>[l]_{\ctr^{\Lambda}_0}\ar@<+6mm>[l]_{\str^{\Lambda}_0}
}
\end{equation*}
such that $\ctr^{\Lambda}_0\itr^{\Lambda}_0 = \Id =
\str^{\Lambda}_0\itr^{\Lambda}_0$. Furthermore, the functor
$\itr^{\Lambda}_0\ctr^{\Lambda}_0$ preserves cofibrations and all
weak-equivalences in~$\dg\Lambda^{\op}_*\Mod$.
\end{prop}

\begin{proof}
The proof of this proposition is straightforward. The construction
of non-unital $\Lambda_*$-modules $\ctr^{\Lambda}_0(M)$,
respectively $\str^{\Lambda}_0(N)$, is similar to the construction
of the connected unital $\Lambda_*$-module $\ctr^1_*(M)$,
respectively $\str^1_*(N)$, given
in~\ref{item:ConnectedAdjointOperads}.

Explicitly, for $\ctr^{\Lambda}_0(M)$ we set simply
$\ctr^{\Lambda}_0(M)(r) = 0$ for $r = 0$ and $\ctr^{\Lambda}_0(M)(r)
= M(r)$ for $r\geq 1$. According to this construction, the functor
$\ctr^{\Lambda}_0$ preserves clearly cofibrations, acyclic
cofibrations and all week-equivalences in~$\dg\Lambda^{\op}_*\Mod$.

For $\str^{\Lambda}_0(N)$ we consider the operation $\eta_0^*:
N(r)\rightarrow N(0)$ associated to the map $\eta_0:
\emptyset\rightarrow\{1,\dots,r\}$. Then, for $r\geq 0$, we set
explicitly $\str^{\Lambda}_0(N)(r) = \ker(\eta_0^*: N(r)\rightarrow
N(0))$. More categorically, the module $\str^{\Lambda}_0(N)(r)$ can
be defined by the pullback diagram
\begin{equation*}
\xymatrix{ \str^{\Lambda}_0(N)(r)\ar@{-->}[r]\ar@{-->}[d] &
N(r)\ar[d]^{\eta_0^*} \\
0\ar[r]^{0} & N(0) }.
\end{equation*}
\end{proof}

Then, assuming theorem~\ref{thm:LambdaModuleModelStructure}, we
obtain:

\begin{prop}\label{prop:NonunitalLambdaModulesModelStructure}
The category of non-unital dg-$\Lambda_*$-modules
$\dg\Lambda^{\op}_*\Mod_0$\linebreak forms a model subcategory
of~$\dg\Lambda^{\op}_*\Mod$ so that a morphism $f: M\rightarrow N$
is a weak-equivalence, respectively a cofibration, a fibration, of
non-unital dg-$\Lambda_*$-modules if and only if $f$ defines a
weak-equivalence, respectively a cofibration, a fibration,
in~$\dg\Lambda^{\op}_*\Mod$.

Furthermore, this category is cofibrantly generated by the morphisms
$\ctr^{\Lambda}_0(i):
\ctr^{\Lambda}_0(A)\rightarrow\ctr^{\Lambda}_0(B)$ associated to a
generating set of cofibrations, respectively acyclic cofibrations,
in $\dg\Lambda^{\op}_*\Mod$.
\end{prop}

\begin{proof}
This assertion is a straightforward consequence of the adjunction
relation and the invariance of~$\ctr^{\Lambda}_0$ with respect to
cofibrations and weak-equivalences.
\end{proof}

The category of augmented unitary $\Lambda^{\op}_*$-modules
$\dg\Lambda^{\op}_*\Mod^1_0/\overline{\C}$, which is a comma
category associated to $\dg\Lambda^{\op}_*\Mod_0$, is equipped with
a canonical induced model structure as usual.

\medskip
In fact, theorem~\ref{thm:LambdaModuleModelStructure} holds for a
category of $\Lambda_*$-objects in any cofibrantly generated model
category and not only in the category of dg-modules. Therefore we
make our arguments as general as possible though we consider
explicitly only dg-modules (and dg-coalgebras in the next section).

In a first stage we make explicit the definition of a fibration in
the category of $\Lambda_*$-modules. Then we prove that the lifting
properties (M4.i-ii) are satisfied, we specify generating
collections of cofibrations and acyclic cofibrations and we deduce
the factorization properties (M5.i-ii) from the classical small
objects argument. The properties (M1-3) are inherited from the
ground model category.

The general idea of our construction is provided by a generalization
of the classical Reedy model structure in the presence of
automorphisms. To be precise, in the usual definition of Reedy
structures one considers functors $F: I\rightarrow\C$ on a fixed
small category $I$ equipped with a grading and a decomposition $I =
\overrightarrow{I}\overleftarrow{I}$, where the direct Reedy
category $\overrightarrow{I}$, respectively the inverse Reedy
category $\overleftarrow{I}$, denote subcategories of $I$ formed by
collections of morphisms that increase, respectively decrease, the
grading. One can observe that classical constructions, given in
Reedy's original article~\cite{Reedy} and in the modern
monographs~\cite{Hirschhorn,HoveyTextBook}, still work for a
category $I$ that contains automorphisms and such that any morphism
$f\in I$ has a unique decomposition $f = \alpha u \beta$ in which
$\alpha\in\overrightarrow{I}$, $\beta\in\overleftarrow{I}$ and $u$
is an automorphism. This extension of Reedy structures does not seem
to occur in the literature in full generality though an instance is
supplied by the strict model category of $\Gamma$-spaces defined
in~\cite{BousfieldFriedlander}. Therefore we give detailed arguments
for the category $I = \Lambda_*^{\op}$. By contravariance, we have
in this case $\overleftarrow{I} = \overrightarrow{\Lambda}_*^{\op}$
for a Reedy direct subcategory of~$\Lambda_*$ and
$\overrightarrow{I}$ is trivial. The latter property simplifies some
constructions.

In the next sections we call the model structure supplied by
theorem~\ref{thm:LambdaModuleModelStructure} the \emph{Reedy model
structure}\index{Reedy!model category!of
$\Lambda_*$-modules}\index{$\Lambda_*$-module!Reedy model category
of $\Lambda_*$-modules} in order to distinguish this model category
from the classical adjoint model structure (the cofibrantly
generated model category of~\cite[Section 11.6]{Hirschhorn}) from
which it differs. To be precise, recall that the classical adjoint
model structure is defined by a transfer of model structures from
dg-modules to $\Lambda_*$-modules through a composite adjunction
\begin{equation*}
\xymatrix{ \dg\Mod^{\N}\ar@<14pt>[r]^{\phi_!} &
\ar@<0pt>[l]_{\phi_!}\dg\Sigma_*\Mod \ar@<14pt>[r]^{\psi_!} &
\ar@<0pt>[l]_{\psi^!}\dg\Lambda^{\op}_*\Mod }.
\end{equation*}
Explicitly, one let a morphism $f: M\rightarrow N$ defines a
weak-equivalence, respectively a fibration, in the adjoint model
category $\dg\Lambda^{\op}_*\Mod$ if $\phi^!\psi^!(f)$ defines a
weak-equivalence, respectively a fibration, in $\dg\Mod^{\N}$.
According to this definition, the Reedy model category and the
adjoint model category have the same weak-equivalences but different
cofibrations and fibrations. In fact, a Reedy fibration defines a
fibration in the adjoint model category and hence the Reedy model
category has less fibrations but more cofibrations than the adjoint
model category. As a byproduct, the identity functor yields a pair
of adjoint derived equivalences between the two homotopy categories.
Thus the Reedy model structure of $\Lambda_*$-modules is different
but Quillen equivalent to the adjoint model structure like a
classical Reedy category (see~\cite[Section 15.6]{Hirschhorn}).

\subsubsection{The category $\Lambda_*$}\glossary{$\Lambda_*$}
Recall that a $\Lambda_*$-module is equivalent to a contravariant
functor on the category $\Lambda_*$ formed by the finite sets
$\underline{r} = \{1,\dots,r\}$ and the injective maps $u:
\{1,\dots,r\}\rightarrow\{1,\dots,s\}$. This equivalence is also a
consequence of the decomposition of $\Lambda_*$ into a direct
subcategory and a subcategory of isomorphisms that we make explicit
in this paragraph.

Namely let $\overrightarrow{\Lambda}_*$ denote the subcategory of
$\Lambda$ whose morphisms
$\alpha\in\overrightarrow{\Lambda}_{r}^{s}$ are the non-decreasing
injections $\alpha: \{1,\dots,r\}\rightarrow\{1,\dots,s\}$. One
observes readily that any morphism $u\in\Lambda_{r}^{s}$ has a
unique $\overrightarrow{\Lambda}_*\Sigma_*$-decomposition $u =
\alpha\sigma$ in which $\sigma\in\Sigma_r$ and
$\alpha\in\overrightarrow{\Lambda}_{r}^{s}$. Moreover, the category
$\overrightarrow{\Lambda}_*$ is the category generated by the
injections $d^i: \{1,\dots,r-1\}\rightarrow\{1,\dots,r\}$ that avoid
$i\in\{1,\dots,r\}$ endowed with the relations $d^j d^i = d^i
d^{j-1}$ for $i<j$. With respect to permutations
$\sigma\in\Sigma_r$, we have also a relation of the form $\sigma d^i
= d^{\sigma(i)}
\partial_i(\sigma)$ which gives the
$\overrightarrow{\Lambda}_*\Sigma_*$-decomposition of $u = \sigma
d^i$. As a consequence, as stated in the
introduction~\ref{subsection:LambdaModulesIntroduction}, the
structure of a $\Lambda_*$-module is indeed equivalent to a
$\Sigma_*$-module equipped with operations $\partial_i:
M(r)\rightarrow M(r-1)$ that satisfy the relations above.

Anyway the category $\overrightarrow{\Lambda}_*$ forms clearly a
direct Reedy category and we have a generalized Reedy decomposition
$\Lambda_* = \overrightarrow{\Lambda}_*\Sigma_*$. As explained
previously, we adapt the classical construction of Reedy model
structures to this context. In the next paragraphs we let
$\Lambda_{*<r}$, respectively $\overrightarrow{\Lambda}_{*<r}$,
denote the comma category of morphisms $\alpha:
\underline{r}'\rightarrow\underline{r}$ in $\Lambda_*$, respectively
$\overrightarrow{\Lambda}_*$, such that $r'<r$. Observe that the
$\overrightarrow{\Lambda}_*\Sigma_*$-decomposition property implies
readily that $\overrightarrow{\Lambda}_{*<r}$ is cofinal in
$\Lambda_{*<r}$.

\subsubsection{Matching modules and fibrations of $\Lambda_*$-modules}\label{item:LambdaModuleFibrations}
In order to obtain an explicit characterization of fibrations we
need to introduce a notion of a \emph{matching object} in the
category of $\Lambda_*$-modules\index{$\Lambda_*$-module!matching
object of a}\index{matching object!of a $\Lambda_*$-module}.
Explicitly, the matching object of a $\Lambda_*$-module is defined
by the collection of dg-modules $\Match{M}(r)$, $r\in\N$, such that
\begin{equation*}
\Match{M}(r) = \mylim_{\alpha:
\underline{r}'\xrightarrow{<}\underline{r}} M(r'),
\end{equation*}
where the limit ranges over the category $\Lambda_{*<r}$ or,
equivalently, over the category $\overrightarrow{\Lambda}_{*<r}$
which is cofinal in $\Lambda_{*<r}$ by the observation of the
previous paragraph. As usual the dg-module $\Match{M}(r)$ can
equivalently be defined by an equalizer
\begin{equation*}
\Match{M}(r) = \ker\bigl(\xymatrix{ *+<2mm>{\prod_{1\leq i\leq r}
M(r-1)}\ar@<2pt>[r]^{d^0}\ar@<-2pt>[r]_{d^1} &
*+<2mm>{\prod_{1\leq i<j\leq r} M(r-2)} }\bigr),
\end{equation*}
where $d^0(x_i)_i = (\partial_i x_j)_{i<j}$ and $d^1(x_i)_i =
(\partial_{j-1} x_i)_{i<j}$.

The matching modules are endowed with natural morphisms $\mu:
M(r)\rightarrow\Match{M}(r)$ induced by the morphisms $\alpha^*:
M(r)\rightarrow M(r')$ on the component of the limit indexed by
$\alpha\in\Lambda_{r'}^{r}$. Explicitly, if we represent an element
of the limit by a collection $(x_\alpha)_{\alpha}$, where
$\alpha\in\Lambda_{r'}^{r}$, then, for $x\in M(r)$, we have
$\mu(x)_\alpha = \alpha^*(x)\in M(r')$. The matching module
$\Match{M}$ is also equipped with a canonical $\Lambda_*$-module
structure such that $\mu: M\rightarrow\Match{M}$ defines a morphism
of $\Lambda_*$-modules. Explicitly, for any morphism $u:
\underline{r}\rightarrow\underline{s}$, if we let $y =
(y_\beta)_\beta$ denote an element of $\Match{M}(s) = \mylim_{\beta}
M(s')$, then we have an associated element $u^*(y) =
(u^*(y)_\alpha)_\alpha$ in $\Match{M}(r) = \mylim_{\alpha} M(r')$
which can be defined by the collection: $u^*(y)_\alpha =
y_{u\alpha}$.

We let a morphism of $\Lambda_*$-modules $p: M\rightarrow N$ be a
fibration if, for all $r\in\N$, the natural morphism
\begin{equation*}
(\mu,p): M(r)\rightarrow\Match{M}(r)\times_{\Match{N}(r)} N(r)
\end{equation*}
defines a fibration in the category of dg-modules and hence in the
category of $\Sigma_r$-modules for the classical model structure of
$\Sigma_r$-modules.\index{Reedy!fibration!of
$\Lambda_*$-modules}\index{fibration!Reedy fibration of
$\Lambda_*$-modules}

\medskip
As stated previously, the properties (M1-3) of a model category are
inherited from the ground category of dg-modules. The first
non-trivial verification is supplied by the following claim:

\begin{claim}\label{claim:LambdaModulesLifting}
The properties (M4.i-ii) hold for the class of weak-equivalences,
cofibrations and fibrations specified in
theorem~\ref{thm:LambdaModuleModelStructure} and
paragraph~\ref{item:LambdaModuleFibrations}. Explicitly, in a
commutative diagram of $\Lambda_*$-modules
\begin{equation*}
\xymatrix{ A\ar@{>->}[]!D-<0pt,4pt>;[d]!U_{i}\ar[r]^{f} & M\ar@{->>}[d]^{p} \\
B\ar[r]_{g}\ar@{-->}[ur]^{h} & N },
\end{equation*}
the lift $h: B\rightarrow M$ exists provided that $p$ is an acyclic
fibration (M4.i), respectively provided that $i$ is an acyclic
cofibration (M4.ii).
\end{claim}

\begin{proof}
In both cases we construct by induction over $r\in\N$ a morphism $h:
B(r)\rightarrow M(r)$ that commutes with the action of morphisms
$u\in\Lambda_{s}^{r}$ such that $s\geq r$. Explicitly, if the
morphisms $h$ are defined for $s<r$, then we have a well defined
induced morphism $h: \Match{B}(r)\rightarrow\Match{M}(r)$ so that we
obtain a commutative diagram
\begin{equation*}
\xymatrix{ A(r)\ar@{>->}[]!D-<0pt,4pt>;[d]_(0.35){i}\ar[r] & M(r)\ar@{->>}[d]^{(\mu,p)} \\
B(r)\ar[r]\ar@{-->}[ur]^{h} & \Match{M}(r)\times_{\Match{N}(r)} N(r)
}
\end{equation*}
in the category of $\Sigma_r$-modules.

The right-hand side morphism is a fibration of $\Sigma_r$-modules by
definition and the left-hand side morphism a cofibration. If $i$ is
an acyclic cofibration, then, by definition, the left-hand side
morphism is obviously an acyclic cofibration as well so that the
lift exists in this case since the lifting property (M4.ii) is known
to be satisfied in the adjoint model category $\Sigma_r$-modules.

If $p$ is an acyclic cofibration, then we claim that the right-hand
side morphism is also an acyclic fibration. In fact, if we forget
the symmetric group action, then the dg-modules $\Match{M}(r)$ can
be identified with the matching modules of the functor underlying
$M$ on the direct Reedy category $\overrightarrow{\Lambda}_*$.
Therefore we deduce from the classical theory that the morphism
$(\mu,p): M(r)\rightarrow\Match{M}(r)\times_{\Match{N}(r)} N(r)$
defines an acyclic fibration in the category of dg-modules (see for
instance~\cite[Proposition 15.3.14]{Hirschhorn}) and hence in the
category of $\Sigma_r$-modules. Finally, the lift exists in this
case as well again since the lifting property (M4.i) is known to be
satisfied in the adjoint model category $\Sigma_r$-modules.

By construction, our morphism $h: M(r)\rightarrow N(r)$ makes the
diagram
\begin{equation*}
\xymatrix{ M(r)\ar[r]^{h}\ar[d] & N(r)\ar[d] \\
\Match{M}(r)\ar[r] & \Match{N}(r) }
\end{equation*}
commute. This property implies immediately that $h$ commutes with
all morphisms $\alpha\in\Lambda_{s}^{r}$ such that $s<r$ and hence
with all morphisms $\alpha\in\Lambda_{s}^{r}$ since $h$ is
$\Sigma_r$-equivariant as well by construction.
\end{proof}

Now we aim to define generating cofibrations and acyclic
cofibrations in the category of $\Lambda_*$-modules so that we can
deduce (M5.i-ii) from the small object argument. For this purpose we
extend the tensor product of dg-modules with simplicial sets
$C\otimes K$ to $\Lambda_*$-sets and we consider for $K$ the obvious
generators of the category of $\Lambda_*$-sets supplied by the
classical Yoneda lemma and the associated latching $\Lambda_*$-sets
that give representatives for the matching space modules. We define
these objects explicitly in the nexts paragraphs.

\subsubsection{Generating $\Lambda_*$-sets}
As usual, for any object $r\in\N$, we have a canonical
$\Lambda_*$-set, denoted by $\Lambda^r$, defined by the functor
corepresented by $r$ and such that $\Lambda^r(s) = \Lambda_{s}^{r}$
for $r\in\N$. By definition, the map $r\mapsto\Lambda^r$ defines a
covariant functor from $\Lambda_*$ to the category of
$\Lambda_*$-sets. Therefore we can form the associated latching
object which are defined explicitly by the colimit
\begin{equation*}
\Latch{\Lambda^r} = \colim_{\alpha:
\underline{r}'\xrightarrow{<}\underline{r}} \Lambda^{r'}.
\end{equation*}
As usual for a functor, this colimit can be determined pointwise so
that $\Latch{\Lambda^r}(s)$ is determined by the equivalent colimit
in the category of sets:
\begin{equation*}
\Latch{\Lambda^r}(s) = \colim_{\alpha:
\underline{r}'\xrightarrow{<}\underline{r}}
\Lambda_{s}^{r'}\quad\text{for all $s\in\N$}.
\end{equation*}
As in the definition of a matching module, we can assume that the
colimit ranges over the category $\Lambda_{*<r}$ or, equivalently,
over the category $\overrightarrow{\Lambda}_{*<r}$ which is cofinal
in $\Lambda_{*<r}$. We have a canonical morphism $\lambda:
\Latch{\Lambda^r}\rightarrow\Lambda^r$ for $r\in\N$. In fact, this
morphism is either trivial or an isomorphism. To be precise, we
obtain:

\begin{obsv}\label{obsv:LambdaSetLatching}\hspace*{2mm}

\begin{enumerate}
\item
For $s<r$, we have $\Latch{\Lambda^r}(s) = \Lambda^r(s) =
\Lambda_{s}^{r}$ and $\lambda:
\Latch{\Lambda^r}(s)\rightarrow\Lambda^r(s)$ is an isomorphism.
\item
For $s=r$, we have $\Latch{\Lambda^r}(r) = \emptyset$ and
$\Lambda^r(r) = \Lambda_{r}^{r} = \Sigma_r$ so that $\lambda:
\Latch{\Lambda^r}(r)\rightarrow\Lambda^r(r)$ is an initial morphism
in the category of sets.
\item
For $s>r$, we have $\Latch{\Lambda^r}(s) = \Lambda^r(s) =
\emptyset$.
\end{enumerate}
\end{obsv}

\begin{proof}
The assertions (b-c) are immediate since $\Lambda_{s}^{r} =
\emptyset$ for $s>r$ and only assertion (a) deserves a verification.
In fact, for $s<r$, the map
\begin{equation*}
\lambda: \colim_{\alpha: \underline{r}'\xrightarrow{<}\underline{r}}
\Lambda_{s}^{r'}\rightarrow\Lambda_{s}^{r}
\end{equation*}
admits a canonical section that maps an element
$u\in\Lambda_{s}^{r}$ to the identity morphism
$\id\in\Lambda_{s}^{s}$ in the component of the colimit indexed by
$u: \underline{s}\rightarrow\underline{r}$. One checks readily that
this section gives also a left-inverse for the map $\lambda$, so
that $\lambda$ is an isomorphism.
\end{proof}

\subsubsection{The tensor product with a $\Lambda_*$-set}
In general the tensor product $C\otimes A$ of an object $C$ with a
$\Lambda_*$-set $A$ is the $\Lambda_*$-module defined by the
coproducts
\begin{equation*}
C\otimes A(r) = \bigoplus_{\alpha\in A(r)} C.
\end{equation*}
In the context of dg-modules the tensor product $C\otimes A(r)$ can
be also be represented by the tensor product of $C$ with the free
module spanned by the set $A(r)$. Hence an element of $C\otimes
A(r)$ can be represented by a tensor $c\otimes\alpha$, where $c\in
C$ and $\alpha\in A(r)$.

For a morphism $i: C\rightarrow D$ in the category of dg-modules, we
form the pushout of $\Lambda_*$-modules
\begin{equation*}
\xymatrix{ C\otimes\Latch{\Lambda^r}\ar[r]\ar[d] & D\otimes\Latch{\Lambda^r}\ar[d] \\
C\otimes\Lambda^r\ar[r] &
C\otimes\Lambda^r\bigoplus_{C\otimes\Latch{\Lambda^r}}
D\otimes\Latch{\Lambda^r}\push{4} }
\end{equation*}
and we consider the natural morphism
\begin{equation*}
C\otimes\Lambda^r\bigoplus_{C\otimes\Latch{\Lambda^r}}
D\otimes\Latch{\Lambda^r}\rightarrow D\otimes\Lambda^r
\end{equation*}
induced by $i\otimes\Lambda^r: C\otimes\Lambda^r\rightarrow
D\otimes\Lambda^r$ and $D\otimes\lambda:
D\otimes\Latch{\Lambda^r}\rightarrow D\otimes\Lambda^r$. We prove
precisely that such morphisms give a set of generating arrows in the
category of~$\Lambda_*$-modules:

\begin{lemm}\label{lemm:LambdaModuleGeneratingCofibrations}
The morphisms
\begin{equation*}
C\otimes\Lambda^r\bigoplus_{C\otimes\Latch{\Lambda^r}}
D\otimes\Latch{\Lambda^r}\rightarrow D\otimes\Lambda^r
\end{equation*}
associated to a generating set of cofibrations, respectively acyclic
cofibrations, of the category of dg-modules $i: C\rightarrow D$
define a generating set of cofibrations, respectively acyclic
cofibrations, in the category of $\Lambda_*$-modules.
\end{lemm}

The proof of this lemma is split up into a sequence of claims and
can be compared with the classical construction of~\cite[Section
15.6]{Hirschhorn}. First, we observe that the morphisms considered
are indeed cofibrations, respectively acyclic cofibrations, in the
category of $\Lambda_*$-modules:

\begin{claim}
If $i: C\rightarrow D$ is a cofibration, respectively an acyclic
cofibration, in the category of dg-modules, then the associated
morphism
\begin{equation*}
C\otimes\Lambda^r\bigoplus_{C\otimes\Latch{\Lambda^r}}
D\otimes\Latch{\Lambda^r}\rightarrow D\otimes\Lambda^r
\end{equation*}
is a cofibration, respectively an acyclic cofibration, in the
category of $\Lambda_*$-modules.
\end{claim}

\begin{proof}
This assertion is a consequence of
observation~\ref{obsv:LambdaSetLatching}. Explicitly, from this
statement one deduces that the morphism
\begin{equation*}
C\otimes\Lambda^r(s)\bigoplus_{C\otimes\Latch{\Lambda^r}(s)}
D\otimes\Latch{\Lambda^r}(s)\rightarrow D\otimes\Lambda^r(s)
\end{equation*}
is an isomorphism for $s<r$, can be identified with
$i\otimes\Sigma_r: C\otimes\Sigma_r\rightarrow D\otimes\Sigma_r$ for
$s = r$ and vanishes for $s>r$. Consequently, this morphism is a
cofibration, respectively an acyclic cofibration, of
$\Sigma_s$-modules in all cases: this assertion is trivial in the
cases $s<r$ and $s>r$ and holds by definition (of cofibrations in
the category of $\Sigma_r$-modules) in the case $s=r$.
\end{proof}

Then we check that the sets of arrows introduced in the lemma detect
acyclic fibrations, respectively fibrations. For this purpose we
observe that the tensor products $C\otimes\Lambda^r$ and
$C\otimes\Latch{\Lambda^r}$ are characterized by natural adjunction
relations. Namely:

\begin{obsv}\label{obsv:FreeLambdaModules}
For a dg-module $C$ and a $\Lambda_*$-module $M$, we have the
adjunction formulas
\begin{multline*}
\Hom_{\dg\Lambda^{\op}_*\Mod}(C\otimes\Lambda^r,M) = \Hom_{\dg\Mod}(C,M(r)) \\
\text{and}\quad\Hom_{\dg\Lambda^{\op}_*\Mod}(C\otimes\Latch{\Lambda^r},M)
= \Hom_{\dg\Mod}(C,\Match{M}(r))
\end{multline*}
and the morphism
\begin{equation*}
\Hom_{\dg\Lambda^{\op}_*\Mod}(C\otimes\Lambda^r,M)\xrightarrow{C\otimes\lambda^*}
\Hom_{\dg\Lambda^{\op}_*\Mod}(C\otimes\Latch{\Lambda^r},M)
\end{equation*}
induced by the latching morphism $\lambda:
\Latch{\Lambda^r}\rightarrow\Lambda^r$ agrees with the morphism
induced by the matching morphism $\mu: M(r)\rightarrow\Match{M}(r)$.
\end{obsv}

\begin{proof}
These adjunction claims are formal consequences of the definition of
the $\Lambda_*$-sets $\Lambda^r$, $\Latch{\Lambda^r}$ and of the
Yoneda lemma.
\end{proof}

As a corollary, we obtain:

\begin{obsv}
For a morphism of $\Lambda_*$-modules $p: M\rightarrow N$, the
lifting problem
\begin{equation*}
\xymatrix{ C\otimes\Lambda^r\bigoplus_{C\otimes\Latch{\Lambda^r}} D\otimes\Latch{\Lambda^r}\ar[r]\ar[d] & M\ar[d]^{p} \\
D\otimes\Lambda^r\ar[r]\ar@{-->}[ur] & N }
\end{equation*}
is equivalent to an adjoint lifting problem
\begin{equation*}
\xymatrix{ C\ar[r]\ar[d] & M(r)\ar[d]^{(\mu,p)} \\
D\ar[r]\ar@{-->}[ur] & \Match{M}(r)\otimes_{\Match{N}(r)} N(r) }
\end{equation*}
in the category of dg-modules.
\end{obsv}

\begin{proof}
This observation is a formal consequence of the adjunction relations
of observation~\ref{obsv:FreeLambdaModules}.
\end{proof}

This observation implies immediately that the arrows introduced in
lemma~\ref{lemm:LambdaModuleGeneratingCofibrations} detect acyclic
fibrations and fibrations and achieve the proof of that
statement.\qed

\medskip
The next claim is a consequence of
lemma~\ref{lemm:LambdaModuleGeneratingCofibrations} and of the
classical small object argument:

\begin{claim}
The properties (M5.i-ii) hold for the class of weak-equivalences,
cofibrations and fibrations specified in
theorem~\ref{thm:LambdaModuleModelStructure} and
paragraph~\ref{item:LambdaModuleFibrations}. Explicitly, any
morphism of $\Lambda_*$-modules $\phi: M\rightarrow N$ admits a
factorization $\phi = p i$ such that $i$ is a cofibration and $p$ an
acyclic fibration (M5.i), respectively such that $i$ is an acyclic
cofibration and $p$ a fibration (M5.ii).
\end{claim}

\begin{proof}
The arguments are classical. We refer to~\cite[Section
10.5]{Hirschhorn}.
\end{proof}

This claim achieves the proof of
theorem~\ref{thm:LambdaModuleModelStructure}.\qed

\subsection{The Reedy model structure for unital operads}\label{subsection:UnitalOperadsReedyModelStructure}
In this section we check that the Reedy model structure of
$\Lambda_*$-modules can be transferred to operads through the
adjunction
\begin{equation*}
\FOp_*:
\dg\Lambda^{\op}_*\Mod^1_0/\overline{C}\rightleftarrows\dg\Op^1_*
:(-)^{-}
\end{equation*}
defined in~\ref{item:FreeUnitalOperads}. Explicitly, we have the
following theorem:

\begin{thm}\label{thm:UnitalOperadModelStructure}
The category of ($\N$ or $\Z$-graded) unital operads $\dg\Op^1_*$ is
equipped with the structure of a cofibrantly generated model
category such that a morphism $\phi: \P\rightarrow\Q$ is a
weak-equivalence, respectively a fibration, if $\phi$ forms a Reedy
weak-equivalence, respectively a Reedy fibration, in the category of
$\Lambda_*$-modules $\dg\Lambda^{\op}_*\Mod^1_0$.
\end{thm}\index{Reedy!model category!of unital operads}

The cofibrations are characterized by the left lifting property with
respect to acyclic fibrations as usual.

Our model structure differs from the adjoint model structure defined
in~\cite{BergerMoerdijkW} like the Reedy model structure of
$\Lambda_*$-modules differs from the adjoint model structure.
Nevertheless one can observe again that the identity functor yields
a pair of adjoint derived equivalences between the two model
categories.

Theorem~\ref{thm:UnitalOperadModelStructure} follows from the
general transfer principle considered in~\emph{loc. cit.} since we
observe in lemma~\ref{lemm:LambdaModuleGeneratingCofibrations} that
the category of $\Lambda_*$-modules is cofibrantly generated and the
scheme of the proof of theorem~\ref{thm:UnitalOperadModelStructure}
is classical. As in~\cite{BergerMoerdijk}, we have essentially to
check that the following property is satisfied:

\begin{claim}\label{claim:AcyclicCellAttachment}
Suppose given a pushout of unital operads
\begin{equation*}
\xymatrix{ \FOp_*(C)\ar@{>->}[]!D-<0pt,4pt>;[d]_(0.35){i}^(0.35){\sim}\ar[r]^{f} & \P\ar[d]^{j} \\
\FOp_*(D)\ar[r]^{g} & \Q\push{4} }
\end{equation*}
such that $i: \FOp_*(C)\rightarrow\FOp_*(D)$ is a morphism of free
operads induced by an acyclic cofibration of $\Lambda_*$-modules $i:
C\wecofib D$. The morphism $j: \P\rightarrow\Q$ is a
weak-equivalence as well.
\end{claim}

\begin{proof}
This property is inherited from the analogous property of non-unital
operads, for which we can refer to~\cite{BergerMoerdijk,Hinich},
since we observe that the functor $\Q\mapsto\overline{\Q}$ from the
category of unital operads to the category of non-unital operads
creates pushouts.
\end{proof}

According to the transfer principle
(see~\cite{BergerMoerdijk,Crans}), this verification achieves the
proof of theorem~\ref{thm:UnitalOperadModelStructure}.\qed

\medskip
Theorem~\ref{thm:UnitalOperadModelStructure} holds for the category
of connected unital dg-operads $\dg\Op^*_*$. In this case we
consider the category $\dg\Lambda^{\op}\Mod^*_0/\overline{\C}$,
introduced in~\ref{item:ConnectedUnitaryOperads}, formed by the
augmented non-unital unitary $\Lambda_*$-modules $M$ such that $M(1)
= \F$. According to observations of this paragraph, we have an
adjunction relation between this category and connected unital
operads:
\begin{equation*}
\FOp_*:
\dg\Lambda^{\op}_*\Mod^*_0/\overline{\C}\rightleftarrows\dg\Op^*_*
:(-)^{-}.
\end{equation*}
Furthermore, the functor $\P\mapsto\overline{\P}$ creates also
pushouts in $\dg\Op^*_*$, so that the argument above remains valid
for~$\dg\Op^*_*$.

On the other hand, one can deduce from the adjunction relations
of~\ref{item:ConnectedAdjointOperads} that the category of connected
unital dg-operads $\dg\Op^*_*$ forms a model subcategory
of~$\dg\Op^1_*$ so that a morphism of connected unital dg-operads
defines a Reedy cofibration, respectively a Reedy fibration, a
weak-equivalence, in $\dg\Op^*_*$ if and only if it defines a Reedy
cofibration, respectively a Reedy fibration, a weak-equivalence,
in~$\dg\Op^1_*$. As a byproduct, we have the following statement:

\begin{fact}
The functor $\str^1_*: \dg\Op^1_*\rightarrow\dg\Op^*_*$, right
adjoint to the category embedding $\itr^1_*:
\dg\Op^1_*\rightarrow\dg\Op^*_*$, preserves fibrations, acyclic
fibrations and all weak-equivalences between fibrant operads.
\end{fact}

This claim could also be deduced from an analogous assertion for
$\Lambda_*$-modules since operad fibrations, respectively
weak-equivalences, are just fibrations, respectively
weak-equivalences, in the $\Lambda_*$-module category. To be
precise, the statement above is implied by the analogous assertion
for the category embedding $\itr^1_*:
\dg\Lambda^{\op}_*\Mod^*_*\hookrightarrow\dg\Lambda^{\op}_*\Mod^1_*$.

\subsection{On Hopf $\Lambda_*$-modules and unital Hopf operads}\label{subsection:HopfOperadsIntroduction}
Recall that the notion of an operad makes sense in any symmetric
monoidal category. In this memoir, if no category is specified, then
an operad refers to an operad in the category of dg-modules
$\dg\Mod$ since this category forms our ground monoidal category. On
the other hand, as explained in the introduction, we consider also
\emph{Hopf operad} structures which are precisely operads in the
category of augmented coassociative dg-coalgebras, denoted by
$\CoAlg^a_+$\glossary{$\CoAlg^a_+$}. The purpose of this section is
to make clear the structure of a Hopf operad and to prove the
analogue of theorem~\ref{thm:UnitalOperadModelStructure} for Hopf
operads. Namely we prove that unital Hopf operads form a model
category.

\subsubsection{On (unital) Hopf operads}\label{item:HopfOperads}
Explicitly, a Hopf operad\index{Hopf operad}\index{operad!Hopf}
consists of a collection of coalgebras $\P(r)\in\CoAlg^a_+$ such
that the operad composition products $\circ_i:
\P(s)\otimes\P(t)\rightarrow\P(s+t-1)$ define morphisms in
$\CoAlg^a_+$ as well as the isomorphisms $w: \P(r)\rightarrow\P(r)$
defined by the action of permutations $w\in\Sigma_r$. As usual for
bialgebra structures, one can equivalently assume that the diagonals
$\Delta: \P(r)\rightarrow\P(r)\otimes\P(r)$ define an operad
morphism, where the tensor product $\P(r)\otimes\P(r)$ is equipped
with place-by-place composition products. The distinguished unital
operation $*\in\P(0)$ and the operad unit $1\in\P(1)$ are supposed
to define group-like elements in~$\P$.

One can observe that the composition products preserve the coalgebra
augmentations $\epsilon: \P(r)\rightarrow\F$ if and only if these
morphisms define a morphism to the operad of commutative algebras
$\C$. Accordingly, any Hopf operad $\P$ is automatically augmented
over the commutative operad $\C$. Furthermore, for a unital Hopf
operad $\P$, we deduce from these observations that the coalgebra
augmentations agree with the augmentation morphisms $\epsilon:
\P(r)\rightarrow\F$ introduced in the introduction of this section
and defined by the operadic composites $\epsilon(p) =
p(*,\dots,*)$\glossary{$\epsilon$}.

Similarly, a \emph{Hopf algebra over $\P$}\index{operad!Hopf!Hopf
algebra over a}\index{Hopf!algebra over a Hopf operad} refers to an
algebra over $\P$ in the category of augmented coassociative
dg-coalgebras. Hence a Hopf algebra consists of a coalgebra
$\Gamma\in\CoAlg^a_+$ equipped with evaluation products
$\P(r)\otimes\Gamma^{\otimes r}\rightarrow\Gamma$ that form
morphisms in $\CoAlg^a_+$. As stated previously for the composition
products of a Hopf operad, one can equivalently assume that the
diagonal of $\Gamma$ defines a morphism of $\P$-algebras $\Delta:
\Gamma\rightarrow\Gamma\otimes\Gamma$, where the tensor product
$\Gamma\otimes\Gamma$ is equipped with a diagonal action of~$\P(r)$.

If $\P$ is a unital Hopf operad\index{Hopf
operad!unital}\index{unital!Hopf operad}\index{operad!unital Hopf},
then a Hopf $\P$-algebra $\Gamma$ comes equipped with a unit
morphism $\eta: \F\rightarrow\Gamma$\index{unit!of a
coalgebra}\index{coalgebra!unit of a} yielded by the unital
operation $*\in\P(0)$. Consequently, the underlying coalgebra of a
Hopf algebra over a unital operad defines an object in the category
$\CoAlg^a_*$\glossary{$\CoAlg^a_*$} of \emph{augmented unitary
coalgebras}\index{unitary!coalgebra}\index{coalgebra!unitary}. By
convention, if a given augmented unitary coalgebra $\Gamma$ is
equipped with the structure of a Hopf $\P$-algebra, where $\P$ is a
unital Hopf operad, then we assume that the unit of $\Gamma$ agrees
with the unital operation.

\subsubsection{Hopf $\Lambda_*$-modules}\label{item:HopfLambdaModules}
Observe that the underlying $\Lambda_*$-module of a unital Hopf
operad defines a $\Lambda_*$-object in the category of augmented
coassociative dg-coalgebras. This structure is called a \emph{Hopf
$\Lambda_*$-module}\index{Hopf
$\Lambda_*$-module}\index{$\Lambda_*$-module!Hopf} according to our
usual conventions. Clearly, a Hopf $\Lambda_*$-module $\Gamma$ is
automatically augmented over the constant $\Lambda_*$-module since
the coalgebra augmentation gives a morphism $\epsilon:
\Gamma\rightarrow\F$ as in the case of unital Hopf operads.

To recapitulate, the forgetful functor of unital operads
$\P\mapsto\overline{\P}$, considered
in~\ref{item:FreeUnitalOperads}, restricts to a functor from the
category of unital Hopf operads to the category of non-unital
unitary Hopf $\Lambda_*$-modules. In the converse direction, if a
unitary $\Lambda_*$-module $\Gamma$ is equipped with a coalgebra
structure, then the associated free unital operad $\FOp_*(\Gamma)$
can be equipped with the structure of a unital Hopf operad so that
$\FOp_*(\Gamma)$\glossary{$\FOp_*(\Gamma)$}\index{free!unital Hopf
operad}\index{Hopf operad!free unital} satisfies the usual universal
property for Hopf $\Lambda_*$-modules. In fact, since the diagonal
of a Hopf operad is supposed to define an operad morphism, the
diagonal of a formal composite $\gamma = (\dots((\gamma_1\circ_{i_2}
\gamma_2)\circ_{i_3}\dots)\circ_{i_l} \gamma_l$ in $\FOp_*(\Gamma)$
can be written
\begin{equation*}
\Delta(\gamma) = \sum_{\gamma_1,\dots,\gamma_r}
(\dots((\gamma'_1\circ_{i_2} \gamma'_2)\circ_{i_3}\dots)\circ_{i_l}
\gamma'_l \otimes (\dots((\gamma''_1\circ_{i_2}
\gamma''_2)\circ_{i_3}\dots)\circ_{i_l} \gamma''_l,
\end{equation*}
where $\Delta(\gamma_i) = \sum_{\gamma_i}
\gamma'_i\otimes\gamma''_i$ denotes the diagonal of $\gamma_i$ in
$\Gamma(r)$.

To conclude, the forgetful functor from unital Hopf operads to
unital operads creates free objects. As a consequence, for unital
Hopf operads, we have an adjunction relation
\begin{equation*}
\FOp_*: \dg\Lambda_*^{\op}\HopfMod^1_0\rightleftarrows\dg\HopfOp^1_*
:(-)^{-},
\end{equation*}\glossary{$\Lambda_*^{\op}\HopfMod^1_0$}\glossary{$\HopfOp^1_*$}
where $\FOp_*$ denotes the free operad functor defined
in~\ref{subsection:LambdaModulesIntroduction}. One can observe in
addition that the forgetful functor creates also all small colimits
in the category of unital Hopf operads since the forgetful functor
from the category of coalgebras to the category of dg-modules has
this property.

\subsubsection{The model category of coalgebras}
According to~\cite{GetzlerJones}, the category of non-negatively
graded dg-coalgebras over a field is equipped with the structure of
a model category such that a morphism $\phi: C\rightarrow D$ is a
weak-equivalence, respectively a cofibration, if $\phi$ defines a
weak-equivalence, respectively a cofibration, in the category of
dg-modules. (Recall that a morphism of dg-modules over a field is a
cofibration if and only if it is injective.) Furthermore, this
category has a set of generating cofibrations, respectively acyclic
cofibrations, defined by the collection of injective morphisms $i:
C\rightarrow D$ such that $D$ is spanned by a finite, respectively
countable, set of homogeneous elements.

As specified in~\ref{subsection:ModelCategoryLambdaModules}, the
Reedy model structure of
theorem~\ref{thm:LambdaModuleModelStructure} can be defined in any
cofibrantly generated ground model category. In particular, for the
$\Lambda_*$-objects in the category of dg-coalgebras we obtain the
following theorem:

\begin{thm}\label{thm:HopfLambdaModuleModelStructure}
The category of $\N$-graded Hopf dg-$\Lambda_*$-modules
$\dg\Lambda^{\op}_*\HopfMod$\glossary{$\dg\Lambda^{\op}_*\HopfMod$}
is equipped with the structure of a cofibrantly generated model
category such that a morphism $\phi: C\rightarrow D$ is a
weak-equivalence, respectively a cofibration, if $\phi$ defines a
weak-equivalence, respectively a cofibration, in the category of
dg-$\Sigma_*$-modules $\dg\Sigma_*\Mod$.
\end{thm}\index{Reedy!model category!of
Hopf $\Lambda_*$-modules}\index{Hopf $\Lambda_*$-module!Reedy model
category of Hopf $\Lambda_*$-modules}

This model structure for Hopf $\Lambda_*$-modules is also
cofibrantly generated. To be more explicit, as in the ground
category of dg-modules, we have the following lemma:

\begin{lemm}
We have a generating set of cofibrations, respectively acyclic
cofibrations, in the category of Hopf $\Lambda_*$-modules defined by
morphisms
\begin{equation*}
C\otimes\Lambda^r\bigoplus_{C\otimes\Latch{\Lambda^r}}
D\otimes\Latch{\Lambda^r}\rightarrow D\otimes\Lambda^r
\end{equation*}
associated to a generating set of cofibrations, respectively acyclic
cofibrations, of the category of dg-coalgebras $i: C\rightarrow
D$.\qed
\end{lemm}

Recall simply that the forgetful functor creates colimits in the
category of coalgebras so that the tensor products $C\otimes K$,
where $K$ is a $\Lambda_*$-set, as well as the coproducts
$C\otimes\Lambda^r\bigoplus_{C\otimes\Latch{\Lambda^r}}
D\otimes\Latch{\Lambda^r}$ have the same realization in the category
of dg-coalgebras as in the category of dg-modules.

\subsubsection{The matching coalgebra of a Hopf $\Lambda_*$-module}
The fibrations can be characterized by the right lifting property as
usual but, as in the context of dg-modules, one can introduced an
appropriate matching object in the category of Hopf
$\Lambda_*$-modules so that a morphism of Hopf $\Lambda_*$-modules
$p: C\rightarrow D$ is a fibration if and only if the morphisms
\begin{equation*}
(\mu,p): C(r)\rightarrow\Match{C}(r)\times_{\Match{D}(r)} D(r)
\end{equation*}
define a fibration in the category of dg-coalgebras for all
$r\in\N$.\index{Reedy!fibration!of Hopf
$\Lambda_*$-modules}\index{fibration!Reedy fibration of Hopf
$\Lambda_*$-modules} The matching object of a Hopf
$\Lambda_*$-module $\Match{C}$\glossary{$\Match{C}$}\index{Hopf
$\Lambda_*$-module!matching object of a}\index{matching object!of a
Hopf $\Lambda_*$-module} is defined by the same limit as in the
category of dg-modules
\begin{equation*}
\Match{C}(r) = \lim_{\alpha:
\underline{r}'\xrightarrow{<}\underline{r}} C(r'),
\end{equation*}
where $\alpha$ ranges over morphisms $\alpha\in\Lambda_{r'}^{r}$
such that $r'<r$, except that we perform this limit in the category
of dg-coalgebras. This matching object $\Match{C}(r)$ can also be
defined by an equalizer as in the category of dg-modules, but an
equalizer of dg-coalgebras. Namely:
\begin{equation*}
\Match{C}(r) = \ker\bigl(\xymatrix{ *+<2mm>{\prod_{1\leq i\leq r}
C(r-1)}\ar@<2pt>[r]^{d^0}\ar@<-2pt>[r]_{d^1} &
*+<2mm>{\prod_{1\leq i<j\leq r} C(r-2)} }\bigr).
\end{equation*}

\medskip
As in the context of $\Lambda^*$-modules, we have an induced model
structure on the subcategory of $\dg\Lambda_*^{\op}\HopfMod$ formed
by non-unital Hopf $\Lambda^*$-modules $\Gamma$. First, we observe
that the adjunction relation of
proposition~\ref{prop:ConnectedAdjoints} can be extended in the
coalgebra context so that we have an adjunction ladder:
\begin{equation*}
\xymatrix{ \dg\Lambda^{\op}_*\HopfMod_0\ar[r]^{\itr^{\Lambda}_0} &
\dg\Lambda^{\op}_*\HopfMod\ar@<-6mm>[l]_{\ctr^{\Lambda}_0}\ar@<+6mm>[l]_{\str^{\Lambda}_0}
}.
\end{equation*}\glossary{$\dg\Lambda^{\op}_*\HopfMod_0$}
In fact, the non-unital object $\ctr^{\Lambda}_0(C)$, respectively
$\str^{\Lambda}_0(D)$, associated to a Hopf $\Lambda_*$-module $C$,
respectively $D$, can be obtained as in the dg-module context except
that in the definition of~$\str^{\Lambda}_0(D)$ we perform the
pullbacks
\begin{equation*}
\xymatrix{ \str^{\Lambda}_0(D)(r)\ar@{-->}[r]\ar@{-->}[d] &
D(r)\ar[d]^{\eta_0^*} \\
0\ar[r]^{0} & D(0) }.
\end{equation*}
in the category of coalgebras. As in the context of
$\Lambda_*$-modules, we observe that the functor
$\itr^{\Lambda}_0\ctr^{\Lambda}_0$ preserves cofibrations and all
weak-equivalences of Hopf $\Lambda_*$-modules. Then we obtain:

\begin{prop}\label{prop:NonunitalHopfLambdaModulesModelStructure}
The category of $\N$-graded Hopf $\Lambda_*$-modules
$\dg\Lambda^{\op}_*\HopfMod_0$ forms a model subcategory
of~$\dg\Lambda^{\op}_*\HopfMod$ so that a morphism $f: C\rightarrow
D$ is a weak-equivalence, respectively a cofibration, a fibration,
in $\dg\Lambda^{\op}_*\HopfMod_0$ if and only if $f$ defines a
weak-equivalence, respectively a cofibration, a fibration,
in~$\dg\Lambda^{\op}_*\HopfMod$.

Furthermore, this category is cofibrantly generated by the morphisms
$\ctr^{\Lambda}_0(i):
\ctr^{\Lambda}_0(A)\rightarrow\ctr^{\Lambda}_0(B)$ associated to a
generating set of cofibrations, respectively acyclic cofibrations,
in $\dg\Lambda^{\op}_*\HopfMod$.\qed
\end{prop}

As in the ground category of dg-modules, we obtain that the functor
\begin{equation*}
\str^{\Lambda}_0:
\dg\Lambda^{\op}_*\HopfMod\rightarrow\dg\Lambda^{\op}_*\HopfMod_0
\end{equation*}
preserves fibrations, acyclic fibrations and weak-equivalences
between fibrant objects by adjunction.

The category of unitary objects $\dg\Lambda^{\op}_*\HopfMod^1_0$ is
also equipped with the canonical model structure of a comma
category.

\medskip
Now we check that the results
of~\ref{subsection:UnitalOperadsReedyModelStructure} can be extended
to the category of unital Hopf operads. Explicitly, the Reedy model
structure of (non-unital) Hopf $\Lambda_*$-modules can be
transferred to unital Hopf operads through the adjunction
\begin{equation*}
\FOp_*: \dg\Lambda_*^{\op}\HopfMod^1_0\rightleftarrows\dg\HopfOp^1_*
:(-)^{-}
\end{equation*}
so that we obtain the following theorem:

\begin{thm}\label{thm:HopfOperadModelStructure}
The category of $\N$-graded unital Hopf dg-operads $\dg\HopfOp^1_*$
is equipped with the structure of a cofibrantly generated model
category such that a morphism $\phi: \P\rightarrow\Q$ is a
weak-equivalence, respectively a fibration, if $\phi$ forms a Reedy
weak-equivalence, respectively a Reedy fibration, in the category of
Hopf $\Lambda_*$-modules $\dg\Lambda^{\op}_*\HopfMod^1_0$.
\end{thm}\index{Reedy!model category!of
unital Hopf operads}\index{Hopf operad!Reedy model category of
unital Hopf operads}

As in the case of unital operads, we check simply that the
assumptions of the transfer principle are satisfied. Essentially, we
check the following property:

\begin{claim}\label{claim:AcyclicHopfCellAttachment}
Suppose given a pushout of unital Hopf operads
\begin{equation*}
\xymatrix{ \FOp_*(C)\ar@{>->}[]!D-<0pt,4pt>;[d]_(0.35){i}^(0.35){\sim}\ar[r]^{f} & \P\ar[d]^{j} \\
\FOp_*(D)\ar[r]^{g} & \Q\push{4} }
\end{equation*}
such that $i: \FOp_*(C)\rightarrow\FOp_*(D)$ is a morphism of free
operads induced by an acyclic cofibration of Hopf
$\Lambda_*$-modules $i: C\wecofib D$. The morphism $j:
\P\rightarrow\Q$ is a weak-equivalence as well.
\end{claim}

\begin{proof}
In fact, this property is clearly inherited from the category of
unital operads since we observe in~\ref{item:HopfLambdaModules} that
the forgetful functor creates the free functor $\FOp_*$ and the
small colimits in the category of unital Hopf operads.
\end{proof}

\subsubsection{Connected unital Hopf
operads}\label{item:ConnectedHopfOperads} As in the context of
operads in dg-modules, we can consider connected unital Hopf operads
$\P$\index{operad!connected}\index{Hopf operad!connected}
characterized by $\P(1) = \F$ and the full category
$\HopfOp^*_*$\glossary{$\HopfOp^*_*$} formed by these objects. The
category embedding $\itr^1_*: \HopfOp^*_*\hookrightarrow\HopfOp^1_*$
admits also a left and a right adjoint $\ctr^1_*,\str^1_*:
\HopfOp^1_*\rightarrow\HopfOp^*_*$ such that $\ctr^1_*\itr^1_* = \Id
=
\str^1_*\itr^1_*$.\glossary{$\itr^1_*$}\glossary{$\ctr^1_*$}\glossary{$\str^1_*$}

In fact, we can consider the underlying categories of unital unitary
Hopf $\Lambda_*$-modules which are associated to these operad
categories and for which we have an adjunction ladder:
\begin{equation*}
\xymatrix{ \Lambda^{\op}_*\HopfMod^*_*\ar[r]^{\itr^1_*} &
\Lambda^{\op}_*\HopfMod^1_*\ar@<-6mm>[l]_{\ctr^1_*}\ar@<+6mm>[l]_{\str^1_*}
}.
\end{equation*}
Then one checks readily that the connected unital unitary Hopf
$\Lambda_*$-module $\str^1_*(\Q)$, respectively $\ctr^1_*(\P)$,
associated to an operad $\Q$, respectively $\P$, forms a Hopf
suboperad of~$\Q$, respectively a quotient operad of~$\P$, and the
maps $\Q\mapsto\str^1_*(\Q)$ and $\P\mapsto \ctr^1_*(\P)$ supply
required adjoint functors at the operad level.

The connected Hopf $\Lambda_*$-module $\str^1_*(N)$ associated to a
unital unitary Hopf $\Lambda_*$-module
$N\in\Lambda^{\op}_*\HopfMod^1_*$ can be defined by the same
pullback diagrams as in the category of dg-modules
\begin{equation*}
\xymatrix{ \str^1_*(N)(r)\ar@{-->}[r]\ar@{-->}[d] &
N(r)\ar[d]^{(\eta_i^*)_i} \\
\F^{\times r}\ar[r]^{1^{\times r}} & N(1)^{\times r} }
\end{equation*}
except that we consider cartesian products in the category of
coalgebras and we perform these pullbacks in the same category. The
other Hopf $\Lambda_*$-module $\ctr^1_*(M)$ is obtained by the same
obvious construction as in the category of dg-modules. Namely we set
$\ctr^1_*(M)(r) = \F$ for $r = 0,1$ and $\ctr^1_*(M)(r) = M(r)$ for
$r\geq 2$. The operations $\partial_i:
\ctr^1_*(M)(r)\rightarrow\ctr^1_*(M)(r-1)$ are given either by the
corresponding operations or by the augmentation of~$M$.

As in the context of operads in dg-modules, we observe that
theorem~\ref{thm:HopfOperadModelStructure} holds for the category of
connected Hopf operads. More precisely, we obtain that the category
of connected unital dg-operads $\dg\HopfOp^*_*$ forms a model
subcategory of~$\dg\HopfOp^1_*$ such that a morphism of connected
unital dg-operads defines a Reedy cofibration, respectively a Reedy
fibration, a weak-equivalence, in $\dg\HopfOp^*_*$ if and only if it
defines a Reedy cofibration, respectively a Reedy fibration, a
weak-equivalence, in~$\dg\HopfOp^1_*$. Furthermore, we have:

\begin{fact}
The functor $\str^1_*: \dg\HopfOp^1_*\rightarrow\dg\HopfOp^*_*$,
right adjoint to the category embedding $\itr^1_*:
\dg\HopfOp^*_*\rightarrow\dg\HopfOp^1_*$, preserves fibrations,
acyclic fibrations and all weak-equivalences between fibrant
operads.
\end{fact}

In fact, since fibrations, respectively weak-equivalences, of unital
Hopf operads are just fibrations, respectively weak-equivalences, of
Hopf $\Lambda_*$-modules, we can also deduce this claim from the
same assertion at the level of Hopf $\Lambda_*$-modules. To be
precise, the statement above is implied by the analogous assertion
for the category embedding $\itr^1_*:
\dg\Lambda^{\op}_*\HopfMod^*_*\hookrightarrow\dg\Lambda^{\op}_*\HopfMod^1_*$.

\subsection{Prospects: cellular operads}
As explained in the memoir introduction, one might be willing to
extend the results of this section to operads equipped with a good
cellular structure. In regard to our construction, a good notion is
supplied by functor operads on an operad $\O$ in Reedy categories.
Explicitly, we assume that $\O$ consists of a sequence of categories
equipped with a Reedy decomposition $\O(r) =
\overrightarrow{\O}(r)\Sigma(r)\overleftarrow{\O}(r)$, where
$\overrightarrow{\O}(r)$ is a direct Reedy category,
$\overleftarrow{\O}(r)$ is an inverse Reedy category and $\Sigma(r)$
consists of isomorphisms, preserved by the operad structure.

Recall that an operad functor in a category $\C$ consists of a
functor $F: \O\rightarrow\C$ equipped with a suitable generalization
of an operad structure (see~\cite{McClureSmithCosimplicial}). The
definitions of this section can be generalized to this context. One
has simply to consider the appropriate generalization of the notion
of a $\Lambda_*$-module (functors $M: \Lambda\O\rightarrow\C$, where
$\Lambda\O$ denotes a semi-direct product of categories) and the
corresponding matching and latching structures.

\section{On Boardman-Vogt' $W$-construction}\label{section:CofibrantHopfOperads}

\subsection{Introduction}
The model category structure implies the existence of Reedy
cofibrant replacements in the category of unital Hopf operads. As
long as we deal with abstract structure results this purely
existence theorem meets our needs. But for effective issues we need
to have explicit cofibrant replacements. For the sake of
completeness (we do not really use the constructions of this section
in the memoir) and for subsequent references, we introduce in this
section a differential graded analogue of the Boardman-Vogt
$W$-construction that fulfils this need. Explicitly, we define a
functor $\P\mapsto W(\P)$, from the category of connected unital
dg-operads to itself, endowed with the following features:

\begin{thm}[Compare with Berger-Moerdijk~\cite{BergerMoerdijkW},
Boardman-Vogt~\cite{BoardmanVogt}]\label{thm:WHopfOperad}\hspace*{2mm}
\begin{enumerate}
\item
The operad $W(\P)$\glossary{$W(\P)$}\index{Boardman-Vogt'
construction}is equipped with a natural morphism $\epsilon:
W(\P)\rightarrow\P$\glossary{$\epsilon$} which is a weak-equivalence
of dg-operads.
\item
The operad $W(\P)$ is quasi-free as a unital operad. To be more
explicit, the operad $W(\P)$ is defined by a free operad
$\FOp_*(W'(\P))$ equipped with a differential
\begin{equation*}
\delta+\partial: \FOp_*(W'(\P))\rightarrow\FOp_*(W'(\P))
\end{equation*}
which is given by the addition of a homogeneous derivation
$\partial: \FOp_*(W'(\P))\rightarrow\FOp_*(W'(\P))$, determined by
the operad structure of $\P$, to the canonical differential $\delta:
\FOp_*(W'(\P))\rightarrow\FOp_*(W'(\P))$, induced by the internal
differential of~$\P$. Furthermore, the morphism $\phi_*:
W(\P)\rightarrow W(\P')$ associated to an operad morphism $\phi:
\P\rightarrow\P'$ is the morphism of quasi-free operads induced by
the morphism of $\Lambda_*$-modules $\phi_*: W'(\P)\rightarrow
W'(\P')$ associated to $\phi$. The morphism $\phi_*:
W(\P)\rightarrow W(\P')$ defines a cofibration in the Reedy model
category of operads provided that $\phi: \P\rightarrow\P'$ forms a
cofibration in the model category of $\Lambda_*$-modules.
\item
If an operad $\P$ is equipped with the structure of a Hopf operad,
then the associated $W$-construction can still be equipped with a
coassociative diagonal functorially in $\P$ so that $\P\mapsto
W(\P)$ defines a functor from the category of connected unital Hopf
operads to itself. Moreover, the augmentation $\epsilon:
W(\P)\rightarrow\P$ defines a weak-equivalence of unital Hopf
operads.\index{Hopf operad!Boardman-Vogt' construction of a}
\end{enumerate}
\end{thm}

We refer to~\ref{item:WQuasiFree} for the general definition of a
quasi-free object in the category of unital operads.

In fact, as observed by Berger-Moerdijk in~\cite{BergerMoerdijkW},
an analogue of the Boardman-Vogt $W$-construction can be defined in
any monoidal model category equipped with a good interval $\I$. Our
differential graded $W$-construction can be identified with an
instance of this general $W$-construction for the category of
dg-modules, respectively for the category of dg-coalgebras. On the
other hand, the classical topological Boardman-Vogt
$W$-construction, introduced in~\cite{BoardmanVogt}, is defined by a
cellular complex. Our differential graded $W$-construction is
precisely the chain complex defined by this cellular object.

In this section, we give a precise account of the definition of
$W(\P)$, though this construction is not original, in order to give
insights into the structure of this operad in the unital context. To
be precise, we have to check that $W(\P)$ that yield cofibrant
resolutions in the Reedy model category as stated in
theorem~\ref{thm:WHopfOperad}. An analogous statement for the
adjoint model structure is proved in~\cite{BergerMoerdijkW} but, for
our purpose, we give other arguments. The Boardman-Vogt construction
is by definition a cellular object in the ground model category.
In~\cite{BergerMoerdijkW,BoardmanVogt} the authors deal only with
this cellular structure. In this memoir, we prove that the
Boardman-Vogt construction $W(\P)$ defines a cellular object in the
category of operads. Explicitly, we prove the following theorem that
implies the cofibrancy claim of theorem~\ref{thm:WHopfOperad}:

\begin{thm}\label{thm:BoardmanVogtOperadicCellularDecomposition}
The Boardman-Vogt construction $W(\P)$ is the colimit of a sequence
of (Hopf) operads
\begin{multline*}
* = W^{-1}(\P)\xrightarrow{j_0}W^{0}(\P)\xrightarrow{j_1}\dots
\xrightarrow{j_d}W^d(\P)\xrightarrow{j_{d+1}}\dots\\
\dots\xrightarrow{}\colim_d W^d(\P) = W(\P)
\end{multline*}
obtained by pushouts
\begin{equation*}
\xymatrix{ \FOp_*(C^d(\P))\ar[r]^{f^d}\ar[d]^{i^d} & W^{d-1}(\P)\ar@{-->}[d]^{j_d} \\
\FOp_*(D^d(\P))\ar@{-->}[r]^{g^d} & W^{d}(\P)\push{4} },
\end{equation*}
where $i^d: \FOp_*(C^d(\P))\rightarrow\FOp_*(D^d(\P))$ is a morphism
of free operads associated to a morphism of unitary (Hopf)
$\Lambda_*$-modules $i^d: C^d(\P)\rightarrow D^d(\P)$.

This decomposition is functorial in $\P$ and in addition the
canonical morphism of unitary (Hopf) $\Lambda_*$-modules
\begin{equation*}
(i^d,\phi): C^d(\P')\bigoplus_{C^d(\P)} D^d(\P)\rightarrow D^d(\P')
\end{equation*}
associated to an operad morphism $\phi: \P\rightarrow\P'$ is a
weak-equivalence, respectively a Reedy cofibration, if $\phi$
defines a weak-equivalence, respectively a Reedy cofibration, in the
category of $\Lambda_*$-modules. In particular, the morphism $i^d:
C^d(\P)\rightarrow D^d(\P)$ is a Reedy cofibration if the operad
$\P$ forms a Reedy cofibrant object in the category of
$\Lambda_*$-modules.
\end{thm}

The decomposition of
theorem~\ref{thm:BoardmanVogtOperadicCellularDecomposition} is
referred to as the operadic cellular decomposition of the
Boardman-Vogt operad $W(\P)$.

Theorem~\ref{thm:BoardmanVogtOperadicCellularDecomposition} holds in
the general framework of~\cite{BergerMoerdijkW}. Nevertheless, for
simplicity, we make our construction explicit and we give
comprehensive proofs only in the context of dg-modules and
dg-coalgebras.

We need to recall conventions and results on trees that give the
structure of the $W$-construction. We devote the next
subsection~\ref{subsection:CellMetricTrees} to this topic. We
achieve the definition of the $W$-construction in the differential
graded context and we prove the assertions of
theorem~\ref{thm:WHopfOperad}
in~\ref{subsection:BoardmanVogtConstruction}. We define the operadic
cellular decomposition of $W(\P)$ and we prove
theorem~\ref{thm:BoardmanVogtOperadicCellularDecomposition}
in~\ref{subsection:BoardmanVogtOperadicCellularDecomposition} of
this section. The figures referred to in the text are displayed in
the appendix subsection~\ref{subsection:figures}.

As explained, the
subsection~\ref{subsection:BoardmanVogtConstruction} does not
provide any original result. Nevertheless our presentation differs
from~\cite{BergerMoerdijkW,BoardmanVogt} at some points since we aim
to introduce another cellular structure on the Boardman-Vogt
construction.

Recall again that we do not really use the constructions of this
section in the memoir. In regard to the needs of the next parts,
this section can be skipped in a first reading. The account of this
section is motivated by the sake of references in view of effective
constructions of operad actions which are postponed to future
articles.

\listoffigures

\subsection{Cell metric trees}\label{subsection:CellMetricTrees}
In the topological context the space $W(\P)(r)$ defined by
Boardman-Vogt in~\cite{BoardmanVogt} is a cubical cellular complex
formed by trees whose vertices are labeled with operations of $\P$
and whose edges are equipped with a length. The authors
of~\cite{BergerMoerdijkW} consider in fact a generalization of this
cellular construction in other categories than topological spaces.

As stated in the section outline, the aim of this subsection is to
recall our conventions for tree structures, borrowed
from~\cite{OperadTextbook}, so that we can make precise the
definition of~$W(\P)$. In addition we prove that the metric trees
considered by Boardman-Vogt form a cellular complex endowed with
good homotopical properties.

Recall that we consider nothing but the chain complex of
Boardman-Vogt' cellular construction. As a byproduct, one can
observe simply that $W(\P)$ is equipped with a diagonal and forms a
Hopf operad if $\P$ is so because $W(\P)$ is defined by a cubical
cellular complex.

\subsubsection{The chain interval}\label{item:ChainInterval} In this memoir we let
$\I$\glossary{$\I$} denote the standard cellular complex of the
topological interval $[0,1]$. Explicitly, the dg-module $\I$ is
generated by elements $\x{0},\x{1}$ in degree $0$, by an element
$\x{01}$ in degree $1$ and is equipped with the differential such
that $\delta(\x{01}) = \x{1}-\x{0}$. One should not be confused by
this notation $\x{0}$ for a basis element because we consider only
basis elements in $\I$. Recall that $\I$ defines a cylinder object
in the category of dg-modules. To be explicit, we have a morphism
$\eta^0: \F\rightarrow\I$\glossary{$\eta^0$}, respectively $\eta^1:
\F\rightarrow\I$\glossary{$\eta^1$}, that maps the ground field to
the submodule $\F\x{0}\subset\I$, respectively $\F\x{1}\subset\I$,
and an augmentation $\epsilon:
\I\xrightarrow{\sim}\F$\glossary{$\epsilon$} that cancels
$\x{01}\in\I$ and such that $\epsilon\cdot\eta^0 =
\epsilon\cdot\eta^1 = \Id$.

In addition we have a morphism of dg-modules $\mu:
\I\otimes\I\rightarrow\I$ induced by the map $(s,t)\mapsto\max(s,t)$
of the topological interval. Explicitly, this morphism satisfies
\begin{equation*}
\mu(\x{0},\x{0}) = \x{0},\qquad\mu(\x{01},\x{0}) = \mu(\x{0},\x{01})
= \x{01},\qquad \mu(\x{0},\x{1}) = \mu(\x{1},\x{0}) =
\mu(\x{1},\x{1}) = \x{1}
\end{equation*}
and vanishes in the other cases.

Observe that $\I$ is equipped with a coassociative diagonal $\Delta:
\I\rightarrow\I\otimes\I$ defined by the classical formulas
\begin{equation*}
\Delta(\x{0}) = \x{0}\otimes\x{0},\quad\Delta(\x{1}) =
\x{1}\otimes\x{1},\quad\text{and}\quad\Delta(\x{01}) =
\x{0}\otimes\x{01} + \x{01}\otimes\x{1}.
\end{equation*}
One checks readily that the morphisms $\eta^0,\eta^1,\epsilon$ and
$\mu$ are morphisms of dg-coalgebras.

Recall that, according to~\cite{BergerMoerdijk}, the background of
the Boardman-Vogt construction is provided by a monoidal model
category equipped with an interval. The dg-module $\I$ is precisely
the standard interval in the category of dg-modules and in the
category of dg-coalgebras.\index{chain interval}

\subsubsection{Tree structures}\label{item:TreeRecalls}
As specified in the introduction, for tree structures, we adopt the
conventions of our article~\cite{OperadTextbook} and we refer to
this reference for more precision.

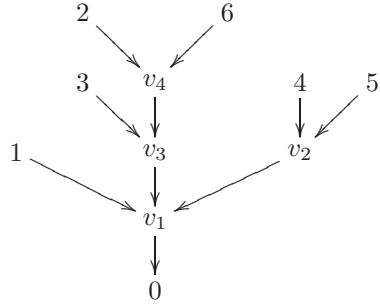
\begin{figure}[b]%
\[\vcenter{\xymatrix@M=0pt@R=5mm@C=5mm{
& *+<2mm>{2}\ar[dr] & & *+<2mm>{6}\ar[dl] & & \\
& *+<2mm>{3}\ar[dr] & *+<2mm>{v_4}\ar[d] & & *+<2mm>{4}\ar[d] & *+<2mm>{5}\ar[dl] \\
*+<2mm>{1}\ar[drr] & & *+<2mm>{v_3}\ar[d] & & *+<2mm>{v_2}\ar[dll] & \\
& & *+<2mm>{v_1}\ar[d] & & & \\
& & *+<2mm>{0} & & & \\ }}\] \caption{The structure of a
tree}\label{figure:AbstractTree}
\end{figure}

Recall that an \emph{$r$-tree}\index{tree!$r$-tree}\index{$r$-tree}
refers to an abstract tree $\tau$ defined by a set of vertices
$V(\tau)$\glossary{$V(\tau)$} and by a set of edges
$E(\tau)$\glossary{$E(\tau)$}, oriented from a source to a target,
equipped with one outgoing edge (the \emph{root}\index{tree!root of
a}\index{root of a tree} of the tree) that targets to an element
denoted by~$0$ and $r$~ingoing edges (the
\emph{leaves}\index{tree!leaves of a}\index{leaves of a tree} of the
tree) in one-to-one correspondence with the elements of
$\{1,\ldots,r\}$ that give the source of these edges (see
figure~\ref{figure:AbstractTree}).

Formally the source of an edge $e\in\E(\tau)$ is specified by an
element $s(e)\in V(\tau)\amalg\{1,\dots,r\}$ and the target by an
element $t(e)\in V(\tau)\amalg\{0\}$. For $x\in
V(\tau)\amalg\{1,\dots,r\}$, we assume that there is one and only
one edge $e\in E(\tau)$ with source $s(e) = x$, so that the
functions $s,t$ define a tree properly, and we assume that there is
one and only one edge such that $t(e) = 0$, so that this edge
represents the root of the tree. The leaves of $\tau$ are the edges
with source $s(e)\in\{1,\dots,r\}$.

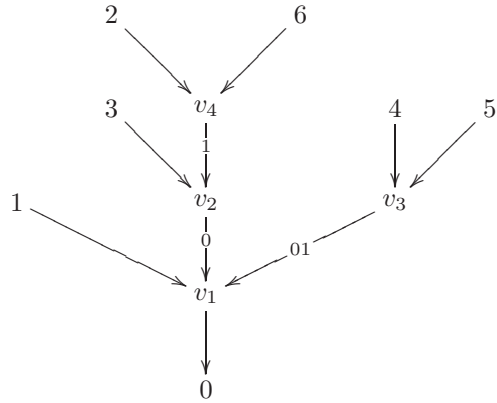
\begin{figure}[b]
\[\vcenter{\xymatrix@M=0pt@!R=4mm@!C=4mm{
& *+<2mm>{2}\ar[dr] & & *+<2mm>{6}\ar[dl] & & \\
& *+<2mm>{3}\ar[dr] & *+<2mm>{v_4}\ar[d]|(0.4)*+<2pt>{\scriptstyle\x{1}} & & *+<2mm>{4}\ar[d] & *+<2mm>{5}\ar[dl] \\
*+<2mm>{1}\ar[drr] & & *+<2mm>{v_2}\ar[d]|(0.4)*+<2pt>{\scriptstyle\x{0}} & & *+<2mm>{v_3}\ar[dll]|*+<2pt>{\scriptstyle\x{01}} & \\
& & *+<2mm>{v_1}\ar[d] & & & \\
& & *+<2mm>{0} & & & \\ }}\] \caption{A cell metric
tree}\label{figure:MetricTree}
\end{figure}

Recall that the tree structure is uniquely determined by a partition
$V(\tau)\amalg\{1,\ldots,r\} = \coprod_{v\in V(\tau)\amalg\{0\}}
I_v$\glossary{$I_v$}, indexed by the set $V(\tau)\amalg\{0\}$,
characterized by $x\in I_y$ if and only if $x$ and $y$ are
respectively the source and the target of an edge of the tree. Hence
the component $I_v$ associated to a vertex $v$ represents the set of
entries of the vertex $v$ in the tree $\tau$. In the example
represented in figure \ref{figure:MetricTree} we have $V(\tau) =
\{v_1,v_2,v_3,v_4\}$, and $I_{0} = \{v_1\}$, $I_{v_1} =
\{1,v_2,v_3\}$, $I_{v_2} = \{3,v_4\}$, $I_{v_3} = \{4,5\}$, $I_{v_4}
= \{2,6\}$.

The set of \emph{internal edges}\index{tree!internal edge of
a}\index{internal edge of a tree} of the tree $\tau$, denoted by
$E'(\tau)$\glossary{$E'(\tau)$}, consists of the edges $e$ which are
neither a leave or the root of $\tau$. Thus an edge $e$ is internal
if and only if $s(e),t(e)\in V(\tau)$.

\subsubsection{Tree morphisms and edge contractions}\label{item:TreeMorphisms}
\begin{figure}[b]%
\[\vcenter{\xymatrix@M=0pt@R=5mm@C=5mm{
& *+<2mm>{2}\ar[dr] & & *+<2mm>{6}\ar[dl] & & \\
& *+<2mm>{3}\ar[dr] & *+<2mm>{v_4}\ar[d] & & *+<2mm>{4}\ar[d] & *+<2mm>{5}\ar[dl] \\
*+<2mm>{1}\ar[drr] & & *+<2mm>{v_3}\ar[d] & & *+<2mm>{v_2}\ar[dll] & \\
& & *+<2mm>{v_1}\ar[d] & & & \\
& & *+<2mm>{0} & & & \\ }} \rightarrow
\vcenter{\xymatrix@M=0pt@R=5mm@C=5mm{
& *+<2mm>{3}\ar[dr] & *+<2mm>{2}\ar[d] & *+<2mm>{6}\ar[dl] & & \\
*+<2mm>{1}\ar[drr] & & *+<2mm>{w_2}\ar[d] & & *+<2mm>{4}\ar[dll] & *+<2mm>{5}\ar[dlll] \\
& & *+<2mm>{w_1}\ar[d] & & & \\
& & *+<2mm>{0} & & & \\ }}\] \caption{An example of a tree
morphism}\label{figure:AbstractTreeMorphism}
\end{figure}
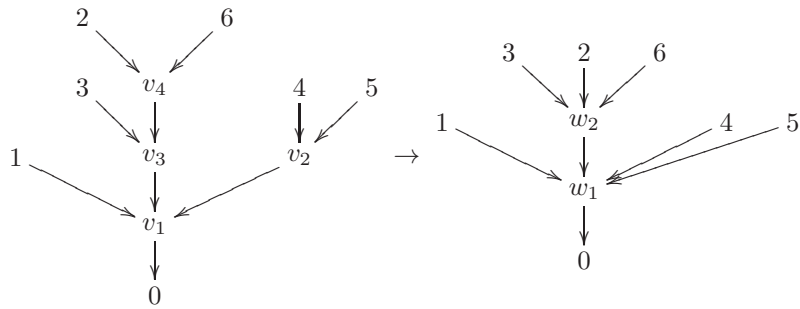

In our construction we consider a category of $r$-trees, denoted by
$\Theta(r)$\glossary{$\Theta(r)$}, in which a morphism $f:
\tau\rightarrow\tau'$ is defined by a contraction of internal edges
$e\in E(\tau)$. Formally, a morphism of $r$-trees $f:
\tau\rightarrow\tau'$\index{tree!morphism of
$r$-trees}\index{morphism of $r$-trees} is defined by a pair of maps
$f_V: V(\tau)\amalg\{0,1,\dots,r\}\rightarrow
V(\tau')\amalg\{0,1,\dots,r\}$ and $f_E: E(\tau)\rightarrow
E(\tau')\amalg V(\tau')$ endowed with the following properties:
\begin{enumerate}
\item\label{item:ExternalEdgeImage}
the map $f_V$ is the identity on $\{0,1,\dots,r\}$;
\item\label{item:ContractedEdge}
for an edge $e\in E(\tau)$ such that $e\in f^{-1}_E V(\tau')$ we
have $f_V(s(e)) = f_V(t(e)) = f_E(e)$;
\item\label{item:UncontractedEdge}
for an edge $e\in E(\tau)$ such that $e\in f^{-1}_E E(\tau')$ we
have $f_V(s(e)) = s(f_E(e))$ and $f_V(t(e)) = t(f_E(e))$.
\end{enumerate}
An example of a tree morphism is represented in
figure~\ref{figure:AbstractTreeMorphism}. In this example the map
$f_V$ is defined by $f_V(v_3) = f_V(v_4) = w_2$ and $f_V(v_1) =
f_V(v_2) = w_1$. Observe that the previous assertions imply that the
map $f_V$ is onto and $f_E$ induces a one-to-one correspondence
between $f^{-1}_E E(\tau')\subset E(\tau)$ and $E(\tau')$. Moreover,
the map $f_E$ is clearly uniquely determined by $f_V$ and
conversely.

Intuitively, the edges $e\in E(\tau)$ such that $e\in f^{-1}_E
V(\tau')$ are contracted to the vertex $v' = f_E(e)$ by the morphism
$f: \tau\rightarrow\tau'$ and the other edges are preserved. Our
assumptions imply that the leaves and the root are fixed by a
morphism of $r$-trees and, as a consequence, only internal edges are
allowed to be contracted. Furthermore, we can observe that the
subsets $f^{-1}_V(v')\subset V(\tau)$ and $f^{-1}_E(v')\subset
E(\tau)$ associated to a vertex $v'\in V(\tau')$ determine a subtree
$\sigma_{v'}$ of~$\tau$. In addition the entry set of $v'$ can be
identified with the image under the map $f_V$ of the entry set of
this subtree $\sigma_{v'}$ (we refer to~\cite{OperadTextbook} for
the formal definition of these notions). In fact, the tree $\tau'$
can be determined by a generalization of the quotient process
of~\cite{OperadTextbook}. Namely, if we perform a sequence of
quotients by the subtrees $\sigma_{v'}$ for $v'\in V(\tau')$, then
we obtain the tree $\tau'$ up to isomorphism.

Clearly, a morphism of $r$-trees is an isomorphism if and only if
the map $f_V$ is one-to-one and $f_E$ defines a one-to-one
correspondence between edge sets.\index{tree!isomorphism of
$r$-trees}

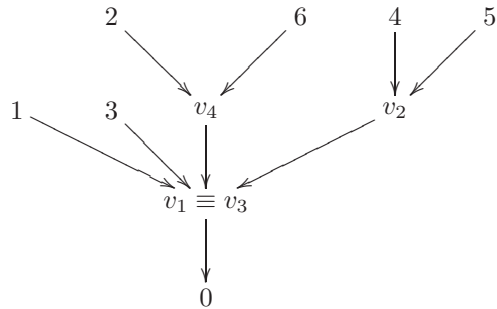
\begin{figure}[b]%
\[\vcenter{\xymatrix@M=0pt@!R=4mm@!C=4mm{
& *+<2mm>{2}\ar[dr] & & *+<2mm>{6}\ar[dl] & *+<2mm>{4}\ar[d] & *+<2mm>{5}\ar[dl] \\
*+<2mm>{1}\ar[drr] & *+<2mm>{3}\ar[dr] & *+<2mm>{v_4}\ar[d] & & *+<2mm>{v_2}\ar[dll] & \\
& & *+<2mm>{v_1\equiv v_3}\ar[d] & & & \\
& & *+<2mm>{0} & & & \\ }}\] \caption{An edge
contraction}\label{figure:EdgeContraction}
\end{figure}

One can also observe that a morphism of $r$-trees $f:
\tau\rightarrow\tau'$ that contracts a single edge $e_0\in E(\tau)$
is equivalent to the contraction process
$\tau\mapsto\tau/e_0$\glossary{$\tau/e_0$} defined
in~\cite{OperadTextbook}. To be precise, for an internal edge
$e_0\in E(\tau)$ one considers the tree $\tau/e_0$ obtained by
identifying the source $s_0 = s(e_0)$ and the target $t_0 = t(e_0)$
of $e_0$ in $\tau$ to a single vertex $s_0\equiv t_0$ whose entries
are defined by the union $I_{t_0}\setminus\{s_0\}\amalg I_{s_0}$ of
the entries of $s_0$ and $t_0$ in $\tau$. Hence we have $V(\tau/e_0)
= V(\tau)/\{s_0\equiv t_0\}$ and $E(\tau/e_0) =
E(\tau)\setminus\{e_0\}$. Furthermore, the quotient map
$(\gamma_{e_0})_V: V(\tau)\rightarrow V(\tau/e_0)$ and the map
$(\gamma_{e_0})_E: E(\tau)\rightarrow E(\tau/e_0)\amalg V(\tau/e_0)$
such that $(\gamma_{e_0})_E(e_0) = \{s_0\equiv t_0\}$ and
$(\gamma_{e_0})_E(e) = e$ for $e\not=e_0$ defines a morphism of
$r$-trees $\gamma_{e_0}: \tau\rightarrow\tau/e_0$ which is called an
\emph{edge contraction}\index{edge contraction}\index{tree!edge
contraction in a}. As an example, the contraction of the edge
$v_3\rightarrow v_1$ in the tree of figure~\ref{figure:AbstractTree}
gives the tree represented in figure~\ref{figure:EdgeContraction}.

Clearly, a morphism of $r$-trees is a composite of tree isomorphisms
and edge contractions.

\subsubsection{Tree categories}\label{item:TreeCategories}
We call a tree \emph{$n$-reduced}\index{tree!$n$-reduced} if all
vertices $v\in V(\tau)$ have more than $n$ entries. For instance, a
tree $\tau$ is called $0$-reduced if it has no \emph{terminal
vertex} (a vertex $v$ is terminal if $I_v = \emptyset$).

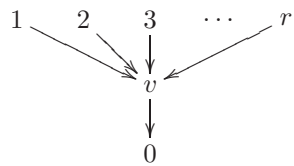
\begin{figure}[b]%
\[\vcenter{\xymatrix@M=0pt@R=5mm@C=5mm{
*+<2mm>{1}\ar[drr] & *+<2mm>{2}\ar[dr] & *+<2mm>{3}\ar[d] & *{\dots} & *+<2mm>{r}\ar[dll] \\
& & *+<2mm>{v}\ar[d] & & \\
& & *+<2mm>{0} & & \\ }}\] \caption{The terminal
$r$-tree}\label{figure:TerminalTree}
\end{figure}

Observe that a $0$-reduced $r$-tree $\tau$ has no automorphism.
Accordingly, if we fix a representative for each isomorphism class
of $0$-reduced tree, then we obtain a category, denoted by
$\Theta'(r)$\glossary{$\Theta'(r)$}, equivalent to the full
subcategory of $\Theta(r)$ generated by $0$-reduced trees in which
all isomorphisms are identities. Hence any morphism in $\Theta'(r)$
is a composite of edge contractions. Moreover, one can observe that
$\Theta'(r)$ defines a poset equipped with a terminal element
represented by the unique $r$-tree $\tau_r$\glossary{$\tau_r$} with
one vertex and no internal edge (see figure
\ref{figure:TerminalTree}).

The category $\Theta'(r)$ is also equipped with a natural grading
$\gr: \Theta'(r)\rightarrow\N$ given by the number of internal edges
and we let $\Theta'_d(r)$\glossary{$\Theta'_d(r)$} denote the
subcategory of $\Theta'(r)$ generated by trees of grading $\leq d$.
Observe that $\Theta'_0(r)$ contains only of the terminal $r$-tree
$\tau_r$. Clearly, any non-identity morphism of $\Theta'(r)$
decreases the grading so that $\Theta'(r)$ defines an inverse Reedy
category.

In fact, in the constructions of this section we consider only the
subcategory of $\Theta'(r)$ formed by $1$-reduced $r$-trees.
Therefore we introduce the notation
$\Theta''(r)$\glossary{$\Theta''(r)$} for this category and we let
also $\Theta''_d(r)$\glossary{$\Theta''(r)$} denote the subcategory
of $\Theta'(r)$ formed by $1$-reduced $r$-trees with no more than
$d$ internal edges.

\subsubsection{Cell metric trees and length tensors}\label{item:MetricTrees}
In the definition of the Boardman-Vogt complex we consider
\emph{cell metric trees}\index{tree!cell metric} $\tau$ equipped
with a \emph{length tensor}\index{length tensor} defined by a tensor
product of interval elements indexed by the internal
edges~of~$\tau$:
\begin{equation*}
\lambda = \bigotimes_{e\in E'(\tau)}\lambda_e\in\bigotimes_{e\in
E'(\tau)}\I.
\end{equation*}
In fact, we identify abusively a length tensor with a basis element
of the tensor product $\bigotimes_{e\in E'(\tau)}\I$, so that a
length tensor is defined by a map $e\mapsto\lambda_e$ which
associates to any internal edge $e\in E'(\tau)$ an element
$\lambda_e\in\I$. In the representation of a tree we decorate the
internal edges $e$ by the corresponding length $\lambda_e$ as in
figure~\ref{figure:MetricTree}.

The module of length tensors\index{length tensor!module of length
tensors} is denoted by
\begin{equation*}
\Cube{\tau} = \bigotimes_{e\in E'(\tau)}\I.
\end{equation*}\glossary{$\Cube{\tau}$}
In fact, this module is nothing but the cellular complex of the cube
build on the internal edges of the tree $\tau$. One observes that
the map $\tau\mapsto\Cube{\tau}$ can be extended to a contravariant
functor from the category of trees to the category of dg-modules.
Explicitly, a tree morphism $f: \tau\rightarrow\tau'$ induces a
dg-module morphism $f^*: \Cube{\tau'}\rightarrow\Cube{\tau}$ that
assigns the length $\lambda_e = \x{0}$ to the edges $e\in E'(\tau)$
which are contracted to a vertex by $f: \tau\rightarrow\tau'$ and
preserves the length of the other edges $e\in E'(\tau)$ which are
mapped to an edge of $\tau'$. Formally, a morphism of trees $f:
\tau\rightarrow\tau'$ yields a partition $E'(\tau) = f^{-1}_E
V(\tau')\amalg f^{-1}_E E'(\tau')$ and $f$ induces a bijection from
$f^{-1}_E E'(\tau')$ to $E'(\tau')$. The dg-module morphism $f^*:
\Cube{\tau'}\rightarrow\Cube{\tau}$ is given by the tensor product
of the canonical isomorphism
\begin{equation*}
\bigotimes_{e'\in E'(\tau')}\I\simeq\bigotimes_{e\in f_E^{-1}
E'(\tau')}\I
\end{equation*}
with the unit map
\begin{equation*}
\bigotimes_{e\in f_E^{-1} V(\tau')}\eta^0: \bigotimes_{e\in f_E^{-1}
V(\tau')}\F\x{0} \rightarrow\bigotimes_{e\in f_E^{-1} V(\tau')}\I.
\end{equation*}
Hence this morphism identifies the module
\begin{equation*}
\Cube{\tau'} = \bigotimes_{e'\in E'(\tau')} \I
\end{equation*}
with the module of length tensors
\begin{equation*}
\Bigl[\bigotimes_{e\in f_E^{-1}
E'(\tau')}\I\Bigr]\otimes\Bigl[\bigotimes_{e\in f_E^{-1} V(\tau')}
\x{0}\Bigr] \subset\bigotimes_{e\in E'(\tau)}\I = \Cube{\tau}
\end{equation*}
in which the edges $e\in f^{-1}_E V(\tau')$, that are contracted to
a vertex in $\tau'$, have length $\lambda_e = \x{0}$. As an example,
an edge contraction $\gamma_{e_0}: \tau\rightarrow\tau/e_0$ induces
a dg-module morphism $(\gamma_{e_0})^*:
\Cube{\tau/e_0}\rightarrow\Cube{\tau}$ that identifies the module
$\Cube{\tau/e_0}$ to the module of length tensors
$\lambda\in\Cube{\tau}$ such that $\lambda_{e_0} = \x{0}$.

Finally, as our chain interval $\I$, equipped with a coassociative
diagonal, defines an interval in the category of dg-coalgebras and
not only to the category of dg-modules, we observe that the
construction of this paragraph gives a functor
$\tau\mapsto\Cube{\tau}$ from the category of trees to the category
of dg-coalgebras as well. Explicitly, the diagonal of~$\Cube{\tau}$
is given by the composite of tensor products of the diagonal of $\I$
with the obvious tensor permutation:
\begin{equation*}
\bigotimes_{e\in E'(\tau)}\I\rightarrow\bigotimes_{e\in E'(\tau)}
(\I\otimes\I) \xrightarrow{\simeq}\Bigl[\bigotimes_{e\in
E'(\tau)}\I\Bigr]\otimes\Bigl[\bigotimes_{e\in E'(\tau)}\I\Bigr].
\end{equation*}
Equivalently, in the definition of $\Cube{\tau} = \bigotimes_{e\in
E'(\tau)}\I$ we consider a tensor product in the monoidal category
of dg-coalgebras and not only in the category of dg-modules. The
morphism $f^*: \Cube{\tau'}\rightarrow\Cube{\tau}$ induced by a tree
morphism defines clearly a morphism of dg-coalgebras as the
definition of $f^*$ can be deduced from the axioms of symmetric
monoidal categories.

\subsubsection{Coends over the category of trees}\label{item:TreeCoends}
In this section we consider coends
\begin{equation*}
\int_{\tau\in\Theta''(r)} \Phi_\tau\otimes\Cube{\tau}
\end{equation*}
associated to covariant functors $\tau\mapsto\Phi_\tau$ from the
category of $1$-reduced trees $\Theta''(r)$
(see~\ref{item:TreeCategories}) to the category of dg-modules,
respectively dg-coalgebras. Recall that the coends in the category
of dg-coalgebras are created in the category of dg-modules like all
colimits.

We aim to determine the homotopy type of these coends. We prove
precisely that $\int_{\tau\in\Theta''(r)}
\Phi_\tau\otimes\Cube{\tau}$ is homotopy equivalent to the object
$\Phi_{\tau_r}$ associated to the terminal $r$-tree $\tau_r$. We
prove in addition that the coend morphism
\begin{equation*}
\phi_*: \int_{\tau\in\Theta''(r)}
\Phi_\tau\otimes\Cube{\tau}\rightarrow\int_{\tau\in\Theta''(r)}
\Psi_\tau\otimes\Cube{\tau}
\end{equation*}
induced by a functor morphism $\phi_\tau:
\Phi_\tau\rightarrow\Psi_\tau$ forms a cofibration, respectively an
acyclic cofibration, if $\phi_\tau: \Phi_\tau\rightarrow\Psi_\tau$
is a pointwise cofibration, respectively a pointwise acyclic
cofibration.

For these purposes we assume that the category of contravariant
functors from the category of trees $\Theta''(r)$ to the category of
dg-modules, respectively to the category of dg-coalgebras, is
equipped with a Reedy model structure. Explicitly, for contravariant
functors $\tau\mapsto C^\tau$ we have a latching object, defined by
\begin{equation*}
\Latch{C}^\tau = \colim_{\tau\xrightarrow{\not=}\tau'} C^{\tau'},
\end{equation*}
and a functor morphism $\phi_\tau: C^\tau\rightarrow D^\tau$ is a
Reedy cofibration, respectively an acyclic Reedy cofibration, if and
only if, for each tree $\tau\in\Theta''(r)$, the morphism
\begin{equation*}
(\phi,\lambda):
C^\tau\bigoplus_{\Latch{C}^\tau}\Latch{D}^\tau\rightarrow D^\tau
\end{equation*}
defines a cofibration, respectively an acyclic cofibration, in the
ground category. Fibrations and weak-equivalences are defined
pointwise. The definitions are dual for covariant functors
$\tau\mapsto\Phi_\tau$ but we do not consider fibrations and
matching objects explicitly in this case.

Recall that a morphism of dg-coalgebras is a cofibration,
respectively an acyclic cofibration, if and only if it defines a
cofibration, respectively an acyclic cofibration, in the category of
dg-modules. The same statement holds for functor morphisms since
colimits and latching objects in the category of dg-coalgebras  are
created in the category of dg-modules. As a byproduct, the
dg-coalgebra case of the next statements is always an immediate
consequence of the dg-module case. Therefore we give proofs in the
context of dg-modules and omit to give more precision in the context
of dg-coalgebras. On the other hand, one can observe that our
assertions and our arguments can be extended to monoidal model
categories including both the category of dg-modules and the
category of dg-coalgebras.

We deduce our results from forthcoming observations about the
functor $\tau\mapsto\Cube{\tau}$ defined by the modules of length
tensors.

\subsubsection{Homotopy properties of the modules of length tensors}
We observe that the functor $\tau\mapsto\Cube{\tau}$ is connected to
the constant functor $\F$ by equivalences of functors on
$\Theta''(r)$. To be more explicit, for each tree $\tau$ we have a
morphism $\eta^0_*: \F\rightarrow\Cube{\tau}$ induced by $\eta^0:
\F\rightarrow\I$. One checks readily that $\eta^0_*$ commutes with
the morphisms $f^*: \Cube{\tau'}\rightarrow\Cube{\tau}$ induced by a
tree morphism $f: \tau\rightarrow\tau'$ so that $\eta^0_*$ defines a
morphism of functors on the category of trees. Similarly, one can
check that the augmentation $\epsilon: \I\rightarrow\F$ induces a
morphism of functors $\epsilon_*: \Cube{\tau}\rightarrow\F$ such
that $\epsilon_*\cdot\eta^0_* = \Id$. Finally, one can observe that
$\eta^0_*: \F\rightarrow\Cube{\tau}$ and $\epsilon_*:
\Cube{\tau}\rightarrow\F$ are coalgebra morphisms because so are
$\eta^0: \F\rightarrow\I$ and $\epsilon: \I\rightarrow\F$. We have
the following assertion:

\begin{obsv}\label{obsv:LengthTensorsHomotopy}
The morphism $\eta^0_*: \F\rightarrow\Cube{\tau}$ is an acyclic
Reedy cofibration in the category of contravariant functors from
$\Theta''(r)$ to the category of dg-modules, respectively to the
category of dg-coalgebras. As a corollary, the morphism $\epsilon_*:
\Cube{\tau}\rightarrow\F$, which is right inverse to $\eta^0_*$,
defines a weak equivalence of functors as well.
\end{obsv}

\begin{proof}
Recall that $\eta^0_*$ defines an acyclic Reedy cofibration if and
only if the morphism
\begin{equation*}
(\eta^0_*,\lambda):
\F\bigoplus_{\Latch{\F}^\tau}\Latch{\Cube{\tau}}\rightarrow\Cube{\tau}
\end{equation*}
defines an acyclic cofibration in the category of dg-modules for all
$\tau\in\Theta''(r)$.

Clearly, for the constant functor we have $\Latch{\F}^{\tau} = 0$ if
$\tau\not=\tau_r$, the terminal $r$-tree, and $\Latch{\F}^{\tau_r} =
\F$. Hence we obtain
\begin{equation*}
\F\bigoplus_{\Latch{\F}^\tau}\Latch{\Cube{\tau}} =
\Latch{\Cube{\tau}}\quad\text{if $\tau\not=\tau_r$}
\qquad\text{and}\qquad\F\bigoplus_{\Latch{\F}^{\tau_r}}\Latch{\Cube{\tau_r}}
= \F = \Cube{\tau_r}.
\end{equation*}
(Recall that $\Cube{\tau_r} = \F$ since $\tau_r$ has no internal
edge.) Consequently, the morphism $\eta^0_*:
\F\rightarrow\Cube{\tau}$ is an acyclic Reedy cofibration if an only
if the canonical map $\Latch{\Cube{\tau}}\rightarrow\Cube{\tau}$ is
an acyclic cofibration of dg-modules (or dg-coalgebras) for
$\tau\not=\tau_r$.

This assertion can be deduced from the axioms of monoidal model
categories and from the general properties of an interval. On the
other hand, in the framework of dg-modules one can observe simply
that the dg-module $\Latch{\Cube{\tau}}$ is identified with the
submodule of $\Cube{\tau}$ spanned by length tensors $\lambda =
\bigotimes_e \lambda_e$ in which some edges have length $\lambda_e =
0$. Consequently, the quotient $\Cube{\tau}/\Latch{\Cube{\tau}}$ can
be identified with a tensor product of dg-modules
$\F\x{01}\xrightarrow{\partial_1}\F\x{1}$ which are clearly acyclic.
\end{proof}

By the way, we observe in this proof that the natural morphism
\begin{equation*}
\Latch{\Cube{\tau}}{\tau} = \colim_{f: \tau'\xrightarrow{\not=}\tau}
\Cube{\tau'}\rightarrow\Cube{\tau}
\end{equation*}
defines an embedding from the latching object $\Latch{\Cube{\tau}}$
to $\Cube{\tau}$. Accordingly, we have the following result:

\begin{obsv}\label{obsv:LengthTensorCofibrancy}
The functor $\tau\mapsto\Cube{\tau}$ defines a Reedy cofibrant
object in the category of contravariant functors from the category
of trees $\Theta''(r)$ to the category of dg-modules, respectively
to the category of dg-coalgebras.\qed
\end{obsv}

We can now prove the results announced in~\ref{item:TreeCoends}.
First, we have the following assertion that arises as a consequence
of observation~\ref{obsv:LengthTensorsHomotopy}:

\begin{claim}\label{claim:CoendCofibrations}
The coend morphism
\begin{equation*}
\eta^0_*: \int_{\tau\in\Theta''(r)}\Phi_\tau\otimes\F
\rightarrow\int_{\tau\in\Theta''(r)}\Phi_\tau\otimes\Cube{\tau}
\end{equation*}
induced by $\eta^0_*: \F\rightarrow\Cube{\tau}$ defines an acyclic
cofibration of dg-modules, respectively dg-coalgebras. As a
consequence, the morphism $\epsilon_*: \Cube{\tau}\rightarrow\F$,
which is right inverse to $\eta^0_*$, defines an inverse weak
equivalence of dg-modules, respectively dg-coalgebras,
\begin{equation*}
\epsilon_*: \int_{\tau\in\Theta''(r)}\Phi_\tau\otimes\Cube{\tau}
\xrightarrow{\sim}\int_{\tau\in\Theta''(r)}\Phi_\tau\otimes\F
\end{equation*}
such that $\epsilon_*\eta^0_* = \Id$.
\end{claim}

\begin{proof}
To be explicit, recall that $\eta^0_*: \F\rightarrow\Cube{\tau}$
defines by observation~\ref{obsv:LengthTensorsHomotopy} an acyclic
Reedy cofibration of functors on $\Theta''(r)$. Accordingly, the
claim is an instance of a general statement that can be deduced from
the categorical definition of a coend and from the axioms of
monoidal model categories. Explicitly, for a fibration of dg-modules
$p: C\fibration D$ the lifting problem
\begin{equation*}
\xymatrix{ \int_{\tau\in\Theta''(r)}\Phi_\tau\otimes\F\ar[r]\ar[d] & C\ar@{->>}[d]^{p} \\
\int_{\tau\in\Theta''(r)}\Phi_\tau\otimes\Cube{\tau}\ar[r]\ar@{-->}[ur]
& D }
\end{equation*}
is equivalent to an adjoint lifting problem
\begin{equation*}
\xymatrix{ \F\ar[r]\ar[d]_{\eta^0_*} & \DGHom(\Phi_\tau,C)\ar[d]^{p_*} \\
\Cube{\tau}\ar[r]\ar@{-->}[ur] & \DGHom(\Phi_\tau,D) }
\end{equation*}
in the category of contravariant functors in $\tau\in\Theta''(r)$.
This lifting problem admits a solution since a fibration of
dg-modules $p: C\rightarrow D$ induces a Reedy fibration on internal
hom-objects $p_*: \DGHom(\Phi_\tau,C)\fibration\DGHom(\Phi_\tau,D)$
and since $\eta^0_*: \F\rightarrow\Cube{\tau}$ is an acyclic Reedy
cofibration. Accordingly, our morphism has the left lifting property
with respect to fibrations of dg-modules and the conclusion follows.
\end{proof}

As the category of trees admits a terminal object $\tau_r$, we
obtain in addition:

\begin{obsv}
We have an isomorphism $\int_{\tau\in\Theta''(r)}
\Phi_\tau\otimes\F\simeq\Phi_{\tau_r}$.\qed
\end{obsv}

Accordingly, we obtain finally:

\begin{lemm}\label{lemm:TreewiseCoends}
Let $\tau\mapsto\Phi_\tau$ denote a covariant functor from the
category of trees to the category of dg-modules, respectively
dg-coalgebras.

The morphism $\epsilon_*: \Cube{\tau}\rightarrow\F$ induces a
weak-equivalence of dg-modules, respectively dg-coalgebras,
\begin{equation*}
\int_{\tau\in\Theta''(r)}\Phi_\tau\otimes\Cube{\tau}
\xrightarrow[\epsilon_*]{\sim}\int_{\tau\in\Theta''(r)}\Phi_\tau\otimes\F
\simeq\Phi_{\tau_r},
\end{equation*}
which has a natural inverse equivalence, such that
$\epsilon_*\eta^0_* = \Id$, induced by $\eta^0_*:
\F\rightarrow\Cube{\tau}$.\qed
\end{lemm}

As stated in~\ref{item:TreeCoends}, we have also the following
assertion which arises as a consequence of
observation~\ref{obsv:LengthTensorCofibrancy}:

\begin{lemm}\label{lemm:TreewiseCoendCofibrations}
Let $\phi_\tau: \Phi_\tau\rightarrow\Psi_\tau$ denote a pointwise
cofibration of covariant functors from the category of trees to the
category of dg-modules, respectively to the category of
dg-coalgebras. The associated coend morphism
\begin{equation*}
\phi_*: \int_{\tau\in\Theta''(r)}\Phi_\tau\otimes\Cube{\tau}
\rightarrow\int_{\tau\in\Theta''(r)}\Psi_\tau\otimes\Cube{\tau}
\end{equation*}
defines a cofibration in the category of dg-modules, respectively in
the category of dg-coalgebras.
\end{lemm}

\begin{proof}
As in the proof of claim~\ref{claim:CoendCofibrations}, we observe
that a lifting problem in the category of dg-modules
\begin{equation*}
\xymatrix{ \int_{\tau\in\Theta''(r)}\Phi_\tau\otimes\Cube{\tau}\ar[r]\ar[d]_{\phi_*} & C\ar@{->>}[d]_{\sim}^{p} \\
\int_{\tau\in\Theta''(r)}\Psi_\tau\otimes\Cube{\tau}\ar[r]\ar@{-->}[ur]
& D },
\end{equation*}
where $p: C\wefib D$ is an acyclic fibration, is equivalent to an
adjoint lifting problem
\begin{equation*}
\xymatrix{ \Phi_\tau\ar[r]\ar[d]_{\phi} & \DGHom(\Cube{\tau},C)\ar[d]^{p_*} \\
\Psi_\tau\ar[r]\ar@{-->}[ur] & \DGHom(\Cube{\tau},D) }
\end{equation*}
in the category of functors in $\tau\in\Theta''(r)$. Recall that the
modules of length tensors $\Cube{\tau}$ define a Reedy cofibrant
object in the category of contravariant functors on~$\Theta''(r)$.
As a consequence, the natural transformation $p_*:
\DGHom(\Cube{\tau},C)\rightarrow\DGHom(\Cube{\tau},D)$, induced by
an acyclic fibration of dg-modules $p: C\wefib D$, defines a Reedy
fibration in the category of covariant functors on~$\Theta''(r)$.
Hence we conclude from our observation that the morphism $\phi_*:
\int_{\tau\in\Theta''(r)}\Phi_\tau\otimes\Cube{\tau}\rightarrow\int_{\tau\in\Theta''(r)}\Psi_\tau\otimes\Cube{\tau}$
satisfies the left lifting property with respect to acyclic
fibrations and defines a cofibration of dg-modules.
\end{proof}

\subsection{The Boardman-Vogt construction}
\label{subsection:BoardmanVogtConstruction} As announced in the
section introduction, the purpose of this subsection is to give a
precise account of the Boardman-Vogt construction in the
differential graded context and to prove the properties of this
construction asserted in theorem~\ref{thm:WHopfOperad} up to the
cofibrancy claim in assertion~(b).

In our construction, we consider only connected unital operads~$\P$
such that $\P(0) = \F$ and $\P(1) = \F$. In this context we set
$W(\P)(0) = \F$ and $W(\P)(1) = \F$ so that $W(\P)$ is still a
connected unital operad. In what follows we give first the
definition of the dg-modules $W(\P)(r)$ for $r\geq 2$. Then we
define the partial composites $\circ_i: W(\P)(s)\otimes
W(\P)(t)\rightarrow W(\P)(s+t-1)$ for operations such that $s,t\geq
2$ and we make a particular case for composites with unital
operations. In view of the construction of~\cite{BergerMoerdijkW},
the connectedness assumption $\P(1) = \F$ is not necessary but
simplifies our construction.

The assertion (c) of theorem~\ref{thm:WHopfOperad}, the Hopf operad
structure, arises from the definition of~$W(\P)$. Only assertions
(a) and (b) require an actual proof supplied by
claim~\ref{claim:WAugmentation} and claim~\ref{claim:WQuasiFree}
respectively.

\subsubsection{Treewise tensor
products}\label{item:TreewiseTensorProducts}\index{treewise tensor
product} For a finite set $I = \{i_1,\dots,i_r\}$, we form the
module $\P(I)$\glossary{$\P(I)$} whose elements represent operations
in $r$ variables indexed by $\{i_1,\dots,i_r\}$ instead of the
integers $\{1,\dots,r\}$. For finite sets $I,J$ and an element $i\in
I$, we have a partial composite operation $\circ_i:
\P(I)\otimes\P(J)\rightarrow\P(I\setminus\{i\}\amalg J)$.

As in~\cite{OperadTextbook}, given a tree $\tau$, we let
$\tau(\P)$\glossary{$\tau(\P)$} denote the module of tensors
\begin{equation*}
\ptensor = \bigotimes_{v\in V(\tau)} p_v\in\bigotimes_{v\in V(\tau)}
\P(I_v)
\end{equation*}
that represent labelings of vertices of $\tau$ by operations
$p_v\in\P(I_v)$ whose variables are in bijection with the entry set
of the associated vertex $I_v$. In fact, we shall consider only
$1$-reduced trees whose vertices have all at least to $2$ entries.
Thus, for the factors $p_v\in\P(I_v)$ of a labeling, the set $I_v$
contains at least two elements.\index{tree!labeling of an $r$-tree}

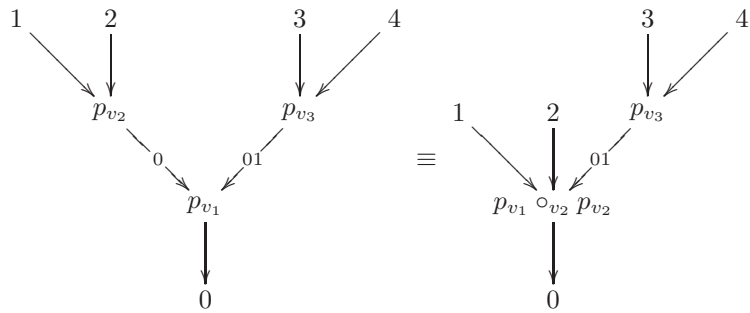
\begin{figure}[b]
\[\vcenter{\xymatrix@M=0pt@!R=4mm@!C=4mm{
*+<2mm>{1}\ar[dr] & *+<2mm>{2}\ar[d] & & *+<2mm>{3}\ar[d] & *+<2mm>{4}\ar[dl] \\
& *+<2mm>{p_{v_2}}\ar[dr]|*+<2pt>{\scriptstyle\x{0}} & & *+<2mm>{p_{v_3}}\ar[dl]|*+<2pt>{\scriptstyle\x{01}} & \\
& & *+<2mm>{p_{v_1}}\ar[d] & & \\
& & *+<2mm>{0} & & \\ }}
\equiv\vcenter{\xymatrix@M=0pt@!R=4mm@!C=4mm{
& & *+<2mm>{3}\ar[d] & *+<2mm>{4}\ar[dl] \\
*+<2mm>{1}\ar[dr] & *+<2mm>{2}\ar[d] & *+<2mm>{p_{v_3}}\ar[dl]|*+<2pt>{\scriptstyle\x{01}} & \\
& *+<2mm>{p_{v_1}\circ_{v_2} p_{v_2}}\ar[d] & & \\
& *+<2mm>{0} & & \\ }}\] \caption{The equivalence relation of the
$W$-construction}\label{figure:LabelledMetricTreeContraction}\end{figure}

One observes that the map $\tau\mapsto\tau(\P)$ can be extended to a
covariant functor from the category of trees to the category of
dg-modules. In particular, for each internal edge $e_0\in E'(\tau)$,
we have a natural morphism
\begin{equation*}
(\gamma_{e_0})_*: \tau(\P)\rightarrow\tau/e_0(\P)
\end{equation*}
from the labellings of the tree $\tau$ to the labellings of the tree
$\tau/e_0$ obtained by the contraction of the edge $e_0$. In general
a morphism $f_*: \tau(\P)\rightarrow\tau'(\P)$ induced by a tree
morphism $f: \tau\rightarrow\tau'$ is a composite of edge
contractions $(\gamma_{e_0})_*$ and isomorphisms. Therefore, for our
purposes, it is sufficient to make the morphisms $(\gamma_{e_0})_*$
explicit. Let $s_0 = s(e_0)$ and $t_0 = t(e_0)$ denote respectively
the source and the target of $e_0$. By definition, the morphism
$\gamma_{e_0}$ preserves the label of vertices $v\not=s_0,t_0$ which
are untouched by the contraction process and labels the collapsed
vertex $s_0\equiv t_0$ in $\tau/e_0$ by the partial composite
$p_{t_0}\circ_{s_0} p_{s_0}\in\P(I_{t_0}\setminus\{s_0\}\amalg
I_{s_0})$ of the labels
$p_{s_0}\in\P(I_{s_0}),p_{t_0}\in\P(I_{t_0})$ of the vertices $s_0$
and $t_0$ in $\tau$ (see figure
\ref{figure:LabelledMetricTreeContraction}).

Clearly, for the terminal $r$-tree $\tau_r$ we have a canonical
isomorphism $\tau_r(\P)\simeq\P(r)$. Accordingly, for each $r$-tree
$\tau$ the terminal morphism $\gamma: \tau\rightarrow\tau_r$ induces
a morphism $\gamma_*: \tau(\P)\rightarrow\P(r)$. Intuitively, a tree
labeling represents a formal operadic composite of operations and
the morphism $\gamma_*: \tau(\P)\rightarrow\P(r)$ is defined by the
evaluation of this composite operation in the operad $\P$.

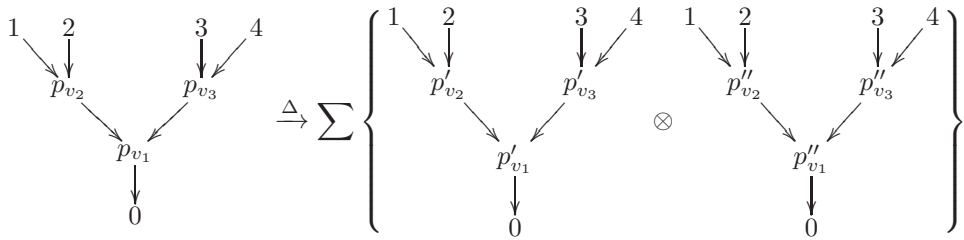
\begin{figure}[b]
\begin{equation*}
\vcenter{\xymatrix@M=0pt@R=5mm@C=3mm{
*+<1mm>{1}\ar[dr] & *+<1mm>{2}\ar[d] & & *+<1mm>{3}\ar[d] & *+<1mm>{4}\ar[dl] \\
& *+<1mm>{p_{v_2}}\ar[dr] & & *+<1mm>{p_{v_3}}\ar[dl] & \\
& & *+<1mm>{p_{v_1}}\ar[d] & & \\
& & *+<1mm>{0} & & }} \xrightarrow{\Delta} \sum\left\{
\vcenter{\xymatrix@M=0pt@R=5mm@C=3mm{
*+<1mm>{1}\ar[dr] & *+<1mm>{2}\ar[d] & & *+<1mm>{3}\ar[d] & *+<1mm>{4}\ar[dl] \\
& *+<1mm>{p'_{v_2}}\ar[dr] & & *+<1mm>{p'_{v_3}}\ar[dl] & \\
& & *+<1mm>{p'_{v_1}}\ar[d] & & \\
& & *+<1mm>{0} & & }} \otimes \vcenter{\xymatrix@M=0pt@R=5mm@C=3mm{
*+<1mm>{1}\ar[dr] & *+<1mm>{2}\ar[d] & & *+<1mm>{3}\ar[d] & *+<1mm>{4}\ar[dl] \\
& *+<1mm>{p''_{v_2}}\ar[dr] & & *+<1mm>{p''_{v_3}}\ar[dl] & \\
& & *+<1mm>{p''_{v_1}}\ar[d] & & \\
& & *+<1mm>{0} & & }}\right\}
\end{equation*}
\caption{The diagonal of a labeled
tree}\label{figure:LabelledTreeDiagonal}
\end{figure}

For a Hopf operad $\P$, the construction of this paragraph gives a
functor $\tau\mapsto\tau(\P)$ from the category of trees to the
category of dg-coalgebras and not only to the category of
dg-modules. Formally, in the definition of $\tau(\P) =
\bigotimes_{v\in V(\tau)} \P(I_v)$ we consider a tensor product in
the monoidal category of dg-coalgebras and not only in the category
of dg-modules (as in the definition of the module of length
tensors). Hence a morphism of $r$-trees $f: \tau\rightarrow\tau'$
induces a coalgebra morphism $f_*: \tau(\P)\rightarrow\tau'(\P)$
since the definition of $f_*$ can be deduced from the axioms of
symmetric monoidal categories and since a Hopf operad refers
precisely to an operad in the category of coalgebras. Explicitly,
for a fixed $r$-tree $\tau$, the diagonal of $\tau(\P) =
\bigotimes_{v\in V(\tau)} \P(I_v)$ is given by the composite of a
tensor product of the internal diagonal of $\P$ with the obvious
tensor permutation. Thus, for an element $\ptensor =
\bigotimes_{v\in V(\tau)} p_v\in\tau(\P)$, we have
\begin{equation*}
\Delta(\ptensor) = \sum \Bigl\{\bigotimes_{v\in V(\tau)}
p'_v\Bigr\}\otimes\Bigl\{\bigotimes_{v\in V(\tau)} p''_v\Bigr\}
\in\tau(\P)\otimes\tau(\P),
\end{equation*}
where we consider the expansion $\Delta(p_v) = \sum p'_v\otimes
p''_v$ of the diagonal of each factor $p_v\in\P(I_v)$. Accordingly,
the diagonal of a labeling is represented by the graphical formula
of figure~\ref{figure:LabelledTreeDiagonal}.

\subsubsection{The $W$-construction}\label{item:WComplex}
In this paragraph we define the chain structure of the Boardman-Vogt
operad $W(\P)$. As specified in the subsection introduction, we set
$W(\P)(0) = \F$ and $W(\P)(1) = \F$ so that $W(\P)$ is still a
connected unital operad. Consequently, in this paragraph and in the
next ones, we consider only the components $W(\P)(r)$ such that
$r\geq 2$ and we define the structure of the reduced operad
associated to $W(\P)$. The definition of the partial composition
products with unital operations $\circ_i: W(\P)(r)\otimes
W(\P)(0)\rightarrow W(\P)(r-1)$ is postponed
to~\ref{item:WUnitComposites}.

For $r\geq 2$, the module $W(\P)(r)$ is defined by the dg-module
coend
\begin{equation*}
W(\P)(r) = \int_{\tau\in\Theta''(r)} \tau(\P)\otimes\Cube{\tau}
\end{equation*}
in which the variable $\tau$ ranges over the category of $1$-reduced
$r$-trees $\Theta''(r)$ introduced in~\ref{item:TreeCategories}. The
map $\P\mapsto W(\P)$ defines clearly a functor on the category of
dg-operads.

Equivalently, the module $W(\P)(r)$ can be defined by the direct sum
\begin{equation*}
W(\P)(r) = \bigoplus_{\tau\in\Theta''(r)}
\tau(\P)\otimes\Cube{\tau}/\equiv
\end{equation*}
together with the relation
\begin{equation*}
\ptensor\otimes f^*\lambda\equiv f_*\ptensor\otimes\lambda
\end{equation*}
for any $r$-tree morphism $f: \tau\rightarrow\tau'$. Clearly, the
relation $\equiv$ is generated by the relations
$\ptensor\otimes(\gamma_{e_0})^*\lambda\equiv(\gamma_{e_0})_*\ptensor\otimes\lambda$
associated to edge contractions $\gamma_{e_0}:
\tau\rightarrow\tau/e_0$ since any morphism of $\Theta''(r)$ is a
composite of edge contractions.

Recall that the morphism $(\gamma_{e_0})^*:
\Cube{\tau/e_0}\rightarrow\Cube{\tau}$ identifies $\Cube{\tau/e_0}$
with the submodule of length tensors $\lambda = \bigotimes_{e\in
E'(\tau)} \lambda_e$ such that $\lambda_{e_0} = \x{0}$. Thus,
intuitively, the relation $\equiv$ identifies an element
$\ptensor\otimes\lambda\in\tau(\P)\otimes\Cube{\tau}$ equipped with
an edge of length $\lambda_{e_0} = \x{0}$ to an element
$\gamma_{e_0}(\ptensor)\otimes\lambda/e_0\in\tau/e_0(\P)\otimes\Cube{\tau/e_0}$
obtained by the contraction of~$e_0$ as represented in
figure~\ref{figure:LabelledMetricTreeContraction}. Accordingly, in
the context of dg-modules, any element of $W(\P)(r)$ has a normal
form $\ptensor\otimes\lambda\in\tau(\P)\otimes\Cube{\tau}$ for a
uniquely determined tree $\tau$ such that $\lambda_e = \x{1}$ or
$\lambda_e = \x{01}$ for all internal edges $e\in E'(\tau)$.

Finally, we observe that $W(\P)(r)$ is a dg-coalgebra if $\P$ is a
Hopf operad. Indeed we observe in this section that the functor
$\tau\mapsto\Cube{\tau}$ targets to the category of coalgebras and
not only to the category of dg-modules as well as the functor
$\tau\mapsto\tau(\P)$ if $\P$ is a Hopf operad. Our claim follows
since the forgetful functors from the category of coalgebras to the
category of dg-modules creates coends. These observations prove also
that, for a Hopf operad, the general definition of $W(\P)(r)$
deduced from monoidal model category axioms agrees with our
elementary construction of~$W(\P)$ within the category of
dg-modules. Explicitly, in order to define the diagonal of an
element $\ptensor\otimes\lambda$ in $W(\P)$, we can simply put
together the diagonal of $\ptensor$, represented in
figure~\ref{figure:LabelledTreeDiagonal}, with the diagonal of
$\lambda$ in the cubical complex $\Cube{\tau}$, as in the
representation of figure~\ref{figure:MetricTreeDiagonal}.

\subsubsection{The cellular differential in the $W$-construction}
\begin{figure}[b]
\begin{multline*}
\vcenter{\xymatrix@M=0pt@R=5mm@C=3mm{*+<1mm>{1}\ar[dr] & *+<1mm>{2}\ar[d] & & *+<1mm>{3}\ar[d] & *+<1mm>{4}\ar[dl] \\
& *+<1mm>{p_{v_2}}\ar[dr]|*+<2pt>{\scriptstyle\x{01}} & & *+<1mm>{p_{v_3}}\ar[dl]|*+<2pt>{\scriptstyle\x{01}} & \\
& & *+<1mm>{p_{v_1}}\ar[d] & & \\
& & *+<1mm>{0} & & \\ }} \xrightarrow{\partial} \pm
\vcenter{\xymatrix@M=0pt@R=5mm@C=3mm{*+<1mm>{1}\ar[dr] & *+<1mm>{2}\ar[d] & & *+<1mm>{3}\ar[d] & *+<1mm>{4}\ar[dl] \\
& *+<1mm>{p_{v_2}}\ar[dr]|*+<2pt>{\scriptstyle\x{1}} & & *+<1mm>{p_{v_3}}\ar[dl]|*+<2pt>{\scriptstyle\x{01}} & \\
& & *+<1mm>{p_{v_1}}\ar[d] & & \\
& & *+<1mm>{0} & & \\ }} \pm
\vcenter{\xymatrix@M=0pt@R=5mm@C=3mm{*+<1mm>{1}\ar[dr] & *+<1mm>{2}\ar[d] & & *+<1mm>{3}\ar[d] & *+<1mm>{4}\ar[dl] \\
& *+<1mm>{p_{v_2}}\ar[dr]|*+<2pt>{\scriptstyle\x{0}} & & *+<1mm>{p_{v_3}}\ar[dl]|*+<2pt>{\scriptstyle\x{01}} & \\
& & *+<1mm>{p_{v_1}}\ar[d] & & \\
& & *+<1mm>{0} & & \\ }}
\\
\shoveright{\pm
\vcenter{\xymatrix@M=0pt@R=5mm@C=3mm{*+<1mm>{1}\ar[dr] & *+<1mm>{2}\ar[d] & & *+<1mm>{3}\ar[d] & *+<1mm>{4}\ar[dl] \\
& *+<1mm>{p_{v_2}}\ar[dr]|*+<2pt>{\scriptstyle\x{01}} & & *+<1mm>{p_{v_3}}\ar[dl]|*+<2pt>{\scriptstyle\x{1}} & \\
& & *+<1mm>{p_{v_1}}\ar[d] & & \\
& & *+<1mm>{0} & & \\ }} \pm
\vcenter{\xymatrix@M=0pt@R=5mm@C=3mm{*+<1mm>{1}\ar[dr] & *+<1mm>{2}\ar[d] & & *+<1mm>{3}\ar[d] & *+<1mm>{4}\ar[dl] \\
& *+<1mm>{p_{v_2}}\ar[dr]|*+<2pt>{\scriptstyle\x{01}} & & *+<1mm>{p_{v_3}}\ar[dl]|*+<2pt>{\scriptstyle\x{0}} & \\
& & *+<1mm>{p_{v_1}}\ar[d] & & \\
& & *+<1mm>{0} & & \\ }}
}\\
\equiv\pm
\vcenter{\xymatrix@M=0pt@R=5mm@C=3mm{*+<1mm>{1}\ar[dr] & *+<1mm>{2}\ar[d] & & *+<1mm>{3}\ar[d] & *+<1mm>{4}\ar[dl] \\
& *+<1mm>{p_{v_2}}\ar[dr]|*+<2pt>{\scriptstyle\x{1}} & & *+<1mm>{p_{v_3}}\ar[dl]|*+<2pt>{\scriptstyle\x{01}} & \\
& & *+<1mm>{p_{v_1}}\ar[d] & & \\
& & *+<1mm>{0} & & \\ }} \pm
\vcenter{\xymatrix@M=0pt@R=5mm@C=3mm{& & & *+<1mm>{3}\ar[d] & *+<1mm>{4}\ar[dl] \\
*+<1mm>{1}\ar[drr] & *+<1mm>{2}\ar[dr] & & *+<1mm>{p_{v_3}}\ar[dl]|*+<2pt>{\scriptstyle\x{01}} & \\
& & *+<1mm>{p_{v_1}\circ_{v_2} p_{v_2}}\ar[d] & & \\
& & *+<1mm>{0} & & \\ }}
\\
\pm
\vcenter{\xymatrix@M=0pt@R=5mm@C=3mm{*+<1mm>{1}\ar[dr] & *+<1mm>{2}\ar[d] & & *+<1mm>{3}\ar[d] & *+<1mm>{4}\ar[dl] \\
& *+<1mm>{p_{v_2}}\ar[dr]|*+<2pt>{\scriptstyle\x{01}} & & *+<1mm>{p_{v_3}}\ar[dl]|*+<2pt>{\scriptstyle\x{1}} & \\
& & *+<1mm>{p_{v_1}}\ar[d] & & \\
& & *+<1mm>{0} & & \\ }} \pm
\vcenter{\xymatrix@M=0pt@R=5mm@C=3mm{*+<1mm>{1}\ar[dr] & *+<1mm>{2}\ar[d] & & & \\
& *+<1mm>{p_{v_2}}\ar[dr]|*+<2pt>{\scriptstyle\x{01}} & & *+<1mm>{3}\ar[dl] &  *+<1mm>{4}\ar[dll] \\
& & *+<1mm>{p_{v_1}\circ_{v_3} p_{v_3}}\ar[d] & & \\
& & *+<1mm>{0} & & \\ }}
\end{multline*}
\caption{The cellular differential of the
$W$-construction}\label{figure:CellBoundaryMetricTree}\end{figure}
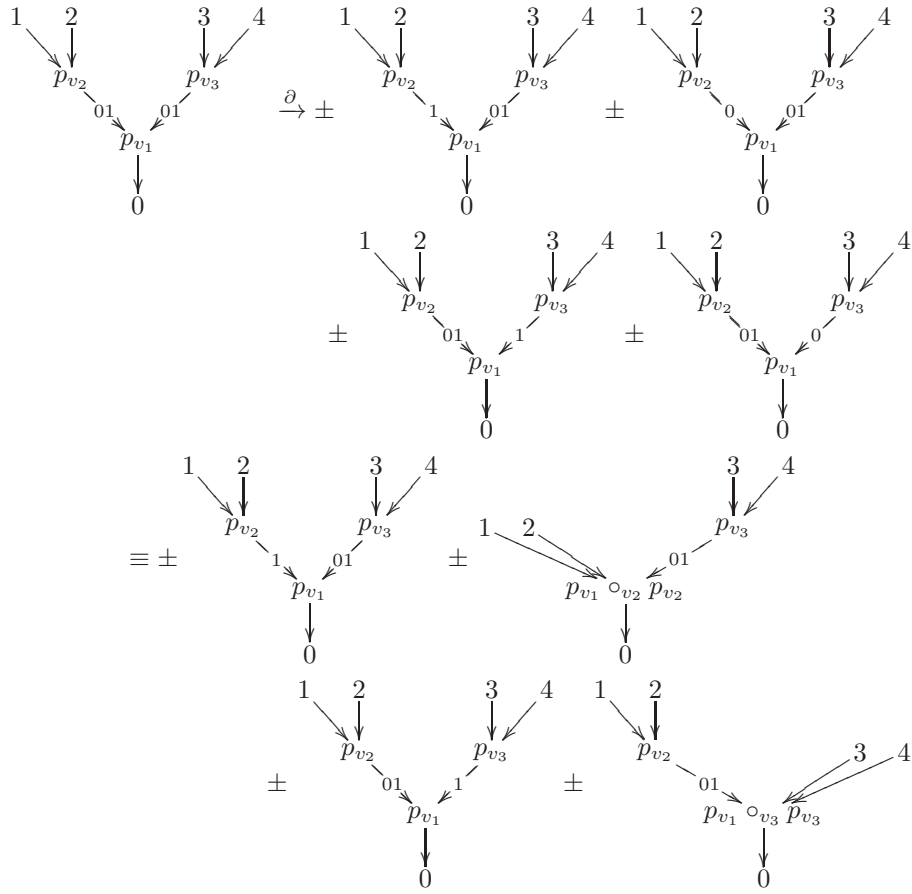

We make explicit the differential of an element of $W(\P)(r)$. By
definition, the module $W(\P)(r)$ is equipped with a natural
differential $\delta: W(\P)(r)\rightarrow W(\P)(r)$ induced by the
internal differential of the operad $\P$ and with a cellular
differential $\partial: W(\P)(r)\rightarrow W(\P)(r)$ induced by the
differential of the chain interval $\I$. Explicitly, this cellular
differential replaces an edge length $\lambda_e = \x{01}$ by a
difference $\partial(\lambda_e) = \x{1} - \x{0}$. The cellular
differential of an explicit element of $W(\P)(r)$ is represented in
figure~\ref{figure:CellBoundaryMetricTree}. We have unspecified
signs determined by an orientation of the cell associated to the
tree $\tau$, equivalent to an ordering of the length tensors
$\bigotimes_{e\in E'(\tau)} \lambda_e$, since the differential
$\partial$ is supposed to operate on the module $\Cube{\tau} =
\bigotimes_{e\in E'(\tau)}\I$ according to the rules of differential
graded algebra. Notice that, according to these definitions, the
dg-module
$(W(\P)(r),\partial)$
forms a cubical complex whose cells are indexed by $1$-reduced
$r$-trees.

\subsubsection{The operad structure of the $W$-construction}\label{item:MetricTreeOperadicComposite}
Observe that $W(\P)(r)$ comes equipped with a natural action of the
symmetric group $\Sigma_r$ given by the leaf reindexing of $r$-trees
as usual. In this paragraph we define the operadic composition
products
\begin{equation*}
\circ_i: W(\P)(s)\otimes W(\P)(t)\rightarrow W(\P)(s+t-1)
\end{equation*}
for $s,t\geq 2$ and this definition gives the structure of the
reduced operad $\overline{W(\P)}$.

The partial composition product $\circ_i$ is deduced from the
grafting process of rooted trees extended to cell metric trees. To
begin with we recall the definition of the tree $\sigma\circ_i\tau$
obtained by grafting the root of a $t$-tree $\tau$ to the $i$th leaf
of an $s$-tree $\sigma$. Formally, the vertex set of
$\sigma\circ_i\tau$ is defined by the union $V(\sigma\circ_i\tau) =
V(\sigma)\amalg V(\tau)$ of the vertices of $\sigma$ and $\tau$ and
the edge set by the union $E(\sigma\circ_i\tau) = E(\sigma)\amalg
E(\tau)/\equiv$ of the edges of $\sigma$ and $\tau$ in which the
root $e_0$ of $\tau$ and the $i$th leaf $e_i$ of $\sigma$ are
identified. For this new edge $\{e_0\equiv e_i\}$, we set
$s(\{e_0\equiv e_i\}) = s(e_0)$ and $t(\{e_0\equiv e_i\}) = t(e_i)$.
The other edges $e\in E(\sigma)\setminus\{e_i\}$ and $f\in
E(\tau)\setminus\{e_0\}$ keep the same source and target in
$\sigma\circ_i\tau$. The leaves of $\sigma\circ_i\tau$ are reindexed
as usual so that our construction produces an $s+t-1$-tree
$\sigma\circ_i\tau$ from an $s$-tree~$\sigma$ and a $t$-tree~$\tau$.

Observe that the internal edges of $\sigma\circ_i\tau$ are either
internal edges of $\sigma$ and $\tau$ or produced by the grafting
$\{e_0\equiv e_i\}$ of the root $e_0$ of $\tau$ with the $i$th leaf
$e_i$ of $\sigma$. Consequently, we have a natural morphism
\begin{equation*}
\circ_i:
\Cube{\sigma}\otimes\Cube{\tau}\rightarrow\Cube{\sigma\circ_i\tau}
\end{equation*}
that assigns the length $\lambda_{\{e_0\equiv e_i\}} = \x{1}$ to the
new internal edge and preserves the length of the other edges of
$\sigma\circ_i\tau$. Formally, we have a partition
$E'(\sigma\circ_i\tau) = E'(\sigma)\amalg E'(\tau)\amalg\{e_0\equiv
e_i\}$ and the morphism $\circ_i$ is defined by
\begin{equation*}
\Bigl[\bigotimes_{e\in
E'(\sigma)}\I\Bigr]\otimes\Bigl[\bigotimes_{e\in E'(\tau)}\I\Bigr]
\xrightarrow{\Id\otimes\eta^1} \Bigl[\bigotimes_{e\in
E'(\sigma)}\I\Bigr]\otimes\Bigl[\bigotimes_{e\in
E'(\tau)}\I\Bigr]\otimes\F\x{1} \hookrightarrow\bigotimes_{e\in
E'(\sigma\circ_i\tau)}\I.
\end{equation*}
One checks readily that this definition gives a morphism of functors
in $\tau\in\Theta''(r)$.

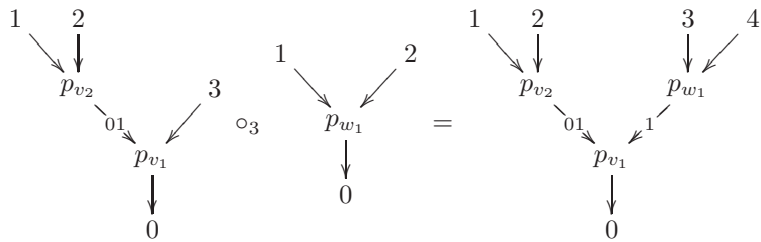
\begin{figure}[b]
\[\vcenter{\xymatrix@M=0pt@R=5mm@C=3mm{
*+<2mm>{1}\ar[dr] & *+<2mm>{2}\ar[d] & & \\
& *+<2mm>{p_{v_2}}\ar[dr]|*+<2pt>{\scriptstyle\x{01}} & & *+<2mm>{3}\ar[dl] \\
& & *+<2mm>{p_{v_1}}\ar[d] & \\
& & *+<2mm>{0} & \\ }}
\circ_{3}\vcenter{\xymatrix@M=0pt@R=5mm@C=3mm{
*+<2mm>{1}\ar[dr] & & *+<2mm>{2}\ar[dl] \\
& *+<2mm>{p_{w_1}}\ar[d] & \\
& *+<2mm>{0} & \\ }} =\vcenter{\xymatrix@M=0pt@R=5mm@C=3mm{
*+<2mm>{1}\ar[dr] & *+<2mm>{2}\ar[d] & & *+<2mm>{3}\ar[d] & *+<2mm>{4}\ar[dl] \\
& *+<2mm>{p_{v_2}}\ar[dr]|*+<2pt>{\scriptstyle\x{01}} & & *+<2mm>{p_{w_1}}\ar[dl]|*+<2pt>{\scriptstyle\x{1}} & \\
& & *+<2mm>{p_{v_1}}\ar[d] & & \\
& & *+<2mm>{0} & & \\ }}\] \caption{The operadic composition product
of the
$W$-construction}\label{figure:GraftingMetricTree}\end{figure}

On the other hand, we have also a natural isomorphism
\begin{equation*}
\circ_i:
\sigma(\P)\otimes\tau(\P)\xrightarrow{\simeq}\sigma\circ_i\tau(\P)
\end{equation*}
since $V(\sigma\circ_i\tau) = V(\sigma)\amalg V(\tau)$. The pairing
of these morphisms gives the required composition product on
$W(P)(r) = \int_{\tau} \tau(\P)\otimes\Cube{\tau}$. The
representation of this grafting process is given by
figure~\ref{figure:GraftingMetricTree}.

One checks readily that the maps $\circ_i$ satisfy the associativity
and commutativity properties of partial composition products of an
operad. Clearly, the only example of a reduced $1$-tree is provided
by the unit tree $\downarrow$ equipped with an empty set of vertices
and only one edge. This object defines a unit for the partial
composition products defined in this paragraph. Accordingly, for $r
= 1$ we can extend the definition of~\ref{item:WComplex} in order to
obtain $W(\P)(1) = \F$, coherently with our conventions.

Observe that the definition of $\circ_i:
\sigma(\P)\otimes\tau(\P)\xrightarrow{\simeq}\sigma\circ_i\tau(\P)$
can be deduced from the axioms of symmetric monoidal categories so
that this morphism defines a morphism of coalgebras if $\P$ is a
Hopf operad. Similarly, the morphism $\circ_i:
\Cube{\sigma}\otimes\Cube{\tau}\rightarrow\Cube{\sigma\circ_i\tau}$
defines a morphism of coalgebras and not only of dg-modules because
$\I$ defines an interval in the category of coalgebras and not only
in the category of dg-modules. From these observations we conclude
that the operadic composition products of the $W$-construction are
morphism of coalgebras so that $W(\P)$ forms a Hopf operad. To
conclude, for a Hopf operad $\P$, our construction returns an
operad~$W(\P)$ in the category of dg-coalgebras. In fact, our
construction agrees with the general construction
of~\cite{BergerMoerdijkW}, defined within the formalism of monoidal
model categories, for the category of dg-coalgebras.

\subsubsection{Unitary operations in the Boardman-Vogt construction}\label{item:WUnitComposites}
We define now the partial composites with unital operations
\begin{equation*}
\circ_i: W(\P)(r)\otimes W(\P)(0)\rightarrow W(\P)(r-1).
\end{equation*}
As specified in the introduction, we set $W(\P)(0) = \P(0) = \F$ so
that $W(\P)$ is still a connected unital operad. For our purpose we
consider an expansion
\begin{equation*}
W(\P)(r) = \bigoplus_{\tau} \tau(\P)\otimes\Cube{\tau}/\equiv
\end{equation*}
that ranges over all $r$-trees, unlike the expansion
of~\ref{item:WComplex}, but which contains more relations so that
this expansion returns the same result for $r\geq 1$ and $W(\P)(0) =
\F$ for $r = 0$. Roughly, we put relations so that any labeling
reduces to the labeling of a reduced tree for $r\geq 1$ or to an
element of $\P(0)$ for $r=0$.

Explicitly, let $\tau$ denote a tree equipped with a vertex $v_0$
such that $I_{v_0}$ is empty and let $e_0$ denote the unique edge
such that $s(e_0) = v_0$. Thus, for a labeling $\ptensor =
\bigotimes_{v\in V(\tau)} p_v\in\tau(\P)$, we have
$p_{v_0}\in\P(0)$. An element
$\ptensor\otimes\lambda\in\tau(\P)\otimes\Cube{\tau}$ is set to be
equivalent to an element of the
summand~$\tau/e_0(\P)\otimes\Cube{\tau/e_0}$ for the tree $\tau/e_0$
obtained by contraction of the edge $e_0$ according to the process
of~\ref{item:TreeMorphisms}. By definition, the terminal vertex
$v_0$ disappears in $\tau/e_0$. Explicitly, we consider the morphism
$(\gamma_{e_0})_*: \tau(\P)\rightarrow\tau/e_0(\P)$ defined
in~\ref{item:WComplex} and we let $(\gamma_{e_0})_*:
\Cube{\tau}\rightarrow\Cube{\tau/e_0}$ denote the dg-module morphism
induced by the augmentation $\epsilon: \I\rightarrow\F$ on the
factor of the tensor product $\Cube{\tau} = \bigotimes_{e\in
E'(\tau)}\I$ indexed by the edge $e_0\in E'(\tau)$. Recall that
$E'(\tau/e_0) = E'(\tau)\setminus\{e_0\}$ so that we have a morphism
\begin{equation*}
\Bigl[\bigotimes_{e\in E'(\tau)\setminus\{e_0\}}\I\Bigr]\otimes\I
\xrightarrow{\Id\otimes\epsilon}\Bigl[\bigotimes_{e\in
E'(\tau)\setminus\{e_0\}}\I\Bigr] \xrightarrow{=}\bigotimes_{e\in
E'(\tau/e_0)}\I.
\end{equation*}
For a tensor $\ptensor\otimes\lambda\in\tau(\P)\otimes\Cube{\tau}$,
we set precisely
\begin{equation*}
\ptensor\otimes\lambda\equiv(\gamma_{e_0})_*(\ptensor)\otimes(\gamma_{e_0})_*(\lambda).
\end{equation*}

Let $\tau$ denote a tree equipped with a vertex $v_0$ such that
$I_{v_0}$ has only one element and let $e_0$, respectively $e_1$,
denote the unique edge such that $s(e_0) = v_0$, respectively
$t(e_1) = v_0$. For a labeling $\ptensor = \bigotimes_{v\in V(\tau)}
p_v$, the element $p_{v_0}\in\P(I_{v_0})$ denotes necessarily the
operad unit $1\in\P(1)$. We consider the tree $\tau\backslash v_0$
in which the vertex $v_0$ is removed. Formally, this tree is defined
by the vertex set $V(\tau\backslash v_0) = V(\tau)\setminus\{v_0\}$
and by the quotient edge set $E(\tau\backslash v_0) =
E(\tau)/\{e_0\equiv e_1\}$ in which the edges $e_0$ and $e_1$ are
identified. An element
$\ptensor\otimes\lambda\in\tau(\P)\otimes\Cube{\tau}$ is set to be
equivalent to an element of the summand~$\tau\backslash
v_0(\P)\otimes\Cube{\tau\backslash v_0}$. Namely we set
\begin{equation*}
\ptensor\otimes\lambda\equiv(\delta_{v_0})_*(\ptensor)\otimes(\delta_{v_0})_*(\lambda)
\end{equation*}
for morphisms $(\delta_{v_0})_*: \tau(\P)\rightarrow\tau\backslash
v_0(\P)$ and $(\delta_{v_0})_*:
\Cube{\tau}\rightarrow\Cube{\tau\backslash v_0}$ obtained as
follows. Since $\P(1) = \F$, we have an isomorphism
\begin{equation*}
\bigotimes_{v\in V(\tau)}
\P(I_v)\xrightarrow{\simeq}\bigotimes_{v\in
V(\tau)\setminus\{v_0\}}\P(I_v)
\end{equation*}
which yields the required morphism $(\delta_{v_0})_*:
\tau(\P)\rightarrow\tau\backslash v_0(\P)$. Equivalently, the tensor
$(\delta_{v_0})_*(\ptensor)$ is obtained simply by removing the unit
element $p_{v_0} = 1$ from $\ptensor = \bigotimes_{v\in V(\tau)}
p_v$. The length tensor $(\delta_{v_0})_*(\lambda)$ is defined by
$(\delta_{v_0})_*(\lambda)_e = \mu(\lambda_{e_0},\lambda_{e_1})$ for
$e = \{e_0\equiv e_1\}$ and $(\delta_{v_0})_*(\lambda)_e =
\lambda_e$ for the other edges of $\tau\backslash v_0$. To be
precise, the edges $e = e_0,e_1$ are not necessarily internal. In
this case we set by convention $\lambda_e = \x{1}$.

\subsubsection{The homotopy type of the $W$-construction -- assertion (a) of theorem~\ref{thm:WHopfOperad}}\label{item:WcomplexHomotopy}
In this paragraph we define the operad equivalence $\epsilon:
W(\P)\rightarrow\P$ specified in the claim (a) of
theorem~\ref{thm:WHopfOperad}. Recall that the module of length
tensors $\Cube{\tau}$ is equipped with natural equivalences
\begin{equation*}
\F\xrightarrow{\eta^0_*}\Cube{\tau}\xrightarrow{\epsilon_*}\F
\end{equation*}
such that $\epsilon_*\eta^0_* = \Id$ and $\eta^0_*$ is an acyclic
Reedy cofibration. Furthermore, by lemma~\ref{lemm:TreewiseCoends},
these morphisms yield inverse equivalences of dg-modules
(respectively dg-coalgebras in the Hopf operad context):
\begin{equation*}
\xymatrix{
*{\displaystyle\tau_r(\P)\simeq\int_{\tau\in\Theta''(r)}\tau(\P)\otimes\F}\ar@<-2pt>[r]_(0.55){\eta^0_*}
&
*{\displaystyle\int_{\tau\in\Theta''(r)}\tau(\P)\otimes\Cube{\tau}}\ar@<-2pt>[l]_(0.45){\epsilon_*}
}.
\end{equation*}
Recall that $\tau_r(\P)\simeq\P(r)$. Consequently, we obtain
morphisms of dg-modules, respectively dg-coalgebras, $\eta:
\P\rightarrow W(\P)$ and $\epsilon: W(\P)\rightarrow\P$ such that
$\epsilon\eta = \Id$ and which are endowed with the following
properties:

\begin{claim}
The morphism $\eta: \P\rightarrow W(\P)$ defines an acyclic
cofibration of dg-modules, respectively dg-coalgebras. As a
consequence, the morphism $\epsilon: W(\P)\rightarrow\P$, which is
right inverse to $\eta$, defines a weak equivalence of dg-modules,
respectively dg-coalgebras, as well.\qed
\end{claim}

We check that $\epsilon: W(\P)\rightarrow\P$ defines indeed an
operad equivalence. In fact, one can obtain a more elementary
definition for $\eta: \P\rightarrow W(\P)$ and $\epsilon:
W(\P)\rightarrow\P$. Namely the morphism $\eta: \P\rightarrow W(\P)$
identifies an operation $p\in\P(r)$ with an element of $\tau_r(\P)$,
the summand of $W(\P)(r) = \bigoplus_{\tau\in\Theta''(r)}
\tau(\P)\otimes\Cube{\tau}/\equiv$ associated to the terminal
$r$-tree $\tau = \tau_r$. In the converse direction, for a tensor
\begin{equation*}
\ptensor\otimes\lambda = \Bigl[\bigotimes_{v\in V(\tau)} p_v\Bigr]
\otimes\Bigl[\bigotimes_{e\in E'(\tau)}\lambda_e\Bigr]
\in\tau(\P)\otimes\Cube{\tau},
\end{equation*}
we have $\epsilon(\ptensor\otimes\lambda) = 0$ if $\lambda_e =
\x{01}$ for an edge $e\in E'(\tau)$. Otherwise, if $\lambda_e =
\x{0}$ or $\x{1}$ for all edges $e\in E'(\tau)$, then the morphism
$\epsilon$ performs the composite of the operations $p_v$ in $\P$.
According to this elementary definition, we obtain:

\begin{obsv}
The morphism $\epsilon: W(\P)\rightarrow\P$ defines a morphism of
unital (Hopf) operads.\qed
\end{obsv}

and finally:

\begin{claim}[Assertion (a) of theorem~\ref{thm:WHopfOperad}]\label{claim:WAugmentation}
The morphism $\epsilon: W(\P)\rightarrow\P$ defines a
weak-equivalence of unital (Hopf) operads.\qed
\end{claim}

Notice that $\eta: \P\rightarrow W(\P)$ does not form an operad
morphism and give only a left inverse of $\epsilon:
W(\P)\rightarrow\P$ in the category of dg-modules, respectively
dg-coalgebras if $\P$ is a Hopf operad.

\subsubsection{The quasi-free operad property -- assertion (b) of theorem~\ref{thm:WHopfOperad}}\label{item:WQuasiFree}
In the differential context, we have a natural notion of a
quasi-free object in the category of unital operads. Explicitly, a
quasi-free unital operad is specified by a pair $\Q =
(\FOp_*(M),\partial)$ and represents a free operad $\FOp_*(M)$
equipped with a non-canonical differential obtained by the addition
of a homogeneous derivation of degree $-1$
\begin{equation*}
\partial: \FOp_*(M)\rightarrow\FOp_*(M)
\end{equation*}
to the natural differential of the free operad $\delta:
\FOp_*(M)\rightarrow\FOp_*(M)$ induced by the internal differential
of~$M$.

In this paragraph we define a graded $\Lambda_*$-module $W'(\P)$
such that $W(\P) = (\FOp_*(W'(\P)),\partial)$ as claimed by
assertion (b) of theorem~\ref{thm:WHopfOperad}. In fact, this
$\Lambda_*$-module $W'(\P)$ is defined by a section of the
indecomposable quotient of $W(\P)$. The construction of this
paragraph makes sense only in the differential graded context and
not in the general framework of~\cite{BergerMoerdijkW} in which a
Boardman-Vogt construction can be defined.

By definition, an element
$\ptensor\otimes\lambda\in\tau(\P)\otimes\Cube{\tau}$ is
decomposable for the operadic composition product of $W(\P)$ if and
only if the length tensor $\lambda$ contains edges of length
$\lambda_e = 1$. Therefore we consider the graded module
$\IntCube{\tau}\subset\Cube{\tau}$ spanned by length tensors
$\lambda = \bigotimes_e \lambda_e$ such that $\lambda_e = \x{0}$ or
$\lambda_e = \x{01}$ for all $e\in E'(\tau)$. This module is not
preserved by the differential of $\Cube{\tau}$ but defines a section
of the quotient dg-module of $\Cube{\tau}$ in which length tensors
$\lambda$ that contain edges of length $\lambda_e = 1$ are canceled.
By an abuse of notation, we identify the section
$\IntCube{\tau}\subset\Cube{\tau}$ with this quotient dg-module so
that we have a natural morphism of dg-modules
$\Cube{\tau}\rightarrow\IntCube{\tau}$.

Then let
\begin{equation*}
W'(\P)(r) = \int_{\tau\in\Theta''(r)}\tau(\P)\otimes\IntCube{\tau}.
\end{equation*}
Clearly, the quotient morphism
$\Cube{\tau}\rightarrow\IntCube{\tau}$ induces a morphism of
dg-modules $W(\P)\rightarrow W'(\P)$ that identifies the module
$W'(\P)$ with the indecomposable quotient of $W(\P)$. On the other
hand, the embedding $\IntCube{\tau}\hookrightarrow\Cube{\tau}$
induces an embedding of graded modules $W'(\P)\hookrightarrow W(\P)$
so that $W'(\P)$ represents a section of the indecomposable quotient
of the operad $W(\P)$. One checks easily that $W'(\P)$ is preserved
by the operadic composites $\partial_i = -\circ_i *$ supplied by the
construction of~\ref{item:WUnitComposites}. Consequently, the module
$W'(\P)$ forms a graded $\Lambda_*$-submodule of $W(\P)$.

\begin{claim}[Assertion (b) of theorem~\ref{thm:WHopfOperad}]\label{claim:WQuasiFree}
The embedding of graded $\Lambda_*$-modules $W'(\P)\hookrightarrow
W(\P)$ induces an isomorphism of graded operads
$\FOp_*(W'(\P))\simeq W(\P)$. Consequently, the operad $W(\P)$ is
quasi-free and we have
\begin{equation*}
W(\P) = (\FOp_*(W'(\P)),\partial)
\end{equation*}
for an operad derivation $\partial:
\FOp_*(W'(\P))\rightarrow\FOp_*(W'(\P))$.
\end{claim}

\begin{proof}
In fact, for an element
$\ptensor\otimes\lambda\in\tau(\P)\otimes\Cube{\tau}$, the edges of
length $\lambda_e = 1$ form the edges of a tree $\sigma$ whose
vertices are the subtrees of $\tau$ cut by these edges precisely.
This process identifies the tensor
$\ptensor\otimes\lambda\in\tau(\P)\otimes\Cube{\tau}$ with an
operadic composite of elements of $W'(\P)$ arranged on a tree
$\sigma$ and hence with an element of the free operad generated by
$W'(\P)$ according to the classical construction of the free operad.
Therefore we have an isomorphism
$\FOp_*(W'(\P))\xrightarrow{\simeq}W(\P)$ as stated.
\end{proof}

To be more precise, notice that the interval differential $\partial$
admits a splitting $\partial = \partial_0+\partial_1$ such that
$\partial_0(\x{01}) = \x{0}$ and $\partial_1(\x{01}) = \x{1}$.
Clearly, the component $\partial_0$ of the interval differential
preserves the module $\IntCube{\tau}\subset\Cube{\tau}$ and defines
the differential of the associated quotient module. As a
consequence, the cellular differential of the $W$-construction is
endowed with a similar splitting $\partial = \partial_0+\partial_1$,
the component $\partial_0$ preserves $W'(\P)$ and $\partial_1$ is
decomposable. Hence the indecomposable quotient of $W(\P)$ can be
identified with the graded module $W'(\P)$ equipped with a cellular
differential defined by $\partial_0$. In fact, the suspended
dg-module $(\Sigma W'(\P),\partial_0)$ can be identified with
$B(\P)$, the operadic bar construction introduced
in~\cite{GetzlerJones} (see also
\cite{OperadTextbook,GinzburgKapranov}), and $\overline{W(\P)} =
(\FOp(W'(\P)),\partial)$ can be identified with $B^c(B(\P))$, the
cobar-bar construction of~$\P$. In the topological context the
suspension of the $\Sigma_*$-space $W'(\P)$ defined by the operadic
indecomposable quotient of $W(\P)$ is homeomorphic to the bar
construction defined in~\cite{Ching} but the relation with
$B^c(B(\P))$ fails in this context.

Clearly, the module $W'(\P)$ is preserved by the morphism $\phi_*:
W(\P)\rightarrow W(\P')$ induced by an operad morphism $\phi:
\P\rightarrow\P'$. Hence the map $\P\mapsto W(\P)$ has the
functoriality property claimed by assertion~(b) of
theorem~\ref{thm:WHopfOperad}. As explained in the introduction, we
deduce the cofibrancy claim from
theorem~\ref{thm:BoardmanVogtOperadicCellularDecomposition}.

\subsection{The cellular decomposition of the $W$-construction}\label{subsection:BoardmanVogtOperadicCellularDecomposition}
In this subsection we define the operadic cellular decomposition of
the operad $W(\P)$ announced by
theorem~\ref{thm:BoardmanVogtOperadicCellularDecomposition}. Namely
we specify suboperads $W^d(\P)\subset W(\P)$ such that $W(\P) =
\colim_d W^d(\P)$ together with $\Lambda_*$-submodules
$C^d(\P)\subset W^{d-1}(\P)$ and $D^d(\P)\subset W^{d}(\P)$ so that
we have an operad pushout
\begin{equation*}
\xymatrix{ \FOp_*(C^d(\P))\ar[r]^{f^d}\ar[d]^{i^d} & W^{d-1}(\P)\ar[d]^{j_d} \\
\FOp_*(D^d(\P))\ar[r]^{g^d} & W^{d}(\P)\push{5} }.
\end{equation*}
For this purpose we consider the quasi-free representation
of~$W(\P)$ defined in~\ref{item:WQuasiFree}. Recall also that the
category of trees is equipped with a grading defined by the number
of internal edges and we let $\Theta''_d(r)$ denote the full
subcategory of~$\Theta''(r)$ formed by trees $\tau\in\Theta''(r)$
such that $\gr(\tau)\leq d$. We consider the graded module
\begin{equation*}
W'{}^d(\P)(r) = \int_{\tau\in\Theta''_d(r)}
\tau(\P)\otimes\IntCube{\tau}
\end{equation*}
that forms clearly a $\Lambda_*$-submodule of~$W'(\P)$. The simplest
definition of~$W^d(\P)$ arises from the following straightforward
observation:

\begin{obsv}
The free unital operad $\FOp_*(W'{}^d(\P))\subset\FOp_*(W'(\P))$ is
preserved by the differential of $W(\P) =
(\FOp_*(W'(\P)),\partial)$. Consequently, the pairs $W^d(\P) =
(\FOp_*(W'{}^d(\P)),\partial)$ define a nested sequence of
quasi-free unital suboperads of $W(\P)$ such that $W(\P) = \colim_d
W^d(\P)$.

If $\P$ is a Hopf operad, then
$\FOp_*(W'{}^d(\P))\subset\FOp_*(W'(\P))$ is preserved by the
diagonal of~$W(\P)$ as well so that $W^d(\P) =
(\FOp_*(W'{}^d(\P)),\partial)$ defines a Hopf suboperad of
$W(\P)$.\qed
\end{obsv}

\subsubsection{Skeletal filtration and decomposition of length tensors}\label{item:DecLengthTensors}
For our purposes we give another equivalent definition of~$W^d(\P)$.
This new definition can also be generalized to the framework
of~\cite{BergerMoerdijkW} unlike the previous construction
of~$W^d(\P)$.

First, we need to introduce a skeletal filtration of~$\Cube{\tau}$
which is defined by the submodules
$\sk_d\Cube{\tau}\subset\Cube{\tau}$ spanned by the length tensors
of degree $\deg(\lambda)\leq d$. In a general framework the module
$\sk_d\Cube{\tau}$ is defined by a colimit of tensor products of the
form
\begin{equation*}
\Bigl\{\bigotimes_e\F\x{0}\Bigr\}
\otimes\Bigl\{\bigotimes_e\F\x{1}\Bigr\}
\otimes\Bigl\{\bigotimes_e\I\Bigr\}\subset\bigotimes_{e\in
E'(\tau)}\I
\end{equation*}
with no more than $d$ factors $\I$.

Then, for an element
$\ptensor\otimes\lambda\in\tau(\P)\otimes\Cube{\tau}$, we observe
that the relation $\ptensor\otimes\lambda\in\FOp_*(W'{}^d(\P))(r)$
is equivalent to a condition on the length tensor only
$\lambda\in\Cube{\tau}$, like the characterization of a decomposable
element given in~\ref{item:WQuasiFree}. Explicitly, as explained in
the proof of claim~\ref{claim:WQuasiFree} the edges of length
$\lambda_e = \x{1}$ decomposes the length tensor
$\lambda\in\Cube{\tau}$ into a composite of elements
$\lambda_u\in\Cube{\tau_u}$ for subtrees $\tau_u\subset\tau$
arranged on a tree $\sigma$. By definition, we have
$\ptensor\otimes\lambda\in\FOp_*(W'{}^d(\P))(r)$ if and only if
these length tensors satisfy $\deg(\lambda_u)\leq d$ for all $u\in
V(\sigma)$.

To make our construction more general, we introduce a submodule
$\dec_d\Cube{\tau}\subset\Cube{\tau}$ spanned by decomposable
tensors such that $\deg(\lambda_u)\leq d$. Formally, this dg-module
$\dec_d\Cube{\tau}$ can be defined by a colimit process. For each
decomposition of a tree $\tau$ into subtrees $\tau_u$ as above, we
consider the tensor product
\begin{equation*}
\bigotimes_{u\in V(\sigma)}\sk_d\Cube{\tau_u}.
\end{equation*}
By definition, we have $E'(\tau) =
E'(\sigma)\amalg\Bigl\{\coprod_{u\in V(\sigma)} E'(\tau_u)\Bigr\}$.
Consequently, we have an embedding
\begin{equation*}
\eta^1_*: \bigotimes_{u\in
V(\sigma)}\sk_d\Cube{\tau_u}\rightarrow\bigotimes_{e\in E'(\tau)}\I
= \Cube{\tau}
\end{equation*}
given by $\eta^1: \F\x{1}\rightarrow\I$ on the factors associated to
an edge $e\in E'(\sigma)$. Furthermore, for each decomposition
refinement, we have a similar morphism
\begin{equation*}
\bigotimes_{u'\in
V(\sigma')}\sk_d\Cube{\tau'_{u'}}\rightarrow\bigotimes_{u\in
V(\sigma)}\sk_d\Cube{\tau'_{u}},
\end{equation*}
that commutes with these embeddings, given by $\eta^1:
\F\x{1}\rightarrow\I$ on the factors associated to an edge
\begin{equation*}
e\in\coprod_{u} E'(\tau_u)\setminus\coprod_{u'} E'(\tau'_{u'}).
\end{equation*}
The module $\dec_d\Cube{\tau}$ is defined by the colimit of these
morphisms for all decompositions into subtrees $\tau_u$. This module
is endowed with a natural embedding
$\dec_d\Cube{\tau}\hookrightarrow\Cube{\tau}$ induced by the
embeddings above. The assumption $\lambda_u\in\sk_d\Cube{\tau_u}$
ensures that a length tensor that arises in this colimit satisfies
$\deg(\lambda_u)\leq d$. The colimit process permits to identify an
element $\lambda_u\in\Cube{\tau_u}$ that contains an edge of length
$(\lambda_u)_e = \x{1}$ to a further decomposition. Hence any
element of $\FOp_*(W'{}^d(\P))(r)$ is equivalent to a representative
$\ptensor\otimes\lambda\in\tau(\P)\otimes\Cube{\tau}$ such that
$\lambda\in\dec_d\Cube{\tau}$.

To conclude, we obtain:

\begin{obsv}\label{obsv:WCellularFiltration}
We have
\begin{equation*}
W^d(\P)(r) = \int_{\tau\in\Theta''(r)}
\tau(\P)\otimes\dec_d\Cube{\tau}.\qed
\end{equation*}
\end{obsv}

Furthermore, observe that we have an isomorphism
$\dec_d(\sigma\circ_i\tau)\simeq\dec_d(\sigma)\otimes\dec_d(\tau)$
for any pair of trees $\sigma\in\Theta''(s)$ and
$\tau\in\Theta''(t)$. Consequently, the composition product
of~$W^d(\P)$ can be obtained by a generalization of the construction
of~\ref{item:MetricTreeOperadicComposite} for the operad $W(\P)$.

\medskip
For our next construction, we record the following assertion which
is a tautological consequence of the definition of the modules
$\dec_d\Cube{\tau}$:

\begin{fact}\label{fact:SkeletonEmbeddings}
For a tree $\tau\in\Theta''(r)$, we have natural embeddings
\begin{equation*}
\xymatrix@R=6mm{\ar@{.}[d] & \ar@{.}[d] & \\
\ar@{^{(}->}[]!D+<0pt,-2pt>;[d] & \ar@{^{(}->}[]!D+<0pt,-2pt>;[d] & \\
\sk_{d-1}\Cube{\tau}\ar@{^{(}->}[]!D+<0pt,-2pt>;[d]\ar@{^{(}->}[]!R+<2pt,0pt>;[r]
&
\dec_{d-1}\Cube{\tau}\ar@{^{(}->}[]!DR+<2pt,-2pt>;[dr]\ar@{^{(}->}[]!D+<0pt,-2pt>;[d] & \\
\sk_d\Cube{\tau}\ar@{^{(}->}[]!R+<2pt,0pt>;[r]\ar@{^{(}->}[]!D+<0pt,-2pt>;[d]
&
\dec_d\Cube{\tau}\ar@{^{(}->}[]!R+<2pt,0pt>;[r]\ar@{^{(}->}[]!D+<0pt,-2pt>;[d] & \Cube{\tau} \\
\ar@{.}[d] & \ar@{.}[d] & \\ && }.
\end{equation*}
\end{fact}

\subsubsection{The cells of the $W$-construction}\label{item:WCells}
In order to define the $\Lambda_*$-modules $C^d(\P)$ and $D^d(\P)$,
we consider the skeletal filtration of the module of length tensors
$\Cube{\tau}$ introduced in the previous paragraph and the
subcategory of $\Theta''(r)$ formed by trees with no more than $d$
internal edges. We set precisely
\begin{equation*}
C^d(\P) = \int_{\tau\in\Theta''_d(r)}
\tau(\P)\otimes\sk_{d-1}\Cube{\tau} \quad\text{and} \quad D^d(\P) =
\int_{\tau\in\Theta''_d(r)} \tau(\P)\otimes\sk_{d}\Cube{\tau}.
\end{equation*}
We have canonical embeddings
\begin{equation*}
\xymatrix{ C^d(\P)\ar[d]_{i^d}\ar[r]^{f^d} &
W^{d-1}(\P)\ar@{^{(}->}[]!DR+<2pt,-2pt>;[dr]\ar@{^{(}->}[]!D+<0pt,-2pt>;[d] & \\
D^d(\P)\ar[r]_{g^d} & W^{d}(\P)\ar@{^{(}->}[]!R+<2pt,0pt>;[r] &
W(\P) }.
\end{equation*}
induced by the embeddings of fact~\ref{fact:SkeletonEmbeddings} and
by the category inclusion $\Theta''_d(r)\subset\Theta''(r)$. One
checks readily that $C^d(\P)$, respectively $D^d(\P)$, is preserved
by the operadic composites with unital operations $\partial_i =
-\circ_i *$ so that the diagram above defines a commutative diagram
of $\Lambda_*$-modules. Furthermore, if $\P$ is a Hopf operad, then
the dg-module $C^d(\P)$, respectively $D^d(\P)$, is equipped with a
natural coalgebra structure and we obtain a commutative diagram of
Hopf $\Lambda_*$-modules. The modules $C^d(\P)$ and $D^d(\P)$ are
also unitary and connected since for all $d\geq 0$ the category
$\Theta''_d(1)$ is reduced to the unit tree.

We consider the commutative square of (Hopf) operads
\begin{equation*}
\xymatrix{ \FOp_*(C^d(\P))\ar[d]_{i^d}\ar[r]^{f^d} & W^{d-1}(\P)\ar@{^{(}->}[]!D+<0pt,-2pt>;[d] \\
\FOp_*(D^d(\P))\ar[r]_{g^d} & W^{d}(\P) }
\end{equation*}
in which $f^d$ and $g^d$ are induced by our morphisms of (Hopf)
$\Lambda_*$-modules $f^d$ and $g^d$.

\begin{claim}[First part of theorem~\ref{thm:BoardmanVogtOperadicCellularDecomposition}]\label{claim:WCellExtensions}
This commutative square defines a pushout in the category of unital
(Hopf) operads.
\end{claim}

\begin{proof}
We use the relation $W^d(\P) = (\FOp_*(W'{}^d(\P)),\partial)$ in our
proof. This argument is valid only in the differential graded
context. Nevertheless one can observe that the claim holds in the
framework of~\cite{BergerMoerdijkW} though the proof becomes more
technical. Therefore we give only a few hints below for a general
proof of the claim.

Recall that the functor $\P\mapsto\overline{\P}$, from the category
of unital operads to the category of non-unital operads, creates
pushouts. Therefore we consider the reduced operads associated to
$W^d(\P)$ and the free non-unital operads $\FOp(C^d(\P))$ and
$\FOp(D^d(\P))$. Recall that the reduced operad of a quasi-free
unital operad is still a quasi-free object in the category of
non-unital operad. For the operad $W^d(\P)$ we obtain precisely
$\overline{W^d(\P)} = (\FOp(W'{}^d(\P)),\partial)$.

Observe that $W'{}^d(\P)(r) = W'{}^{d-1}(\P)(r)\oplus E^d(\P)(r)$,
where $E^d(\P)(r)$ denotes the submodule of~$W'{}^d(\P)(r) =
\int_{\tau\in\Theta''(r)} \tau(\P)\otimes\Cube{\tau}$ spanned by
tensors $\ptensor\otimes\lambda\in\tau(\P)\otimes\Cube{\tau}$ such
that $\gr(\tau) = d$ and $\lambda_e = \x{01}$ for all edges $e\in
E'(\tau)$. Equivalently, the module $E^d(\P)(r)$ consists of tensors
$\ptensor\otimes\lambda\in\tau(\P)\otimes\Cube{\tau}$ such that the
length tensor $\lambda$ has a top degree $\deg(\lambda) = d$.
Clearly, the modules $E^d(\P)(r)$ form a graded $\Sigma_*$-modules
and the splitting $W'{}^d(\P) = W'{}^{d-1}(\P)\oplus E^d(\P)(r)$
holds in the category of graded $\Sigma_*$-modules. Observe that we
have a similar splitting $D^d(\P)(r) = C^d(\P)(r)\oplus E^d(\P)(r)$
for the same $\Sigma_*$-module $E^d(\P)$. If we forget about the
differential of the quasi-free operads $W^d(\P)$, then these
observations imply that the diagram of free graded operads
\begin{equation*}
\xymatrix{ \FOp(C^d(\P))\ar[d]_{i^d}\ar[r]^{f^d} & \FOp(W'{}^{d-1}(\P))\ar[d] \\
\FOp(D^d(\P))\ar[r]_{g^d} & \FOp(W'{}^d(\P)) }
\end{equation*}
form a pushout. Our claim follows immediately from this assertion.
\end{proof}

For a proof of claim~\ref{claim:WCellExtensions} in the general
framework of~\cite{BergerMoerdijkW}, we use the coend representation
of $W^d(\P)$ supplied by observation~\ref{obsv:WCellularFiltration}.
One can observe further that the free operad $\FOp_*(D^d(\P))$ is
given by a coend of tensor products
\begin{equation*}
\bigotimes_u \tau_u(\P)\otimes\sk_d\Cube{\tau_u},
\end{equation*}
where $\tau_u$ ranges over collection of trees in $\Theta''_d(r)$
arranged on a tree $\sigma$. This assertion follows from a
straightforward interchange of colimits. In addition the map $g^d:
\FOp_*(D^d(\P))\rightarrow W^d(\P)$ can be identified with a coend
morphism induced by canonical morphisms
\begin{equation*}
\bigotimes_u
\tau_u(\P)\otimes\sk_d\Cube{\tau_u}\rightarrow\rho(\P)\otimes\dec_d\Cube{\rho},
\end{equation*}
where $\rho$ is a tree formed by a composite of the trees $\tau_u$
along $\sigma$. The free operad $\FOp_*(C^d(\P))$ and the attaching
map $f^d: \FOp_*(C^d(\P))\rightarrow W^{d-1}(\P)$ have a similar
representation. Then one can use a treewise representation of
operadic composites in order to check by hand that $W^d(\E)$
satisfies the universal property of an operad pushout.

\begin{claim}[Second part of theorem~\ref{thm:BoardmanVogtOperadicCellularDecomposition}]\label{claim:CellInclusions}
The morphism of unitary (Hopf) $\Lambda_*$-modules
\begin{equation*}
(i^d,\phi): C^d(\P')\bigoplus_{C^d(\P)} D^d(\P)\rightarrow D^d(\P')
\end{equation*}
associated to an operad morphism $\phi: \P\rightarrow\P'$ is a
weak-equivalence, respectively a Reedy cofibration, if $\phi$
defines a weak-equivalence, respectively a Reedy cofibration, in the
category of (Hopf) $\Lambda_*$-modules.
\end{claim}

\begin{proof}
This statement is valid in the framework of~\cite{BergerMoerdijkW}
like claim~\ref{claim:WCellExtensions} but for simplicity we give
arguments that make sense only in the differential graded context.

Consider the module $E^d(\P)(r)\subset D^d(\P)(r)$, introduced in
the proof of claim~\ref{claim:WCellExtensions}, spanned by tensors
$\ptensor\otimes\lambda\in\tau(\P)\otimes\Cube{\tau}$ such that
$\deg(\tau) = d$ and $\lambda_e = \x{01}$ for all edges $e\in
E'(\tau)$. Recall that these modules $E^d(\P)(r)$ form a
$\Sigma_*$-module and we have a natural splitting $D^d(\P) =
C^d(\P)\oplus E^d(\P)$. In fact, the $\Sigma_*$-module $E^d(\P)$,
equipped with a differential $\delta$ induced by the internal
differential of~$\P$, can be identified with a quotient object of
$D^d(\P)$ so that we have a short exact sequence
\begin{equation*}
\xymatrix{ C^d(\P)\ar@{>->}[]!R+<4pt,0pt>;[r] & D^d(\P)\ar@{->>}[r]
& E^d(\P) }.
\end{equation*}

Observe that the coproduct $C^d(\P')\bigoplus_{C^d(\P')} D^d(\P)$
admits a similar splitting
\begin{equation*}
C^d(\P')\bigoplus_{C^d(\P')} D^d(\P) = C^d(\P')\oplus E^d(\P)
\end{equation*}
so that the morphism of the claim $(i^d,\phi)$ fits a diagram of
short exact sequences of dg-modules
\begin{equation*}
\xymatrix{ C^d(\P')\ar@{>->}[]!R+<4pt,0pt>;[r]\ar[d]^{=} &
C^d(\P')\bigoplus_{C^d(\P')} D^d(\P)\ar@{->>}[r]\ar[d]^{(i^d,\phi)}
& E^d(\P)\ar[d]^{\phi} \\
C^d(\P')\ar@{>->}[]!R+<4pt,0pt>;[r] & D^d(\P')\ar@{->>}[r] &
E^d(\P') },
\end{equation*}
where $\phi: E^d(\P)\rightarrow E^d(\P')$ denotes the natural
morphism induced by $\phi: \P\rightarrow\P'$.

By definition, we have
\begin{equation*}
E^d(\P)(r) = \bigoplus_{\gr(\tau) =
d}\tau(\P)\otimes\Bigl[\bigotimes_{e\in E'(\tau)}\x{01}\Bigr]
\simeq\bigoplus_{\gr(\tau) = d}\Sigma^d\tau(\P),
\end{equation*}
where the sum ranges over isomorphism classes of reduced $r$-trees
with $d$ internal edges. Recall that a reduced $r$-tree has no
automorphism (see~\ref{item:TreeCategories}). As a consequence, no
relation occurs in the expansion above. Hence if $\phi:
\P\rightarrow\P'$ is a weak-equivalence of operads, then the induced
morphism $\phi: E^d(\P)\rightarrow E^d(\P')$ is a weak-equivalence
of dg-modules and we deduce from the short exact sequence that the
morphism $(\phi,i^d)$ is a weak-equivalence as well, as stated by
the claim.

Recall that a morphism $\phi: \P\rightarrow\P'$ is a Reedy
cofibration in the category of (Hopf) $\Lambda_*$-modules if $\phi:
\P\rightarrow\P'$ defines a cofibration in the category of
$\Sigma_*$-modules. Furthermore, in this case the morphism $\phi:
\P\rightarrow\P'$ can be decomposed into a sequence of
$\Sigma_*$-modules embeddings
\begin{equation*}
\P = M_{-1}\hookrightarrow M_0\hookrightarrow\dots\hookrightarrow
M_n\hookrightarrow\dots \hookrightarrow\colim_n M_n = \P'
\end{equation*}
such that $\delta(M_n)\subset M_{n-1}$ and where
$M_{n-1}\hookrightarrow M_n$ is a split injective morphism of
$\Sigma_*$-modules with a projective cokernel. The treewise tensor
products of~\ref{item:TreewiseTensorProducts} can be defined for
elements of a $\Sigma_*$-modules (see~\cite{OperadTextbook}) so that
we have a sequence of $\Sigma_*$-modules $E^d(M_n)$ defined by
\begin{equation*}
E^d(M_n)(r) = \bigoplus_{\gr(\tau) = d}\Sigma^d\tau(M_n).
\end{equation*}
One can check that the morphisms $E^d(M_{n-1})\hookrightarrow
E^d(M_n)$ are split injective and have a projective cokernel as well
because the symmetric group operates freely on the treewise tensor
products of a free $\Sigma_*$-module. Furthermore, the cellular
differential of $D^d(\P')$ satisfies $\partial(E^d(M_n))\subset
C^d(\P')$ and we have $\delta(E^d(M_n))\subset E^d(M_{n-1})$ for the
differential induced by the internal differential of~$\P'$.
Therefore the sequence of $\Sigma_*$-module inclusions
\begin{multline*}
C^d(\P')\oplus E^d(\P) = C^d(\P')\oplus E^d(M_{-1})\hookrightarrow\dots\\
\dots\hookrightarrow C^d(\P')\oplus E^d(M_n)\hookrightarrow\dots\\
\dots\hookrightarrow\colim_n C^d(\P')\oplus E^d(M_n) =
C^d(\P')\oplus E^d(\P')
\end{multline*}
fulfils our requirements for a cofibration decomposition and show
that the morphism
\begin{equation*}
(i^d,\phi): C^d(\P')\bigoplus_{C^d(\P)} D^d(\P)\rightarrow D^d(\P')
\end{equation*}
defines a cofibration in the category of~$\Sigma_*$-modules. This
conclusion achieves the proof of claim~\ref{claim:CellInclusions}.
\end{proof}

This claim achieves the proof
theorem~\ref{thm:BoardmanVogtOperadicCellularDecomposition} and, as
a byproduct, of the cofibrancy claim in assertion~(b) of
theorem~\ref{thm:WHopfOperad}.\qed

\clearpage
\subsection{Appendix: figures}

\begin{figure}[b]
\begin{multline*}
\vcenter{\xymatrix@M=0pt@R=5mm@C=3mm{
*+<1mm>{1}\ar[dr] & *+<1mm>{2}\ar[d] & & *+<1mm>{3}\ar[d] & *+<1mm>{4}\ar[dl] \\
& *+<1mm>{p_{v_2}}\ar[dr]|*+<2pt>{\scriptstyle\x{01}} & & *+<1mm>{p_{v_3}}\ar[dl]|*+<2pt>{\scriptstyle\x{01}} & \\
& & *+<1mm>{p_{v_1}}\ar[d] & & \\
& & *+<1mm>{0} & & }} \xrightarrow{\Delta} \pm
\vcenter{\xymatrix@M=0pt@R=5mm@C=3mm{
*+<1mm>{1}\ar[dr] & *+<1mm>{2}\ar[d] & & *+<1mm>{3}\ar[d] & *+<1mm>{4}\ar[dl] \\
& *+<1mm>{p'_{v_2}}\ar[dr]|*+<2pt>{\scriptstyle\x{0}} & & *+<1mm>{p'_{v_3}}\ar[dl]|*+<2pt>{\scriptstyle\x{0}} & \\
& & *+<1mm>{p'_{v_1}}\ar[d] & & \\
& & *+<1mm>{0} & & }} \otimes \vcenter{\xymatrix@M=0pt@R=5mm@C=3mm{
*+<1mm>{1}\ar[dr] & *+<1mm>{2}\ar[d] & & *+<1mm>{3}\ar[d] & *+<1mm>{4}\ar[dl] \\
& *+<1mm>{p''_{v_2}}\ar[dr]|*+<2pt>{\scriptstyle\x{01}} & & *+<1mm>{p''_{v_3}}\ar[dl]|*+<2pt>{\scriptstyle\x{01}} & \\
& & *+<1mm>{p''_{v_1}}\ar[d] & & \\
& & *+<1mm>{0} & & }}
\\
\pm \vcenter{\xymatrix@M=0pt@R=5mm@C=3mm{
*+<1mm>{1}\ar[dr] & *+<1mm>{2}\ar[d] & & *+<1mm>{3}\ar[d] & *+<1mm>{4}\ar[dl] \\
& *+<1mm>{p'_{v_2}}\ar[dr]|*+<2pt>{\scriptstyle\x{0}} & & *+<1mm>{p'_{v_3}}\ar[dl]|*+<2pt>{\scriptstyle\x{01}} & \\
& & *+<1mm>{p'_{v_1}}\ar[d] & & \\
& & *+<1mm>{0} & & }} \otimes \vcenter{\xymatrix@M=0pt@R=5mm@C=3mm{
*+<1mm>{1}\ar[dr] & *+<1mm>{2}\ar[d] & & *+<1mm>{3}\ar[d] & *+<1mm>{4}\ar[dl] \\
& *+<1mm>{p''_{v_2}}\ar[dr]|*+<2pt>{\scriptstyle\x{01}} & & *+<1mm>{p''_{v_3}}\ar[dl]|*+<2pt>{\scriptstyle\x{1}} & \\
& & *+<1mm>{p''_{v_1}}\ar[d] & & \\
& & *+<1mm>{0} & & }}
\\
\pm \vcenter{\xymatrix@M=0pt@R=5mm@C=3mm{
*+<1mm>{1}\ar[dr] & *+<1mm>{2}\ar[d] & & *+<1mm>{3}\ar[d] & *+<1mm>{4}\ar[dl] \\
& *+<1mm>{p'_{v_2}}\ar[dr]|*+<2pt>{\scriptstyle\x{01}} & & *+<1mm>{p'_{v_3}}\ar[dl]|*+<2pt>{\scriptstyle\x{0}} & \\
& & *+<1mm>{p'_{v_1}}\ar[d] & & \\
& & *+<1mm>{0} & & }} \otimes \vcenter{\xymatrix@M=0pt@R=5mm@C=3mm{
*+<1mm>{1}\ar[dr] & *+<1mm>{2}\ar[d] & & *+<1mm>{3}\ar[d] & *+<1mm>{4}\ar[dl] \\
& *+<1mm>{p''_{v_2}}\ar[dr]|*+<2pt>{\scriptstyle\x{1}} & & *+<1mm>{p''_{v_3}}\ar[dl]|*+<2pt>{\scriptstyle\x{01}} & \\
& & *+<1mm>{p''_{v_1}}\ar[d] & & \\
& & *+<1mm>{0} & & }}
\\
\pm \vcenter{\xymatrix@M=0pt@R=5mm@C=3mm{
*+<1mm>{1}\ar[dr] & *+<1mm>{2}\ar[d] & & *+<1mm>{3}\ar[d] & *+<1mm>{4}\ar[dl] \\
& *+<1mm>{p'_{v_2}}\ar[dr]|*+<2pt>{\scriptstyle\x{01}} & & *+<1mm>{p'_{v_3}}\ar[dl]|*+<2pt>{\scriptstyle\x{01}} & \\
& & *+<1mm>{p'_{v_1}}\ar[d] & & \\
& & *+<1mm>{0} & & }} \otimes \vcenter{\xymatrix@M=0pt@R=5mm@C=3mm{
*+<1mm>{1}\ar[dr] & *+<1mm>{2}\ar[d] & & *+<1mm>{3}\ar[d] & *+<1mm>{4}\ar[dl] \\
& *+<1mm>{p''_{v_2}}\ar[dr]|*+<2pt>{\scriptstyle\x{1}} & & *+<1mm>{p''_{v_3}}\ar[dl]|*+<2pt>{\scriptstyle\x{1}} & \\
& & *+<1mm>{p''_{v_1}}\ar[d] & & \\
& & *+<1mm>{0} & & }}
\end{multline*}\caption{The diagonal of the $W$-construction}\label{figure:MetricTreeDiagonal}\end{figure}
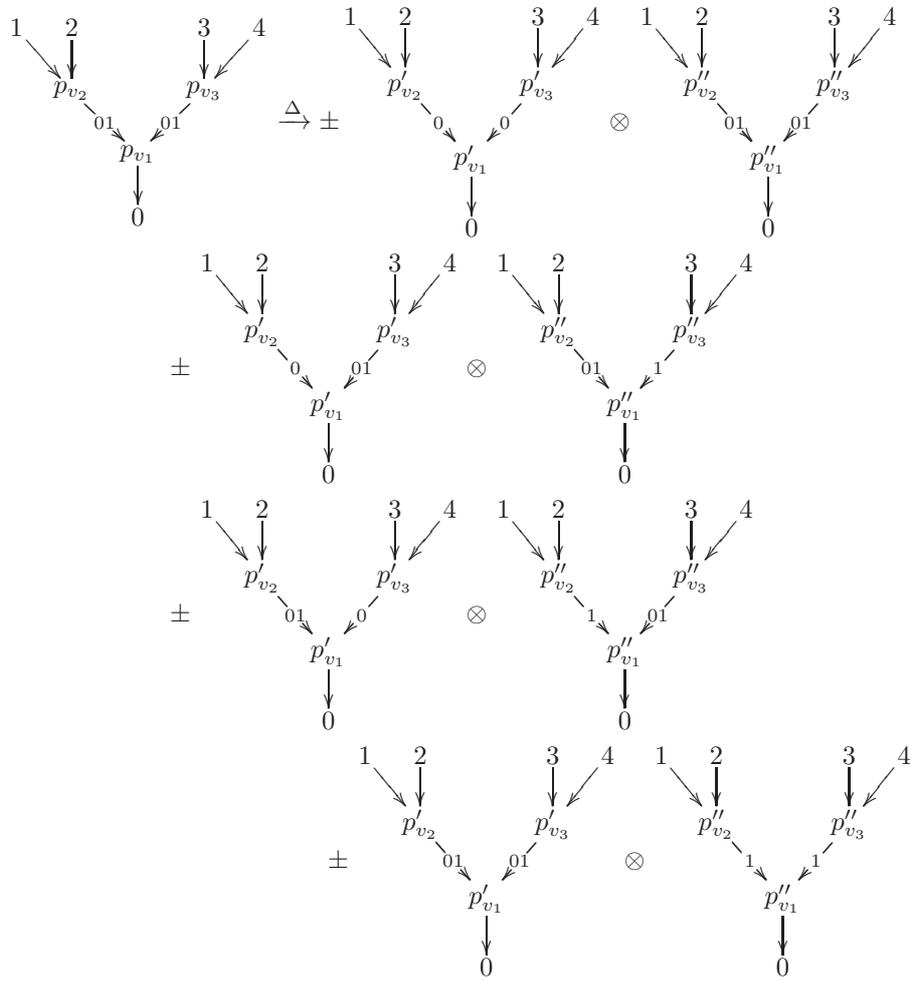

\clearpage

\part*{Construction of Hopf operad actions}

\section{Cofree coalgebras and quasi-cofree Hopf $\Lambda_*$-modules}\label{section:CocellularCoalgebras}

\subsection{Introduction}
The goal of this section is to introduce a suitable notion of
\emph{cocellular complex} in order to obtain an effective class of
fibrations in the category of Hopf operads. Our cocellular objects
do not clearly generate the class of fibrations in the category of
Hopf algebras but the \emph{Hopf operad of bar operations}
$\HopfOp_B^{\P}$ defined in the next section (our main device for
the construction of operad actions) has such a cocellular structure.
Hence the framework introduced in this section enables us to prove
that the functor $\P\mapsto\HopfOp_B^{\P}$ maps fibrations,
respectively acyclic fibrations, of dg-operads to fibrations,
respectively acyclic fibrations, of unital Hopf operads. This
assertion allows us to deduce the existence of operad morphisms
$\Q\rightarrow\HopfOp_B^{\P}$ from model category arguments.

For simplicity, we perform our constructions in the category of
coalgebras first and we extend our definitions to Hopf
$\Lambda_*$-modules next.

In our constructions, we consider $\Z$-graded objects which are
outside the coalgebra model category considered
in~\ref{section:UnitaryHopfOperads}. Thus, for our purposes, we have
to extend the notion of a fibration and of an acyclic fibration to
this context. In fact, we deal with lifting problems
\begin{equation*}
\xymatrix{ A\ar[r]\ar[d]_{i} & C\ar[d]^{q} \\
B\ar[r]\ar@{-->}[ur] & D }
\end{equation*}
such that $i: A\rightarrow B$ is a morphism of non-negatively graded
objects. Therefore we distinguish morphisms of (possibly
$\Z$-graded) dg-coalgebras $q: C\rightarrow D$ such that this
lifting problem has a solution for any acyclic cofibration,
respectively cofibration, of $\N$-graded dg-coalgebras $i:
A\rightarrow B$. Clearly, we can also characterize such morphisms of
$\Z$-graded coalgebras $q: C\rightarrow D$ by the right lifting
property with respect to a generating set of acyclic cofibrations,
respectively cofibrations, of $\N$-graded dg-coalgebras. By
convention we refer abusively to these class of morphisms $q:
C\rightarrow D$ as fibrations, respectively acyclic fibrations, of
dg-coalgebras.

Equivalently, let $\dg_{\N}\CoAlg$, respectively $\dg_{\Z}\CoAlg$,
denote the category of $\N$-graded, respectively $\Z$-graded
dg-coalgebras. One can observe that the category embedding
$i^{\dg}_+: \dg_{\N}\CoAlg\rightarrow\dg_{\Z}\CoAlg$ has a right
adjoint $\str^{\dg}_+: \dg_{\Z}\CoAlg\rightarrow\dg_{\N}\CoAlg$.
Clearly, a morphism of coalgebras $q: C\rightarrow D$ forms a
fibration, respectively an acyclic fibration, according to the
convention above if the associated morphism $\str^{\dg}_+(q):
\str^{\dg}_+(C)\rightarrow\str^{\dg}_+(D)$ defines a fibration,
respectively an acyclic fibration in the model category of
$\N$-graded dg-coalgebras.\index{fibration!of a $\Z$-graded
coalgebra}\index{fibration!of a $\Z$-graded Hopf object}

One observes that the category of augmented coalgebras is equipped
with cofree objects characterized by the usual universal property.
Namely, for any dg-module $V$, we have a coalgebra, denoted by
$\Gamma(V)$, equipped with a dg-module morphism $\pi:
\Gamma(V)\rightarrow V$ such that any dg-module morphism $f:
\Gamma\rightarrow\Gamma(V)$ where $\Gamma$ is a dg-coalgebra admits
one and only one factorization $f = \pi\cdot\nabla_f$ where
$\nabla_f: \Gamma\rightarrow\Gamma(V)$ is a coalgebra morphism. For
our needs we give an explicit construction of this
object~$\Gamma(V)$ in~\ref{subsection:CofreeCoalgebras}. The cofree
coalgebra is equipped with a natural differential $\delta:
\Gamma(V)\rightarrow\Gamma(V)$ induced by the internal differential
of~$V$. For our purposes we consider \emph{quasi-cofree coalgebras}
which are specified as usual by a cofree coalgebra $\Gamma =
\Gamma(V)$ equipped with a coderivation $\partial:
\Gamma(V)\rightarrow\Gamma(V)$ such that $\delta+\partial$ defines
the differential of $\Gamma$.

Our goal is to give sufficient conditions for a morphism of
quasi-cofree coalgebras
\begin{equation*}
\nabla_f:
(\Gamma(V),\partial_\alpha)\rightarrow(\Gamma(W),\partial_\beta)
\end{equation*}
to be a fibration, respectively an acyclic fibration. Explicitly, we
observe that the coderivation of a quasi-cofree coalgebra
$(\Gamma(V),\partial_\alpha)$ is equivalent to a morphism of
dg-coalgebras $\nabla_\alpha:
(\Gamma(V),\partial_\alpha)\rightarrow\Gamma(\Delta^1\wedge V)$,
where $\Delta^1\wedge V$ is the cone of~$V$ in the category of
dg-modules. For our purpose we determine the structure of
morphisms~$\nabla_f$ that fit a pullback diagram of the form
\begin{equation*}
\xymatrix{ (\Gamma(V),\partial_\alpha)\ar[r]\ar[d]
& \Gamma(\Delta^1\wedge V)\ar[d] \\
(\Gamma(W),\partial_\beta)\ar[r] & \Gamma(\Delta^1\wedge
W)\times_{\Gamma(S^1\wedge W)}\Gamma(S^1\wedge V) },
\end{equation*}
where we consider natural morphisms of cofree coalgebras
\begin{equation*}
\xymatrix{ \Gamma(\Delta^1\wedge V)\ar[r]\ar[d] & \Gamma(S^1\wedge V)\ar[d] \\
\Gamma(\Delta^1\wedge W)\ar[r] & \Gamma(S^1\wedge W) }
\end{equation*}
induced by a morphism of dg-modules $f: V\rightarrow W$ and a
canonical morphism $\sigma\wedge-: \Delta^1\wedge-\rightarrow
S^1\wedge-$. By immediate categorical arguments, we obtain that
$\nabla_f$ defines a fibration, respectively an acyclic fibration,
of dg-coalgebras if $f: V\rightarrow W$ is a fibration, respectively
an acyclic fibration, of dg-modules. This program is carried out
in~\ref{subsection:CocellularCoalgebras}. This construction admits a
natural and straightforward extension to the category of Hopf
$\Lambda_*$-modules. For the sake of completeness, we state
explicitly the results that we obtain in this context
in~\ref{subsection:CocellularHopfLambdaModules}.

For our needs, we consider morphisms of quasi-cofree coalgebras
$\nabla_f:
(\Gamma(V),\partial_\alpha)\rightarrow(\Gamma(W),\partial_\beta)$
obtained as the limit of a tower of morphisms
\begin{multline*}
(\Gamma(V),\partial_\alpha) = \lim_m(\Gamma(\ckcell_m
V),\partial_\alpha)
\rightarrow\dots\\
\dots\rightarrow(\Gamma(\ckcell_m V),\partial_\alpha)
\xrightarrow{\pk_m}(\Gamma(\ckcell_{m-1} V),\partial_\alpha)
\rightarrow\dots\\
\dots\rightarrow(\Gamma(\ckcell_{-1} V),\partial_\alpha) = \F
\end{multline*}
such that the maps $\pk_m: (\Gamma(\ckcell_m
V),\partial_\alpha)\rightarrow(\Gamma(\ckcell_{m-1}
V),\partial_\alpha)$ fit coextension diagrams as above.

One can observe that any morphism of quasi-cofree coalgebras where
$V$ and $W$ are non-negatively graded has a such a decomposition
that arises from the canonical degreewise filtration of dg-modules.
Moreover this argument can also be applied to the truncation
$\str^{\dg}_+\nabla_f$ of a morphisms of quasi-cofree coalgebras
$\nabla_f$ because the truncation functor preserves quasi-cofree
objects. In the memoir, we do not give this general construction. In
fact, in the next section, we define only a specific decomposition
\begin{multline*}
\HopfOp_B^{\P} = \lim_m\ckcell_m\HopfOp_B^{\P}
\rightarrow\dots\\
\dots\rightarrow\ckcell_{m}\HopfOp_B^{\P}
\rightarrow\ckcell_{m-1}\HopfOp_B^{\P}
\rightarrow\dots\\
\dots\rightarrow\ckcell_0\HopfOp_B^{\P} = \C,
\end{multline*}
obtained similarly as the general construction, but more natural in
regard to the Hopf operad of universal bar operations
$\HopfOp_B^{\P}$. Then we use this construction to prove that the
functor $\P\mapsto\HopfOp_B^{\P}$ maps an operad fibration,
respectively acyclic fibration, to a fibration, respectively an
acyclic fibration, of unital Hopf operads.

Here is the plan of this section.
In~\ref{subsection:CofreeCoalgebras} we recall categorical
properties of coalgebras: existence of small limits, of cofree
objects and adjoint functors. This subsection does not contain any
original result but the explicit constructions recalled there allows
us to introduce a notion of a quasi-cofree object in the category of
coalgebras. In the next
subsection~\ref{subsection:CocellularCoalgebras} we study the
structure of quasi-cofree coalgebras with the aim to define
coalgebra fibrations. In the final
subsection~\ref{subsection:CocellularHopfLambdaModules} we extend
these results to Hopf $\Lambda_*$-modules.

\subsection{Cofree coalgebras}\label{subsection:CofreeCoalgebras}
As stated, the aim of this subsection is to recall categorical
properties of augmented coassociative coalgebras. For our needs the
main issue is to give an effective construction of cofree objects.
For the sake of precision, we state the definition of these objects
in a proposition:

\begin{prop}\label{prop:CofreeCoalgebras}\index{cofree!coalgebra}\index{coalgebra!cofree} Any dg-module $V$ has an
associated cofree augmented coassociative coalgebra
$\Gamma(V)$\glossary{$\Gamma(V)$} equipped with a morphism of
dg-modules $\pi: \Gamma(V)\rightarrow V$\glossary{$\pi$}
characterized by the classical universal property. Namely any
morphism of dg-modules $f: \Gamma\rightarrow V$ where $\Gamma$ is a
coalgebra has a unique factorization $f = \pi\cdot\nabla_f$ such
that $\nabla_f: \Gamma\rightarrow\Gamma(V)$\glossary{$\nabla_f$} is
a coalgebra morphism.
\end{prop}

We refer to~\cite{Sweedler} for a proof of this result in the case
of non-graded coalgebras, to~\cite{GetzlerGoerss} for the case of
non-negatively graded coalgebras and to~\cite{SmithJustin} for a
generalization in the context of coalgebras over an operad. In this
subsection we give an explicit realization of $\Gamma(V)$ which is a
special instance of the construction of the latter reference in the
case of coassociative coalgebras. As such this subsection does not
contain any original idea and our account is only motivated by the
applications of the next subsection.

As explained in the introduction, we have to deal with $\Z$-graded
coalgebras. Therefore we prove
proposition~\ref{prop:CofreeCoalgebras} in this framework. On the
other hand, our assertion holds in both the category of $\Z$-graded
coalgebras and the category of $\N$-graded coalgebras. In fact, if
$V$ is an $\N$-graded dg-module, then the associated cofree object
in the category $\Z$-graded coalgebras turns out to be $\N$-graded
and gives also a realization of the cofree object cogenerated by $V$
in the category of $\N$-graded coalgebras.

\subsubsection{An inductive construction}
We define a nested sequence of dg-modules
\begin{equation*}
\xymatrix{ *{\dots}\ar@{^{(}->}[]!R+<4pt,0pt>;[r]!L-<2pt,0pt> &
*{\Gamma_{r+1}(V)}\ar@{^{(}->}[]!R+<4pt,0pt>;[r]!L-<2pt,0pt> &
*{\Gamma_{r}(V)}\ar@{^{(}->}[]!R+<4pt,0pt>;[r]!L-<2pt,0pt> &
*{\dots}\ar@{^{(}->}[]!R+<4pt,0pt>;[r]!L-<2pt,0pt> &
*{\Gamma_{1}(V) = \prod_{n=0}^{\infty} V^{\otimes n}} }
\end{equation*}
and we prove that the module $\Gamma_\infty(V) =
\bigcap_{r=1}^{\infty} \Gamma_r(V)$ is equipped with a coalgebra
structure and represents the cofree coalgebra cogenerated by $V$.

The modules $\Gamma_r(V)$ are defined by induction. We consider the
natural map
\begin{equation*}
\xymatrix{ \{\prod_{m} V^{\otimes m}\}\otimes\{\prod_{n} V^{\otimes
n}\}\ar[r]^{\nabla_\Pi} & \prod_{m,n}\{V^{\otimes m}\otimes
V^{\otimes n}\} }
\end{equation*}
and the composite
\begin{equation*}
\xymatrix{ \prod_{N} V^{\otimes
N}\ar[r]^(0.3){\{\Delta_N\}}\ar`u[r]`/4pt[rr]^{\Delta_\Pi}`_d[rr]+U
& \prod_{N}\{\prod_{m+n = N} V^{\otimes m}\otimes V^{\otimes
n}\}\ar[r]^(0.55){\simeq} & \prod_{m,n}\{V^{\otimes m}\otimes
V^{\otimes n}\} },
\end{equation*}
where $\Delta_N: V^{\otimes N}\rightarrow\prod_{m+n = N} V^{\otimes
m}\otimes V^{\otimes n}$ denotes the deconcatenation of tensors. In
general, the natural map $A\otimes\prod_j B_j\rightarrow\prod_j
A\otimes B_j$ is an embedding provided that the module $A$ is free
over the ground ring. Accordingly, in our context, the map
$\nabla_\Pi$ is always injective since our ground ring $\F$ is
supposed to be a field. Notice that the map $\Delta_\Pi$ is
injective as well.

For $r = 1$, we set $\Gamma_1(V) = \prod_{n=0}^{\infty} V^{\otimes
n}$. By induction, we have a module $\Gamma_{r}(V)$ equipped with an
embedding $\iota_r: \Gamma_{r}(V)\hookrightarrow\prod_{n=0}^{\infty}
V^{\otimes n}$. The next module $\Gamma_{r+1}(V)$ is defined by the
fiber product
\begin{equation*}
\xymatrix{ \Gamma_{r+1}(V)\pull{6} 
\ar@{-->}[d]\ar@{-->}[r] &
\Gamma_{r}(V)\otimes\Gamma_{r}(V)\ar[d]^{\nabla_\Pi\cdot\iota_r\otimes\iota_r} \\
\Gamma_{r}(V)\ar[r]^(0.35){\Delta_\Pi\cdot\iota_r} &
\prod_{m,n}\{V^{\otimes m}\otimes V^{\otimes n}\} }.
\end{equation*}
The composite map $\nabla_\Pi\cdot\iota_r\otimes\iota_r$ is
injective according to the observations above. Therefore the map
$\Gamma_{r+1}(V)\rightarrow\Gamma_{r}(V)$, which is defined by a
base extension of $\nabla_\Pi\cdot\iota_r\otimes\iota_r$, is
injective as well.

As announced, we consider the dg-module $\Gamma_\infty(V) =
\bigcap_{r=0}^{\infty} \Gamma_r(V)$ equipped with the embedding
$\iota_\infty: \Gamma_{\infty}(V)\hookrightarrow\prod_{n=0}^{\infty}
V^{\otimes n}$.\glossary{$\iota_\infty$}

\begin{claim}\label{claim:CoFreeDiagonal}
The module $\Gamma_{\infty}(V)$ is equipped with a diagonal
$\Delta_{\infty}$ such that we have a commutative diagram
\begin{equation*}
\xymatrix{
\Gamma_{\infty}(V)\ar@{^{(}->}[]!D-<0pt,4pt>;[d]^(0.4){\iota_\infty}\ar@{-->}[r]^(0.4){\Delta_{\infty}}
&
\Gamma_{\infty}(V)\otimes\Gamma_{\infty}(V)\ar@{^{(}->}[]!D-<0pt,4pt>;[d]^(0.4){\nabla_\Pi\cdot\iota_\infty\otimes\iota_\infty} \\
\prod_N V^{\otimes N}\ar[r]^(0.4){\Delta_\Pi} & \prod_{m,n}
\{V^{\otimes m}\otimes V^{\otimes n}\} }.
\end{equation*}
\end{claim}

\begin{proof}
By definition of $\Gamma_{r+1}(V)$, for $r\in\N$, we have maps
$\Delta_r: \Gamma_{r+1}(V)\rightarrow\Gamma_r(V)\otimes\Gamma_r(V)$
which fit a commutative diagram as in the claim statement.
Accordingly, these maps restrict to a map from the intersection
$\Gamma_{\infty}(V) = \bigcap_{r=1}^{\infty}\Gamma_r(V)$ to the
module $\bigcap_{r=1}^{\infty}\{\Gamma_r(V)\otimes\Gamma_r(V)\}$ and
the claim is a consequence of the next assertion.
\end{proof}

\begin{claim}
The natural map
\begin{equation*}
\Gamma_\infty(V)\otimes\Gamma_\infty(V) = \{\bigcap_{r=1}^{\infty}
\Gamma_r(V)\}\otimes\{\bigcap_{s=1}^{\infty} \Gamma_s(V)\}
\rightarrow\bigcap_{r=1}^{\infty}\{\Gamma_r(V)\otimes\Gamma_r(V)\}
\end{equation*}
is an isomorphism.
\end{claim}

\begin{proof}
This map connects submodules of $\prod_{m,n} V^{\otimes m}\otimes
V^{\otimes n}$ and hence represents a submodule inclusion. We prove
that any element
$\omega\in\bigcap_{r=1}^{\infty}\{\Gamma_r(V)\otimes\Gamma_r(V)\}$
belongs to
$\bigcap_{r=1}^{\infty}\{\Gamma_r(V)\}\otimes\bigcap_{s=1}^{\infty}\{\Gamma_r(V)\}$.

There is a finitely generated module $\Omega_1\subset\Gamma_1(V)$
such that $\omega\in\Omega_1\otimes\Omega_1$ inside
$\Gamma_1(V)\otimes\Gamma_1(V)$. If we let $\Omega_r =
\Omega_1\cap\Gamma_r(V)$, then we have
$\{\Omega_1\otimes\Omega_1\}\cap\{\Gamma_r(V)\otimes\Gamma_r(V)\} =
\Omega_r\otimes\Omega_r$ and hence
$\omega\in\bigcap_{r=1}^{\infty}\{\Omega_r\otimes\Omega_r\}$. On the
other hand, since $\Omega_1$ is finitely generated, the sequence
$\Omega_r$ is necessary stationary: we have $\Omega_{r_0} =
\Omega_{\infty}\subset\Gamma_{\infty}(V)$ for some $r_0<\infty$.
Hence the relation $\omega\in\Omega_{r_0}\otimes\Omega_{r_0}$
implies $\omega\in\Gamma_\infty(V)\otimes\Gamma_\infty(V)$.
\end{proof}

\begin{lemm}\label{lemm:CoFreeCoAlg}
The module $\Gamma(V) = \Gamma_{\infty}(V)$ defines a realization of
the cofree augmented coassociative coalgebra cogenerated by $V$.

To be precise, the coproduct of $\Gamma(V)$ is given by the map
$\Delta_{\infty}: \Gamma_{\infty}(V)\rightarrow\Gamma_{\infty}(V)$
of claim~\ref{claim:CoFreeDiagonal}, the augmentation $\epsilon:
\Gamma(V)\rightarrow\F$ is defined by the composite of the embedding
$\iota_{\infty}: \Gamma_{\infty}(V)\hookrightarrow\prod_n V^{\otimes
n}$ with the projection onto the component $n=0$ of the product
$\prod_n V^{\otimes n}$ and the universal morphism $\pi:
\Gamma(V)\rightarrow V$ is defined by the composite of
$\iota_{\infty}$ with the projection onto the component $n=1$.
\end{lemm}

\begin{proof}
Observe that the diagonal $\Delta_\infty$ specified in
claim~\ref{claim:CoFreeDiagonal} is coassociative simply because the
deconcatenation of tensors is componentwise coassociative. One
checks similarly that the composite of the embedding
$\iota_{\infty}: \Gamma_{\infty}(V)\hookrightarrow\prod_n V^{\otimes
n}$ with the projection onto the component $n=0$ of the product
$\prod_n V^{\otimes n}$ defines an augmentation for the diagonal
$\Delta_\infty$.

We check the universal property of a cofree coalgebra. Let $\Gamma$
be a coalgebra equipped with a morphism of dg-modules $f:
\Gamma\rightarrow V$. Consider the map $\widehat{\nabla}_f:
\Gamma\rightarrow\prod_n V^{\otimes n}$ which maps an element
$\gamma\in\Gamma$ to the collection of tensors $\{f^{\otimes
n}\cdot\Delta^n(\gamma)\}$, where
$\Delta^n(\gamma)\in\Gamma^{\otimes n}$ denotes the $n$-fold
diagonal of $\gamma\in\Gamma$.

We claim that this map $\widehat{\nabla}_f$ admits a sequence of
factorizations
\begin{equation*}
\xymatrix{
\Gamma\ar[]!D-<0pt,2pt>;[d]!U+<0pt,2pt>^(0.5){\nabla_{\infty}}
\ar@/^5pt/[]!R+<2pt,0pt>;[drr]!U+<0pt,2pt>^(0.92){\nabla_{r+1}}
\ar@/^5pt/[]!R+<3pt,0pt>;[drrr]!U+<0pt,2pt>^(0.95){\nabla_{r}}
\ar@/^5pt/[]!R+<4pt,0pt>;[drrrrr]!UL+<8pt,2pt>^(0.80){\nabla_{1}=\widehat{\nabla}_f} &&&& \\
*{\Gamma_{\infty}(V)}\ar@{^{(}->}[]!R+<4pt,0pt>;[r]!L-<2pt,0pt> &
*{\dots}\ar@{^{(}->}[]!R+<4pt,0pt>;[r]!L-<2pt,0pt> &
*{\Gamma_{r+1}(V)}\ar@{^{(}->}[]!R+<4pt,0pt>;[r]!L-<2pt,0pt> &
*{\Gamma_{r}(V)}\ar@{^{(}->}[]!R+<4pt,0pt>;[r]!L-<2pt,0pt> &
*{\dots}\ar@{^{(}->}[]!R+<4pt,0pt>;[r]!L-<2pt,0pt> &
*{\Gamma_{1}(V) = \prod_n V^{\otimes n}} }.
\end{equation*}
Indeed, by induction, we are given a map $\nabla_r:
\Gamma\rightarrow\Gamma_r(V)$ such that $\iota_r\cdot\nabla_r =
\widehat{\nabla}_f$. Observe that the map $\widehat{\nabla}_f$ makes
the diagram
\begin{equation*}
\xymatrix{ \Gamma\ar[d]^{\widehat{\nabla}_f}\ar[r]^(0.4){\Delta} &
\Gamma\otimes\Gamma\ar[d]^{\nabla_\Pi\cdot\widehat{\nabla}_f\otimes\widehat{\nabla}_f} \\
\prod_N V^{\otimes N}\ar[r]^(0.4){\Delta_\Pi} & \prod_{m,n}
\{V^{\otimes m}\otimes V^{\otimes n}\} }
\end{equation*}
commute and, as a consequence, our map $\nabla_r:
\Gamma\rightarrow\Gamma_r(V)$ fits a commutative diagram
\begin{equation*}
\xymatrix{
\Gamma\ar[r]^{\Delta}\ar@/_10pt/[ddr]_{\nabla_r}\ar@{-->}[dr] &
\Gamma\otimes\Gamma\ar@/^/[dr]^{\nabla_r\otimes\nabla_r} & \\
& \Gamma_{r+1}(V)\pull{6} 
\ar@{-->}[d]\ar@{-->}[r] & \Gamma_{r}(V)\otimes\Gamma_{r}(V)\ar[d]^{\nabla_\Pi\cdot\iota_r\otimes\iota_r} \\
& \Gamma_{r}(V)\ar[r]^(0.35){\Delta_\Pi\cdot\iota_r} &
\prod_{m,n}\{V^{\otimes m}\otimes V^{\otimes n}\} }.
\end{equation*}
Accordingly, the existence of $\nabla_{r+1}:
\Gamma\rightarrow\Gamma_{r+1}(V)$ follows from the fiber product
definition of~$\Gamma_{r+1}(V)$.

According to this construction, the resulting map
\begin{equation*}
\Gamma\xrightarrow{\nabla_\infty}\bigcap_{r=1}^{\infty}\Gamma_r(V) =
\Gamma_\infty(V)
\end{equation*}
commutes with the diagonal of~$\Gamma_\infty(V)$. Hence we conclude
that the map $f: \Gamma\rightarrow V$ lifts to a coalgebra morphism
$\nabla_\infty: \Gamma\rightarrow\Gamma_\infty(V)$. The next
assertion implies any coalgebra morphism $\nabla:
\Gamma\rightarrow\Gamma_\infty(V)$ that lifts the map $f:
\Gamma\rightarrow V$ defines a factorization of the map
$\widehat{\nabla}_f: \Gamma\rightarrow\prod_n V^{\otimes n}$ through
$\Gamma_\infty(V)$. Therefore the uniqueness property of the
coalgebra lifting follows from the injectivity of the canonical
embedding $\iota_\infty: \Gamma_\infty(V)\hookrightarrow\prod_n
V^{\otimes n}$. This observation achieves the proof of lemma
\ref{lemm:CoFreeCoAlg}.
\end{proof}

\begin{obsv}\label{obsv:CoFreeEmbedding}
The embedding $\iota_{\infty}:
\Gamma_{\infty}(V)\hookrightarrow\prod_n V^{\otimes
n}$\glossary{$\iota_\infty$} can be identified with the composite
\begin{equation*}
\Gamma_{\infty}(V)
\xrightarrow{\{\Delta_\infty^n\}}\prod_n\Gamma_{\infty}(V)^{\otimes
n}\xrightarrow{\{\pi^{\otimes n}\}}\prod_n V^{\otimes n},
\end{equation*}
where $\Delta_\infty^n$ denotes the $n$-fold diagonal of the cofree
coalgebra.
\end{obsv}

\begin{proof}
This assertion is an immediate consequence of the definition of the
diagonal given in claim~\ref{claim:CoFreeDiagonal}.
\end{proof}

Recall that the category of augmented coassociative coalgebras is
denoted by $\CoAlg^a_+$. Beside the construction of cofree
coalgebras, we recall that~$\CoAlg^a_+$ has all small limits. Our
constructions are standard for a category of coalgebras over a
comonad. First, we have the following classical result:

\begin{lemm}
The forgetful functor from $\CoAlg^a_+$ to the category of
dg-modules creates the equalizers which are reflexive in the
category of dg-modules.
\end{lemm}

\begin{proof}
Explicitly, we consider a pair of coalgebra morphisms
\begin{equation*}
d^0,d^1: \Gamma^0\rightarrow\Gamma^1
\end{equation*}
together with a map $s^0: \Gamma^1\rightarrow\Gamma^0$ such that
$s^0 d^0 = s^0 d^1 = \Id$.

One checks readily that $\ker(d^0-d^1)\otimes\ker(d^0-d^1)$ is the
equalizer of the morphisms
\begin{equation*}
(d^0\otimes d^0,d^1\otimes d^1:
\Gamma^0\otimes\Gamma^0\rightarrow\Gamma^1\otimes\Gamma^1
\end{equation*}
in the category of dg-modules. Indeed, if the ground ring is a
field, then we have $\ker(d^0-d^1)\otimes\ker(d^0-d^1) =
\ker((d^0-d^1)\otimes\Id)\cap\ker(\Id\otimes(d^0-d^1))$. On the
other hand, the relation $(d^0\otimes d^0-d^1\otimes d^1)(\gamma) =
0$ implies $(\Id\otimes s^0)(d^0\otimes d^0-d^1\otimes d^1)(\gamma)
= (d^0\otimes\Id-d^1\otimes\Id)(\gamma) = 0$ and symmetrically
$(\Id\otimes d^0-\Id\otimes d^1)(\gamma) = 0$. Hence, for a
reflexive pair, we have
\begin{equation*}
\ker(d^0\otimes d^0-d^1\otimes d^1) =
\ker(d^0-d^1)\otimes\ker(d^0-d^1).
\end{equation*}

One deduces from this observation that $\ker(d^0-d^1)$ forms a
subcoalgebra of $\Gamma^0$ and the lemma follows.
\end{proof}

Then we obtain:

\begin{lemm}
The category of augmented coalgebras $\CoAlg^a_+$ has small
products.
\end{lemm}

\begin{proof}
This assertion is classical for coalgebras over a comonad. In fact,
one can observe that a product $\prod_\alpha X_\alpha$ can be
defined by a reflexive coequalizer of cofree objects. Namely:
\begin{equation*}
\prod_\alpha X_\alpha = \ker\bigl(\xymatrix{ \Gamma(\prod_{\alpha}
X_\alpha)\ar@<+2pt>[r]^{d^0}\ar@<-2pt>[r]_{d^1}
  & \Gamma(\prod_{\alpha} \Gamma(X_\alpha))\ar@/_6mm/[l]_{s^0}}\bigr),
\end{equation*}
where $d^0$ is induced by the coalgebra structure morphisms
$\rho_\alpha: X_\alpha\rightarrow\Gamma(X_\alpha)$ and $d^1$ is the
composite of the comonad coproduct $\nu:
\Gamma(X)\rightarrow\Gamma(\Gamma(X))$ of the cofree coalgebra
functor with the morphism $(p_\alpha)_*:
\Gamma(\Gamma(\prod_{\alpha}
X_\alpha))\rightarrow\Gamma(\prod_{\alpha} \Gamma(X_\alpha))$
induced by the canonical projections.
\end{proof}

Observe also that $\CoAlg^a_+$ comes equipped with a final object
which is defined by the ground field $* = \F$ since any coalgebra
$\Gamma\in\CoAlg^a_+$ is supposed to be augmented over $\F$.

By standard categorical constructions, the existence of all small
limits can be deduced from these particular cases (reflexive
equalizers, small products and the final object). Hence we obtain
the expected proposition:

\begin{prop}
The category of augmented coalgebras $\CoAlg^a_+$ has all small
limits.\qed
\end{prop}

As recalled in the introduction of this section, the model structure
used in this memoir has been defined for $\N$-graded coalgebras
only; we have to consider $\Z$-graded coalgebras but we shall deal
only with lifting problems
\begin{equation*}
\xymatrix{ A\ar[r]\ar[d]_{i} & C\ar[d]^{q} \\
B\ar[r]\ar@{-->}[ur] & D }
\end{equation*}
such that $i: A\rightarrow B$ is a morphism of non-negatively graded
objects. In fact, one can use the following proposition in order to
put such problems into an $\N$-graded framework:

\begin{prop}
The category embedding $i^{\dg}_+:
\dg_{\N}\CoAlg^a_+\hookrightarrow\dg_{\Z}\CoAlg^a_+$ has a right
adjoint $\str^{\dg}_+:
\dg_{\Z}\CoAlg^a_+\rightarrow\dg_{\N}\CoAlg^a_+$.
\end{prop}\glossary{$i^{\dg}_+$}\glossary{$\str^{\dg}_+$}

In fact, the construction of this functor $\str^{\dg}_+:
\dg_{\Z}\CoAlg^a_+\rightarrow\dg_{\N}\CoAlg^a_+$ is not essential
for our purposes and we do not use this functor explicitly.
Therefore we can skip the detailed verification of this proposition.
On the other hand, the proposition is a straightforward consequence
of the special adjoint functor theorem. Recall simply that colimits
in a coalgebra category are created in the ground category of
dg-modules. As a consequence, the functor $i^{\dg}_+$ preserves
colimits. The category of dg-coalgebras has also a set of
generators. One can observe more precisely that any $\Z$-graded
coalgebra is a colimit of finite dimensional coalgebras as in the
non-graded framework, for which we refer to the classical
book~\cite{Sweedler}, or as in the $\N$-graded framework, for which
we refer to the article~\cite{GetzlerGoerss}. In fact, the proof
given in these references can be extended to the $\Z$-graded
context.

One can also adapt the classical arguments used for another adjoint
functor in~\ref{subsection:MorphismCoalgebras} in order to obtain an
explicit realization of the functor~$\str^{\dg}_+$. Namely recall
first that we have a standard truncation functor on dg-modules.
Then, for a cofree coalgebra $\Gamma(V)$, we are clearly forced to
set $\str^{\dg}_+\Gamma(V) = \Gamma(\str^{\dg}_+ V)$. One checks
that the map $V\mapsto\Gamma(\str^{\dg}_+ V)$ extends appropriately
to a functor on the full subcategory of $\CoAlg^a_+$ formed by
cofree coalgebras. Finally, in the general case, the truncation
$\str^{\dg}_+ K$ of a coalgebra $K$ is obtained by an equalizer of
cofree coalgebra truncations since any coalgebra $K$ is the
equalizer of a natural pair of cofree coalgebra morphisms associated
to $K$.

As mentioned in the introduction, one can also check that the
truncation functor preserves quasi-cofree coalgebras, the
generalization of cofree coalgebras that we define in the next
subsection.

\subsection{Quasi-cofree coalgebras and coextension diagrams}\label{subsection:CocellularCoalgebras} In the
differential graded context, the cofree coalgebra $\Gamma(V)$ is
equipped with a natural differential $\delta:
\Gamma(V)\rightarrow\Gamma(V)$ induced by the internal differential
of $V$. As explained in the introduction of this section, we
consider \emph{quasi-cofree
coalgebra}\index{coalgebra!quasi-cofree}\index{quasi-cofree!coalgebra}
structures defined by a cofree coalgebra $\Gamma = \Gamma(V)$
equipped with a coderivation $\partial:
\Gamma(V)\rightarrow\Gamma(V)$ such that $\delta+\partial$ defines
the differential of $\Gamma$.

The following useful assertion generalizes a classical result for
the tensor coalgebra:

\begin{lemm}\label{lemm:QuasiCofreeCoalgebraStructure}
For any homogeneous morphism $\alpha: \Gamma(V)\rightarrow V$, there
is a unique coderivation $\partial_\alpha:
\Gamma(V)\rightarrow\Gamma(V)$ such that $\alpha =
\pi\partial_\alpha$.\glossary{$\partial_\alpha$}

Assume that the morphism $\alpha: \Gamma(V)\rightarrow V$ is
homogeneous of degree~$-1$. The sum $\delta+\partial_\alpha$ defines
a differential on $\Gamma(V)$ so that $\Gamma =
(\Gamma(V),\partial_\alpha)$ defines a quasi-cofree coalgebra if and
only if we have the relation
\begin{equation*}
\delta(\alpha) + \alpha\partial_\alpha = 0
\end{equation*}
in $\DGHom(\Gamma(V),V)$. Furthermore, the morphism of graded
coalgebras $\nabla_f: K\rightarrow\Gamma(V)$\glossary{$\nabla_f$}
induced by a homogeneous morphism $f: K\rightarrow V$ of degree~$0$
defines a morphism of differential graded coalgebras $\nabla_f:
K\rightarrow\Gamma$, where $\Gamma = (\Gamma(V),\partial_\alpha)$,
if and only if we have the relation
\begin{equation*}
\delta(f) + \alpha\nabla_f = 0
\end{equation*}
in $\DGHom(K,V)$.
\end{lemm}

\begin{proof}
The construction of $\partial_\alpha$ is similar to the construction
of the coalgebra morphism $\nabla_f: K\rightarrow\Gamma(V)$
associated to a morphism of dg-modules $f: K\rightarrow\Gamma(V)$.
In particular, we deduce the existence of $\partial_\alpha$ from our
realization of the cofree coalgebra $\Gamma(V) = \Gamma_\infty(V)$.
Explicitly, we consider the morphism $\widehat{\partial}_\alpha:
\Gamma(V)\rightarrow\prod_n V^{\otimes n}$ which maps an element
$\gamma\in\Gamma(V)$ to the collection
\begin{equation*}
\{\pi^{\otimes i-1}\otimes \alpha\otimes\pi^{\otimes
n-i}\cdot\Delta^n(\gamma)\}.
\end{equation*}
One checks readily that this morphism fits a commutative diagram
\begin{equation*}
\xymatrix{ \Gamma(V)\ar[r]\ar[d]^{\widehat{\partial}_\alpha} &
\Gamma(V)\otimes\Gamma(V)\ar[d]^{\nabla_\Pi\cdot(\widehat{\partial}_\alpha\otimes 1+1\otimes\widehat{\partial}_\alpha)} \\
\prod_N V^{\otimes N}\ar[r]^(0.4){\Delta_\Pi} & \prod_{m,n}
\{V^{\otimes m}\otimes V^{\otimes n}\} }
\end{equation*}
and, as in the proof of lemma \ref{lemm:CoFreeCoAlg}, we deduce from
this property that $\widehat{\partial}_\alpha$ restricts to
morphisms
\begin{equation*}
\xymatrix{
\Gamma(V)\ar[]!D-<0pt,2pt>;[d]!U+<0pt,2pt>^(0.5){\partial_{\infty}}
\ar@/^5pt/[]!R+<2pt,0pt>;[drr]!U+<0pt,2pt>^(0.92){\partial_{r+1}}
\ar@/^5pt/[]!R+<3pt,0pt>;[drrr]!U+<0pt,2pt>^(0.95){\partial_{r}}
\ar@/^5pt/[]!R+<4pt,0pt>;[drrrrr]!UL+<8pt,2pt>^(0.80){\partial_1 = \widehat{\partial}_\alpha} &&&& \\
*{\Gamma_{\infty}(V)}\ar@{^{(}->}[]!R+<4pt,0pt>;[r]!L-<2pt,0pt> &
*{\dots}\ar@{^{(}->}[]!R+<4pt,0pt>;[r]!L-<2pt,0pt> &
*{\Gamma_{r+1}(V)}\ar@{^{(}->}[]!R+<4pt,0pt>;[r]!L-<2pt,0pt> &
*{\Gamma_{r}(V)}\ar@{^{(}->}[]!R+<4pt,0pt>;[r]!L-<2pt,0pt> &
*{\dots}\ar@{^{(}->}[]!R+<4pt,0pt>;[r]!L-<2pt,0pt> &
*{\Gamma_{1}(V) = \prod_n V^{\otimes n}} }
\end{equation*}
such that $\partial_\alpha = \partial_\infty$ defines a coderivation
of $\Gamma(V) = \Gamma_\infty(V)$. In fact, for a coderivation, the
relation $\alpha = \pi\partial_\alpha$ implies that the composite of
$\partial_\alpha$ with the embedding $\iota_\infty:
\Gamma_\infty(V)\hookrightarrow\prod_n V^{\otimes n}$ agrees with
the map $\widehat{\partial}_\alpha$. Therefore the coderivation
$\partial_\alpha$ is uniquely characterized by this relation $\alpha
= \pi\partial_\alpha$. The verification of the other assertions of
the lemma is similar and straighforward.
\end{proof}

As explained in the introduction of this section, we aim to
determine the structure of morphism of quasi-cofree coalgebras
$\nabla_f:
(\Gamma(V),\partial_\alpha)\rightarrow(\Gamma(W),\partial_\beta)$
that fit a pullback diagram of the form
\begin{equation*}
\xymatrix{ (\Gamma(V),\partial_\alpha)\ar[r]\ar[d]
& \Gamma(\Delta^1\wedge V)\ar[d] \\
(\Gamma(W),\partial_\beta)\ar[r] & \Gamma(\Delta^1\wedge
W)\times_{\Gamma(S^1\wedge W)}\Gamma(S^1\wedge V) }.
\end{equation*}
First, we define precisely the dg-modules $\Delta^1\wedge E$ and
$S^1\wedge E$ that occur in this construction. In fact, we consider
nothing but the classical cone and suspension functors in the
category of dg-modules.

\subsubsection{The cone sequence of a dg-module}\label{item:ClassicalCocells}
Explicitly, for a dg-module $E$ (possibly $\Z$-graded), we let
$\Delta^1\wedge E$\glossary{$\Delta^1\wedge E$}, respectively
$S^1\wedge E$\glossary{$S^1\wedge E$}, denote the quotient
$\I\otimes E/\x{0}\otimes E$, respectively $\I\otimes E/\x{0}\otimes
E\oplus\x{1}\otimes E$, of the tensor product of~$E$ with the
classical interval~$\I$ of the category of dg-modules. Recall that
this dg-module~$\I$ is spanned by homogeneous elements
$\x{0},\x{1},\x{01}$ of degree $\deg(\x{0}) = \deg(\x{1}) = 0$ and
$\deg(\x{01}) = 1$ respectively and the differential $\partial:
\I\rightarrow\I$ is defined by $\partial(\x{01}) = \x{1} - \x{0}$
(see~\ref{item:ChainInterval}).

Equivalently, the dg-module $\Delta^1\wedge E$ can be defined by
\begin{equation*}
\Delta^1\wedge E = \x{01}\otimes E\oplus\x{1}\otimes E.
\end{equation*}
The differential of $\Delta^1\wedge E$ can be decomposed into a
natural differential $\delta: \Delta^1\wedge
E\rightarrow\Delta^1\wedge E$ induced by the internal differential
of $E$ and an extra term $\partial: \Delta^1\wedge
E\rightarrow\Delta^1\wedge E$ induced by the differential of $\I$.
By definition, this differential $\partial$ maps a tensor
$\x{01}\otimes x\in\x{01}\otimes E$ to a corresponding element
$\x{1}\otimes x\in\x{1}\otimes E$ and vanishes on the other
component of $\Delta^1\wedge E$. We have similarly
\begin{equation*}
S^1\wedge E = \x{01}\otimes E
\end{equation*}
so that $S^1\wedge E$ can be identified with the suspension of $E$.

Clearly, we have a natural morphism of dg-modules $\sigma\wedge E:
\Delta^1\wedge E\rightarrow S^1\wedge E$\glossary{$\sigma\wedge E$}
which can be identified with the projection onto the component
$\x{01}\otimes E$ of $\Delta^1\wedge E = \x{01}\otimes
E\oplus\x{1}\otimes E$. We have also a natural embedding $d^1:
E\hookrightarrow\Delta^1\wedge E$ which identifies the dg-module $E$
with the component $\x{1}\otimes E$ of~$\Delta^1\wedge E$.

\medskip
We begin our constructions with the following simple observation:

\begin{obsv}
Any quasi-cofree coalgebra $\Gamma = (\Gamma(V),\partial_\alpha)$ is
endowed with a morphism of dg-modules $\pi_\alpha:
\Gamma\rightarrow\Delta^1\wedge V$ such that $\pi_\alpha(\gamma) =
\x{01}\otimes\alpha(\gamma)+\x{1}\otimes\sigma(\gamma)$, for
$\gamma\in\Gamma(V)$.\glossary{$\pi_\alpha$}
\end{obsv}

\begin{proof}
One can observe precisely that the commutation of $\alpha$ with
differentials, given by the commutativity of the square
\begin{equation*}
\xymatrix{
\Gamma(V)\ar[r]^{\pi_{\alpha}}\ar[d]_{\delta+\partial_\alpha} &
\Delta^1\wedge V\ar[d]_{\delta+\partial} \\
\Gamma(V)\ar[r]^{\pi_{\alpha}} & \Delta^1\wedge V },
\end{equation*}
is equivalent to the relation $\delta(\alpha)+\alpha\partial_\alpha
= 0$ of lemma~\ref{lemm:QuasiCofreeCoalgebraStructure}.
\end{proof}

In fact, we have the following general assertion:

\begin{fact}\label{fact:ConeMorphism}
A map $\tilde{e}: U\rightarrow\Delta^1\wedge E$, where $U$ is a
dg-module, defines a dg-module morphism if and only if we have
$\tilde{e}(u) = -\x{01}\otimes\delta(e)(u)+\x{1}\otimes e(u)$, for a
homogeneous map $e: U\rightarrow E$ of degree $0$.
\end{fact}

In the previous observation the differential of $U =
(\Gamma(V),\partial_\alpha)$ includes the coderivation
$\partial_\alpha$. Hence the differential of the homogeneous map
$\pi_\alpha: \Gamma(V)\rightarrow V$ is given by
$\delta\pi_\alpha-\pi_\alpha\delta-\pi_\alpha\partial_\alpha =
0-\alpha = -\alpha$.

These assertions yield also to the following useful observation:

\begin{obsv}\label{obsv:ConeCoalgebraMorphism}
Suppose given a dg-module morphism $\tilde{v}:
K\rightarrow\Delta^1\wedge V$, where $K$ is a dg-coalgebra, and
consider the equivalent homogeneous map $v: K\rightarrow V$. The
coalgebra map $\nabla_v: K\rightarrow\Gamma(V)$ induced by $v:
K\rightarrow V$ defines a morphism to the quasi-cofree coalgebra
$\Gamma = (\Gamma(V),\partial_\alpha)$ if and only if it makes
commute the diagram
\begin{equation*}
\xymatrix{ K\ar@{-->}[d]_{\nabla_v}\ar[dr]^{\tilde{v}} & \\
\Gamma(V)\ar[r]_{\pi_\alpha} & \Delta^1\wedge V }.
\end{equation*}
In addition this assertion holds as soon as the composites
of~$\tilde{v}$ and~$\pi_\alpha\nabla_v$ with~$\sigma\wedge V:
\Delta^1\wedge V\rightarrow S^1\wedge V$ agree.
\end{obsv}

\begin{proof}
These claims are immediate: we have by definition
$\pi_\alpha\nabla_v(x) = \x{01}\otimes\alpha\nabla_v(x) +
\x{1}\otimes v(x)$ so that the relation $\pi_\alpha\nabla_v =
\tilde{v}$ is equivalent to the relation $\alpha\nabla_v = -
\delta(v)$ of lemma~\ref{lemm:QuasiCofreeCoalgebraStructure}.
\end{proof}

\subsubsection{Coalgebra coextensions}\label{item:CofreeCoextensions}
We consider now a morphism of quasi-cofree coalgebras
\begin{equation*}
\nabla_f:
(\Gamma(V),\partial_\alpha)\rightarrow(\Gamma(W),\partial_\beta)
\end{equation*}
induced by a morphism of dg-modules $f: V\rightarrow W$. As this
morphism $f$ is supposed to commute with the internal differentials
of $V$ and $W$, the relation of
lemma~\ref{lemm:QuasiCofreeCoalgebraStructure}, that gives the
commutation of $\nabla_f$ with quasi-cofree coalgebra differentials,
is equivalent to $f\alpha = \beta\nabla_f$.

Assume that the diagram
\begin{equation*}
\xymatrix{
(\Gamma(V),\partial_\beta)\ar[r]^{\pi_{\alpha}}\ar[d]_{\nabla_f} &
\Delta^1\wedge V\ar[r]^{\sigma\wedge V} & S^1\wedge V\ar[d]^{S^1\wedge f} \\
(\Gamma(W),\partial_\beta)\ar[r]_{\pi_{\beta}}\ar@{-->}[urr]^{\rho}
& \Delta^1\wedge W\ar[r]_{\sigma\wedge W} & S^1\wedge W }
\end{equation*}
admits a lifting $\rho: (\Gamma(W),\partial_\beta)\rightarrow
S^1\wedge V$ in the category of dg-modules. Observe that the
composite $\rho_\alpha = \sigma\wedge V\cdot\pi_\alpha$ is given by
$\rho_\alpha(\gamma) = \x{01}\otimes\alpha(\gamma)$, for all
$\gamma\in\Gamma(V)$, and similarly for $\rho_\beta = \sigma\wedge
W\cdot\pi_\beta$. As a consequence, if $f$ is epi, then the
existence of~$\rho$ is equivalent to the existence of a lifting map
$\beta': \Gamma(V)\rightarrow W$ in the diagram
\begin{equation*}
\xymatrix{ \Gamma(V)\ar[r]^{\alpha}\ar[d]_{\nabla_f} & V\ar[d]^{f} \\
\Gamma(W)\ar[r]_{\beta}\ar@{-->}[ur] & W }
\end{equation*}
since one can set $\rho(\gamma) = \x{01}\otimes\beta'(\gamma)$ and
the epimorphism assumption implies that $\rho$ commutes
automatically with differentials.

Anyway, in this situation, we have a diagram of dg-module morphisms
\begin{equation*}
\xymatrix{
(\Gamma(V),\partial_\beta)\ar[r]^{\pi_{\alpha}}\ar[d]_{\nabla_f} &
\Delta^1\wedge V\ar[d]^{(\Delta^1\wedge f,\sigma\wedge V)} \\
(\Gamma(W),\partial_\beta)\ar[r]_(0.38){(\pi_{\beta},\rho)} &
\Delta^1\wedge W\times_{S^1\wedge W} S^1\wedge V }
\end{equation*}
and we can consider the associated diagram of coalgebra morphisms
\begin{equation*}
\xymatrix@C=16mm{
(\Gamma(V),\partial_\alpha)\ar[r]^{\nabla_{\pi_{\alpha}}}\ar[d]_{\nabla_f}
&
\Gamma(\Delta^1\wedge V)\ar[d]^{(\nabla_{\Delta^1\wedge f},\nabla_{\sigma\wedge V})} \\
(\Gamma(W),\partial_\beta)\ar[r]_(0.35){(\nabla_{\pi_{\beta}},\nabla_{\rho})}
& \Gamma(\Delta^1\wedge W)\times_{\Gamma(S^1\wedge W)}
\Gamma(S^1\wedge V) }.
\end{equation*}
We obtain:

\begin{claim}\label{claim:CofreeCoextension}
The diagram above forms a pullback in the category of dg-coalgebras.
\end{claim}

\begin{proof}
We suppose given a pair of coalgebra morphisms $\nabla_{\tilde{v}}:
K\rightarrow\Gamma(\Delta^1\wedge V)$ and $\nabla_w:
K\rightarrow(\Gamma(W),\partial_W)$ that fit the commutative diagram
\begin{equation*}
\xymatrix@C=16mm{
K\ar@{-->}[dr]^{\nabla_v}\ar@/_/[ddr]!L_{\nabla_w}\ar@/^/[rrd]^{\nabla_{\tilde{v}}}
& &
\\ & (\Gamma(V),\partial_\alpha)\ar[r]^{\nabla_{\pi_{\alpha}}}\ar[d]_{\nabla_f} &
\Gamma(\Delta^1\wedge V)\ar[d]^{(\nabla_{\Delta^1\wedge f},\nabla_{\sigma\wedge V})} \\
&
(\Gamma(W),\partial_\beta)\ar[r]_(0.35){(\nabla_{\pi_{\beta}},\nabla_{\rho})}
& \Gamma(\Delta^1\wedge W)\times_{\Gamma(S^1\wedge W)}
\Gamma(S^1\wedge V) }
\end{equation*}
and we check the existence of a unique filling morphism $\nabla_v:
K\rightarrow(\Gamma(V),\partial_\alpha)$ that fits this diagram.
According to the fact~\ref{fact:ConeMorphism}, the coalgebra
morphism $\nabla_{\tilde{v}}: K\rightarrow\Gamma(\Delta^1\wedge V)$
is equivalent to a morphism of dg-modules $\tilde{v}:
K\rightarrow\Delta^1\wedge V$ such that $\tilde{v}(x) =
-\x{01}\otimes\delta(v)(x) + \x{1}\otimes v(x)$ for a homogeneous
map $v: K\rightarrow V$ of degree $0$. We check that this map $v$
induces a morphism of coalgebras $\nabla_v:
K\rightarrow(\Gamma(V),\partial_\alpha)$ that provides a filling
morphism in our diagram.

The commutativity of the diagram
\begin{equation*}
\xymatrix{
K\ar@/_/[ddr]!L_{\nabla_w}\ar@/^2mm/[rrd]^{\nabla_{\tilde{v}}}\ar@/^6mm/[rrrd]^{\tilde{v}}
& & & \\ & (\Gamma(V),\partial_\alpha)\ar[r]\ar[d]^{\nabla_f} &
\Gamma(\Delta^1\wedge V)\ar[d]^{\nabla_{\Delta^1\wedge f}}\ar[r] & \Delta^1\wedge V\ar[d]^{\Delta^1\wedge f} \\
& (\Gamma(W),\partial_\beta)\ar[r]\ar@/_6mm/[rr]_{\pi_\beta} &
\Gamma(\Delta^1\wedge W)\ar[r] & \Delta^1\wedge W }
\end{equation*}
gives $\Delta^1\wedge f\cdot\tilde{v} = \pi_\beta\cdot\nabla_w$ and
from this relation we deduce the identity $f v = w$. As a
consequence, for the induced coalgebra maps, we obtain
$\nabla_f\nabla_v = \nabla_w$. Then the commutativity of the diagram
\begin{equation*}
\xymatrix{
K\ar@/_/[ddr]!L_{\nabla_w}\ar@/^2mm/[rrd]^{\nabla_{\tilde{v}}}\ar@/^6mm/[rrrd]^{\tilde{v}} & & & \\
& (\Gamma(V),\partial_\alpha)\ar[r]\ar[d]^{\nabla_f} &
\Gamma(\Delta^1\wedge V)\ar[d]^{\nabla_{\sigma\wedge V}}\ar[r] & \Delta^1\wedge V\ar[d]^{\sigma\wedge V} \\
& (\Gamma(W),\partial_\beta)\ar[r]\ar@/_6mm/[rr]_{\rho} &
\Gamma(S^1\wedge V)\ar[r] & S^1\wedge V }
\end{equation*}
gives $\sigma\wedge V\cdot\tilde{v} = \rho\cdot\nabla_w =
\rho\cdot\nabla_f\cdot\nabla_v$. By the very definition of $\rho$,
we have $\rho\cdot\nabla_f = \sigma\wedge V\cdot\pi_\alpha$. Hence
we obtain the relation $\sigma\wedge V\cdot\tilde{v} = \sigma\wedge
V\cdot\pi_\alpha\cdot\nabla_v$. By
observation~\ref{obsv:ConeCoalgebraMorphism}, this assertion implies
that $\nabla_v$ defines a coalgebra morphism to the quasi-cofree
coalgebra $(\Gamma(V),\partial_\alpha)$ and we have in addition
$\nabla_{\tilde{v}} = \nabla_{\pi_\alpha}\nabla_v$. As we obtain
also $\nabla_v\nabla_f = \nabla_w$, this observation completes the
proof of the existence of a filling morphism $\nabla_v$. As this
filling morphism is clearly unique, this achieves the proof of
claim~\ref{claim:CofreeCoextension}.
\end{proof}

Observe that we have the identity
\begin{equation*}
\Gamma(\Delta^1\wedge W)\times_{\Gamma(S^1\wedge W)}
\Gamma(S^1\wedge V) = \Gamma(\Delta^1\wedge W\times_{S^1\wedge W}
S^1\wedge V)
\end{equation*}
in the category of coalgebras. Moreover, the morphism
$(\nabla_{\Delta^1\wedge f},\nabla_{\sigma\wedge V})$ involved in
the construction of~\ref{item:CofreeCoextensions} can be identified
with the morphism of cofree coalgebras induced by the dg-module
morphism
\begin{equation*}
(\Delta^1\wedge f,\sigma\wedge V): \Delta^1\wedge
V\rightarrow\Delta^1\wedge W\times_{S^1\wedge W} S^1\wedge V.
\end{equation*}
For our purposes we give a more explicit form to this morphism.
Namely we have the following simple assertion:

\begin{obsv}\label{obsv:ConeSuspensionFiberedProduct}
The dg-module $\Delta^1\wedge W\times_{S^1\wedge W} S^1\wedge V$ can
be identified with the direct sum $\x{01}\otimes V\oplus\x{1}\otimes
W$ equipped with a differential given by the sum of the internal
differentials of~$V$ and~$W$ with an extra term $\partial$ such that
$\partial(\x{01}\otimes v) = \x{1}\otimes v$, for all $v\in V$.

The morphism $(\Delta^1\wedge f,\sigma\wedge V): \Delta^1\wedge
V\rightarrow\Delta^1\wedge W\times_{S^1\wedge W} S^1\wedge V$ is
also identified with the direct sum
\begin{equation*}
\x{01}\otimes V\oplus\x{1}\otimes V\xrightarrow{(\x{01}\otimes\Id,
\x{1}\otimes f)}\x{01}\otimes V\oplus\x{1}\otimes W
\end{equation*}
where we consider the identity map $\x{01}\otimes\Id: \x{01}\otimes
V\rightarrow\x{01}\otimes V$ and the map $\x{1}\otimes f:
\x{1}\otimes V\rightarrow\x{1}\otimes W$ on the components
of~$\Delta^1\wedge V$.\qed
\end{obsv}

From this observation we obtain immediately:

\begin{claim}
If $f: V\rightarrow W$ is a fibration, respectively an acyclic
fibration, of dg-modules, then so is $(\Delta^1\wedge f,\sigma\wedge
V): \Delta^1\wedge V\rightarrow\Delta^1\wedge W\times_{S^1\wedge W}
S^1\wedge V$.\qed
\end{claim}

By adjunction, the morphism of cofree coalgebras
\begin{equation*}
\nabla_{(\Delta^1\wedge f,\sigma\wedge V)}: \Gamma(\Delta^1\wedge
V)\rightarrow\Gamma(\Delta^1\wedge W\times_{S^1\wedge W} S^1\wedge
V)
\end{equation*}
associated to $(\Delta^1\wedge f,\sigma\wedge V)$ defines a
fibration, respectively an acyclic fibration, in the category of
dg-coalgebras if $f$ is so. As a corollary, by standard arguments,
we obtain:

\begin{lemm}\label{lemm:CofreeCoextensionFibrations}
In the situation of~\ref{item:CofreeCoextensions}, if $f$ is a
fibration, respectively an acyclic fibration in the category of
dg-modules, then $\nabla_f$ defines a fibration, respectively an
acyclic fibration in the category of dg-coalgebras.\qed
\end{lemm}

\subsubsection{Functorial coextensions}\label{item:DiagramCofreeCoextensions}
For our purposes we need a relative version of the result of
lemma~\ref{lemm:CofreeCoextensionFibrations}. To be precise, we
consider now a commutative square of quasi-cofree coalgebra
morphisms
\begin{equation*}
\xymatrix{
(\Gamma(V),\partial_{\alpha})\ar[d]_{\nabla_{f}}\ar[r]^{\nabla_v} &
(\Gamma(V'),\partial_{\alpha'})\ar[d]_{\nabla_{f'}} \\
(\Gamma(W),\partial_{\beta})\ar[r]^{\nabla_w} &
(\Gamma(W'),\partial_{\beta'}) }
\end{equation*}
yielded by a commutative square of dg-module morphisms
\begin{equation*}
\xymatrix{ V\ar[d]_{f}\ar[r]^{v} &
V'\ar[d]_{f'} \\
W\ar[r]^{w} & W' }.
\end{equation*}
In~\ref{item:CofreeCoextensions} we observed that $\nabla_f$
satisfies the commutation relation $\beta\nabla_f = f\alpha$ and so
does $\nabla_{f'}$ since $f$, respectively $f'$, is supposed to
commute with internal differentials of dg-modules. Similarly, as
$\nabla_v$ and $\nabla_w$ are supposed to be morphisms of
dg-coalgebras induced by morphisms of dg-modules $v$ and $w$, we
obtain a commutative cube:
\begin{equation*}
\xymatrix@!C=4mm@!R=4mm{
\Gamma(V)\ar[rr]^(0.4){\alpha}\ar[dd]\ar[dr] && V\ar'[d][dd]\ar[dr] & \\
& \Gamma(V')\ar[rr]^(0.4){\alpha'}\ar[dd] && V'\ar[dd] \\
\Gamma(W)\ar'[r]^(0.8){\beta}[rr]\ar[dr] && W\ar[dr] & \\
& \Gamma(W')\ar[rr]^(0.4){\beta'} && W' }
\end{equation*}

As in~\ref{item:CofreeCoextensions}, we assume the existence of
liftings $\rho$ and $\rho'$ of the coderivation maps. Naturally we
assume in addition that these liftings commute with the morphisms
$\nabla_v$ and $\nabla_w$ and hence fit the commutative cube above.
Clearly, if $f$ and $g$ are epi, then this functoriality requirement
is automatically satisfied. Anyway, in this situation, the maps $v$
and $w$ yield a morphism between the coalgebra pullbacks associated
to $\nabla_f$ and $\nabla_{f'}$. Explicitly, these coalgebra
pullbacks form the back and front square of a commutative cube. By
pulling back the front corners of this cube, we obtain a commutative
square
\begin{equation*}
\xymatrix{
(\Gamma(V),\partial_\alpha)\ar[d]_{(\nabla_f,\nabla_v)}\ar[r] &
\Gamma(\Delta^1\wedge V)\ar[d] \\
(\Gamma(W),\partial_\beta)\times_{(\Gamma(W'),\partial_{\beta'})}
(\Gamma(V'),\partial_{\alpha'}) \ar[r] & \Gamma(\Delta^1\wedge
W\times_{?} S^1\wedge V\times_{?}\Delta^1\wedge V') }
\end{equation*}
in which the bottom right hand-side corner can be defined precisely
by the limit of the corner diagram of cofree coalgebras
\begin{equation*}
\xymatrix@!C=8mm@!R=8mm{ && \Gamma(S^1\wedge V)\ar'[d][dd]\ar[dr] & \\
& \Gamma(\Delta^1\wedge V')\ar[dd]\ar[rr] && \Gamma(S^1\wedge V')\ar[dd] \\
\Gamma(\Delta^1\wedge W)\ar'[r][rr]\ar[dr] && \Gamma(S^1\wedge W)\ar[dr] & \\
& \Gamma(\Delta^1\wedge W')\ar[rr] && \Gamma(S^1\wedge W') }.
\end{equation*}

\begin{obsv}
The cartesian product
$(\Gamma(W),\partial_\beta)\times_{(\Gamma(W'),\partial_{\beta'})}
(\Gamma(V'),\partial_{\alpha'})$ can be identified with the
quasi-cofree coalgebra $(\Gamma(X),\partial_\gamma)$ such that $X =
W\times_{W'} V'$ and where the coderivation $\partial_\gamma$ is
induced by the map $(\beta,\alpha'): \Gamma(W\times_{W'}
V')\rightarrow W\times_{W'} V'$.

Furthermore, our commutative square can be identified with the
pullback diagram
\begin{equation*}
\xymatrix@C=16mm{
(\Gamma(V),\partial_\alpha)\ar[r]^{\nabla_{\pi_{\alpha}}}\ar[d]_{\nabla_{(f,v)}}
&
\Gamma(\Delta^1\wedge V)\ar[d]^{(\nabla_{\Delta^1\wedge(f,v)},\nabla_{\sigma\wedge V})} \\
(\Gamma(X),\partial_\gamma)\ar[r]_(0.35){(\nabla_{\pi_{\gamma}},\nabla_{\sigma})}
& \Gamma(\Delta^1\wedge X)\times_{\Gamma(S^1\wedge X)}
\Gamma(S^1\wedge V) }
\end{equation*}
associated to the morphism $(f,v): V\rightarrow W\times_{W'} V'$ and
where $\sigma: (\Gamma(X),\partial_\gamma)\rightarrow S^1\wedge V$
is given by the composite
\begin{equation*}
(\Gamma(X),\partial_\gamma)\xrightarrow{\pr_1}(\Gamma(W),\partial_\beta)
\xrightarrow{\rho}S^1\wedge V.
\end{equation*}
\end{obsv}

\begin{proof}
This assertion follows from straightforward verifications. Observe
simply that we have $\Delta^1\wedge(W\times_{W'} V') =
\Delta^1\wedge W\times_{\Delta^1\wedge W'}\Delta^1\wedge V'$ in view
of our definition of the functor $\Delta^1\wedge-$ and similarly for
$S^1\wedge(W\times_{W'} V')$.
\end{proof}

\begin{lemm}\label{lemm:DiagramCofreeCoextensionFibrations}
In the situation of~\ref{item:DiagramCofreeCoextensions}, if the
morphism
\begin{equation*}
\nabla_w:
(\Gamma(W),\partial_\beta)\rightarrow(\Gamma(W'),\partial_{\beta'})
\end{equation*}
is a fibration of dg-coalgebras and $(f,v): V\rightarrow
W\times_{W'} V'$ is a fibration of dg-modules, then
\begin{equation*}
\nabla_v:
(\Gamma(V),\partial_\alpha)\rightarrow(\Gamma(V'),\partial_{\alpha'})
\end{equation*}
is fibration of dg-coalgebras as well.

If we assume furthermore that $\nabla_w$ and $(f,v)$ are
weak-equivalences so that $\nabla_w$ forms an acyclic fibration of
dg-coalgebras and $(f,v)$ forms an acyclic fibration of dg-modules,
then $\nabla_v$ forms an acyclic fibration of dg-coalgebras as well.
\end{lemm}

For $(f,v)$, one can observe that the weak-equivalence property is
satisfied as soon as we assume that $w$ forms an acyclic fibration
and $v$ is a weak-equivalence. In fact, consider the diagram
\begin{equation*}
\xymatrix{ V\ar@/_12mm/[dd]^(0.75){\sim}_(0.75){v}\ar[d]^{(f,v)} & \\
W\times_{W'} V'\ar[r]\ar[d] & W\ar@{->>}[d]^{\sim}_{w} \\
V'\ar[r] & W' }.
\end{equation*}
The pullback of $w: W\rightarrow W'$ is an acyclic fibration by
standard model category arguments. If we assume that $v$ is a
weak-equivalence, then, by the two out of three axiom, we conclude
that $(f,v)$ is a weak-equivalence as well, as claimed.

\begin{proof}
According to lemma~\ref{lemm:CofreeCoextensionFibrations}, the
morphism $\nabla_{(f,v)}$ forms a fibration, respectively an acyclic
fibration, if $(f,v)$ is so. On the other hand, by the general model
category argument recalled just above this proof, in the pullback
diagram
\begin{equation*}
\xymatrix{ (\Gamma(X),\partial_\gamma)\ar[r]\ar[d] & (\Gamma(W),\partial_\beta)\ar[d]^{\nabla_w} \\
(\Gamma(V'),\partial_{\alpha'})\ar[r] &
(\Gamma(W'),\partial_{\beta'}) }
\end{equation*}
the left-hand side morphism forms a fibration, respectively an
acyclic fibration, if $\nabla_w$ is so. As a conclusion, under the
assumptions of the lemma, the morphism $\nabla_v$ can be decomposed
into a composite of fibrations, respectively acyclic fibrations, and
hence forms a fibration, respectively an acyclic fibration, as well.
\end{proof}

\subsection{Quasi-cofree Hopf $\Lambda_*$-modules}\label{subsection:CocellularHopfLambdaModules}
The purpose of this subsection is to extend the constructions of the
previous section to the category of Hopf $\Lambda_*$-modules. First,
observe that the category of Hopf $\Lambda_*$-modules is endowed
with a cofree object functor $\Gamma:
\Lambda_*^{\op}\Mod\rightarrow\Lambda_*^{\op}\HopfMod$. In fact, the
\emph{cofree Hopf $\Lambda_*$-module}\index{cofree!Hopf
$\Lambda_*$-module}\index{Hopf $\Lambda_*$-module!cofree}
$\Gamma(M)$\glossary{$\Gamma(M)$} associated to a $\Lambda_*$-module
$M$ can be defined explicitly by the cofree coalgebras
\begin{equation*}
\Gamma(M)(r) = \Gamma(M(r))
\end{equation*}
on the dg-modules $M(r)$. Furthermore, the operations $\partial_i:
\Gamma(M)(r)\rightarrow\Gamma(M)(r-1)$ can be identified with the
morphisms of cofree coalgebras $\nabla_{\partial_i}:
\Gamma(M(r))\rightarrow\Gamma(M(r-1))$ induced by the operations
$\partial_i: M(r)\rightarrow M(r-1)$ on $M(r)$.

For our purposes we extend the notion of a quasi-cofree object to
Hopf $\Lambda_*$-modules. Then we observe that the categorical
results of the previous section hold as well for Hopf
$\Lambda_*$-modules. Finally, we obtain useful sufficient condition
for a map to be a fibrations and acyclic fibrations in that
category.

\subsubsection{Quasi-cofree Hopf $\Lambda_*$-modules}\label{item:QuasiCofreeHopfLambdaModules}
As usual, in the differential graded framework, the cofree Hopf
$\Lambda_*$-module $\Gamma(M)$ is equipped with a natural
differential $\delta: \Gamma(M)\rightarrow\Gamma(M)$ induced by the
internal differential of~$M$.

The structure of a \emph{quasi-cofree Hopf
$\Lambda_*$-module}\index{quasi-cofree!Hopf
$\Lambda_*$-module}\index{Hopf $\Lambda_*$-module!quasi-cofree} is
defined by a pair
\begin{equation*}
\Gamma = (\Gamma(M),\partial_\alpha)
\end{equation*}
where $M$ is a $\Lambda_*$-module, $\Gamma(M)$ denotes the
associated cofree Hopf $\Lambda_*$-module and $\partial_\alpha:
\Gamma(M)\rightarrow\Gamma(M)$ denotes a homogeneous morphism of
$\Lambda_*$-modules which forms a coderivation with respect to the
coalgebra structure so that the sum $\delta+\partial$ defines the
differential of the Hopf $\Lambda_*$-module $\Gamma$.

The coderivation $\partial_\alpha:
\Gamma(M)\rightarrow\Gamma(M)$\glossary{$\partial_\alpha$} is
defined by a collection of coderivations of cofree coalgebras
$\partial_\alpha: \Gamma(M(r))\rightarrow\Gamma(M(r))$ which in turn
can be specified by a collection of homogeneous maps $\alpha:
\Gamma(M(r))\rightarrow M(r)$. As specified in the definition, the
coderivation $\partial_\alpha$ is supposed to define a homogeneous
morphism of $\Lambda_*$-modules $\partial_\alpha:
\Gamma(M)\rightarrow\Gamma(M)$. Thus we assume explicitly that
$\partial_\alpha$ commutes with the action of $\Sigma_r$ and with
the operations $\partial_i: \Gamma(M(r))\rightarrow\Gamma(M(r-1))$.
One can assume equivalently that the permutations $w\in\Sigma_r$
induce morphisms of quasi-cofree coalgebras $w:
(\Gamma(M(r)),\partial_\alpha)\rightarrow(\Gamma(M(r)),\partial_\alpha)$
and similarly for the operations $\partial_i: M(r)\rightarrow
M(r-1)$.

Observe that the matching objects of a quasi-cofree Hopf
$\Lambda_*$-module $\Gamma =
(\Gamma(M),\partial_\alpha)$\index{matching object!of a quasi-cofree
Hopf $\Lambda_*$-module} can be defined by the quasi-cofree
coalgebras
\begin{equation*}
\Match{\Gamma}(r) = (\Gamma(\Match{M}(r)),\partial_\alpha),
\end{equation*}
where $\partial_\alpha$ is induced by the homogeneous map $\alpha:
\Gamma(\Match{M}(r))\rightarrow\Match{M}(r)$ that fits the
commutative diagrams
\begin{equation*}
\xymatrix{
\Gamma(\Match{M}(r))\ar[r]^{\alpha}\ar[d]_{\nabla_{\partial_i}} &
\Match{M}(r)\ar[d]^{\partial_i} \\
\Gamma(\Match{M}(r-1))\ar[r]^{\alpha} & \Match{M}(r) },
\end{equation*}
for $i = 1,\dots,r$. Similarly, if a morphism of quasi-cofree Hopf
$\Lambda_*$-modules $\nabla_f: \Gamma\rightarrow\Gamma'$ where
$\Gamma = (\Gamma(M),\partial_\alpha)$ and $\Gamma' =
(\Gamma(M'),\partial_{\alpha'})$ is induced by a morphism of
$\Lambda_*$-modules $f: M\rightarrow M'$, then we have
\begin{equation*}
\Match{\Gamma}(r)\times_{\Match{\Gamma'}(r)}\Gamma'(r) =
(\Gamma(\Match{M}(r)\times_{\Match{M'}(r)}
M'(r)),\partial_{(\alpha,\alpha')}).
\end{equation*}
Moreover, the morphism $(\mu,\nabla_f):
\Gamma(r)\rightarrow\Match{\Gamma}(r)\times_{\Match{\Gamma'}(r)}\Gamma'(r)$
can be identified with the morphism of quasi-cofree coalgebras
induced by the morphism of dg-modules $(\mu,f):
M(r)\rightarrow\Match{M}(r)\times_{\Match{M'}(r)} M'(r)$.

\subsubsection{Coalgebra coextensions}\label{item:CofreeLambdaModuleCoextensions}
As in the case of coalgebras, we consider a morphism of quasi-cofree
Hopf $\Lambda_*$-modules
\begin{equation*}
\nabla_f:
(\Gamma(M),\partial_\alpha)\rightarrow(\Gamma(N),\partial_\beta)\glossary{$\nabla_f$}
\end{equation*}
induced by a morphism of dg-modules $f: M\rightarrow N$.
Furthermore, we assume the existence of a lifting $\rho:
(\Gamma(M),\partial_\beta)\rightarrow S^1\wedge N$ in a diagram of
the form
\begin{equation*}
\xymatrix{
(\Gamma(M),\partial_\beta)\ar[r]^{\pi_{\alpha}}\ar[d]_{\nabla_f} &
\Delta^1\wedge M\ar[r]^{\sigma\wedge M} & S^1\wedge N\ar[d]^{S^1\wedge f} \\
(\Gamma(M),\partial_\beta)\ar[r]_{\pi_{\beta}}\ar@{-->}[urr]^{\rho}
& \Delta^1\wedge N\ar[r]_{\sigma\wedge N} & S^1\wedge N },
\end{equation*}
where the $\Lambda_*$-modules $\Delta^1\wedge E$ are defined by
$(\Delta^1\wedge E)(r) = \Delta^1\wedge
E(r)$\glossary{$\Delta^1\wedge E$} and similarly for $S^1\wedge
E$\glossary{$S^1\wedge E$}. To be precise, we assume that $\rho$
defines a lifting in the category of dg-modules and hence commutes
with the action of permutations and with the operations
$\partial_i$.

Clearly, if $f$ is epi, then the existence of~$\rho$ is equivalent
to the existence of a lifting map $\beta': \Gamma(M(r))\rightarrow
N(r)$ in the diagrams
\begin{equation*}
\xymatrix{ \Gamma(M(r))\ar[r]^{\alpha}\ar[d]_{\nabla_f} & M(r)\ar[d]^{f} \\
\Gamma(N(r))\ar[r]_{\beta}\ar@{-->}[ur] & N(r) },
\end{equation*}
for $r\in\N$. The epimorphism assumption implies that $\rho$
commutes automatically with differentials and with
$\Lambda_*$-module operations.

Anyway, as in the case of coalgebras, we have a diagram of Hopf
$\Lambda_*$-modules
\begin{equation*}
\xymatrix@C=16mm{
(\Gamma(M),\partial_\alpha)\ar[r]^{\nabla_{\pi_{\alpha}}}\ar[d]_{\nabla_f}
&
\Gamma(\Delta^1\wedge M)\ar[d]^{(\nabla_{\Delta^1\wedge f},\nabla_{\sigma\wedge M})} \\
(\Gamma(N),\partial_\beta)\ar[r]_(0.35){(\nabla_{\pi_{\beta}},\nabla_{\rho})}
& \Gamma(\Delta^1\wedge N)\times_{\Gamma(S^1\wedge N)}
\Gamma(S^1\wedge M) }
\end{equation*}
and, furthermore, we obtain:

\begin{claim}
The diagram above forms a pullback in the category of Hopf
$\Lambda_*$-modules.\qed
\end{claim}

Then:

\begin{claim}
If $f: M\rightarrow N$ is a Reedy fibration, respectively an acyclic
Reedy fibration, of $\Lambda_*$-modules, then so is $(\Delta^1\wedge
f,\sigma\wedge M): \Delta^1\wedge M\rightarrow\Delta^1\wedge
N\times_{S^1\wedge N} S^1\wedge M$.
\end{claim}

\begin{proof}
This assertion can be deduced readily from the componentwise
expansion of~$(\Delta^1\wedge f,\sigma\wedge M)$ given in
observation~\ref{obsv:ConeSuspensionFiberedProduct}. Observe simply
that the matching functor commutes with the cone construction
$\Delta^1\wedge-$ so that
\begin{equation*}
\Match{(\Delta^1\wedge M)}(r) = \Delta^1\wedge\Match{M}(r) =
\x{01}\otimes\Match{M}(r)\oplus\x{1}\otimes\Match{M}(r)
\end{equation*}
and similarly for~$\Match{(S^1\wedge M)}$.
\end{proof}

And, by the usual categorical arguments, we obtain:

\begin{lemm}
In the situation of~\ref{item:CofreeLambdaModuleCoextensions}, if
$f$ is a Reedy fibration, respectively an acyclic Reedy fibration,
in the category of $\Lambda_*$-modules, then $\nabla_f$ defines a
Reedy fibration, respectively an acyclic Reedy fibration, in the
category of Hopf $\Lambda_*$-modules.\qed
\end{lemm}

\subsubsection{Functorial coextensions}\label{item:DiagramCofreeLambdaModuleCoextensions}
Finally, we give also a relative version of this result in the
context of Hopf $\Lambda_*$-modules. We consider a commutative
square of morphisms of quasi-cofree Hopf $\Lambda_*$-modules
\begin{equation*}
\xymatrix{
(\Gamma(M),\partial_{\alpha})\ar[d]_{\nabla_{f}}\ar[r]^{\nabla_v} &
(\Gamma(M'),\partial_{\alpha'})\ar[d]_{\nabla_{f'}} \\
(\Gamma(N),\partial_{\beta})\ar[r]^{\nabla_w} &
(\Gamma(N'),\partial_{\beta'}) }
\end{equation*}
yielded by a commutative square of $\Lambda_*$-module morphisms
\begin{equation*}
\xymatrix{ M\ar[d]_{f}\ar[r]^{v} &
M'\ar[d]_{f'} \\
N\ar[r]^{w} & N' }.
\end{equation*}
As in the case of coalgebras, we assume the existence of functorial
liftings $\rho$ and $\rho'$ of the coderivation maps. Naturally we
still assume that $\rho$ and $\rho'$ form $\Lambda_*$-module
morphisms. Furthermore, if $f$ and $g$ are epi, then these
functoriality requirements are automatically satisfied.

\begin{lemm}\label{lemm:DiagramCofreeCoextensionLambdaModuleFibrations}
In the situation
of~\ref{item:DiagramCofreeLambdaModuleCoextensions}, if the morphism
\begin{equation*}
\nabla_w:
(\Gamma(N),\partial_\beta)\rightarrow(\Gamma(N'),\partial_{\beta'})
\end{equation*}
is a Reedy fibration of Hopf $\Lambda_*$-modules and $(f,v):
M\rightarrow N\times_{N'} M'$ is a Reedy fibration of
$\Lambda_*$-modules, then
\begin{equation*}
\nabla_v:
(\Gamma(M),\partial_\alpha)\rightarrow(\Gamma(M'),\partial_{\alpha'})
\end{equation*}
is also a Reedy fibration of Hopf $\Lambda_*$-modules.

If we assume furthermore that $\nabla_w$ and $(f,v)$ are
weak-equivalences so that $\nabla_w$ forms an acyclic Reedy
fibration of Hopf $\Lambda_*$-modules and $(f,v)$ forms an acyclic
Reedy fibration of $\Lambda_*$-modules, then $\nabla_v$ forms an
acyclic Reedy fibration of Hopf $\Lambda_*$-modules as well.\qed
\end{lemm}

As in the case of coalgebras, one can observe that the
weak-equivalence property for $(f,v)$ is satisfied as soon as we
assume that $w$ forms an acyclic Reedy fibration of
$\Lambda_*$-modules and $v$ is a weak-equivalence.

\begin{proof}
One can prove this lemma by a tedious but straightforward
generalization of the arguments involved in the case of coalgebras.
On the other hand, by definition, the morphism $\nabla_v:
(\Gamma(M),\partial_\alpha)\rightarrow(\Gamma(M'),\partial_{\alpha'})$
forms a Reedy fibration, respectively an acyclic Reedy fibration, in
the category of Hopf $\Lambda_*$-modules if and only if the morphism
\begin{equation*}
\nabla_{(\mu,v)}:
(\Gamma(M(r)),\partial_\alpha)\rightarrow(\Gamma(\Match{M}(r)\times_{\Match{M'}(r)}
M'(r)),\partial_{(\alpha,\alpha')})
\end{equation*}
forms a fibration in the category of dg-coalgebras. Accordingly, our
claims can also be deduced from the statement of
lemma~\ref{lemm:DiagramCofreeCoextensionFibrations} applied to the
coalgebra morphisms
\begin{equation*}
\xymatrix{
(\Gamma(M(r)),\partial_{\alpha})\ar[d]_{\nabla_{f}}\ar[r]^(0.3){\nabla_{(\mu,v)}}
& (\Gamma(\Match{M}(r)\times_{\Match{M'}(r)}
M'(r)),\partial_{(\alpha,\alpha')})\ar[d]_{\nabla_{f'}} \\
(\Gamma(N(r)),\partial_{\beta})\ar[r]^(0.3){\nabla_{(\mu,w)}} &
(\Gamma(\Match{N}(r)\times_{\Match{N'}(r)}
N'(r)),\partial_{(\beta,\beta')}) }.
\end{equation*}
This observation achieves the proof of
lemma~\ref{lemm:DiagramCofreeCoextensionLambdaModuleFibrations}.\end{proof}

\section{Hopf endomorphism operads}\label{section:OperadActionConstruction}

\subsection{Introduction and sketch of the section results}
The main results of this memoir are proved in this section. For that
reason we give in this introduction a detailed summary of the
content of this section. Thus we give precisions on the memoir
objectives recalled from the memoir introduction.

As in~\cite{Bar1}, we let $\K$\glossary{$\K$} denote the chain
$A_\infty$-operad of Stasheff's associahedra. We consider operads
$\P$ equipped with a fixed operad morphism $\K\rightarrow\P$ so that
we can extend the bar construction $A\mapsto B(A)$ to the category
of $\P$-algebras. To be more precise, as explained in the
introduction of~\ref{section:UnitaryHopfOperads}, it is natural to
assume that $A$ is a non-unitary algebra and hence to consider the
non-unital version of the $A_\infty$-operad $\K$ and a non-unital
operad~$\P$ under~$\K$. Otherwise we have to replace the algebra $A$
by its augmentation ideal $\overline{A}$ and the operad $\P$ by the
associated reduced operad $\overline{\P}$.

Recall that $B(A) = T^c(\Sigma A)$ the \emph{tensor coalgebra}
cogenerated by the suspended dg-module $\Sigma A$. The differential
of $B(A)$ is the sum of the internal differential of $A$ with an
additional component $\partial: T^c(\Sigma A)\rightarrow T^c(\Sigma
A)$ defined by a coalgebra coderivation such that
\begin{equation*}
\partial(a_1\otimes\dots\otimes a_n)
= \sum_{r=2}^n\sum_{k=1}^{n-r+1} \pm
a_1\otimes\dots\otimes\mu_r(a_k,\dots,a_{k+r-1})\otimes\dots\otimes
a_n,
\end{equation*}
where the $\mu_r$'s are the standard generators of Stasheff's operad
$\K$. These operations can abusively be identified with their image
under the operad morphism $\K\rightarrow\P$. The bar complex $B(A)$
is equipped with a canonical unit $\F\rightarrow B(A)$ and forms a
connected augmented unitary coalgebra.

\subsubsection*{The Hopf endomorphism operad of the bar complex}
The first goal of this section is to extend the endomorphism operad
construction of~\cite{Bar1} in the context of Hopf operads. More
explicitly, we prove the following theorem:

\begin{thm}\label{thm:BarHopfEndomorphismOperadRecall}
Let $\P$ denote a (non-unital) operad in dg-modules equipped with an
operad morphism $\K\rightarrow\P$, where $\K$ denotes Staheff's
chain operad. There is a universal unital Hopf operad $\Q =
\HopfEnd_B^{\P}$ such that the bar complex of a $\P$-algebra $B(A)$
is equipped with the structure of a Hopf algebra over $\Q$,
functorially in $A\in\P\Alg$.

More precisely, the Hopf operad $\HopfEnd_B^{\P}$ operates
functorially on the coalgebra $B(A)$ and so that the unital
operation $*: \F\rightarrow B(A)$ agrees with the unit of $B(A)$.
Furthermore we have a one-to-one correspondance between such Hopf
operad actions and morphisms of unital Hopf operads $\rho:
\Q\rightarrow\HopfEnd_B^{\P}$.
\end{thm}

The proof of this statement relies on the existence of internal
hom-objects in the category of coalgebras. To be precise, as we
consider on one hand the bar complex, which forms an augmented
unitary coalgebra, and on the other hand the underlying coalgebras
of a unital Hopf operad, which are supposed to form augmented but
non-unitary coalgebras, we have to adapt our definitions.
Explicitly, to any pair of augmented unitary coalgebras $L$ and $M$,
we associate an augmented coalgebra $\HopfHom(L,M)$ such that a
morphism $\phi^{\sharp}: K\rightarrow\HopfHom(L,M)$ is equivalent to
a morphism of augmented coalgebras $\phi: K\otimes L\rightarrow M$
that makes the diagram
\begin{equation*}
\xymatrix{ K\otimes\F\ar[d]_{K\otimes\eta}\ar[r]^{\epsilon} &
\F\ar[d]_{\eta} \\ K\otimes
  L\ar[r]^{\phi} & M }
\end{equation*}
commute. Then, for any coalgebra $\Gamma$, the objects
$\HopfEnd_\Gamma(r) = \HopfHom(\Gamma^{\otimes r},\Gamma)$ form
clearly a unital Hopf operad and a morphism of unital Hopf operads
$\rho: \Q\rightarrow\HopfEnd_\Gamma$ is equivalent to a Hopf operad
action of $\Q$ on $\Gamma$ such that the unital operation $*:
\F\rightarrow\Gamma$ agrees with the internal unit of the
coalgebra~$\Gamma$.

Formally, the endomorphism operad $\HopfEnd_B^{\P}$ is defined by
the end
\begin{equation*}
\HopfEnd_B^{\P}(r) = \int^{A\in\P\Alg}\HopfHom(B(A)^{\otimes
r},B(A)),
\end{equation*}
also denoted by
\begin{equation*}
\HopfEnd_B^{\P}(r) = \HopfHom_{A\in\P\Alg}(B(A)^{\otimes r},B(A)),
\end{equation*}
and where $A$ ranges over the category of $\P$-algebras. Roughly,
the idea is to replace morphisms by natural transformations. The end
construction permits simply to extend this process to internal
hom-objects. According to this construction, we have a morphism
\begin{equation*}
\HopfEnd_B^{\P}(r)\rightarrow\HopfHom(B(A)^{\otimes r},B(A)),
\end{equation*}
for all $A\in\P\Alg$, so that $\HopfEnd_B^{\P}$ operates on the bar
complex $B(A)$ functorially in $A\in\P\Alg$. Thus the existence
assertion of theorem~\ref{thm:BarHopfEndomorphismOperadRecall} is a
formal consequence of categorical properties of the category of
coalgebras.

\subsubsection*{Fibration properties and contruction of operad actions}
The second goal of this section is to give more insights into the
structure of $\HopfEnd_B^{\P}$ so that we can prove the existence
and uniqueness results stated in the introduction of this memoir.
Let us outline our constructions more precisely.

First, we observe that $\HopfEnd_B^{\P}$ forms a quasi-cofree object
in the category of Hopf $\Lambda_*$-modules. Recall that a
quasi-cofree Hopf $\Lambda_*$-module $\Gamma = (\Gamma(M),\partial)$
is defined explicitly by a collection of cofree coalgebras
$\Gamma(M(r))$ associated to a $\Lambda_*$-module $M$ and equipped
with coderivations $\partial: \Gamma(M(r))\rightarrow\Gamma(M(r))$,
which preserve the $\Lambda_*$-module structure, such that the
differential of $\Gamma$ is given by the sum $\delta+\partial$,
where $\delta$ is the natural differential of the cofree coalgebra
induced by the internal differential of~$M$. For the Hopf
endomorphism operad $\HopfEnd_B^{\P}$, we obtain precisely
\begin{equation*}
\HopfEnd_B^{\P} = (\Gamma(\PrimEnd_B^{\P}),\partial)
\end{equation*}
for the $\Lambda_*$-module
\begin{equation*}
\PrimEnd_B^{\P}(r) = \DGHom_{A\in\P\Alg} \Hom(T^c(\Sigma A)^{\otimes
  r},\Sigma A)
\end{equation*}
formed by all homogeneous morphisms $\theta_A: T^c(\Sigma
A)^{\otimes r}\rightarrow\Sigma A$ which are functorial in $A$.

Then, as in~\cite{Bar1}, we introduce a smaller operad
$\HopfOp_B^{\P}$, the \emph{Hopf operad of universal bar
operations}, that behaves better than the endomorphism operad
$\HopfEnd_B^{\P}$. As in the context of dg-modules, this operad is
equipped with a split injective morphism $\nabla_\Theta:
\HopfOp_B^{\P}\hookrightarrow\HopfEnd_B^{\P}$ which becomes an
isomorphism if the operad $\P$ is $\Sigma_*$-cofibrant or if the
ground field $\F$ is infinite. Explicitly, one observes that the
homogeneous morphisms $p_A: \Sigma A^{\otimes
n_1}\otimes\dots\otimes\Sigma A^{\otimes n_r}\rightarrow\Sigma A$
associated to operations $p\in\P(n_1+\dots+n_r)$ span a submodule
$\PrimOp_B^{\P}(r)$ of $\PrimEnd_B^{\P}(r)$. The operad
$\HopfOp_B^{\P}$ consists of quasi-cofree subcoalgebras of
$\HopfEnd_B^{\P}(r)$ such that
\begin{equation*}
\HopfOp_B^{\P}(r) = (\Gamma(\PrimOp_B^{\P}(r)),\partial).
\end{equation*}
One checks also that $\PrimOp_B^{\P}$ forms also a
$\Lambda_*$-submodule of $\PrimEnd_B^{\P}$ so that $\HopfOp_B^{\P}$
forms a quasi-cofree Hopf $\Lambda_*$-module.

The endomorphism operad $\HopfEnd_B^{\P}$, as well as
$\HopfOp_B^{\P}$, is composed of $\Z$-graded coalgebras unlike the
Hopf operads considered in~\ref{section:UnitaryHopfOperads} since a
dg-module of homogeneous morphisms is naturally a $\Z$-graded
object. Similarly, one can observe that these operads are unital but
non-connected. On the other hand, we consider only morphisms
$\Q\rightarrow\HopfOp_B^{\P}$ where $\Q$ is a non-negatively graded
unital Hopf operad. One can check that the coalgebra truncation
functor $\str^{\dg}_+$\glossary{$\str^{\dg}_+$} defined
in~\ref{subsection:CofreeCoalgebras} induces a truncation functor on
operad categories so that a universal $\N$-graded unital Hopf operad
is associated to any $\Z$-graded unital Hopf operad. Accordingly,
any morphism as above admits a factorization
\begin{equation*}
\xymatrix@!C=6mm{ \Q\ar@{-->}[dr]\ar[rr] && \HopfOp_B^{\P}
\\ & \str^{\dg}_+(\HopfOp_B^{\P})\ar[ur] & }
\end{equation*}
and in the applications one can replace $\HopfOp_B^{\P}$ by the
non-negatively graded operad
$\str^{\dg}_+(\HopfOp_B^{\P})$\glossary{$\str^{\dg}_+$}. If $\Q$ is
connected, then we can consider a further factorization by a
morphism $\Q\rightarrow\str^1_*\str^{\dg}_+(\HopfOp_B^{\P})$ where
$\str^1_*\str^{\dg}_+(\HopfOp_B^{\P})$ denotes the connected unital
Hopf operad associated to $\str^{\dg}_+(\HopfOp_B^{\P})$ by the
construction of~\ref{item:ConnectedHopfOperads}. For the sake of
completness, recall that this functor $\str^1_*$ preserves
fibrations, acyclic fibrations and all weak-equivalences between
fibrant objects.

As in the previous section, we let a morphism of $\Z$-graded unital
Hopf operads $\phi: \P\rightarrow\Q$ be a Reedy fibration,
respectively an acyclic Reedy fibration, if the associated morphism
$\str^{\dg}_+(\phi): \str^{\dg}_+(\P)\rightarrow\str^{\dg}_+(\Q)$ is
so in the model category of non-negatively graded unital Hopf
operads. Equivalently, the morphism $\phi: \P\rightarrow\Q$ is a
Reedy fibration, respectively an acyclic Reedy fibration, if it has
the left lifting property with respect to acyclic Reedy
cofibrations, respectively Reedy cofibrations, of non-negatively
graded unital Hopf operads.\index{Reedy!fibration!of $\Z$-graded
Hopf operads}\index{fibration!of $\Z$-graded Hopf operads}

We prove that the Hopf operad of bar operations $\HopfOp_B^{\P}$ is
endowed with the following property:

\begin{thm}\label{thm:HopfEndomorphismOperadFibration}
The functor $\P\mapsto\HopfOp_B^{\P}$ maps a fibration, respectively
an acyclic fibration, of non-unital operads under $\K$ to a Reedy
fibration, respectively an acyclic Reedy fibration, of unital Hopf
operads. In particular, the Hopf operad $\HopfOp_B^{\P}$ defines a
fibrant object in the category of Hopf operads, for any operad $\P$
under $\K$.
\end{thm}

One can observe that the Hopf endomorphism operad $\HopfEnd_B^{\P}$
preserve fibrations like $\HopfOp_B^{\P}$ but not acyclic
fibrations.

\medskip
Recall that the commutative operad $\C$ forms the final object of
the category of unital Hopf operads. In the introduction of this
memoir we mention that the classical shuffle product of tensors
corresponds to a morphism $\nabla_c: \C\rightarrow\HopfEnd_B^{\C}$.
In the final part of this section we check that this morphism admits
a factorization
\begin{equation*}
\xymatrix{ \C\ar@{-->}[dr]!UL_{\nabla_\gamma}\ar[r]^(0.35){\nabla_c} & \HopfEnd_B^{\C} \\
& \HopfOp_B^{\C}\ar@{^{(}->}[]!U+<0pt,4pt>;[u]_(0.35){\nabla_\Theta}
& }.
\end{equation*}
In addition we prove the following result:

\begin{thm}\label{thm:HopfOperadLifting}
Any morphism of unital Hopf operads $\nabla_\rho:
\Q\rightarrow\HopfOp_B^{\C}$, where $\Q$ is connected and
non-negatively graded, makes commute the diagram
\begin{equation*}
\xymatrix@!C=6mm{ \Q\ar[dr]_{\epsilon}\ar[rr]^{\nabla_\rho} &&
\HopfOp_B^{\C} \\ & \C\ar[ur]_{\nabla_\gamma} & }.
\end{equation*}
\end{thm}

Accordingly, the morphism $\nabla_\gamma$ induces an isomorphism
\begin{equation*}
\nabla_\gamma:
\C\xrightarrow{\simeq}\str^1_*\str^{\dg}_+(\HopfOp_B^{\C})
\end{equation*}
that identifies $\str^1_*\str^{\dg}_+(\HopfOp_B^{\C})$ with the
commutative operad.

Finally, our main existence and uniqueness result,
theorem~\ref{thm:HopfExistenceUniqueness} of the introduction,
arises as a formal consequence of
theorems~\ref{thm:HopfEndomorphismOperadFibration}
and~\ref{thm:HopfOperadLifting}. We recall our statement for the
sake of completeness. In fact, the claim of
theorem~\ref{thm:HopfExistenceUniqueness} is valid for any unital
Hopf operad (not necessarily $E_\infty$):

\begin{thm}\label{thm:HopfExistenceUniquenessRecall}
Let $\E$ be a $\Sigma_*$-cofibrant non-unital $E_\infty$-operad. Let
$\Q$ be a unital Hopf operad.

\begin{enumerate}
\item
The bar complex of an $\E$-algebra $B(A)$ can be equipped with the
structure of a Hopf $\Q$-algebra, functorially in $A$, and so that
the unital operation $\Q(0)\rightarrow B(A)$ agrees with the natural
unit of the bar complex $\F\rightarrow B(A)$ provided that $\Q$ is a
Reedy cofibrant object in the category of unital Hopf operads.
\item
Any such $\Q$-algebra structure where $\Q$ is connected and
non-negatively graded satisfies the requirement of the uniqueness
theorem of~\cite{Bar1}. More explicitly, if the unit condition of
claim (a) is satisfied and the operad $\Q$ is connected and
non-negatively graded, then, for a commutative algebra $A$, the
$\Q$-algebra structure of $B(A)$ reduces automatically to the
classical commutative algebra structure of $B(A)$, the one given by
the shuffle product of tensors.
\end{enumerate}
\end{thm}

We give the proof of these claims assuming
theorems~\ref{thm:HopfEndomorphismOperadFibration}
and~\ref{thm:HopfOperadLifting}. As explained in the introduction of
this memoir, for an $E_\infty$-operad $\E$, we consider the lifting
problem
\begin{equation*}
\xymatrix{ & & \HopfOp_B^{\E}\ar[r]\ar@{->>}[d]^{\sim} & \HopfEnd_B^{\E}\ar[d] \\
\Q\ar@{-->}^{\exists\nabla_\rho}[urr]\ar[r] &
\C\ar[r]^(0.33){\nabla_\gamma} & \HopfOp_B^{\C}\ar[r] &
\HopfEnd_B^{\C} }
\end{equation*}
which has automatically a solution if $\Q$ is cofibrant since, by
theorem~\ref{thm:HopfEndomorphismOperadFibration}, the augmentation
of an $E_\infty$-operad induces an acyclic fibration $\epsilon_*:
\HopfOp_B^{\E}\wefib\HopfOp_B^{\C}$. Then the composite morphism
\begin{equation*}
\Q\xrightarrow{\nabla_\rho}\HopfOp_B^{\E}\xrightarrow{\nabla_\Theta}\HopfEnd_B^{\E}
\end{equation*}
provides an operad action on the bar complex that fulfils the
existence assertions of
theorem~\ref{thm:HopfExistenceUniquenessRecall}.

Conversely, by the universal definition of the Hopf operad
$\HopfEnd_B^{\E}$, any operad action on $B(A)$ that satisfies our
unit requirement is determined by an operad morphism $\nabla:
\Q\rightarrow\HopfEnd_B^{\E}$. If $\E$ is $\Sigma_*$-cofibrant, then
any such morphism factors through $\HopfOp_B^{\E}$ since the
embedding $\nabla_\Theta:
\HopfOp_B^{\E}\hookrightarrow\HopfEnd_B^{\E}$ is an isomorphism.
Then we deduce from theorem~\ref{thm:HopfOperadLifting} that the
diagram
\begin{equation*}
\xymatrix{ \Q\ar[r]^(0.4){\nabla}\ar[d] & \HopfOp_B^{\E}\ar[d] \\
\C\ar[r]^(0.4){\nabla_\gamma} & \HopfOp_B^{\C} }
\end{equation*}
commutes automatically. For a commutative algebra $A$, it follows
that the action of $\Q$ on $B(A)$ reduces to the classical
commutative operad action as claimed by
theorem~\ref{thm:HopfExistenceUniquenessRecall}. \qed

\medskip
The lifting process can be simplified if we assume that $\Q$ is a
connected operad. Namely, if we apply the truncation functors, then,
by theorem~\ref{thm:HopfOperadLifting}, our lifting problem becomes
equivalent to
\begin{equation*}
\xymatrix{ & \str^1_*\str^{\dg}_+(\HopfOp_B^{\E})\ar@{->>}[d]^{\sim} \\
\Q\ar@{-->}[ur]^{\exists\nabla_\rho}\ar[r] & \C }.
\end{equation*}
Recall also that the truncation functors preserve acyclic
fibrations. Therefore we obtain that the augmentation of an
$E_\infty$-operad induces an acyclic fibration $\epsilon_*:
\str^1_*\str^{\dg}_+(\HopfOp_B^{\E})\wefib\C$ since this morphism
represents the truncation of the morphism $\epsilon_*:
\HopfOp_B^{\E}\wefib\HopfOp_B^{\C}$ which forms an acyclic fibration
by theorem~\ref{thm:HopfEndomorphismOperadFibration}. Finally, we
obtain that any morphism $\nabla_\rho: \Q\rightarrow\HopfOp_B^{\E}$
where $\Q$ is a non-negatively graded connected unital Hopf operad
fits this lifting diagram since $\C$ forms the terminal object in
the category of unital Hopf operads. As a corollary, by usual model
category arguments, we obtain:

\begin{thm}\label{thm:HopfOperadMorphismHomotopyUniqueness}
Let $\E$ be a non-unital $E_\infty$-operad. Any pair of morphisms of
unital Hopf operads $\nabla_0,\nabla_1: \Q\rightarrow\HopfOp_B^{\E}$
where $\Q$ is a connected and non-negatively graded unital Hopf
operads are left-homotopic.
\end{thm}

As in~\cite{Bar1}, one could give an interpretation of this
uniqueness assertion at the algebra level. Namely suppose given a
pair of morphisms of unital Hopf operads
\begin{equation*}
\nabla_0,\nabla_1: \Q\rightarrow\HopfOp_B^{\E}
\end{equation*}
which provide the chain complex $B(A)$ with the structure of a Hopf
$\Q$-algebra. Then the Hopf $\Q$-algebras $(B(A),\nabla_0)$ and
$(B(A),\nabla_1)$ can be connected by weak-equivalences of Hopf
$\Q$-algebras
\begin{equation*}
(B(A),\nabla_0)\xleftarrow{\sim}\cdot\xrightarrow{\sim}(B(A),\nabla_1).
\end{equation*}
To be precise, we have not checked this claim in full generality.
But if we assume that $A$ is a non-negatively graded dg-algebras,
then the bar complex $B(A)$ belongs to the category of
non-negatively graded dg-coalgebra for which we have a model
structure by~\cite{GetzlerGoerss}. In this context the claim can be
deduced from the results of~\cite{Rezk} extended to the ground model
category of dg-coalgebras.

\subsubsection*{Section outline}
Here is the plan of this section:
in~\ref{subsection:MorphismCoalgebras} we prove the existence of an
internal hom in the category of coalgebras; in
subsections~\ref{subsection:BarHopfEndomorphismOperad}-\ref{subsection:BarHopfOperations}
we define the Hopf endomorphism operad $\HopfEnd_B^{\P}$, the
related Hopf operad of bar operations $\HopfOp_B^{\P}$ and we make
explicit the internal structure of these operads; we prove the
fibration properties asserted by
theorem~\ref{thm:HopfEndomorphismOperadFibration}
in~\ref{subsection:FibrationProperties} and we give the proof of
theorem~\ref{thm:HopfOperadLifting}
in~\ref{subsection:CommutativeHopfBarAction}.

\subsection{Morphism coalgebras}\label{subsection:MorphismCoalgebras}
The purpose of this subsection is to construct \emph{morphism
coalgebras} $\HopfHom(L,M)$ that satisfy the adjunction property
specified in the section introduction. Before we shall give
precisions on the coalgebra categories that occur in the definition
of~$\HopfHom(L,M)$.

\subsubsection{Augmented and unitary coalgebras}
In general we work within the category of \emph{augmented
coassociative coalgebras} denoted by
$\CoAlg^a_+$\glossary{$\CoAlg^a_+$}. But, as explained above, we
consider also coalgebras in the category of \emph{non-augmented
coalgebras} $\CoAlg^a$\glossary{$\CoAlg^a$}, in the category of
\emph{augmented unitary coalgebras}
$\CoAlg^a_*$\glossary{$\CoAlg^a_*$} and in the category of
\emph{connected coalgebras} $\CoAlg^a_0$\glossary{$\CoAlg^a_0$}.

To be precise, an object $K\in\CoAlg^a_+$\index{coalgebra!augmented}
denotes a coassociative coalgebra equipped with an augmentation
defined by a morphism of coalgebras $\epsilon:
K\rightarrow\F$\glossary{$\epsilon$} ; an object
$L\in\CoAlg^a_*$\index{coalgebra!augmented unitary} denotes a
coassociative coalgebra equipped with an augmentation and a
coalgebra unit defined by a morphism of augmented coalgebra $\eta:
\F\rightarrow L$. The unit cokernel of a unitary coalgebra
$\overline{L} = \coker(\eta: \F\rightarrow
L)$\glossary{$\overline{L}$} defines an object of $\CoAlg^a$.
Clearly, the map $L\mapsto\overline{L}$ defines an equivalence
between the category of augmented unitary coalgebras $\CoAlg^a_*$
and the category of non-augmented coalgebras
$\CoAlg^a$\index{coalgebra!non-augmented} since any augmented
unitary coalgebra $L$ has a natural decomposition $L =
\F\oplus\overline{L}$.

\subsubsection{Cofree unitary coalgebras}\label{item:CofreeUnitalCoalgebras}
The cofree coalgebra $\Gamma(V)$ defined in the previous subsection
is characterized by the adjunction relation
\begin{equation*}
\Hom_{\dg\Mod}(K,V) = \Hom_{\CoAlg^a_+}(K,\Gamma(V)),
\end{equation*}
for any augmented coalgebra $K\in\CoAlg^a_+$. Notice that
$\Gamma(V)$ is equipped with a canonical unit $\eta:
\F\rightarrow\Gamma(V)$ induced by the null morphism $0:
0\rightarrow V$. The unit cokernel of $\Gamma(V)$ is denoted by
$\overline{\Gamma}(V)$\glossary{$\overline{\Gamma}(V)$}. One
observes readily that the following adjunction relations hold
\begin{equation*}
\Hom_{\dg\Mod}(\overline{L},V) =
\Hom_{\CoAlg^a}(\overline{L},\overline{\Gamma}(V)) =
\Hom_{\CoAlg^a_*}(L,\Gamma(V)),
\end{equation*}
for any augmented unitary coalgebra $L\in\CoAlg^a_*$. Accordingly,
the cofree coalgebra $\Gamma(V)$ forms also a cofree object in the
category of augmented unitary coalgebras and its unit cokernel
$\overline{\Gamma}(V)$ forms a cofree object in $\CoAlg^a$. Clearly,
a morphism of cofree coalgebras $\nabla_f:
\Gamma(V)\rightarrow\Gamma(W)$ preserves units if and only if the
associated map $f: \Gamma(V)\rightarrow W$ cancels the unit of
$\Gamma(V)$ and hence is equivalent to a map $f:
\overline{\Gamma}(V)\rightarrow W$.

We consider also quasi-cofree objects $\Gamma =
(\Gamma(V),\partial_\alpha)$ in the category of augmented unitary
coalgebras $\CoAlg^a_*$. We assume in this case that the
coderivation $\partial_\alpha$ cancels the unit of $\Gamma(V)$ so
that the unit morphism of the cofree coalgebra $\eta:
\F\rightarrow\Gamma(V)$ defines a morphism of dg-coalgebras $\eta:
\F\rightarrow\Gamma$. As in the context of coalgebra morphisms, it
is equivalent to assume that the coderivation $\partial_\alpha$ is
induced by a map $\alpha: \overline{\Gamma}(V)\rightarrow V$.

\medskip
Recall that a tensor product of coassociative coalgebras $K\otimes
L$ is equipped with a natural coalgebra structure so that the
coalgebra categories considered in this memoir are symmetric
monoidal. In the language of monoidal categories, our internal
morphism coalgebra $\HopfHom(L,M)$\index{morphism
coalgebra}\glossary{$\HopfHom(L,M)$} can be characterized by a
closure property. In fact, we define $\HopfHom(L,M)$ by the
following assertion:

\begin{prop}\label{prop:CoalgebraHom}
We have a bifunctor
\begin{equation*}
\HopfHom: \CoAlg^{a\op}_*\times\CoAlg^a_*\rightarrow\CoAlg^a_+
\end{equation*}
that satisfies the adjunction relation
\begin{equation*}
\Hom_{\CoAlg^a}(K\otimes\overline{L},\overline{M}) =
\Hom_{\CoAlg^a_+}(K,\HopfHom(L,M)),
\end{equation*}
for all augmented coalgebras $K\in\CoAlg^a_+$ and all augmented
unitary coalgebras $L,M\in\CoAlg^a_*$.
\end{prop}

Recall that the forgetful functor from a coalgebra category to the
category of dg-modules creates colimits. This assertion implies
immediately that the functor $K\mapsto K\otimes\overline{L}$
preserves colimits. On the other hand, we mention
in~\ref{subsection:CofreeCoalgebras} that the category of coalgebras
is equipped with a set of generating objects. Consequently, the
proposition can be deduced from the special adjoint functor theorem.

Nevertheless we prefer to give another proof of
proposition~\ref{prop:CoalgebraHom} so that we can obtain an
explicit construction of~$\HopfHom(L,M)$ in the case of a
quasi-cofree coalgebra $M = (\Gamma(V),\partial)$. Before observe
that the adjunction relation of the proposition is equivalent to the
adjunction property specified in the section introduction. More
explicitly, we have the following assertion:

\begin{fact}\label{fact:CoalgebraAdjunction}
Let $K\in\CoAlg^a_+$ denote an augmented coalgebra. Let
$L,M\in\CoAlg^a_*$ be augmented unitary coalgebras as in
proposition~\ref{prop:CoalgebraHom}. We have a one-to-one
correspondance bewteen morphisms of augmented coalgebras $\phi:
K\otimes L\rightarrow M$ that make the diagram
\begin{equation*}
\xymatrix{ K\otimes\F\ar[r]^{\epsilon}\ar[d]_{K\otimes\eta} & \F\ar[d]_{\eta} \\
K\otimes L\ar[r]^{\phi} & M }
\end{equation*}
commute and morphisms of non-augmented coalgebras $\overline{\phi}:
K\otimes\overline{L}\rightarrow\overline{M}$ such that
$\overline{\phi}$ is defined by the restriction of $\phi: K\otimes
L\rightarrow M$ to the summand $K\otimes\overline{L}$ of the tensor
product $K\otimes L$.

Hence a morphism of augmented coalgebras $\phi^{\sharp}:
K\rightarrow\HopfHom(L,M)$ is equivalent to a morphism of augmented
coalgebras $\phi: K\otimes L\rightarrow M$ that makes the diagram
above commute.
\end{fact}

As explained previously, we aim to prove the existence of morphism
coalgebras $\HopfHom(L,M)$ by an effective construction. In fact,
our arguments rely on a classical proof of the existence of adjoint
functors in a category of coalgebras over a comonad.

To begin with, we have the following immediate observation:

\begin{obsv}\label{obsv:CoFreeCoalgebraHom}
For a cofree coalgebra $M = \Gamma(V)$, the required morphism
coalgebra is given by the cofree coalgebra
\begin{equation*}
\HopfHom(L,\Gamma(V)) = \Gamma(\DGHom(\overline{L},V))
\end{equation*}
since we have adjunction relations
\begin{multline*}
\Hom_{\CoAlg^a}(K\otimes\overline{L},\overline{\Gamma}(V)) \\
= \Hom_{\dg\Mod}(K\otimes\overline{L},V) = \Hom_{\dg\Mod}(K,\DGHom(\overline{L},V)) \\
= \Hom_{\CoAlg^a_+}(K,\Gamma(\DGHom(\overline{L},V))).\qed
\end{multline*}
\end{obsv}\index{morphism coalgebra!of a cofree coalgebra}

These adjunction relations are functorial in $K\in\CoAlg^a_+$ and
$L\in\CoAlg^a_*$ but the middle terms are not functors in $M =
\Gamma(V)\in\CoAlg^a_*$. Therefore we prove directly that
$\HopfHom(L,\Gamma(V)) = \Gamma(\DGHom(\overline{L},V))$ defines a
functor on the full subcategory of $\CoAlg^a_*$ generated by cofree
coalgebras $M = \Gamma(V)$. For this purpose we consider the
augmentation morphism of the adjunction of
observation~\ref{obsv:CoFreeCoalgebraHom}:
\begin{equation*}
\ev_{\Gamma}:
\Gamma(\DGHom(\overline{L},V))\otimes\overline{L}\rightarrow\Gamma(V).\glossary{$\ev_\Gamma$}
\end{equation*}
An explicit definition of this morphism can be obtained by going
through our adjunction relations. We obtain precisely:

\begin{obsv}\label{obsv:CoFreeCoalgebraEvaluation}
In the adjunction of observation~\ref{obsv:CoFreeCoalgebraHom}
\begin{equation*}
\Hom_{\CoAlg^a}(K\otimes\overline{L},\overline{\Gamma}(V)) =
\Hom_{\CoAlg^a_+}(K,\Gamma(\DGHom(\overline{L},V))),
\end{equation*}
the augmentation morphism can be identified with the coalgebra
morphism
\begin{equation*}
\ev_{\Gamma}:
\Gamma(\DGHom(\overline{L},V))\otimes\overline{L}\rightarrow\overline{\Gamma}(V)
\end{equation*}
induced by the composite morphism of dg-modules
\begin{equation*}
\Gamma(\DGHom(\overline{L},V))\otimes\overline{L}
\rightarrow\DGHom(\overline{L},V)\otimes\overline{L}
\xrightarrow{\ev}V,
\end{equation*}
where we consider the universal morphism of the cofree coalgebra
\begin{equation*}
\Gamma(\DGHom(\overline{L},V))\rightarrow\DGHom(\overline{L},V),
\end{equation*}
and the augmentation morphism
\begin{equation*}
\ev: \DGHom(\overline{L},V)\otimes\overline{L}\rightarrow
V\glossary{$\ev$}
\end{equation*}
of the adjunction relation
\begin{equation*}
\Hom_{\dg\Mod}(K\otimes\overline{L},V) =
\Hom_{\dg\Mod}(K,\DGHom(\overline{L},V))
\end{equation*}
in the category of dg-modules.\qed
\end{obsv}

Then let
\begin{equation*}
\ev_\Gamma^\sharp:
\Gamma(\DGHom(\overline{L},V))\rightarrow\DGHom(\overline{L},\overline{\Gamma}(V)).
\end{equation*}
denote the adjoint morphism of $\ev_{\Gamma}:
\Gamma(\DGHom(\overline{L},V))\otimes\overline{L}\rightarrow\overline{\Gamma}(V)$
in the category of dg-modules. We obtain:

\begin{claim}\label{claim:CoFreeCoalgebraHom}
We suppose given a morphism of cofree coalgebras $\nabla_f:
\Gamma(V)\rightarrow\Gamma(W)$ induced by a morphism of dg-modules
$f: \overline{\Gamma}(V)\rightarrow W$ so that $\nabla_f$ preserves
the unit of the cofree coalgebra. We consider the morphism of cofree
coalgebras
\begin{equation*}
\nabla_{f_*\ev_\Gamma^\sharp}:
\Gamma(\DGHom(\overline{L},V))\rightarrow\Gamma(\DGHom(\overline{L},W)).
\end{equation*}
induced by the composite
\begin{equation*}
\Gamma(\DGHom(\overline{L},V))\xrightarrow{\ev_\Gamma^\sharp}\DGHom(\overline{L},\overline{\Gamma}(V))
\xrightarrow{f_*}\DGHom(\overline{L},W).
\end{equation*}
We claim that this morphism fits a commutative diagram
\begin{equation*}
\xymatrix{
\Hom_{\CoAlg^a}(K\otimes\overline{L},\overline{\Gamma}(V))\ar[d]^{(\nabla_f)_*}\ar[r]^(0.45){=}
&
\Hom_{\CoAlg^a_+}(K,\Gamma(\DGHom(\overline{L},V)))\ar[d]^{(\nabla_{f_*\ev_\Gamma^\sharp})_*} \\
\Hom_{\CoAlg^a}(K\otimes\overline{L},\overline{\Gamma}(W))\ar[r]^(0.45){=}
& \Hom_{\CoAlg^a_+}(K,\Gamma(\DGHom(\overline{L},W))) }
\end{equation*}
and hence makes the adjunction relation of
observation~\ref{obsv:CoFreeCoalgebraHom} functorial with respect to
the coalgebra morphism $\nabla_f: \Gamma(V)\rightarrow\Gamma(W)$.
\end{claim}

\begin{proof}
One can check that $\nabla_{f_*\ev_\Gamma^\sharp}$ makes the diagram
\begin{equation*}
\xymatrix{
\Gamma(\DGHom(\overline{L},V))\otimes\overline{L}\ar[r]^(0.7){\ev_\Gamma}\ar[d]_{\nabla_{f_*\ev_\Gamma^\sharp}\otimes\overline{L}}
& \overline{\Gamma}(V)\ar[d]^{\nabla_f} \\
\Gamma(\DGHom(\overline{L},W))\otimes\overline{L}\ar[r]^(0.7){\ev_\Gamma}
& \overline{\Gamma}(W) }
\end{equation*}
commute. For this purpose it suffices to compare the composite
\begin{equation*}
\Gamma(\DGHom(\overline{L},V))\otimes\overline{L}
\xrightarrow{\nabla_{f_*\ev_\Gamma^\sharp}\otimes\overline{L}}\Gamma(\DGHom(\overline{L},W))\otimes\overline{L}
\xrightarrow{\ev_\Gamma}\overline{\Gamma}(W)\xrightarrow{\pi}W
\end{equation*}
with $\pi\cdot\nabla_f\cdot\ev_\Gamma = f\cdot\ev_\Gamma$. The
identity of these morphisms follows from formal verifications
involving essentially adjunctions in the category of dg-modules.

According to observation~\ref{obsv:CoFreeCoalgebraHom}, the adjoint
of a coalgebra morphism $\nabla_u:
K\rightarrow\Gamma(\DGHom(\overline{L},V))$ is obtained by the
composite
\begin{equation*}
K\otimes\overline{L}
\xrightarrow{\nabla_u\otimes\overline{L}}\Gamma(\DGHom(\overline{L},V))\otimes\overline{L}
\xrightarrow{\ev_\Gamma}\overline{\Gamma}(V).
\end{equation*}
Therefore our adjunction claim follows from the commutativity of the
diagram
\begin{equation*}
\xymatrix{
K\otimes\overline{L}\ar[r]^(0.35){\nabla_u\otimes\overline{L}} &
\Gamma(\DGHom(\overline{L},V))\otimes\overline{L}\ar[r]^(0.7){\ev_\Gamma}\ar[d]_{\nabla_{f_*\ev_\Gamma^\sharp}}
& \overline{\Gamma}(V)\ar[d]_{\nabla_f} \\
& \DGHom(\overline{L},W))\otimes\overline{L}\ar[r]^(0.7){\ev_\Gamma}
& \overline{\Gamma}(W) }
\end{equation*}
which is proved above.
\end{proof}

Finally, claim~\ref{claim:CoFreeCoalgebraHom} gives the following
result:

\begin{lemm}\label{lemm:CoFreeCoalgebraHom}
The map $V\mapsto\Gamma(\DGHom(\overline{L},V))$ extends to a
functor on the full subcategory of $\CoAlg^a_*$ formed by cofree
coalgebras $\Gamma(V)\in\CoAlg^a_*$ so that the adjunction relation
\begin{equation*}
\Hom_{\CoAlg^a}(K\otimes\overline{L},\overline{\Gamma}(V)) =
\Hom_{\CoAlg^a_+}(K,\Gamma(\DGHom(\overline{L},V)))
\end{equation*}
is functorial in $K$, $L$ and $\Gamma(V)\in\CoAlg^a_*$.
\end{lemm}

The crux of our construction is supplied by this statement. In fact,
according to general categorical constructions, any coalgebra $M$ is
the equalizer of a natural pair of cofree coalgebra morphisms
\begin{equation*}
\xymatrix{ M\ar[r] &
\Gamma(V^0)\ar@<+1mm>[r]^{d^0}\ar@<-1mm>[r]_{d^1} & \Gamma(V^1) }.
\end{equation*}
As a consequence:

\begin{fact}
The coalgebra $\HopfHom(L,M)$ can be defined by the equalizer
diagram
\begin{equation*}
\xymatrix{ \HopfHom(L,M)\ar[r] &
\Gamma(\DGHom(\overline{L},V^0))\ar@<+1mm>[r]^{d^0}\ar@<-1mm>[r]_{d^1}
& \Gamma(\DGHom(\overline{L},V^1)) },
\end{equation*}
where the morphisms $d^0,d^1$ are deduced from
claim~\ref{claim:CoFreeCoalgebraHom}. The adjunction relation
\begin{equation*}
\Hom_{\CoAlg^a}(K\otimes\overline{L},\overline{M}) =
\Hom_{\CoAlg^a_+}(K,\HopfHom(L,M))
\end{equation*}
follows then from the case of cofree coalgebras stated in
lemma~\ref{lemm:CoFreeCoalgebraHom} by an immediate and classical
exactness argument.
\end{fact}

This assertion achieves the proof of
proposition~\ref{prop:CoalgebraHom}.\qed
\medskip

As usual for internal hom-objects, we have a composition product
\begin{equation*}
\HopfHom(M,N)\otimes\HopfHom(L,M)\xrightarrow{\circ}\HopfHom(L,N)
\end{equation*}
equivalent to the composite evaluation morphism
\begin{equation*}
\HopfHom(M,N)\otimes\HopfHom(L,M)\otimes
L\xrightarrow{\Id\otimes\ev}\HopfHom(M,N)\otimes
M\xrightarrow{\ev}N.
\end{equation*}
For cofree coalgebras $M = \Gamma(V)$ and $N = \Gamma(W)$, we obtain
readily:

\begin{obsv}\label{obsv:CoalgebraComposition}\index{morphism
coalgebra!composites in a} For cofree coalgebras $M = \Gamma(V)$ and
$N = \Gamma(W)$, the composition product
\begin{equation*}
\HopfHom(\Gamma(V),\Gamma(W))\otimes\HopfHom(L,\Gamma(V))\xrightarrow{\circ}\HopfHom(L,\Gamma(W))
\end{equation*}
can be identified with the coalgebra morphism
\begin{equation*}
\Gamma(\DGHom(\overline{\Gamma}(V),W))\otimes\Gamma(\DGHom(\overline{L},V))
\xrightarrow{\circ_\Gamma}\Gamma(\DGHom(\overline{L},W))
\end{equation*}
induced by the composite
\begin{multline*}
\Gamma(\DGHom(\overline{\Gamma}(V),W))\otimes\Gamma(\DGHom(\overline{L},V))
\\
\xrightarrow{\pi_*\otimes\ev^\sharp_\Gamma}
\DGHom(\overline{\Gamma}(V),W)\otimes\DGHom(\overline{L},\overline{\Gamma}(V))
\xrightarrow{\circ}\DGHom(\overline{L},W)
\end{multline*}
where we consider the universal projection
\begin{equation*}
\Gamma(\DGHom(\overline{\Gamma}(V),W))
\xrightarrow{\pi}\DGHom(\overline{\Gamma}(V),W),
\end{equation*}
the morphism
\begin{equation*}
\Gamma(\DGHom(\overline{L},V))\xrightarrow{\ev^\sharp_\Gamma}\DGHom(\overline{L},\overline{\Gamma}(V))
\end{equation*}
and the composition product of homogeneous maps of dg-modules
\begin{equation*}
\DGHom(\overline{\Gamma}(V),W)\otimes\DGHom(\overline{L},\overline{\Gamma}(V))
\xrightarrow{\circ}\DGHom(\overline{L},W).\qed
\end{equation*}
\end{obsv}

We observe that a morphism of augmented unitary coalgebras $\phi:
L\rightarrow M$ is equivalent to a group-like element
$\phi\in\HopfHom(L,M)$. Formally, the module $\F[X]$ spanned by a
set $X$ is usually equipped with the structure of an augmented
coalgebra in which the basis elements $x\in X$ are group-like. This
coalgebra satisfies the adjunction relation
$\Hom_{\CoAlg^a_+}(\F[X],\Gamma) = \Hom_{\Set}(X,\Gr(\Gamma))$,
where $\Gr(\Gamma)$ denotes the set of group-like elements in an
augmented coalgebra $\Gamma\in\CoAlg^a_+$. One checks readily that
the evaluation of morphisms yields a coalgebra morphism
\begin{equation*}
\F[\Hom_{\CoAlg^a_*}(L,M)]\otimes L\rightarrow M
\end{equation*}
that satisfies the requirement of
fact~\ref{fact:CoalgebraAdjunction}. Accordingly, this morphism is
equivalent to a morphism of augmented coalgebras
\begin{equation*}
\F[\Hom_{\CoAlg^a_*}(L,M)]\rightarrow\HopfHom(L,M).
\end{equation*}
We obtain:

\begin{claim}\label{claim:GroupLikeCoalgebraMorphisms}
The coalgebra morphism above yields a canonical bijection
\begin{equation*}
\Hom_{\CoAlg^a_*}(L,M)\xrightarrow{\simeq}\Gr(\HopfHom(L,M)).
\end{equation*}
\end{claim}

\begin{proof}
This bijection property is stated is a remark. Therefore we just
sketch the proof of this statement. By left exactness, it is
sufficient to check the assertion for a cofree coalgebra $M =
\Gamma(V)$ for which we have $\HopfHom(L,\Gamma(V)) =
\Gamma(\DGHom(\overline{L},V))$. In this case, by adjunction, we
have natural bijections
\begin{equation*}
\Hom_{\CoAlg^a_*}(L,\Gamma(V))
\xrightarrow{\simeq}\Hom_{\dg\Mod}(L,\Gamma(V))
\xrightarrow{\simeq}\Gr(\Gamma(\DGHom(\overline{L},V)))
\end{equation*}
whose composite can be identified with the map of the claim.
\end{proof}

The last assertion of the proof can be deduced from the following
observation:

\begin{obsv}
For a cofree coalgebra $M = \Gamma(V)$, the map
\begin{equation*}
\F[\Hom_{\CoAlg^a_*}(L,\Gamma(V))]\rightarrow\HopfHom(L,\Gamma(V)).
\end{equation*}
can be identified with the coalgebra morphism
\begin{equation*}
\F[\Hom_{\CoAlg^a_*}(L,\Gamma(V))]\rightarrow\Gamma(\DGHom(\overline{L},V))
\end{equation*}
induced by the canonical morphism of dg-modules
\begin{multline*}
\F[\Hom_{\CoAlg^a_*}(L,\Gamma(V))]\\
\xrightarrow{\simeq}\F[\Hom_{\dg\Mod}(\overline{L},V)]
\rightarrow\Hom_{\dg\Mod}(\overline{L},V)\\
\hookrightarrow\DGHom(\overline{L},V).\qed
\end{multline*}
\end{obsv}

This assertion arises as a straightforward consequence of our
constructions. To conclude this set of observations, we can form an
enriched category of augmented unitary coalgebras in which morphism
objects are given by the morphism coalgebras $\HopfHom(L,M)$. The
bijection
$\Hom_{\CoAlg^a_+}(L,M)\xrightarrow{\simeq}\Gr(\HopfHom(L,M))$
yields an embedding from the ground category to the enriched
category of augmented unitary coalgebras. The morphism coalgebra
$\HopfHom(L,M)$ extends clearly to a bifunctor on this enriched
category. Observe also that for a group-like element
$\nabla_f\in\HopfHom(\Gamma(V),\Gamma(W)$, equivalent to a morphism
$\nabla_f: \Gamma(V)\rightarrow\Gamma(W)$, the composition process
of observation~\ref{obsv:CoalgebraComposition} extends the
definition of claim~\ref{claim:CoFreeCoalgebraHom}.

\medskip We claim that our construction permits to obtain an
explicit realization of the morphism coalgebra $\HopfHom(L,M)$ for a
quasi-cofree coalgebra $M = (\Gamma(V),\partial_\alpha)$. We obtain
precisely the following result:

\begin{lemm}\label{lemm:QuasiCofreeHopfHom}\index{morphism
coalgebra!of a quasi-cofree coalgebra} If $M$ is a quasi-cofree
coalgebra, then so is the morphism coalgebra $\HopfHom(L,M)$. To be
more explicit, suppose given a quasi-cofree unitary coalgebra $M =
(\Gamma(V),\partial_\alpha)$, where the coderivation
$\partial_\alpha: \Gamma(V)\rightarrow\Gamma(V)$ is induced by a
homogeneous map $\alpha: \overline{\Gamma}(V)\rightarrow V$. Then we
have:
\begin{equation*}
\HopfHom(L,M) =
(\Gamma(\DGHom(\overline{L},V)),\partial_{\alpha_*\ev_\Gamma^\sharp}),
\end{equation*}
where the coderivation $\partial_{\alpha_*\ev_\Gamma^\sharp}:
\Gamma(\DGHom(\overline{L},V))\rightarrow\Gamma(\DGHom(\overline{L},V))$
is induced by the composite map
\begin{equation*}
\Gamma(\DGHom(\overline{L},V))\xrightarrow{\ev_\Gamma^\sharp}\DGHom(\overline{L},\overline{\Gamma}(V))
\xrightarrow{\alpha_*}\DGHom(\overline{L},V).
\end{equation*}
Furthermore, the adjunction augmentation
\begin{equation*}
\ev: \HopfHom(L,M)\otimes\overline{L}\rightarrow\overline{M}
\end{equation*}
can be identified with a coalgebra morphism
\begin{equation*}
\ev_\Gamma:
(\Gamma(\DGHom(\overline{L},V)),\partial_{\alpha_*\ev_\Gamma^\sharp})\otimes\overline{L}
\rightarrow(\overline{\Gamma}(V),\partial_\alpha)
\end{equation*}
supplied by observation~\ref{obsv:CoFreeCoalgebraEvaluation}.
\end{lemm}

Roughly, we check that the equations of
lemma~\ref{lemm:QuasiCofreeCoalgebraStructure} involved in
quasi-cofree coalgebra structures hold for the pair
$(\Gamma(\DGHom(\overline{L},V)),\partial_{\alpha_*\ev_\Gamma^\sharp})$
and that the adjunction relation of
observation~\ref{obsv:CoFreeCoalgebraHom} extends to quasi-cofree
coalgebras. Basically, our results are consequences of the following
commutation relation between the coderivations and the adjunction
augmentation:

\begin{claim}\label{claim:AdjunctionAugmentationCoderivation}
The coderivation $\partial_{\alpha_*\ev_\Gamma^\sharp}$ makes the
diagram
\begin{equation*}
\xymatrix{\Gamma(\DGHom(\overline{L},V))\otimes\overline{L}\ar[r]^(0.7){\ev_\Gamma}
\ar[d]_{\partial_{\alpha_*\ev_\Gamma^\sharp}\otimes\overline{L}}
& \overline{\Gamma}(V)\ar[d]^{\partial_{\alpha}} & \\
\Gamma(\DGHom(\overline{L},V))\otimes\overline{L}\ar[r]^(0.7){\ev_\Gamma}
& \overline{\Gamma}(V) }
\end{equation*}
commute.
\end{claim}

\begin{proof}
In the diagram both composites define coderivations $\partial:
L'\rightarrow\overline{\Gamma}(V)$, where $L' =
\Gamma(\DGHom(\overline{L},V))$. One can extend the correspondence
of lemma~\ref{lemm:QuasiCofreeCoalgebraStructure} to this relative
context. Hence these coderivations agree if and only if their
composite with the universal morphism $\pi: \Gamma(V)\rightarrow V$
agree. Thus, as in the proof of
claim~\ref{claim:CoFreeCoalgebraHom}, it suffices to compare the
composite
\begin{equation*}
\Gamma(\DGHom(\overline{L},V))\otimes\overline{L}
\xrightarrow{\partial_{\alpha_*\ev_\Gamma^\sharp}\otimes\overline{L}}\Gamma(\DGHom(\overline{L},V))\otimes\overline{L}
\xrightarrow{\ev_\Gamma}\overline{\Gamma}(V)\xrightarrow{\pi}V
\end{equation*}
with $\pi\cdot\alpha\cdot\ev_\Gamma = \alpha\cdot\ev_\Gamma$. The
identity of these morphisms follows also from formal verifications.
\end{proof}

As a corollary, we obtain:

\begin{claim}
The relation $\delta(\alpha) + \alpha\partial_\alpha = 0$ implies
the same relation $\delta(\alpha_*\ev_\Gamma^\sharp) +
\alpha_*\ev_\Gamma^\sharp\cdot\partial_{\alpha_*\ev_\Gamma^\sharp} =
0$ for the coderivation $\partial_{\alpha_*\ev_\Gamma^\sharp}$.
Hence the pair
\begin{equation*}
(\Gamma(\DGHom(\overline{L},V)),\partial_{\alpha_*\ev_\Gamma^\sharp})
\end{equation*}
defines a quasi-cofree coalgebra.
\end{claim}

\begin{proof}
By adjunction, we deduce from
claim~\ref{claim:AdjunctionAugmentationCoderivation} that the
diagram
\begin{equation*}
\xymatrix{\Gamma(\DGHom(\overline{L},V))\ar[r]^{\ev_\Gamma^\sharp}\ar[d]_{\partial_{\alpha_*\ev_\Gamma^\sharp}}
& \DGHom(\overline{L},\Gamma(V))\ar[d]^{(\partial_{\alpha})_*} \\
\Gamma(\DGHom(\overline{L},V))\ar[r]^{\ev_\Gamma^\sharp} &
\DGHom(\overline{L},\Gamma(V))\ar[d]^{\alpha_*} \\ &
\DGHom(\overline{L},V) }
\end{equation*}
commute. Hence the relation $\delta(\alpha) + \alpha\partial_\alpha
= 0$ implies $\delta(\alpha)_*\ev_\Gamma^\sharp +
\alpha\partial_\alpha\ev_\Gamma^\sharp =
\delta(\alpha_*\ev_\Gamma^\sharp) +
\alpha_*\ev_\Gamma^\sharp\cdot\partial_{\alpha_*\ev_\Gamma^\sharp} =
0$.
\end{proof}

Then lemma~\ref{lemm:QuasiCofreeHopfHom} is a consequence of the
following claim:

\begin{claim}
A dg-module map $u: K\otimes\overline{L}\rightarrow V$ induces a
morphism of dg-coalgebras
\begin{equation*}
\nabla_u:
K\otimes\overline{L}\rightarrow(\overline{\Gamma}(V),\partial_\alpha)
\end{equation*}
to the quasi-cofree coalgebra $M =
(\overline{\Gamma}(V),\partial_\alpha)$ if and only if the adjoint
map $u^\sharp: K\rightarrow\DGHom(\overline{L},V)$ induces a
morphism of dg-coalgebras
\begin{equation*}
\nabla_{u^\sharp}:
K\rightarrow(\Gamma(\DGHom(\overline{L},V)),\partial_{\alpha_*\ev_\Gamma^\sharp})
\end{equation*}
to the quasi-cofree coalgebra
$(\Gamma(\DGHom(\overline{L},V)),\partial_{\alpha_*\ev_\Gamma^\sharp})$.

Accordingly, the adjunction relation
$\nabla_u\mapsto\nabla_{u^\sharp}$ yields an adjunction relation
\begin{equation*}
\Hom_{\CoAlg^a}(K\otimes\overline{L},(\overline{\Gamma}(V),\partial_\alpha))
=
\Hom_{\CoAlg^a_+}(K,(\Gamma(\DGHom(\overline{L},V)),\partial_{\alpha_*\ev_\Gamma^\sharp}))
\end{equation*}
for the quasi-cofree coalgebra $M =
(\overline{\Gamma}(V),\partial_\alpha)$.
\end{claim}

\begin{proof}
The map $u: K\otimes\overline{L}\rightarrow V$ is not assumed to
commute with internal differentials. Accordingly, the induced
coalgebra morphism $\nabla_u:
K\otimes\overline{L}\rightarrow\overline{\Gamma}(V)$ forms only a
morphism of graded coalgebras. Hence, in the correspondence
$\nabla_u\mapsto\nabla_{u^\sharp}$, we assume implicitly that we
extend the adjunction relations of
observation~\ref{obsv:CoFreeCoalgebraHom} to morphisms of graded
coalgebras. The morphism $\nabla_u$ can also be obtained from
$\nabla_{u^\sharp}$ by the composite
\begin{equation*}
K\otimes\overline{L}
\xrightarrow{\nabla_{u^\sharp}\otimes\overline{L}}\Gamma(\DGHom(\overline{L},V))\otimes\overline{L}
\xrightarrow{\ev_\Gamma}\overline{\Gamma}(V)
\end{equation*}
as in the case of an actual morphism of dg-coalgebras.

Explicitly, we check that the coalgebra map $\nabla_u:
K\otimes\overline{L}\rightarrow\overline{\Gamma}(V)$ induced by $u:
K\otimes\overline{L}\rightarrow V$ satisfies the commutation
relation $(\delta+\partial_{\alpha})\nabla_u = \nabla_u\delta$ for
the differential $\delta+\partial_\alpha$ of the quasi-cofree
coalgebra $M = (\overline{\Gamma}(V),\partial_\alpha)$ if and only
if the adjoint coalgebra map $\nabla_{u^\sharp}:
K\rightarrow\Gamma(\DGHom(\overline{L},V))$ satisfies the
commutation relation
$(\delta+\partial_{\alpha_*\ev_\Gamma^\sharp})\nabla_{u^\sharp} =
\nabla_{u^\sharp}\delta$. This claim is a formal consequence of the
observation above and of the commutation assertion of
claim~\ref{claim:AdjunctionAugmentationCoderivation}.
\end{proof}

This claim achieves the proof of
lemma~\ref{lemm:QuasiCofreeHopfHom}.\qed

\medskip
One can observe further that the morphism
\begin{equation*}
(\nabla_f)_*:
\HopfHom(L,(\Gamma(V),\partial_\alpha))\rightarrow\HopfHom(L,(\Gamma(W),\partial_\beta))
\end{equation*}
induced by a morphism of quasi-cofree coalgebras $\nabla_f:
(\Gamma(V),\partial_\alpha)\rightarrow(\Gamma(W),\partial_\beta)$ is
obtained by the same construction as in the case of cofree
coalgebras. Namely this morphism can be identified with a morphism
of quasi-cofree coalgebras
\begin{equation*}
\nabla_{f_*\ev^\sharp_\Gamma}:
(\Gamma(\DGHom(\overline{L},V)),\partial_{\alpha_*\ev_\Gamma^\sharp})
\rightarrow(\Gamma(\DGHom(\overline{L},W)),\partial_{\beta_*\ev^\sharp_\Gamma})
\end{equation*}
which is induced by the composite map
\begin{equation*}
\Gamma(\DGHom(\overline{L},V))
\xrightarrow{\ev_\Gamma^\sharp}\DGHom(\overline{L},\overline{\Gamma}(V))
\xrightarrow{f_*}\DGHom(\overline{L},W).
\end{equation*}
This observation can be extended to the composition product
\begin{multline*}
\HopfHom((\Gamma(V),\partial_\alpha),(\Gamma(W),\partial_\beta))
\otimes\HopfHom(L,(\Gamma(V),\partial_\alpha))\\
\xrightarrow{\circ}\HopfHom(L,(\Gamma(W),\partial_\beta))
\end{multline*}
which is also given by the construction of
observation~\ref{obsv:CoalgebraComposition}.

\medskip
In the next subsections, we consider morphism coalgebras
$\HopfHom(L,M)$ for connected coalgebras $L$ and $M$. Therefore, for
our needs, we state connected variants of our previous assertions.

\subsubsection{Connected coalgebras}
First, recall that a unitary coalgebra $L$ is
\emph{connected}\index{connected coalgebra} if the iterated
coproduct of the associated non-unitary coalgebra $\Delta^n:
\overline{L}\rightarrow\overline{L}^{\otimes n}$ vanishes for $n$
sufficiently large. The full subcategory of $\CoAlg^a_*$ formed by
connected coalgebras is denoted by
$\CoAlg^a_0$\glossary{$\CoAlg^a_0$}. This category is equipped with
cofree object like $\CoAlg^a_*$. One observes precisely that the
cofree object cogenerated by $V$ in $\CoAlg^a_0$ is realized by the
\emph{tensor coalgebra} $T^c(V)$\glossary{$T^c(V)$} defined by the
direct sum
\begin{equation*}
T^c(V) = \bigoplus_{n=0}^{\infty} V^{\otimes n}
\end{equation*}
and equipped with the diagonal $\Delta_\amalg: T^c(V)\rightarrow
T^c(V)\otimes T^c(V)$ induced by the deconcatenation of tensors. In
fact, we have clearly an adjunction relation
\begin{equation*}
\Hom_{\dg\Mod}(\overline{L},V) = \Hom_{\CoAlg^a_0}(L,T^c(V)),
\end{equation*}
for any connected coalgebra $L\in\CoAlg^a_0$. The unit cokernel of
the tensor coalgebra is also denoted by
$\overline{T}^c(V)$\glossary{$\overline{T}^c(V)$}.

\medskip
Our results on the morphism coalgebra $\HopfHom(L,M)$ associated to
connected coalgebras $L$ and $M$ are consequences of the following
observation:

\begin{obsv}\label{obsv:ConnectedProduct}
If $L$ is a connected coalgebra, then the iterated coproduct
$\Delta^n:
K\otimes\overline{L}\rightarrow(K\otimes\overline{L})^{\otimes n}$
vanishes for $n$ large, for any coalgebra $K$. Accordingly, the
tensor product $K\otimes\overline{L}$ forms also a connected (but
non-augmented) coalgebra.\qed
\end{obsv}

Accordingly, for a connected coalgebra $L$, The map $K\mapsto
K\otimes\overline{L}$ defines a functor from the category of
augmented coalgebras to the category of connected coalgebras. We
have clearly:

\begin{fact}\label{fact:ConnectedCoalgebraAdjunction}\index{morphism
coalgebra!of a connected coalgebra}\index{connected
coalgebra!morphism coalgebra of a} The restriction of the functor
$M\mapsto\HopfHom(L,M)$ to connected coalgebras $M$ defines a right
adjoint of the functor $K\mapsto K\otimes\overline{L}$ considered
above.
\end{fact}

On the other hand, one can deduce from
observation~\ref{obsv:ConnectedProduct} that the relations of
observation~\ref{obsv:CoFreeCoalgebraHom} hold for a connected
cofree coalgebra $M = T^c(V)$. Explicitly, we have the following
assertion:

\begin{obsv}\label{obsv:ConnectedCoFreeHom}\index{morphism
coalgebra!of a tensor coalgebra}\index{tensor coalgebra!morphism
coalgebra of a} If $L$ is a connected coalgebra, then we have
adjunction relations
\begin{multline*}
\Hom_{\CoAlg^a}(K\otimes\overline{L},\overline{T}^c(V)) \\
= \Hom_{\dg\Mod}(K\otimes\overline{L},V) = \Hom_{\dg\Mod}(K,\DGHom(\overline{L},V)) \\
= \Hom_{\CoAlg^a_+}(K,\Gamma(\DGHom(\overline{L},V)))
\end{multline*}
Accordingly, for a connected coalgebra $L$ and a cofree connected
coalgebra $M = T^c(V)$, we obtain:
\begin{equation*}
\HopfHom(L,T^c(V)) = \Gamma(\DGHom(\overline{L},V)).\qed
\end{equation*}
\end{obsv}

In fact, one can adapt the previous constructions in order to obtain
explicit realizations of the morphism coalgebras $\HopfHom(L,M)$
associated to connected coalgebras as in the non-connected context.
First, we have a tractable realization of the adjunction
augmentation
\begin{equation*}
\ev_T:
\Gamma(\DGHom(\overline{L},V))\otimes\overline{L}\rightarrow\overline{T}^c(V).\glossary{$\ev_T$}
\end{equation*}
Namely:

\begin{obsv}\label{obsv:TensorCoalgebraEvaluation}
In the adjunction of observation~\ref{obsv:ConnectedCoFreeHom}
\begin{equation*}
\Hom_{\CoAlg^a}(K\otimes\overline{L},\overline{T}^c(V)) =
\Hom_{\CoAlg^a_+}(K,\Gamma(\DGHom(\overline{L},V))),
\end{equation*}
the augmentation morphism can be identified with the morphism of
connected coalgebras
\begin{equation*}
\ev_T:
\Gamma(\DGHom(\overline{L},V))\otimes\overline{L}\rightarrow\overline{T}^c(V)
\end{equation*}
induced by the composite morphism of dg-modules
\begin{equation*}
\Gamma(\DGHom(\overline{L},V))\otimes\overline{L}
\rightarrow\DGHom(\overline{L},V)\otimes\overline{L}
\xrightarrow{\ev}V
\end{equation*}
as in the case of non-connected coalgebras.\qed
\end{obsv}

Then, for quasi-cofree connected coalgebras, we obtain:

\begin{lemm}\label{lemm:ConnectedQuasiCofreeHopfHom}\index{morphism
coalgebra!of a quasi-cofree connected
coalgebra}\index{quasi-cofree!connected coalgebra!morphism coalgebra
of a} If $L$ is a connected coalgebra and $M$ is a quasi-cofree
connected coalgebra, then the morphism coalgebra $\HopfHom(L,M)$ is
a quasi-cofree (non-connected) coalgebra. More precisely, if $M =
(T^c(V),\partial_\alpha)$, for a coderivation $\partial_\alpha:
T^c(V)\rightarrow T^c(V)$ induced by a homogeneous map $\alpha:
\overline{T}^c(V)\rightarrow V$, we have:
\begin{equation*}
\HopfHom(L,M) =
(\Gamma(\DGHom(\overline{L},V)),\partial_{\alpha_*\ev_T^\sharp}),
\end{equation*}
where the coderivation $\partial_{\alpha_*\ev_T^\sharp}:
\Gamma(\DGHom(\overline{L},V))\rightarrow\Gamma(\DGHom(\overline{L},V))$
is induced by the composite map
\begin{equation*}
\Gamma(\DGHom(\overline{L},V))\xrightarrow{\ev_T^\sharp}\DGHom(\overline{L},\overline{T}^c(V))
\xrightarrow{\alpha_*}\DGHom(\overline{L},V).
\end{equation*}
Furthermore, the adjunction augmentation
\begin{equation*}
\ev: \HopfHom(L,M)\otimes\overline{L}\rightarrow\overline{M}
\end{equation*}
can be identified with the coalgebra morphism
\begin{equation*}
\ev_T: (\Gamma(\DGHom(\overline{L},V)),\partial)\otimes\overline{L}
\rightarrow(\overline{T}^c(V),\partial)
\end{equation*}
supplied by observation~\ref{obsv:TensorCoalgebraEvaluation}.
\end{lemm}

\begin{proof}
The proof relies on the same verifications as in the case of
non-connected quasi-cofree coalgebras stated in
lemma~\ref{lemm:QuasiCofreeHopfHom}.
\end{proof}

For our needs we give the construction of the composition product
\begin{equation*}
\HopfHom(M,N)\otimes\HopfHom(L,M)\xrightarrow{\circ}\HopfHom(L,N)
\end{equation*}
in the case of quasi-cofree connected coalgebras $M =
(T^c(V),\partial_\alpha)$ and $N = (T^c(W),\partial_\beta)$. By
definition, we have simply to determine explicitly the adjoint map
of the composite evaluation morphism
\begin{equation*}
\HopfHom(M,N)\otimes\HopfHom(L,M)\otimes
L\xrightarrow{\Id\otimes\ev}\HopfHom(M,N)\otimes
M\xrightarrow{\ev}N.
\end{equation*}
From the assertions of lemma~\ref{lemm:ConnectedQuasiCofreeHopfHom}
we obtain:

\begin{obsv}\label{obsv:ConnectedCoalgebraComposition}\index{morphism
coalgebra!composites in a} For quasi-cofree connected coalgebras $M
= (T^c(V),\partial_\alpha)$ and $N = (T^c(W),\partial_\beta)$, the
composition product
\begin{equation*}
\HopfHom(M,N)\otimes\HopfHom(L,M)\xrightarrow{\circ}\HopfHom(L,N)
\end{equation*}
can be identified with the morphism of quasi-cofree coalgebras
\begin{equation*}
\bigl(\Gamma(\DGHom(\overline{T}^c(V),W)),\partial\bigr)
\otimes\bigl(\Gamma(\DGHom(\overline{L},V)),\partial\bigr)
\xrightarrow{\circ_\Gamma}\bigl(\Gamma(\DGHom(\overline{L},W)),\partial\bigr)
\end{equation*}
induced by the composite
\begin{multline*}
\Gamma(\DGHom(\overline{T}^c(V),W))\otimes\Gamma(\DGHom(\overline{L},V))
\\
\xrightarrow{\pi\otimes\ev^\sharp_T}
\DGHom(\overline{T}^c(V),W)\otimes\DGHom(\overline{L},\overline{T}^c(V))
\xrightarrow{\circ}\DGHom(\overline{L},W)
\end{multline*}
where we consider the universal projection
\begin{equation*}
\Gamma(\DGHom(\overline{T}^c(V),W))
\xrightarrow{\pi}\DGHom(\overline{T}^c(V),W),
\end{equation*}
the morphism
\begin{equation*}
\Gamma(\DGHom(\overline{L},V))\xrightarrow{\ev^\sharp_T}\DGHom(\overline{L},\overline{T}^c(V))
\end{equation*}
and the composition product of homogeneous maps of dg-modules
\begin{equation*}
\DGHom(\overline{T}^c(V),W)\otimes\DGHom(\overline{L},\overline{T}^c(V))
\xrightarrow{\circ}\DGHom(\overline{L},W).\qed
\end{equation*}
\end{obsv}

As explained before for cofree coalgebras, this construction covers
the definition of the coalgebra morphism
\begin{equation*}
(\nabla_f)_*:
\HopfHom(L,(T^c(V),\partial_\alpha))\rightarrow\HopfHom(L,(T^c(W),\partial_\beta))
\end{equation*}
induced by a morphism of quasi-cofree connected coalgebras
$\nabla_f:
(T^c(V),\partial_\alpha)\rightarrow(T^c(W),\partial_\beta)$ since
$\nabla_f$ can be identified with a group-like element in the
morphism coalgebra
$\HopfHom((T^c(V),\partial_\alpha),(T^c(W),\partial_\beta))$.

The relation $\HopfHom(L,M) =
(\Gamma(\DGHom(\overline{L},V)),\partial)$ is also trivially
functorial in $L\in\CoAlg^a_0$. To be explicit, if $\phi:
K\rightarrow L$ is a morphism of connected coalgebras, then the
associated morphism
\begin{equation*}
\phi^*:
\HopfHom(L,(T^c(V),\partial_\alpha))\rightarrow\HopfHom(K,(T^c(V),\partial_\alpha))
\end{equation*}
can be identified with the morphism of quasi-cofree coalgebras
\begin{equation*}
\phi^*:
\bigl(\Gamma(\DGHom(\overline{L},V)),\partial_{\alpha_*\ev_T^\sharp}\bigr)
\rightarrow\bigl(\Gamma(\DGHom(\overline{K},V)),\partial_{\alpha_*\ev_T^\sharp}\bigr)
\end{equation*}
induced by the map
\begin{equation*}
\DGHom(\overline{L},V)\xrightarrow{\phi^*}\DGHom(\overline{K},V).
\end{equation*}

\medskip
Observe that the adjunction augmentations $\ev_\Gamma$ and $\ev_T$
have also an explicit description with respect to our construction
of the cofree coalgebra. Namely we have the following
straightforward assertion:

\begin{obsv}\label{obsv:CoalgebraEvaluation}
The adjunction augmentation
\begin{equation*}
\ev_\Gamma:
\Gamma(\DGHom(\overline{L},V))\otimes\overline{L}\rightarrow\overline{\Gamma}(V)
\end{equation*}
defined in observation~\ref{obsv:CoFreeCoalgebraEvaluation} can also
be obtained by a restriction of the natural morphism
\begin{equation*}
\ev_\Pi: \Bigl\{\prod_n \DGHom(\overline{L},V)^{\otimes
n}\Bigr\}\otimes\overline{L}\rightarrow\prod_{n>0} V^{\otimes n}
\end{equation*}
defined componentwise by the composite
\begin{multline*}
\Bigl\{\prod_n \DGHom(\overline{L},V)^{\otimes
n}\Bigr\}\otimes\overline{L} \\
\xrightarrow{\pr_n\otimes\Id}\DGHom(\overline{L},V)^{\otimes
n}\otimes\overline{L}
\xrightarrow{\Id\otimes\Delta^n}\DGHom(\overline{L},V)^{\otimes
n}\otimes\overline{L}{}^{\otimes n}\xrightarrow{\ev^{\otimes
n}}V^{\otimes n}.
\end{multline*}
If the coalgebra $L$ is connected, then this morphism admits a
further restriction to the adjunction augmentation
\begin{equation*}
\ev_T:
\Gamma(\DGHom(\overline{L},V))\otimes\overline{L}\rightarrow\overline{T}^c(V)
\end{equation*}
of observation~\ref{obsv:TensorCoalgebraEvaluation}.

To summarize, whenever it makes sense, we have a commutative
diagram:
\begin{equation*}
\xymatrix{
\Gamma(\DGHom(\overline{L},V))\otimes\overline{L}\ar@{-->}[d]_{\ev_T}\ar[r]^{=}
&
\Gamma(\DGHom(\overline{L},V))\otimes\overline{L}\ar@{-->}[d]_{\ev_\Gamma}\ar@{^{(}->}[]!R+<4pt,0pt>;[r]
& \Bigl\{\prod_n \DGHom(\overline{L},V)^{\otimes
n}\Bigr\}\otimes\overline{L}\ar[d]_{\ev_\Pi} \\
\overline{T}^c(V)\ar@{^{(}->}[]!R+<4pt,0pt>;[r] &
\overline{\Gamma}(V)\ar@{^{(}->}[]!R+<4pt,0pt>;[r] & \prod_{n>0}
V^{\otimes n}
}.\glossary{$\ev_\Gamma$}\glossary{$\ev_\Pi$}\glossary{$\ev_T$}\qed
\end{equation*}
\end{obsv}

Recall simply that the cofree coalgebra is endowed with a natural
embedding
\begin{equation*}
\Gamma(\DGHom(\overline{L},M))\hookrightarrow\prod_n
\DGHom(\overline{L},M)^{\otimes n}
\end{equation*}
given by the composite
\begin{equation*}
\Gamma(\DGHom(\overline{L},M))
\xrightarrow{\{\Delta^n\}}\prod_n\Gamma(\DGHom(\overline{L},M))^{\otimes
n}\xrightarrow{\{\pi^{\otimes n}\}}\prod_n
\DGHom(\overline{L},M)^{\otimes n},
\end{equation*}
where we consider the $n$-fold coproducts $\Delta^n$ and the
canonical projection $\pi$ of the cofree coalgebra
$\Gamma(\DGHom(\overline{L},M))$ (see
observation~\ref{obsv:CoFreeEmbedding}).

\subsection{The Hopf endomorphism operad of the bar complex}\label{subsection:BarHopfEndomorphismOperad}
Let $\Gamma$ be an augmented unitary coalgebra. The \emph{Hopf
endomorphism operad}\index{Hopf endomorphism operad!of a
coalgebra}\index{operad!Hopf endomorphism operad!of a coalgebra} of
$\Gamma$ is the unital Hopf operad
$\HopfEnd_\Gamma$\glossary{$\HopfEnd_\Gamma$} defined by the
morphism coalgebras $\HopfEnd_\Gamma(r) = \HopfHom(\Gamma^{\otimes
r},\Gamma)$. Observe that $\HopfEnd_\Gamma(0) = \F$. This operad
satisfies the following expected property (announced in the section
introduction):

\begin{prop}
The Hopf operad $\HopfEnd_\Gamma$ operates on the coalgebra $\Gamma$
so that $\Gamma$ forms a Hopf algebra over $\HopfEnd_\Gamma$ and the
unital operation $\HopfEnd_\Gamma(0)\rightarrow\Gamma$ agrees with
the unit of $\Gamma$.

Furthermore, the Hopf endomorphism operad $\HopfEnd_\Gamma$ is the
universal Hopf operad with this property. To be more explicit, we
have a one-to-one correspondance between Hopf algebra structures as
above and morphisms of unital Hopf operads
$\P\rightarrow\HopfEnd_\Gamma$.
\end{prop}

\begin{proof}
This proposition is a direct consequence of the adjunction relation
for the internal hom $\HopfHom(\Gamma^{\otimes r},\Gamma)$. Observe
simply that a Hopf operad action such that $\P(0)\rightarrow\Gamma$
agrees with the unit of $\Gamma$ gives rise to a commutative diagram
of coalgebras
\begin{equation*}
\xymatrix{\P(r)\ar[d]\ar[r]^{\epsilon} & \F\ar[d] \\
\P(r)\otimes\Gamma^{\otimes r}\ar[r] & \Gamma }
\end{equation*}
as in fact~\ref{fact:CoalgebraAdjunction}. Indeed recall that the
composition product $\P(r)\otimes\P(0)^{\otimes r}\rightarrow\P(0)$
is supposed to commute with coalgebra augmentations. Consequently,
the augmentation $\epsilon: \P(r)\rightarrow\F$ is equivalent to a
composite with unital operations. More precisely, we have a
commutative diagram
\begin{equation*}
\xymatrix{ \P(r)\ar[d]^{\simeq}\ar[r]^{\epsilon} & \F\ar[d]^{=} \\
\P(r)\otimes\P(0)^{\otimes r}\ar[r] & \P(0) }
\end{equation*}
from which we deduce the commutativity of the considered diagram.
\end{proof}

\subsubsection{The Hopf endomorphism operad of the bar complex}
The \emph{Hopf endomorphism operad of the bar construction in the
category of $\P$-algebras}\index{Hopf endomorphism operad!of the bar
construction}\index{operad!Hopf endomorphism operad!of the bar
construction} is defined by the formula
\begin{equation*}
\HopfEnd_B^{\P}(r) = \HopfHom_{A\in\P\Alg}(B(A)^{\otimes
r},B(A)),\glossary{$\HopfEnd_B^{\P}(r)$}
\end{equation*}
where, according to the conventions of~\cite{Bar1}, the right
hand-side denotes the end of the bifunctor $\HopfHom(B(A)^{\otimes
r},B(A))$ over the category of $\P$-algebras. As explained in the
section introduction, this operad $\HopfEnd_B^{\P}$ satisfies the
following feature by construction:

\begin{fact}[theorem~\ref{thm:BarHopfEndomorphismOperadRecall}]
The Hopf endomorphism operad $\HopfEnd_B^{\P}$ operates functorially
on the coalgebra $B(A)$ so that $B(A)$ forms a Hopf algebra over
$\HopfEnd_B^{\P}$ and the unital operation
$\HopfEnd_B^{\P}(0)\rightarrow B(A)$ is given by the unit of the bar
construction $B(A)$.

Furthermore, the Hopf endomorphism operad $\HopfEnd_B^{\P}$ is the
universal unital Hopf operad with this property. To be more
explicit, we have a one-to-one correspondance between such Hopf
actions and morphisms of Hopf operads $\Q\rightarrow\HopfEnd_B^{\P}$
that preserve unital operations.
\end{fact}

The aim of this subsection is to give more insights into the
structure of the Hopf endomorphism operad $\HopfEnd_B^{\P}$. For
this purpose, we give first structure results for the Hopf
endomorphism operad of the bar complex of a fixed $\P$-algebra
$\HopfEnd_{B(A)}$. Then we use the relation
\begin{equation*}
\HopfEnd_B^{\P}(r) = \int^{A\in\P\Alg} \HopfEnd_{B(A)}(r)
\end{equation*}
to extend our results to the Hopf endomorphism operad
$\HopfEnd_B^{\P}$.

\subsubsection{On the coalgebra structure of the Endomorphism operad of
the bar complex}\label{item:BarHopfEndomorphismOperad} First,
observe that the morphism coalgebra $\HopfHom(L,B(A))$ is
quasi-cofree, for any connected coalgebra $L$, since the bar complex
$B(A)$ is defined by a quasi-cofree connected coalgebra $B(A) =
(T^c(\Sigma A),\partial)$. More precisely, according to results of
the previous section, we have
\begin{equation*}
\HopfHom(L,B(A)) = (\Gamma(\DGHom(\overline{L},\Sigma A)),\partial),
\end{equation*}
for a coderivation $\partial: \Gamma(\DGHom(\overline{L},\Sigma
A))\rightarrow\Gamma(\DGHom(\overline{L},\Sigma A))$ determined
functorially by the differential of~$B(A)$. For $L = B(A)^{\otimes
r}$, we have identifications
\begin{multline*}
\overline{L} = \bigoplus_{m_1+\dots+m_r>0} \Sigma A^{\otimes
m_1+\dots+m_r}\\
\text{and}\quad\DGHom(\overline{L},\Sigma A) =
\prod_{m_1+\dots+m_r>0} \DGHom(\Sigma A^{\otimes
m_1+\dots+m_r},\Sigma A).
\end{multline*}
Consequently, we obtain
\begin{multline*}
\HopfEnd_{B(A)}(r) = (\Gamma(\PrimEnd_{B(A)}(r)),\partial)\\
\text{where}\quad\PrimEnd_{B(A)}(r) = \prod_{m_1+\dots+m_r>0}
\DGHom(\Sigma A^{\otimes m_1+\dots+m_r},\Sigma
A).\glossary{$\PrimEnd_{B(A)}$}
\end{multline*}

To be precise, the differential of~$L = B(A)^{\otimes r}$ is given
by the sum of the internal differential of~$\Sigma A$ together with
extra terms $\partial^h_i$ given by the bar coderivation $\partial:
B(A)\rightarrow B(A)$ on the $i$th factor of $L$. By convention we
assume that the coderivation $\partial:
\Gamma(\PrimEnd_{B(A)}(r))\rightarrow\Gamma(\PrimEnd_{B(A)}(r))$ of
$\HopfEnd_{B(A)}(r)$ includes these bar coderivations on the source
$\partial^h_i$ as well as a coderivation $\partial^v$ determined by
the bar coderivation on the target. Accordingly, we assume that the
differential $\delta:
\Gamma(\PrimEnd_{B(A)}(r))\rightarrow\Gamma(\PrimEnd_{B(A)}(r))$ is
determined by the internal differential of~$A$
only.\glossary{$\partial^h_i$}\glossary{$\partial^v$}

Observe that the action of permutations $w\in\Sigma_r$ on
$\HopfHom(B(A)^{\otimes r},B(A))$ as well as the partial composites
with unital operations $\partial_i = -\circ_i
*$ are induced by coalgebra morphisms on the source $B(A)^{\otimes r}$.
But, according to the constructions of the previous section, any
operation on the source of~$\HopfHom(L,B(A))$ can be identified with
the morphism of quasi-cofree coalgebras induced the same operation
on the source of~$\DGHom(\overline{L},B(A))$. Therefore we obtain
immediately:

\begin{prop}\label{prop:LambdaModuleEndomorphismOperad}
The modules $\PrimEnd_{B(A)}(r)$ are equipped with a
$\Lambda_*$-module structure so that the coalgebras
\begin{equation*}
\HopfEnd_{B(A)}(r) = \bigl(\Gamma(\PrimEnd_{B(A)}(r)),\partial\bigr)
\end{equation*}
form a quasi-cofree Hopf $\Lambda_*$-module.
\end{prop}

The purpose of the next paragraphs is to make explicit the structure
of the endomorphism operad $\HopfEnd_B^{\P}$ with respect to the
representation supplied by the assertion above. More specifically,
we give an explicit definition of the coderivations $\partial^h_i$
and $\partial^v$ of~$\HopfEnd_{B(A)}(r)$
in~\ref{item:BarHopfEndomorphismDifferentials}, of the
$\Lambda^*$-module structure
in~\ref{item:QuasiCofreeBarEndomorphisms} and of the operad
composition products in~\ref{item:BarHopfEndomorphismComposites}. We
obtain our results simply by going through the definition of the
endomorphism operad and the constructions
of~\ref{subsection:MorphismCoalgebras}. More specifically, for the
coderivations, we deduce the construction
of~\ref{item:BarHopfEndomorphismDifferentials} from
lemma~\ref{lemm:ConnectedQuasiCofreeHopfHom}; for $\Lambda_*$-module
operations and composition products, we deduce the constructions of
paragraphs~\ref{item:QuasiCofreeBarEndomorphisms}-\ref{item:BarHopfEndomorphismComposites}
from observation~\ref{obsv:ConnectedCoalgebraComposition} and
remarks below this statement.

But first, in the next paragraph, we give an explicit definition of
the operad evaluation product $\HopfEnd_{B(A)}(r)\otimes
B(A)^{\otimes r}\rightarrow B(A)$.

\subsubsection{On the evaluation product}\label{item:BarHopfEndomorphismEvaluation}
By definition, the evaluation product of the Hopf endomorphism
operad $\HopfEnd_{B(A)}$ is an instance of the universal evaluation
product of a morphism coalgebra
\begin{equation*}
\HopfHom(L,M)\otimes L\xrightarrow{\ev}M
\end{equation*}
in the case $L = B(A)^{\otimes r}$ and $M = B(A) = (T^c(\Sigma
A),\partial)$. We apply the explicit constructions of the previous
subsection in order to obtain the expansion of an operation
$\ev(\gamma): B(A)^{\otimes n}\rightarrow B(A)$ associated to an
element $\gamma\in\HopfEnd_{B(A)}(r)$.

First, let $\theta = \pi(\gamma)$ denotes the image of~$\gamma$
under the canonical projection $\pi:
\HopfEnd_{B(A)}(r)\rightarrow\PrimEnd_{B(A)}(r)$. Recall that
$\theta$ represents a natural transformation $\theta:
\overline{L}\rightarrow\Sigma A$ equivalent to a collection of
natural maps
\begin{equation*}
\theta = \{\theta_{m_*}\}\in\prod_{m_*} \DGHom_{A\in\P\Alg}(\Sigma
A^{\otimes m_1+\dots+m_r},\Sigma A)
\end{equation*}
that represent the components of the map~$\theta$. According to
lemma~\ref{lemm:ConnectedQuasiCofreeHopfHom} and
observation~\ref{obsv:ConnectedCoFreeHom}, the natural
transformation $\theta$ determines the composite of the
operation~$\ev(\gamma): B(A)^{\otimes n}\rightarrow B(A)$ with the
natural projection $\pi: B(A)\rightarrow\Sigma A$. Equivalently, for
each component $\Sigma A^{\otimes m_1+\dots+m_r}\subset
B(A)^{\otimes r}$, we have a commutative diagram:
\begin{equation*}
\xymatrix@C=24mm{ \Sigma A^{\otimes
m_1+\dots+m_r}\ar[r]^(0.6){\theta_{m_*}}\ar@{^{(}->}[]!D-<0pt,4pt>;[d]
& \Sigma A \\
B(A)^{\otimes r}\ar[r]^(0.6){\ev(\gamma)} & B(A)\ar[u]_{\pi} },
\end{equation*}
where we consider the natural projection $\pi: B(A)\rightarrow\Sigma
A$ of the bar complex.

Then consider the tensors
\begin{equation*}
\pi^{\otimes n}\cdot\Delta^n(\gamma) =
\sum\theta^1\otimes\dots\otimes\theta^n\in\PrimEnd_{B(A)}(r)^{\otimes
n}
\end{equation*}
associated to the $n$-fold coproduct of~$\gamma$ in
$\HopfEnd_{B(A)}(r)$. According to
observation~\ref{obsv:CoalgebraEvaluation}, the operation
$\ev(\gamma): B(A)^{\otimes r}\rightarrow B(A)$ has an expansion of
the form
\begin{equation*}
\ev(\gamma)(\alpha_1,\dots,\alpha_r) = \sum_{n=1}^{\infty}
\Bigl\{\sum\Bigl[
\underbrace{\theta^1(\alpha_1^1,\dots,\alpha_r^1)\otimes\dots
\otimes\theta^n(\alpha_1^n,\dots,\alpha_r^n)}_{\in\Sigma A^{\otimes
n}}\Bigr]\Bigr\}
\end{equation*}
for all elements in the bar complex $\alpha_1,\dots,\alpha_r\in
B(A)$, where we consider the $n$-fold diagonals $\sum
\alpha_i^1\otimes\dots\otimes\alpha_i^n\in B(A)^{\otimes n}$ of the
tensors $\alpha_i\in B(A)$.

\subsubsection{On differentials}\label{item:BarHopfEndomorphismDifferentials}\glossary{$\partial^h_i$}\glossary{$\partial^v$}
Recall that
\begin{equation*}
\HopfHom(B(A)^{\otimes r},B(A))
\xrightarrow{\partial^h_i}\HopfHom(B(A)^{\otimes r},B(A))
\end{equation*}
denotes the coderivation of~$\HopfHom(B(A)^{\otimes r},B(A))$
induced by the bar coderivation on the $i$th factor of the tensor
product $B(A)^{\otimes r}$. By construction, these coderivations
$\partial^h_i = \partial_{\beta^h_i}$ are induced by maps
\begin{equation*}
\PrimOp_{B(A)}(r)\xrightarrow{\beta^h_i}\PrimOp_{B(A)}(r).
\end{equation*}
Explicitly, for a map $\theta: \Sigma A^{\otimes
m_1+\dots+m_r}\rightarrow\Sigma A$, the components
of~$\beta^h_i(\theta)$ are given by the composite of~$\theta$ with
the natural transformations
\begin{multline*}
\Sigma A^{\otimes m_1}\otimes\dots\otimes\Sigma A^{\otimes
m_i+n-1}\otimes
\dots\otimes\Sigma A^{\otimes m_r} \\
\xrightarrow{\partial^h_i}\Sigma A^{\otimes
m_1}\otimes\dots\otimes\Sigma A^{\otimes m_i}\otimes
\dots\otimes\Sigma A^{\otimes m_r}
\end{multline*}
induced by the operations $\mu_n\in\K(n)$ on the factors
$a_k\otimes\dots\otimes a_{k+n-1}$ such that
$m_1+\dots+m_{i-1}+1\leq k<k+n-1\leq m_1+\dots+m_{i-1}+m_i+n-1$.
Accordingly, the map $\beta^h_i$ maps the component $\DGHom(\Sigma
A^{\otimes m_1+\dots+m_r},\Sigma A)$ of $\PrimEnd_{B(A)}(r)$ to the
components $\DGHom(\Sigma A^{\otimes m_1+\dots+(m_i+n-1)+m_r},\Sigma
A)$ such that $n = 2,3,\dots$.

The other coderivation of~$\HopfEnd_{B(A)}(r)$, denoted by
\begin{equation*}
\HopfHom(B(A)^{\otimes
r},B(A))\xrightarrow{\partial^v}\HopfHom(B(A)^{\otimes r},B(A)),
\end{equation*}
is yielded by the bar coderivation of the target. According to
lemma~\ref{lemm:ConnectedQuasiCofreeHopfHom} and
observation~\ref{obsv:CoalgebraEvaluation}, this coderivation can be
identified with a coderivation $\partial^v =
\partial_{\beta^v}$ of the cofree coalgebra
$\Gamma(\PrimEnd_{B(A)}(r))$ induced by a homogeneous morphism
\begin{equation*}
\xymatrix{
\Gamma(\PrimEnd_{B(A)}(r))\ar@{^{(}->}[]!R+<4pt,0pt>;[r]\ar`d[r]`/4pt[rr]_{\beta^v}`_u[rr]+D
& \prod_n\PrimEnd_{B(A)}(r)^{\otimes n}\ar[r]^{\mu_*} &
\PrimEnd_{B(A)}(r) },
\end{equation*}
where $\mu_*$ has a component
\begin{equation*}
\bigotimes_{j=1}^n \DGHom(\Sigma A^{\otimes
m^j_1+\dots+m^j_r},\Sigma A) \xrightarrow{\mu_*}\DGHom(\Sigma
A^{\otimes m_1+\dots+m_r},\Sigma A),
\end{equation*}
for all collections such that $m_i = m^1_i+\dots+m^n_i$, for $i =
1,\dots,r$. Explicitly, the morphism $\mu_*$ maps a tensor product
of homogeneous morphisms $\theta_j: \Sigma A^{\otimes
m^j_1+\dots+m^j_r}\rightarrow\Sigma A$ to the composite
\begin{multline*}
\Sigma A^{\otimes m_1+\dots+m_r}
\xrightarrow{\shuffle^*} \\
\Sigma A^{\otimes m^1_1+\dots+m^1_r}\otimes\dots\otimes\Sigma
A^{\otimes m^n_1+\dots+m^n_r}
\xrightarrow{\theta_1\otimes\dots\otimes\theta_n} \Sigma
A\otimes\dots\otimes\Sigma A \xrightarrow{\mu_n} \Sigma A
\end{multline*}
where $\shuffle^* =
\shuffle(m^j_i)\in\Sigma_{m_1+\dots+m_r}$\glossary{$\shuffle(m^j_i)$}
denotes the bloc permutation which shuffles the $rn$ tensor
groupings $\Sigma A^{\otimes m^j_i}$ according to the permutation
$\shuffle\in\Sigma_{rn}$ such that
\begin{equation*}
\shuffle((j-1)r+i) = (i-1)n+j,\quad\text{for $i = 1,\dots,r$ and $j
= 1,\dots,n$}.
\end{equation*}

\subsubsection{On $\Lambda_*$-module structures}\label{item:QuasiCofreeBarEndomorphisms} By definition, the
action of the symmetric group $\Sigma_r$ on $\HopfEnd_{B(A)}(r) =
\HopfHom(B(A)^{\otimes r},B(A))$ is given by tensor permutations on
the source. Then we observe in~\ref{item:BarHopfEndomorphismOperad}
that this action is induced by an action of $\Sigma_r$ on the
dg-module $\PrimEnd_{B(A)}(r)$. In fact, a permutation
$w\in\Sigma_r$ gives rise to a natural isomorphism
\begin{equation*}
\Sigma A^{\otimes
m_{w^{-1}(1)}+\dots+m_{w^{-1}(r)}}\xrightarrow{\simeq}\Sigma
A^{\otimes m_1+\dots+m_r}
\end{equation*}
given by the permutation of the $r$ blocs of tensors $A^{\otimes
m_i}$ specified by $w$. The composite of a homogeneous map $\theta:
\Sigma A^{m_1+\dots+m_r}\rightarrow\Sigma A$ with this isomorphism
yields a map
\begin{equation*}
w\theta: \Sigma
A^{m_{w^{-1}(1)}+\dots+m_{w^{-1}(r)}}\rightarrow\Sigma A
\end{equation*}
and this process defines the action of~$\Sigma_r$ on the dg-module
\begin{equation*}
\PrimEnd_{B(A)}(r) = \prod_{m_1+\dots+m_r>0}\DGHom(\Sigma A^{\otimes
m_1+\dots+m_r},\Sigma A).
\end{equation*}

Similarly, the morphisms $\partial_i:
\HopfEnd_{B(A)}(r)\rightarrow\HopfEnd_{B(A)}(r-1)$ are induced by
the canonical morphisms
\begin{equation*}
B(A)^{\otimes r-1} \simeq
B(A)\otimes\dots\otimes\F\otimes\dots\otimes B(A) \hookrightarrow
B(A)^{\otimes r}
\end{equation*}
given by the insertion of a unit at the $i$th position of the tensor
product. Consequently, this operation on $\HopfEnd_{B(A)}(r)$ is
induced by an operation on $\PrimEnd_{B(A)}(r)$:
\begin{equation*}
\partial_i: \HopfEnd_{B(A)}(r)\rightarrow\HopfEnd_{B(A)}(r-1).
\end{equation*}
Explicitly, this operation is given by the projection onto the
components
\begin{equation*}
\DGHom(\Sigma A^{\otimes m_1+\dots+0+\dots+m_r},\Sigma A)
\end{equation*}
for which $m_i=0$. Observe that these components of
$\PrimEnd_{B(A)}(r)$ are naturally identified with components
of~$\PrimEnd_{B(A)}(r-1)$.

\subsubsection{On composition products}\label{item:BarHopfEndomorphismComposites}
By definition, and according to the result of
observation~\ref{obsv:ConnectedCoalgebraComposition}, the partial
composition product
\begin{equation*}
\circ_i:
\HopfEnd_{B(A)}(s)\otimes\HopfEnd_{B(A)}(t)\rightarrow\HopfEnd_{B(A)}(s+t-1)
\end{equation*}
can be identified with the morphism of cofree coalgebras induced by
a morphism of dg-modules
\begin{multline*}
\PrimEnd_{B(A)}(s)\otimes\Gamma(\PrimEnd_{B(A)}(t)) \\
\hookrightarrow\PrimEnd_{B(A)}(s)\otimes\prod_m\PrimEnd_{B(A)}(t)^{\otimes
m} \xrightarrow{\gamma_i}\PrimEnd_{B(A)}(s+t-1)
\end{multline*}
which has components
\begin{multline*}
\DGHom(\Sigma A^{\otimes m_1+\dots+m_s},\Sigma A)
\otimes\Bigl\{\bigotimes_{k=1}^{m}\DGHom(\Sigma A^{n^k_1+\dots+n^k_t},\Sigma A)\Bigr\} \\
\xrightarrow{\gamma_i}\DGHom(\Sigma
A^{m_1+\dots+(n_1+\dots+n_t)+\dots+m_s},\Sigma A),
\end{multline*}
for all collections such that $m = m_i$ and $n_j =
n^1_j+\dots+n^m_j$. Explicitly, for homogeneous maps $\phi: \Sigma
A^{\otimes m_1+\dots+m_s}\rightarrow\Sigma A$ and $\psi_k: \Sigma
A^{\otimes n^k_1+\dots+n^k_t}\rightarrow\Sigma A$, $k = 1,\dots,m$,
the map $\gamma_i(\phi\otimes\{\psi_1\otimes\dots\otimes\psi_m\})$
is defined by the composite
\begin{multline*}
\Sigma A^{\otimes m_1+\dots+n_1+\dots+n_t+\dots+m_s}\\
\xrightarrow{\shuffle_i^*} \Sigma A^{\otimes m_1}\otimes\dots
\otimes\Bigl\{\bigotimes_{k=1}^{m} \Sigma A^{\otimes
n^k_1+\dots+n^k_t}\Bigr\}
\otimes\dots\otimes\Sigma A^{\otimes m_s}\\
\xrightarrow{\Id\otimes\dots\otimes\{\bigotimes_{k=1}^{m}\psi_k\}\otimes\dots\otimes\Id}
\shoveright{\Sigma A^{\otimes m_1}\otimes\dots\otimes\Sigma A^{\otimes m}\otimes\dots\otimes\Sigma A^{\otimes m_s}}\\
\xrightarrow{\phi}\Sigma A,
\end{multline*}
where $\shuffle_i^* = \shuffle_i(n^k_j)\in\Sigma_{n_1+\dots+n_t}$
denotes the bloc permutation which shuffles the $t m$ tensor
groupings $\Sigma A^{\otimes n^k_j}$ according to the permutation
$\shuffle\in\Sigma_{t m}$ defined
in~\ref{item:BarHopfEndomorphismDifferentials}.

\subsubsection{On operad units}\label{item:BarHopfEndomorphismUnits}
The operad unit of the Hopf endomorphism operad is represented by
the collection
\begin{equation*}
\{\Id^{\otimes n}\}\in\prod_n\DGHom(\Sigma A,\Sigma A)^{\otimes n},
\end{equation*}
where $\Id\in\DGHom(\Sigma A,\Sigma A)$ is the identity morphism. In
fact, this collection represents the image of the identity morphism
$\Id: B(A)\rightarrow B(A)$ under the bijection
\begin{equation*}
\Hom_{\CoAlg^a_*}(B(A),B(A))\xrightarrow{\simeq}\Gr(\HopfHom(B(A),B(A)))
\end{equation*}
defined by claim~\ref{claim:GroupLikeCoalgebraMorphisms}.

Observe that the identity morphism $\Id\in\DGHom(\Sigma A,\Sigma A)$
specifies also a unit element in $\PrimEnd_{B(A)}(1)$. Accordingly,
the module $\PrimEnd_{B(A)}$ forms a unitary $\Lambda_*$-module.

\medskip
As announced, we deduce structure results on the Hopf endomorphism
operad $\HopfEnd_B^{\P}$ from the assertions obtained in the
previous paragraphs on $\HopfEnd_{B(A)}$. This pre-statement follows
from the following formal assertion:

\begin{fact}
Consider quasi-cofree coalgebras of the form
\begin{equation*}
(\Gamma(G(X,Y),\partial_{X,Y}),
\end{equation*}
where $G: \C^{\op}\times\C\rightarrow\dg\Mod$ is a bifunctor on an
essentially small category and such that the coderivation
$\partial_{X,Y}$ is natural in $X,Y\in\C$. Let $G_{X\in\C}(X,X) =
\int^{X\in\C} G(X,X)$. The morphisms of quasi-cofree coalgebras
\begin{equation*}
(\Gamma(G_{X\in\C}(X,X),\partial_{X,X})\rightarrow
(\Gamma(G(X,X),\partial_{X,X})
\end{equation*}
induced by the natural morphism of dg-modules
$G_{X\in\C}(X,X)\rightarrow G(X,X)$ yields an end isomorphism
\begin{equation*}
(\Gamma(G_{X\in\C}(X,X),\partial_{X,X})\xrightarrow{=} \int^{X\in\C}
(\Gamma(G(X,X),\partial_{X,X}).
\end{equation*}
\end{fact}

Then we deduce from
proposition~\ref{prop:LambdaModuleEndomorphismOperad}:

\begin{prop}\index{Hopf endomorphism operad!of the bar
construction}\index{operad!Hopf endomorphism operad!of the bar
construction} The Hopf endomorphism operad $\HopfEnd_B^{\P}$ forms a
quasi-cofree Hopf $\Lambda_*$-module such that
\begin{equation*}
\HopfEnd_B^{\P}(r) =
\bigl(\Gamma(\PrimEnd_B^{\P}(r)),\partial\bigr),
\end{equation*}
where
\begin{multline*}
\PrimEnd_B^{\P}(r) = \int^{A\in\P\Alg} \PrimEnd_{B(A)}(r)\\ =
\prod_{m_1+\dots+m_r>0} \DGHom_{A\in\P\Alg}(\Sigma A^{\otimes
m_1+\dots+m_r},\Sigma A).\glossary{$\PrimEnd_B^{\P}$}
\end{multline*}
\end{prop}

The differential of~$\HopfEnd_B^{\P}$, the $\Lambda_*$-module
structure and the operad structure can also be deduced from the
constructions of
paragraphs~\ref{item:BarHopfEndomorphismDifferentials}-\ref{item:BarHopfEndomorphismComposites}
extended to natural transformations.

\subsection{The Hopf operad of universal bar operations}\label{subsection:BarHopfOperations}
\index{Hopf operad!of universal bar operations}\index{operad!Hopf
operad of universal bar operations} We use ideas of~\cite[Section
1.2]{Bar1} in order to reduce the structure of the operad
$\HopfEnd_B^{\P}$. Precisely, for an operad $\P$, we consider the
natural morphism
\begin{equation*}
\P(m)\rightarrow\DGHom_{A\in\P\Alg}(A^{\otimes m},A)
\end{equation*}
which identifies an operad element $p\in\P(m)$ with the associated
operation $p: A^{\otimes m}\rightarrow A$ defined for $A$ a
$\P$-algebra.

\begin{fact}[See \emph{loc. cit.}]
The morphism
\begin{equation*}
\P(m)\rightarrow\Hom_{A\in\P\Alg}(A^{\otimes m},A)
\end{equation*}
is split injective in general and defines an isomorphism if the
operad $\P$ is $\Sigma_*$-projective or if the ground field $\F$ is
infinite.
\end{fact}

For our purpose we consider the module
\begin{equation*}
\PrimOp_B^{\P}(r) = \prod_{m_1+\dots+m_r>0}
\Lambda\P(m_1+\dots+m_r),\glossary{$\PrimOp_B^{\P}$}
\end{equation*}
where $\Lambda\P$ denotes the operadic suspension of $\P$
(see~\cite{GetzlerJones}), and the associated cofree coalgebra
\begin{equation*}
\HopfOp_B^{\P}(r) = \Gamma(\PrimOp_B^{\P}(r)).
\end{equation*}
Recall that the suspension $\Lambda\P$\glossary{$\Lambda\P$} of an
operad $\P$ is the operad whose algebras are suspensions $\Sigma A$
of $\P$-algebras $A$. This operad satisfies the relation
$\Lambda\P(r) = \Sigma^{1-r}\P(r)\otimes\sgn(r)$, where $\sgn(r)$
denotes the signature representation of $\Sigma_r$.

We have canonical embeddings
\begin{equation*}
\Theta_{m_*}:
\Lambda\P(m_1+\dots+m_r)\hookrightarrow\DGHom_{A\in\P\Alg}(\Sigma
A^{\otimes m_1+\dots+m_r},\Sigma A)
\end{equation*}
which induce embeddings of dg-modules $\Theta:
\PrimOp_B^{\P}(r)\rightarrow\PrimEnd_B^{\P}(r)$ and embeddings of
coalgebras $\nabla_\Theta:
\HopfOp_B^{\P}(r)\rightarrow\HopfEnd_B^{\P}(r)$. These embeddings
form an isomorphism if the operad $\P$ is $\Sigma_*$-projective or
if the ground field is infinite. As a consequence, the module
$\HopfOp_B^{\P}$ is equipped with the structure of a differential
graded operad that reflects the structure of the Hopf endomorphism
operad of the bar construction. Precisely, we have the following
assertion:

\begin{prop}\label{prop:BarHopfOperationDefinition}
The coalgebras $\HopfOp_B^{\P}(r)$ are equipped with the structure
of a differential graded Hopf operad so that $\HopfOp_B^{\P}$ forms
a quasi-cofree Hopf $\Lambda_*$-module and the canonical embeddings
\begin{equation*}
\nabla_\Theta: \HopfOp_B^{\P}(r)\rightarrow\HopfEnd_B^{\P}(r)
\end{equation*}
form a morphism of differential graded Hopf operads.
\end{prop}

We check this proposition simply by going through the constructions
of~\ref{subsection:BarHopfEndomorphismOperad}. In fact, we perform
analogous constructions for $\HopfOp_B^{\P}$ so that we provide
$\HopfOp_B^{\P}$ with the structure of a suboperad of
$\HopfEnd_B^{\P}$. Then we obtain explicit definitions for the
structure of the Hopf operad $\HopfOp_B^{\P}$.

First, for the differential, we have:

\begin{claim}\label{claim:BarHopfOperationDifferentials}\hspace*{2mm}
\begin{enumerate}\renewcommand{\labelenumi}{(\alph{enumi})}
\item
The dg-coalgebra $\HopfOp_B^{\P}(r)$ can be equipped with
coderivations $\partial^h_i$ that correspond to the terms
$\partial^h_i$ of the differential of~$\HopfEnd_B^{\P}(r)$ under the
embedding $\HopfOp_B^{\P}(r)\hookrightarrow\HopfEnd_B^{\P}(r)$.
\item
The dg-coalgebra $\HopfOp_B^{\P}(r)$ can be equipped with a
coderivation $\partial^v$ that corresponds to the term $\partial^v$
of the differential of~$\HopfEnd_B^{\P}(r)$ under the embedding
$\HopfOp_B^{\P}(r)\hookrightarrow\HopfEnd_B^{\P}(r)$.
\end{enumerate}

Consequently, the dg-module $\HopfOp_B^{\P}(r)$ forms a quasi-cofree
subcoalgebra of $\HopfEnd_B^{\P}(r)$.
\end{claim}

\begin{proof}\glossary{$\partial^h_i$}\glossary{$\partial^v$}
The construction of the differentials $\partial^h_i$ and
$\partial^v$ of $\HopfOp_B^{\P}(r)$ are immediate consequences of
the assertions of~\ref{item:BarHopfEndomorphismDifferentials}.

For $i = 1,\dots,r$, the coderivation $\partial^h_i:
\HopfOp_B^{\P}(r)\rightarrow\HopfOp_B^{\P}(r)$ is induced by a map
$\beta^h_i: \PrimOp_B^{\P}(r)\rightarrow\PrimOp_B^{\P}(r)$ given
componentwise by maps
\begin{equation*}
\Lambda\P(m_1+\dots+m_i+\dots+m_r)\xrightarrow{\beta^h_i}\Lambda\P(m_1+\dots+(m_i+n-1)+\dots+m_r)
\end{equation*}
such that $\beta^h_i(p) = \sum_t p\circ_t\mu_n$, where the summation
ranges over the interval $t =
m_1+\dots+m_{i-1}+1,\dots,m_1+\dots+m_{i-1}+m_i$.

The term $\partial^v: \HopfOp_B^{\P}(r)\rightarrow\HopfOp_B^{\P}(r)$
of the differential is defined by a coderivation $\partial^v =
\partial_{\beta^v}$ of the cofree coalgebra $\HopfOp_B^{\P}(r) =
\Gamma(\PrimOp_B^{\P}(r))$ determined by a homogeneous morphism
\begin{equation*}
\xymatrix{
\Gamma(\PrimOp_B^{\P}(r))\ar@{^{(}->}[]!R+<4pt,0pt>;[r]\ar`d[r]`/4pt[rr]_{\beta^v}`_u[rr]+D
& \prod_n\PrimOp_B^{\P}(r)^{\otimes n}\ar[r]^{\mu_*} &
\PrimOp_B^{\P}(r) },
\end{equation*}
where $\mu_*$ has a component
\begin{equation*}
\bigotimes_{j=1}^n
\Lambda\P(m^j_1+\dots+m^j_r)\xrightarrow{\mu_*}\Lambda\P(m_1+\dots+m_r),
\end{equation*}
for all collections such that $m_i = m^1_i+\dots+m^n_i$. This
morphism $\mu_*$ is defined explicitly by the formula
$\mu_*(p_1\otimes\dots\otimes p_n) =
\shuffle(m^j_i)\cdot\mu_n(p_1,\dots,p_n)$, where
$\shuffle(m^j_i)\in\Sigma_{m_1+\dots+m_r}$ denotes the bloc
permutation introduced
in~\ref{item:BarHopfEndomorphismDifferentials}.
\end{proof}

Then for the $\Lambda_*$-module structure:

\begin{claim}\label{claim:BarOperationsUnitaryComposites}
The modules $\PrimOp_B^{\P}(r)$ can be equipped with an action of
the symmetric group $\Sigma_r$ and with operations $\partial_i:
\PrimOp_B^{\P}(r)\rightarrow\PrimOp_B^{\P}(r-1)$ so that
$\PrimOp_B^{\P}$ forms a $\Lambda_*$-submodule of~$\PrimEnd_B^{\P}$.

The module $\PrimOp_B^{\P}(1)$ is also equipped with a unit element
$1\in\PrimOp_B^{\P}(1)$ that corresponds to the unit element
of~$\PrimEnd_B^{\P}(1)$ (see~\ref{item:BarHopfEndomorphismUnits}).
Accordingly, the module $\PrimOp_B^{\P}$ forms a unitary
$\Lambda_*$-submodule of~$\PrimEnd_B^{\P}$.
\end{claim}

\begin{proof}
This result is also an immediate consequence of the explicit
definition of the $\Lambda_*$-module structure of $\PrimEnd_B^{\P}$
given in~\ref{item:QuasiCofreeBarEndomorphisms}

Explicitly, the action of a permutation $w\in\Sigma_r$ on
\begin{equation*}
\PrimOp_B^{\P}(r) = \prod_{m_1+\dots+m_r>0} \Lambda\P(m_1+\dots+m_r)
\end{equation*}
is given componentwise by the action of the bloc permutations
\begin{equation*}
\Lambda\P(m_1+\dots+m_r)
\xrightarrow{w(m_{w^{-1}(1)},\dots,m_{w^{-1}(r)})}\Lambda\P(m_{w^{-1}(1)}+\dots+m_{w^{-1}(r)})
\end{equation*}
and the operation
\begin{equation*}
\partial_i: \PrimOp_B^{\P}(r)\rightarrow\PrimOp_B^{\P}(r-1)
\end{equation*}
is given by the projection onto the components
of~$\PrimOp_B^{\P}(r)$ such that $m_i=0$ which can be identified
with components of $\PrimOp_B^{\P}(r-1)$ as in the case of
$\PrimEnd_B^{\P}(r)$.

The unit operation $1\in\Lambda\P(1)$ corresponds tautologically to
the identity morphism $\Id: \Sigma A\rightarrow\Sigma A$ that
defines the unit element of the module $\PrimEnd_B^{\P}(1)$ and
hence specifies an appropriate unit element in $\PrimOp_B^{\P}(1)$.
\end{proof}

Finally, for the operad structure:

\begin{claim}\label{claim:BarHopfOperationComposites}
The dg-coalgebras $\HopfOp_B^{\P}(r) = \Gamma(\PrimOp_B^{\P}(r))$
are equipped with an action of the symmetric group induced by the
action of $\Sigma_r$ on $\PrimOp_B^{\P}(r)$. We have also
composition products
\begin{equation*}
\circ_i:
\HopfOp_B^{\P}(s)\otimes\HopfOp_B^{\P}(t)\rightarrow\HopfOp_B^{\P}(s+t-1)
\end{equation*}
so that $\HopfOp_B^{\P}$ forms a Hopf suboperad of
$\HopfEnd_B^{\P}$.
\end{claim}

\begin{proof}
As in the case of the Hopf endomorphism operad $\HopfEnd_B^{\P}$,
the partial composition product
\begin{equation*}
\circ_i:
\HopfOp_B^{\P}(s)\otimes\HopfOp_B^{\P}(t)\rightarrow\HopfOp_B^{\P}(s+t-1)
\end{equation*}
is the morphism of cofree coalgebras induced by a morphism of
dg-modules
\begin{multline*}
\PrimOp_B^{\P}(s)\otimes\Gamma(\PrimOp_B^{\P}(t)) \\
\hookrightarrow\PrimOp_B^{\P}(s)\otimes\prod_m\PrimOp_B^{\P}(t)^{\otimes
m} \xrightarrow{\gamma_i}\PrimOp_B^{\P}(s+t-1)
\end{multline*}
which has components
\begin{multline*}
\Lambda\P(m_1+\dots+m_s)\otimes\Bigl\{\bigotimes_{k=1}^{m}\Lambda\P(n^k_1+\dots+n^k_t)\Bigr\} \\
\xrightarrow{\gamma_i}\Lambda\P(m_1+\dots+n_1+\dots+n_t+\dots+m_s),
\end{multline*}
for all collections such that $m = m_i$ and $n_j =
n^1_j+\dots+n^m_j$. We have explicitly
$\gamma_i(p\otimes\{q_1\otimes\dots\otimes q_m\}) =
\shuffle_i(n^k_j)\cdot p(1,\dots,q_1,\dots,q_m,\dots,1)$, where the
operation $q_k$ is substituted to the entry $k =
m_1+\dots+m_{i-1}+k$ of $p$ and where $\shuffle_i(n^k_j)$ denotes
the bloc permutation introduced
in~\ref{item:BarHopfEndomorphismComposites}.

Observe also that the operad unit of~$\HopfEnd_B^{\P}$ corresponds
to an element of~$\HopfOp_B^{\P}$. Namely this element can be
represented by the collection
\begin{equation*}
\{1^{\otimes n}\}\in\prod_n\PrimOp_B^{\P}(1)^{\otimes n},
\end{equation*}
where $1\in\Lambda\P(1)$ is the unit of~$\P$.
\end{proof}

This claim achieves the proof of
proposition~\ref{prop:BarHopfOperationDefinition}. To recapitulate,
we have the following result:

\begin{lemm}\label{lemm:BarHopfOperationDefinition}
The dg-coalgebras $\HopfOp_B^{\P}(r)$ equipped with the
differentials supplied by
claim~\ref{claim:BarHopfOperationDifferentials} and the structure
specified by claim~\ref{claim:BarHopfOperationComposites} form a
differential graded unital Hopf operad. Moreover, the canonical
embeddings
\begin{equation*}
\nabla_\Theta: \HopfOp_B^{\P}(r)\hookrightarrow\HopfEnd_B^{\P}(r)
\end{equation*}
define a natural morphism of unital Hopf operads which is an
isomorphism if the operad $\P$ is $\Sigma_*$-projective or if the
ground field $\F$ is infinite.
\end{lemm}
In addition we obtain the following assertion:

\begin{obsv}
The Hopf operad $\HopfOp_B^{\P}(r)$ forms a quasi-cofree Hopf
$\Lambda_*$-module such that
\begin{equation*}
\HopfOp_B^{\P}(r) = \bigl(\Gamma(\PrimOp_B^{\P}(r)),\partial\bigr),
\end{equation*}
where $\partial$ is composed of the coderivations $\partial^v$ and
$\partial^h_i$ specified
in~\ref{claim:BarHopfOperationDifferentials}.\qed
\end{obsv}

\subsection{Fibration properties}\label{subsection:FibrationProperties}
One can deduce from formal properties of monoidal model categories
that a Hopf endomorphism operad $\HopfEnd_\Gamma$ is a Reedy fibrant
object if $\Gamma$ is a fibrant unitary dg-coalgebras. One can
extend this assertion to connected unitary coalgebras because the
adjoint definition of morphism coalgebras holds in the connected
context as well by fact~\ref{fact:ConnectedCoalgebraAdjunction}. As
any quasi-cofree connected coalgebra $\Gamma =
(T^c(V),\partial_\alpha)$ forms a fibrant object in the category of
connected coalgebras, we obtain that the Hopf endomorphism operad of
the bar complex $\HopfEnd_{B(A)}$ is Reedy fibrant.

In the next paragraphs, we give another more effective proof of this
assertion in order to extend our results to the Hopf operad of bar
operations $\HopfOp_B^{\P}$. Explicitly, as mentioned in the
introduction of~\ref{section:CocellularCoalgebras}, we prove that
the augmentation morphism of the Hopf endomorphism operad
$\HopfEnd_{B(A)}$ splits up into a sequence
\begin{multline*}
\HopfEnd_{B(A)} = \lim_m\ckcell_m\HopfEnd_{B(A)}
\rightarrow\dots\\
\dots\rightarrow\ckcell_{m}\HopfEnd_{B(A)}
\rightarrow\ckcell_{m-1}\HopfEnd_{B(A)}\rightarrow\dots\\
\dots\rightarrow\ckcell_0\HopfEnd_{B(A)} = \C,
\end{multline*}
where
$\ckcell_{m}\HopfEnd_{B(A)}\rightarrow\ckcell_{m-1}\HopfEnd_{B(A)}$
is obtained by a coextension of a cofree Hopf $\Lambda_*$-modules
morphism. Then we extend this result to the Hopf endomorphism operad
of the bar complex $\HopfEnd_B^{\P}$ and to the Hopf operad of bar
operations~$\HopfOp_B^{\P}$.

\subsubsection{On $\Lambda_*$-module structures and the canonical filtration of the bar complex}
Recall that the $\Lambda_*$-module structure of a Hopf endomorphism
operad $\HopfEnd_\Gamma(r) = \HopfHom(\Gamma^{\otimes r},\Gamma)$ is
deduced from operations on the source. To be more precise, one can
observe that the tensor powers of a unitary coalgebra
$\Gamma^{\otimes r}$ form a left Hopf $\Lambda^*$-module. The
symmetric group operates on $\Gamma^{\otimes r}$ by tensor
permutations and we have operations $\partial^i: \Gamma^{\otimes
r-1}\rightarrow\Gamma^{\otimes r}$ given by the insertion of a
coalgebra unit $*: \F\rightarrow\Gamma$ at the $i$th place. The
induced operations $w^*: \HopfHom(\Gamma^{\otimes
r},\Gamma)\rightarrow\HopfHom(\Gamma^{\otimes r},\Gamma)$ and
$(\partial^i)^*: \HopfHom(\Gamma^{\otimes
r},\Gamma)\rightarrow\HopfHom(\Gamma^{\otimes r-1},\Gamma)$
determine the $\Lambda^*$-module structure
of~$\HopfHom(\Gamma^{\otimes r},\Gamma)$.

The quotient Hopf $\Lambda_*$-modules $\ckcell_{m}\HopfEnd_{B(A)}$
are associated to a sequence of Hopf $\Lambda_*$-submodules
of~$B(A)^{\otimes r}$. Explicitly, for $m\in\N$, we let
$\skcell_m(B(A)^{\otimes r})$ denote the submodule of~$B(A)^{\otimes
r}$ such that:
\begin{equation*}
\skcell_m(B(A)^{\otimes r}) = \bigoplus_{m\geq m_1+\dots+m_r\geq 0}
\Sigma A^{\otimes m_1+\dots+m_r}.
\end{equation*}
Recall that $B(A)^{\otimes r}$ is equipped with coderivations
denoted by $\partial^h_i$ and given by the bar coderivation of the
$i$th factor of~$B(A)^{\otimes r}$, for $i = 1,\dots,r$. Clearly,
the submodule $\skcell_m B(A)^{\otimes r}$ is preserved by the bar
coderivations $\partial^h_i$ and hence form a dg-submodule of
$B(A)$. In fact, we have more precisely:

\begin{obsv}\label{obsv:DifferentialBarFiltration}
We have $\partial^h_i(\skcell_m(B(A)^{\otimes
r}))\subset\skcell_{m-1}(B(A)^{\otimes r})$.
\end{obsv}

\begin{proof}
By definition, the bar coderivation $\partial^h_i$ maps the module
$\Sigma A^{\otimes m_i}\subset B(A)$ into components $\Sigma
A^{\otimes m_i-n+1}\subset B(A)$ such that $n\geq 2$. Hence the
assertion is immediate.
\end{proof}

Then, as expected, we have clearly:

\begin{fact}
The dg-modules
\begin{equation*}
\skcell_m(B(A)^{\otimes r})\subset B(A)^{\otimes r}
\end{equation*}
are preserved by the diagonal and by the differential
of~$B(A)^{\otimes r}$, by the action of the symmetric group and by
the operations $\partial^i: B(A)^{\otimes r}\rightarrow
B(A)^{\otimes r-1}$. Hence these dg-modules define a nested sequence
of Hopf $\Lambda_*$-submodules of~$B(A)^{\otimes r}$ such that
\begin{equation*}
B(A)^{\otimes r} = \colim_m\skcell_m(B(A)^{\otimes r}).
\end{equation*}
\end{fact}

As a corollary we obtain:

\begin{fact}
The coalgebras
\begin{equation*}
\ckcell_m\HopfEnd_{B(A)}(r) = \HopfHom(\skcell_m B(A)^{\otimes
r},B(A))
\end{equation*}
define a tower of quotient Hopf $\Lambda_*$-modules
of~$\HopfEnd_{B(A)}$ such that
\begin{equation*}
\HopfEnd_{B(A)} = \mylim_m\ckcell_m\HopfEnd_{B(A)}.
\end{equation*}
Furthermore, we have
\begin{equation*}
\ckcell_m\HopfEnd_{B(A)} =
(\Gamma(\ckcell_m\PrimEnd_{B(A)}),\partial),
\end{equation*}
where $\ckcell_m\PrimEnd_{B(A)}$ denotes the quotient
$\Lambda_*$-module of~$\PrimEnd_{B(A)}$ defined by
\begin{equation*}
\ckcell_m\PrimEnd_{B(A)}(r) = \prod_{m\geq
m_1+\dots+m_r>0}\DGHom(\Sigma A^{\otimes m_1+\dots+m_r},\Sigma A).
\end{equation*}
\end{fact}

Observe that $\ckcell_m\HopfEnd_{B(A)}$ forms also a unital unitary
quotient Hopf $\Lambda_*$-module of~$\HopfEnd_{B(A)}$. Explicitly,
we have clearly $\ckcell_m\HopfEnd_{B(A)}(0) = \F$. Furthermore, the
collection
\begin{equation*}
\{\Id^{\otimes n}\}\in\prod_n\DGHom(\Sigma A,\Sigma A)^{\otimes n}
\end{equation*}
that represents the unit element of~$\HopfEnd_{B(A)}$ specifies
clearly a unit element in~$\ckcell_m\HopfEnd_{B(A)}$ for all $m\geq
1$ and projects tautologically to the unit element of the
commutative operad for $m = 0$.

\medskip
Clearly, the Hopf $\Lambda_*$-module
\begin{equation*}
\HopfEnd_B^{\P}(r) = \HopfHom_{A\in\P\Alg}(B(A)^{\otimes r},B(A))
\end{equation*}
is equipped with a similar decomposition induced by the
decomposition of~$\HopfEnd_{B(A)}$. Explicitly, we have
$\HopfEnd_B^{\P} = \mylim_{m\in\N}\ckcell_m\HopfEnd_B^{\P}$ for
unital unitary Hopf $\Lambda_*$-modules such that:
\begin{equation*}
\ckcell_m\HopfEnd_B^{\P}(r) = \HopfHom_{A\in\P\Alg}(\skcell_m
B(A)^{\otimes r},B(A)).
\end{equation*}
Furthermore, we have
\begin{equation*}
\ckcell_m\HopfEnd_B^{\P} =
(\Gamma(\ckcell_m\PrimEnd_B^{\P}),\partial),
\end{equation*}
for unitary quotient $\Lambda_*$-modules of $\PrimEnd_B^{\P}$, that
can be defined by
\begin{equation*}
\ckcell_m\PrimEnd_{B(A)}(r) = \prod_{m\geq
m_1+\dots+m_r>0}\DGHom_{A\in\P\Alg}(\Sigma A^{\otimes
m_1+\dots+m_r},\Sigma A).
\end{equation*}

\medskip
We check that the requirements
of~\ref{item:CofreeLambdaModuleCoextensions} are satisfied for the
Hopf $\Lambda_*$-modules
$(\Gamma(\ckcell_m\PrimEnd_{B(A)}),\partial)$ so that the morphisms
\begin{equation*}
\pk_m:
(\Gamma(\ckcell_{m}\PrimEnd_{B(A)}),\partial)\rightarrow(\Gamma(\ckcell_{m-1}\PrimEnd_{B(A)}),\partial)
\end{equation*}
have the structure specified
in~\ref{subsection:CocellularHopfLambdaModules} and similarly for
the functorial Hopf endomorphism operad $\HopfEnd_B^{\P}$. First,
for~$\HopfEnd_{B(A)}$, we have explicitly:

\begin{claim}\label{claim:HopfEndomorphismsDecomposition}
Let $K = \PrimEnd_{B(A)}$. The homogeneous maps $\beta =
\beta^h_i,\beta^v: \Gamma(K)\rightarrow K$ that determine the
coderivations $\partial = \partial^h_i,\partial^v$ of
$\HopfEnd_{B(A)}$ admit factorizations
\begin{equation*}
\xymatrix{ \Gamma(K)\ar[d]\ar[r]^{\beta} &
\Pi\ar[d] \\
\Gamma(\ckcell_m K)\ar@{-->}[r]^{\beta} & \ckcell_{m-1} K }.
\end{equation*}
As a corollary, the projection morphisms
\begin{equation*}
\pk_m: (\Gamma(\ckcell_m
K),\partial)\rightarrow(\Gamma(\ckcell_{m-1} K),\partial)
\end{equation*}
fit coextension diagrams of the form
\begin{equation*}
\xymatrix{ (\Gamma(\ckcell_m K),\partial)\ar[r]\ar[d]
& \Gamma(\Delta^1\wedge\ckcell_m K)\ar[d] \\
(\Gamma(\ckcell_{m-1} K),\partial)\ar[r] & \Gamma(\Delta^1\wedge
\ckcell_{m-1} K\times_{S^1\wedge\ckcell_{m-1} K} S^1\wedge\ckcell_m
K) }.
\end{equation*}
\end{claim}

\begin{proof}
For the coderivations $\partial = \partial^h_i$, the assertion is a
corollary of observation~\ref{obsv:DifferentialBarFiltration} since
these coderivations are induced by the components $\partial^h_i$ of
the differentials of~$B(A)^{\otimes r}$. Explicitly, by
observation~\ref{obsv:DifferentialBarFiltration}, the
maps~$\partial^h_i$ have factorizations
\begin{equation*}
\xymatrix{ B(A)^{\otimes r}\ar[r]^{\partial^h_i} & B(A)^{\otimes r} \\
\skcell_{m} B(A)^{\otimes r}\ar@{-->}[r]\ar[u] & \skcell_{m-1}
B(A)^{\otimes r}\ar[u] }
\end{equation*}
which yield a factorization at the level of~$\PrimEnd_{B(A)}$:
\begin{equation*}
\xymatrix{ \PrimEnd_{B(A)}\ar[r]^{\beta^h_i} & \PrimEnd_{B(A)} \\
\ckcell_{m}\PrimEnd_{B(A)}\ar@{-->}[r]\ar[u] &
\ckcell_{m-1}\PrimEnd_{B(A)}\ar[u] }.
\end{equation*}

Recall that the other differential $\partial = \partial^v$ is
induced by homogeneous morphisms
\begin{equation*}
\Gamma(\PrimEnd_{B(A)}(r))\hookrightarrow\prod_n
\PrimEnd_{B(A)}(r)^{\otimes n}\xrightarrow{\mu_*}\PrimEnd_{B(A)}(r)
\end{equation*}
which admit a component
\begin{equation*}
\bigotimes_{j=1}^n \DGHom(\Sigma A^{\otimes
m^j_1+\dots+m^j_r},\Sigma A)\xrightarrow{\mu_*}\DGHom(\Sigma
A^{\otimes m_1+\dots+m_r},\Sigma A)
\end{equation*}
for all $n\geq 2$ and for all collections $(m^j_i)$ such that $m_i =
m^1_i+\dots+m^n_i$. Clearly, if $m^j_1+\dots+n^j_r>m-1$ for some
$j$, as we assume $m^k_1+\dots+n^k_r>0$ for all $k$, we obtain
$m_1+\dots+m_r = \sum_{i,j} m^j_i>m^j_1+\dots+n^j_r>m-1$.

Accordingly, for a given tensor $u_1\otimes\dots\otimes
u_n\in\PrimEnd_{B(A)}(r)^{\otimes n}$, if we have
$u_j\in\prod_{m^j_1+\dots+m^j_r>m-1} \DGHom(\Sigma
A^{m^j_1+\dots+m^j_r},\Sigma A)$ for some $j$, then we obtain
$\mu_*(u_1\otimes\dots\otimes u_n)\in\prod_{m_1+\dots+m_r>m}
\DGHom(\Sigma A^{m^j_1+\dots+m^j_r},\Sigma A)$ so that $\mu_*$
admits a factorization
\begin{equation*}
\xymatrix{\prod_n \PrimEnd_{B(A)}(r)^{\otimes n}\ar[r]^{\mu_*}\ar[d] & \PrimEnd_{B(A)}(r)\ar[d] \\
\ckcell_{m-1}\PrimEnd_{B(A)}(r)^{\otimes n}\ar@{-->}[r]^{\mu_*} &
\ckcell_m\PrimEnd_{B(A)}(r) }
\end{equation*}
and this assertion implies our claim for the differential
$\partial^v$.
\end{proof}

Observe that the lifting construction of the claim is natural
in~$A$. To give a more proper assertion we should extend this
structure result to Hopf $\Lambda_*$-modules of morphisms
$\HopfEnd_{B(A),B(A')}$, which are defined by
\begin{equation*}
\HopfEnd_{B(A),B(A')}(r) = \HopfHom(B(A)^{\otimes r},B(A')).
\end{equation*}
One checks precisely that the coderivation liftings
\begin{equation*}
\xymatrix{
\Gamma(\ckcell_m\PrimEnd_{B(A),B(A')})\ar@{-->}[r]^{\beta} &
\ckcell_{m-1}\PrimEnd_{B(A),B(A')} }
\end{equation*}
are functorial in $A$ and $A'$. In fact, this assertion holds simply
because the modules $\ckcell_m\PrimEnd_{B(A),B(A')}$ are quotient
of~$\PrimEnd_{B(A),B(A')}$. As a consequence, we obtain:

\begin{fact}
The coderivation liftings of
claim~\ref{claim:HopfEndomorphismsDecomposition} induce a
coderivation lifting
\begin{equation*}
\xymatrix{ \Gamma(\ckcell_m\PrimEnd_B^{\P})\ar@{-->}[r]^{\beta} &
\ckcell_{m-1}\PrimEnd_B^{\P} }
\end{equation*}
on the end of the $\Lambda_*$-modules $\ckcell_m\PrimEnd_{B(A)}$. As
a consequence, the results of
claim~\ref{claim:HopfEndomorphismsDecomposition} hold for the
$\Lambda_*$-module $K = \PrimEnd_B^{\P}$ and the quotient Hopf
$\Lambda_*$-modules
\begin{equation*}
\ckcell_m\HopfEnd_B^{\P} =
(\Gamma(\ckcell_m\PrimEnd_B^{\P}),\partial)
\end{equation*}
of the Hopf endomorphism operad $\HopfEnd_B^{\P}$.
\end{fact}

We prove now that the Hopf operad of universal operations
$\HopfOp_B^{\P}$ is equipped with the same decomposition as the Hopf
endomorphism operad $\HopfEnd_B^{\P}$. We consider the dg-modules
\begin{equation*}
\ckcell_m\PrimOp_B^{\P}(r) = \prod_{m\geq m_1+\dots+m_r>0}
\Lambda\P(m_1+\dots+m_r)
\end{equation*}
which come equipped with a canonical embedding
\begin{equation*}
\ckcell_m\Theta:
\ckcell_m\PrimOp_B^{\P}(r)\hookrightarrow\ckcell_m\PrimEnd_B^{\P}(r).
\end{equation*}
These modules equipped with the canonical projections
\begin{equation*}
\pk_m:
\ckcell_m\PrimOp_B^{\P}(r)\rightarrow\ckcell_{m-1}\PrimOp_B^{\P}(r)
\end{equation*}
form clearly a subtower of~$\ckcell_m\PrimEnd_B^{\P}(r)$ such that
\begin{equation*}
\PrimOp_B^{\P}(r) = \mylim_m\ckcell_m\PrimOp_B^{\P}(r).
\end{equation*}
We have in addition:

\begin{fact}
The dg-modules $\ckcell_m\PrimOp_B^{\P}(r)$ defined above can be
identified with the image of~$\PrimOp_B^{\P}(r)$ under the composite
map
\begin{equation*}
\PrimOp_B^{\P}(r)\hookrightarrow\PrimEnd_B^{\P}(r)\rightarrow\ckcell_m\PrimEnd_B^{\P}(r)
\end{equation*}
\end{fact}

As a corollary, we obtain the following assertions:

\begin{fact}\hspace*{2mm}

\begin{enumerate}
\item
The dg-modules $\ckcell_m\PrimOp_B^{\P}(r)$ form unitary quotient
$\Lambda_*$-modules of~$\PrimOp_B^{\P}$ so that the canonical
embedding $\Theta: \PrimOp_B^{\P}\hookrightarrow\PrimEnd_B^{\P}$
splits up into an embedding of $\Lambda_*$-module towers
\begin{equation*}
\{\ckcell_m\PrimOp_B^{\P}\}_m\hookrightarrow\{\ckcell_m\PrimEnd_B^{\P}\}_m.
\end{equation*}
\item
The coalgebras $\Gamma(\ckcell_m\PrimOp_B^{\P}(r))$ are preserved by
the differential of $\HopfOp_B^{\P}$ and define unital unitary
quotient Hopf $\Lambda_*$-modules of~$\HopfOp_B^{\P}$. Accordingly,
the canonical morphism $\nabla_\Theta:
\HopfOp_B^{\P}\hookrightarrow\HopfEnd_B^{\P}$ splits up into a
morphism of Hopf $\Lambda_*$-module towers
\begin{equation*}
\{\ckcell_m\HopfOp_B^{\P}\}_m\hookrightarrow\{\ckcell_m\HopfEnd_B^{\P}\}_m,
\end{equation*}
where $\ckcell_m\HopfOp_B^{\P} =
(\Gamma(\ckcell_m\PrimOp_B^{\P}),\partial)$.
\end{enumerate}
\end{fact}

We have in addition:

\begin{fact}
The assertion of claim~\ref{claim:HopfEndomorphismsDecomposition}
holds for $K = \PrimOp_B^{\P}$. To be more precise, let $K' =
\PrimEnd_B^{\P}$. The factorization of the homogeneous maps $\beta =
\beta^h_i,\beta^v: \Gamma(K')\rightarrow K'$ that determine the
coderivations $\partial =
\partial^h_i,\partial^v$ of the quasi-cofree Hopf
$\Lambda_*$-module~$\HopfEnd_B^{\P} = (\Gamma(K'),\partial)$ admit a
restriction to $\Gamma(\ckcell_m K)$.

As a corollary, the projection morphisms
\begin{equation*}
\pk_m:
(\Gamma(\ckcell_m\PrimOp_B^{\P}),\partial)\rightarrow(\Gamma(\ckcell_{m-1}\PrimOp_B^{\P}),\partial)
\end{equation*}
can be obtained by coextensions of the form specified
in~\ref{subsection:CocellularHopfLambdaModules} like the projection
morphism of Hopf endomorphism operads $\HopfEnd_B^{\P}$ and
$\HopfEnd_{B(A)}$.
\end{fact}

This assertion can either be deduced from the relationship between
the Hopf operad of bar operations $\HopfOp_B^{\P}$ and the Hopf
endomorphism operad $\HopfEnd_B^{\P}$ or can be checked directly as
in the proof of claim~\ref{claim:HopfEndomorphismsDecomposition} for
the Hopf endomorphism operad $\HopfEnd_{B(A)}$.

We use the decomposition $\HopfOp_B^{\P} =
\mylim_m\ckcell_m\HopfOp_B^{\P}$ and the results
of~\ref{section:CocellularCoalgebras} in order to prove the
fibration properties asserted by
theorem~\ref{thm:HopfEndomorphismOperadFibration}. Namely we prove:

\begin{claim}[theorem~\ref{thm:HopfEndomorphismOperadFibration}]\label{claim:HopfBarOperationsFibrationRecall}
The morphism
\begin{equation*}
\phi_*: \HopfOp_B^{\P}\rightarrow\HopfOp_B^{\P'}
\end{equation*}
induced by a fibration, respectively an acyclic fibration, of
non-unital operads $\phi: \P\rightarrow\P'$ forms a Reedy fibration,
respectively an acyclic Reedy fibration, of Hopf
$\Lambda_*$-modules.
\end{claim}

For this purpose we check first the following statement:

\begin{obsv}\label{obsv:LambdaModuleFibrations}
If $\phi: \P\rightarrow\P'$ is a fibration, respectively an acyclic
fibration, of non-unital operads, then the map
\begin{equation*}
\ckcell_m\PrimOp_B^{\P} \xrightarrow{(\pk_m,\ckcell_m\phi_*)}
\ckcell_{m-1}\PrimOp_B^{\P}\times_{\ckcell_{m-1}\PrimOp_B^{\P'}}\ckcell_m\PrimOp_B^{\P'}
\end{equation*}
is a Reedy fibration, respectively an acyclic Reedy fibration, of
$\Lambda_*$-modules.
\end{obsv}

\begin{proof}
Recall that a morphism of $\Lambda_*$-modules $f: M\rightarrow N$
forms a Reedy fibration, respectively an acyclic Reedy fibration, if
the morphisms $(\mu,f):
M(r)\rightarrow\Match{M}(r)\times_{\Match{N}(r)} N(r)$ are
fibrations, respectively acyclic fibrations, of dg-modules, for all
$r\in\N$.

The matching object of the $\Lambda_*$-module $M =
\ckcell_m\PrimOp_B^{\P}$ can clearly be identified with the
restricted product
\begin{equation*}
\Match{\PrimOp_B^{\P}}(r) = {\prod_{m_*}}' \Lambda\P(m_1+\dots+m_r)
\end{equation*}
which ranges over the collections $m\geq m_1+\dots+m_r>0$ such that
$m_i = 0$ for some~$i$. For $M = \ckcell_m\PrimOp_B^{\P}$ and
\begin{equation*}
N =
\ckcell_{m-1}\PrimOp_B^{\P}\times_{\ckcell_{m-1}\PrimOp_B^{\P'}}\linebreak\ckcell_m\PrimOp_B^{\P'},
\end{equation*}
the relative matching object $\Match{M}(r)\times_{\Match{N}(r)}
N(r)$ can be identified with the direct product
\begin{equation*}
{\prod_{m_*}}'\Lambda\P(m_1+\dots+m_r)\times{\prod_{m_*}}''\Lambda\P'(m_1+\dots+m_r),
\end{equation*}
where $\prod'_{m_*}$ ranges over all collections $m\geq
m_1+\dots+m_r>0$ such that $m_i = 0$ for some $i$ or $m-1\geq
m_1+\dots+m_r>0$ and $\prod''_{m_*}$ ranges over all collections
$m\geq m_1+\dots+m_r>0$ such that $m_i>0$ for $i = 1,\dots,r$.
Moreover, the relative matching morphism $(\phi,\mu)$ is given
componentwise by the identity of $\Lambda\P(m_1+\dots+m_r)$ or by
the morphism $\phi:
\Lambda\P(m_1+\dots+m_r)\rightarrow\Lambda\P'(m_1+\dots+m_r)$.
Accordingly, this morphism forms clearly a fibration, respectively
an acyclic fibration, if $\phi$ is so.
\end{proof}

Then we deduce an inductive proof of the properties of
claim~\ref{claim:HopfBarOperationsFibrationRecall} from the results
of~\ref{section:CocellularCoalgebras}. Explicitly, we prove the
following statement:

\begin{claim}
Under the assumptions of
claim~\ref{claim:HopfBarOperationsFibrationRecall}, the morphism
\begin{equation*}
\ckcell_m\phi_*:
\ckcell_m\HopfOp_B^{\P}\rightarrow\ckcell_m\HopfOp_B^{\P'}
\end{equation*}
is a Reedy fibration, respectively an acyclic Reedy fibration of
$\Lambda_*$-modules, for all $m\in\N$.
\end{claim}

\begin{proof}
We apply
lemma~\ref{lemm:DiagramCofreeCoextensionLambdaModuleFibrations} to
the commutative square
\begin{equation*}
\xymatrix@C=18mm{
\ckcell_{m}\HopfOp_B^{\P}\ar[d]^{\pk_m}\ar[r]^{\ckcell_m\phi_*} &
\ckcell_{m}\HopfOp_B^{\P'}\ar[d]^{\pk_m} \\
\ckcell_{m-1}\HopfOp_B^{\P}\ar[r]^{\ckcell_{m-1}\phi_*} &
\ckcell_{m-1}\HopfOp_B^{\P'} }.
\end{equation*}
Under the assertion of
observation~\ref{obsv:LambdaModuleFibrations}, we obtain that
$\ckcell_{m}\phi_*$ forms a fibration, respectively an acyclic
fibration if $\ckcell_{m-1}\phi_*$ is so. Hence the claim follows by
induction.
\end{proof}

This result achieves the proof of
claim~\ref{claim:HopfBarOperationsFibrationRecall} and
theorem~\ref{thm:HopfEndomorphismOperadFibration}.

\subsection{Universal bar operations for the commutative operad}\label{subsection:CommutativeHopfBarAction}
As recalled in the memoir and section introductions, the bar complex
$B(A)$ of a commutative algebra $A$ is equipped with an associative
and commutative product $\smile: B(A)\otimes B(A)\rightarrow B(A)$,
the shuffle product of tensors, that provides $B(A)$ with the
structure of an associative and commutative differential graded Hopf
algebra. This product is equivalent to a morphism of unital Hopf
operads $\nabla_c: \C\rightarrow\HopfEnd_B^{\C}$.

In the introduction we claim also that:

\begin{lemm}\label{lemm:ShuffleStructure}
The morphism $\nabla_c: \C\rightarrow\HopfEnd_B^{\C}$ admits a
factorization
\begin{equation*}
\xymatrix{ \C\ar@{-->}[dr]!UL_{\nabla_\gamma}\ar[r]^(0.35){\nabla_c} & \HopfEnd_B^{\C} \\
& \HopfOp_B^{\C}\ar@{^{(}->}[]!U+<0pt,4pt>;[u]_(0.35){\nabla_\Theta}
& }.
\end{equation*}
\end{lemm}

Here is the proof of this lemma.

\begin{proof}
The $r$-fold shuffle product yields the operad evaluation product
\begin{equation*}
\ev: \C(r)\otimes B(A)^{\otimes r}\rightarrow B(A).
\end{equation*}
To make explicit the associated morphism $\nabla_c:
\C(r)\rightarrow\HopfEnd_B^{\C}(r)$, we consider the composite of
this evaluation product with the projection $B(A)\rightarrow\Sigma
A$ and the associated adjoint morphism which gives a map
\begin{equation*}
c: \C(r)\rightarrow\PrimEnd_B^{\C}(r).
\end{equation*}
According to the construction of
lemma~\ref{lemm:ConnectedQuasiCofreeHopfHom}, the morphism
$\nabla_c$ is the coalgebra morphism
\begin{equation*}
\nabla_c: \C(r)\rightarrow(\Gamma(\PrimEnd_B^{\C}(r)),\partial)
\end{equation*}
induced by this map $c$.

We apply this construction. First, the components
of~$c(1_r)\in\PrimEnd_B^{\C}(r)$, where $1_r\in\C(r)$ denotes the
generator of~$\C(r)$, are given by composites
\begin{equation*}
\Sigma A^{\otimes m_1}\otimes\dots\otimes\Sigma A^{\otimes m_r}
\hookrightarrow B(A)^{\otimes
r}\xrightarrow{\smile}B(A)\rightarrow\Sigma A.
\end{equation*}
By definition of the shuffle product, these composites are given by
the identical morphism of~$\Sigma A$ on components such that
\begin{equation*}
m_k = \begin{cases} 0 & \text{for $k\not=i$}, \\ 1 & \text{for
$k=i$} \end{cases}
\end{equation*}
and vanish otherwise. As a consequence, the maps
\begin{equation*}
c(1_r)_{m_*}\in\DGHom_{A\in\C\Alg}(\Sigma A^{\otimes
m_1+\dots+m_r},\Sigma A)
\end{equation*}
are identified with the natural transformations
$\Theta(\gamma_{m_*})$ associated to elements
$\gamma_{m_*}\in\Lambda\C(m_1+\dots+m_r)$ given by operad units
$1\in\Lambda\C(0+\dots+1+\dots+0)$. Accordingly, the maps $c:
\C(r)\rightarrow\PrimEnd_B^{\C}(r)$ admit factorizations
\begin{equation*}
\xymatrix{ \C(r)\ar@{-->}[dr]!UL_{\nabla_\gamma}\ar[r]^(0.4){\nabla_c} & \PrimEnd_B^{\C}(r) \\
& \PrimOp_B^{\C}(r)\ar@{^{(}->}[]!U+<0pt,4pt>;[u]_(0.35){\Theta} &
}.
\end{equation*}
and the induced coalgebra morphisms
\begin{equation*}
\xymatrix{ \C(r)\ar@{-->}[dr]!UL_{\nabla_\gamma}\ar[r]^(0.35){\nabla_c} & (\Gamma(\PrimEnd_B^{\C}(r)),\partial) \\
&
(\Gamma(\PrimOp_B^{\C}(r)),\partial)\ar@{^{(}->}[]!U+<0pt,4pt>;[u]_(0.35){\Theta}
}.
\end{equation*}
return a factorization $\nabla_\gamma: \C\rightarrow\HopfOp_B^{\C}$
of the given operad morphism $\nabla_c$.

Observe that this factorization defines automatically an operad
morphism since $\nabla_\Theta:
\HopfOp_B^{\C}\hookrightarrow\HopfEnd_B^{\C}$ is an injection. This
assertion can also be checked directly from the explicit form of the
map $\nabla_\gamma$.
\end{proof}

\medskip
The goal of the next paragraphs is to prove that any morphism of
unital Hopf operads
\begin{equation*}
\nabla_\rho: \Q\rightarrow\HopfOp_B^{\C},
\end{equation*}
where $\Q$ is connected and non-negatively graded, matches the same
construction. This gives the result of
theorem~\ref{thm:HopfOperadLifting} stated in the introduction of
this section. In fact, our arguments uses only the unital unitary
Hopf $\Lambda_*$-module structure of a unital Hopf operad. Therefore
our claim holds naturally for connected unital unitary Hopf
$\Lambda_*$-modules and not only for connected unital Hopf operads.
Finally, we check the following statement:

\begin{claim}[Compare with theorem~\ref{thm:HopfOperadLifting}]\label{claim:HopfOperadLifting}
Any morphism of unital unitary Hopf $\Lambda_*$-modules
$\nabla_\rho: M\rightarrow\HopfOp_B^{\C}$, where $M$ is connected
and non-negatively graded, makes commute the diagram
\begin{equation*}
\xymatrix@!C=6mm{ M\ar[dr]_{\epsilon}\ar[rr]^{\nabla_\rho} &&
\HopfOp_B^{\C} \\ & \C\ar[ur]_{\nabla_\gamma} & }.
\end{equation*}
\end{claim}

Recall that a unital unitary $\Lambda_*$-module refers to a
$\Lambda_*$-module equipped with a distinguished element $*\in M(0)$
that spans $M(0)$ and with a unit element $1\in M(1)$ such that
$\partial_1(1) = *$. Furthermore, a unital unitary
$\Lambda_*$-module is connected if $M(1) = \F$ as in the case of a
unital operad. These objects are equipped with a canonical
augmentation $\epsilon: M\rightarrow\C$ given by the
$\Lambda_*$-module operation $\eta^*_0: M(r)\rightarrow M(0) = \F$
associated to the initial map $\eta_0:
\emptyset\rightarrow\{1,\dots,r\}$. In the case of an operad, the
augmentation can be identified with the operadic composite with
unital operations $\epsilon(p) = p(*,\dots,*)$.

Any morphism of Hopf $\Lambda_*$-modules $\nabla_\rho:
M\rightarrow\HopfOp_B^{\C}$ is determined by a collection of
dg-module maps $\rho: M(r)\rightarrow\PrimOp_B^{\C}(r)$ as in the
case of the commutative operad $M = \C$ since $\HopfOp_B^{\C}$
consists of quasi-cofree Hopf coalgebras. We check that these maps
agree with the maps $\gamma: \C(r)\rightarrow\PrimOp_B^{\C}(r)$
considered in the proof of lemma~\ref{lemm:ShuffleStructure}.
Accordingly, we have automatically $\rho = \gamma\epsilon$ so that
$\nabla_\rho = \nabla_\gamma\cdot\epsilon$.

For our purpose we determine the components of~$\PrimOp_B^{\C}(r)$
of degree $*\geq 0$. To begin with, we have the following easy
observation:

\begin{obsv}\label{obsv:CommutativeOperadSuspension}
In degree $*>0$, we have $\Lambda\C(m)_* = 0$, for all $m>0$.

In degree $*=0$, we have $\Lambda\C(m)_0 = \F$ if $m = 1$ and
$\Lambda\C(m)_0 = 0$ for all $m>1$.\qed
\end{obsv}

\begin{proof}
By definition, we have
\begin{equation*}
\Lambda\C(m)_* = \Sigma^{1-m}\C(m)_* = \C(m)_{*+m-1}
\end{equation*}
(up to signs). Consequently, the module $\Lambda\C(m)_*$ is
concentrated in degree $*=1-m$. The observation follows.
\end{proof}

As a consequence, we obtain immediately:

\begin{obsv}\label{obsv:CommutativeBarPrimaryOperations}
In degree $*>0$, we have $\PrimOp_B^{\C}(r)_* = 0$. In degree $*=0$,
we have $\PrimOp_B^{\C}(1)_0 = \F$ and $\PrimOp_B^{\C}(r)_0 =
\prod_{m_*} \Lambda\C(m_1+\dots+m_r)_0$ is reduced to components of
the form $\Lambda\C(0+\dots+1+\dots+0)_0 = \F$, where $m_k = 0$
except for one index $k = i$ for which we have $m_i = 1$.\qed
\end{obsv}

In order to determine the components of a morphism
$\Q(r)\rightarrow\PrimOp_B^{\C}(r)$, we consider again the
$\Lambda_*$-module operation $\eta_i^*: M(r)\rightarrow M(1)$
associated to the map $\eta_i: \{1\}\rightarrow\{1,\dots,r\}$ such
that $\eta_i(1) = i$. In fact, we observe that these operations
isolate the degree $0$ components of $\PrimOp_B^{\C}$. More
formally, we have the following assertion:

\begin{obsv}\label{obsv:CommutativeBarOperations}
In degree $0$, we obtain $\PrimOp_B^{\C}(1)_0 = \Lambda\C(1) = \F$
and the operations
\begin{equation*}
\eta_i^*: \PrimOp_B^{\C}(r)\rightarrow\PrimOp_B^{\C}(1)
\end{equation*}
yield an isomorphism
\begin{equation*}
(\eta_i^*)_i: \PrimOp_B^{\C}(r)_0\xrightarrow{\simeq}\F^{\times r}.
\end{equation*}
\end{obsv}

\begin{proof}
By definition of the $\Lambda_*$-module structure
of~$\PrimOp_B^{\C}$ (see
claim~\ref{claim:BarOperationsUnitaryComposites}) the morphism
\begin{equation*}
\eta_i^*: \prod_{m_*} \Lambda\C(m_1+\dots+m_r)_0\rightarrow\F
\end{equation*}
is an identical morphism on the component
$\Lambda\C(m_1+\dots+m_r)_0$ such that
\begin{equation*}
m_k = \begin{cases} 0 & \text{for $k\not=i$}, \\ 1 & \text{for
$k=i$} \end{cases}
\end{equation*}
and vanishes otherwise.

The claim follows immediately since the module $\PrimOp_B^{\C}(r)_0$
is the product of these components $\Lambda\C(0+\dots+1+\dots+0)_0$
according to the previous
observation~\ref{obsv:CommutativeBarPrimaryOperations}.
\end{proof}

Then claim~\ref{claim:HopfOperadLifting} arises as a corollary of
the following result:

\begin{obsv}
Suppose given a morphism of unital unitary Hopf $\Lambda_*$-modules
$\nabla_\rho: M\rightarrow\HopfOp_B^{\C}$ induced by homogeneous
maps of degree $0$:
\begin{equation*}
\rho: M(r)\rightarrow\PrimOp_B^{\C}(r).
\end{equation*}
If we assume that $M$ is connected, so that $M(1) = \F$, then the
maps $\rho$ are determined in degree $*=0$ by the commutative
triangle
\begin{equation*}
\xymatrix{ M(r)_0\ar[dr]_{(\eta_i^*)_i}\ar[rr] &&
\PrimOp_B^{\C}(r)_0\ar[dl]_{\simeq}^{(\eta_i^*)_i} \\ & \F^{\times
r} & }.
\end{equation*}
If we assume that $M$ is non-negatively graded, then the other
components of~$\rho$ vanish since we observe that
$\PrimOp_B^{\C}(r)_* = 0$ in degree $*>0$.
\end{obsv}

\begin{proof}
Recall that the maps
\begin{equation*}
\rho: M(r)\rightarrow\PrimOp_B^{\C}(r)
\end{equation*}
are given by the composite of the associated coalgebra morphisms
\begin{equation*}
\nabla_\rho: M(r)\rightarrow(\Gamma(\PrimOp_B^{\C}(r)),\partial)
\end{equation*}
with the projection
$\Gamma(\PrimOp_B^{\C}(r))\rightarrow\PrimOp_B^{\C}(r)$. These
projections commute with $\Lambda_*$-module operations by
construction since $\HopfOp_B^{\C} =
(\Gamma(\PrimOp_B^{\C}),\partial)$ is defined as a quasi-cofree Hopf
$\Lambda_*$-module. Accordingly, so do the maps $\rho:
M(r)\rightarrow\PrimOp_B^{\C}(r)$.

Recall also that the operad unit of $\HopfOp_B^{\C}$ is represented
by the collection
\begin{equation*}
1 = \{1^{\otimes n}\}\in\prod_n\PrimOp_B^{\C}(1).
\end{equation*}
Consequently, the projection
$\HopfOp_B^{\C}(1)\rightarrow\PrimOp_B^{\C}(1)$ maps this unit to
the element $1\in\Lambda\C(1)$ that generates the degree $0$
component of~$\PrimOp_B^{\C}(1)$. Since a morphism of unitary
$\Lambda_*$-modules is supposed to preserve units, we have
$\nabla_\rho(1) = \{1^{\otimes n}\}$. As a consequence, if $M$ is
connected, then the component $\rho:
M(1)\rightarrow\PrimOp_B^{\C}(1)$ of the map~$\rho$ is given by the
identical morphism of~$M(1) =\linebreak \PrimOp_B^{\C}(1)_0 = \F$.

Then our claim follows from the commutativity of the diagrams
\begin{equation*}
\xymatrix{ M(r)\ar[d]^{\eta_i^*}\ar[r] & \PrimOp_B^{\C}(r)\ar[d]^{\eta_i^*} \\
M(1)\ar[r] & \PrimOp_B^{\C}(1) },
\end{equation*}
for $r\geq 1$ and $i = 1,\dots,r$, and from the previous
observation.
\end{proof}

These verifications achieve the proof of
claim~\ref{claim:HopfOperadLifting} and hence of
theorem~\ref{thm:HopfOperadLifting}.\qed

\medskip
As asserted in the section introduction,
claim~\ref{claim:HopfOperadLifting} implies that the morphism
$\nabla_\gamma: \C\rightarrow\HopfOp_B^{\C}$ induces an isomorphism
\begin{equation*}
\nabla_\gamma:
\C\xrightarrow{\simeq}\str^1_*\str^{\dg}_+(\HopfOp_B^{\C})
\end{equation*}
that identifies the image of~$\HopfOp_B^{\C}$ under the truncation
functors with the commutative operad. As a consequence, we obtain
that the augmentation morphism of an $E_\infty$-operad induces an
acyclic fibration of connected unital Hopf operads
\begin{equation*}
\epsilon_*: \str^1_*\str^{\dg}_+(\HopfOp_B^{\E})\wefib\C
\end{equation*}
since truncation functors preserve fibrations and acyclic
fibrations. As a corollary, for a connected unital Hopf operad~$\Q$,
the lifting problem considered in the section introduction is
equivalent to adjoint lifting problems
\begin{equation*}
\vcenter{\xymatrix{ & \str^1_*\str^{\dg}_+(\HopfOp_B^{\E})\ar@{->>}[d]^{\sim} \\
\Q\ar@{-->}[ur]^{\exists\nabla_\rho}\ar[r] & \C }}
\quad\Leftrightarrow\quad\vcenter{\xymatrix{ & \HopfOp_B^{\E}\ar[d] \\
\Q\ar@{-->}[ur]^{\exists\nabla_\rho}\ar[r] & \C }}
\end{equation*}
which have automatically a solution if $\Q$ is cofibrant.

\subsection{Prospects: actions of cellular operads}
Thinking about it, in our construction, we consider a natural
cocellular decomposition of the Hopf operad of universal bar
operations $\HopfOp_B^{\P}$ which arises from the degreewise
filtration of the bar complex $B(A)$. In subsequent work, we plan to
study operadic cellular structures giving rise to refinements of
this cocellular decomposition. As alluded to in the memoir
introduction, this might shed light on the algebraic structure of
the bar complex for subclasses of particular algebras. This prospect
motivates in part our presentation choices and the detailed accounts
of~\ref{subsection:BarHopfEndomorphismOperad}-\ref{subsection:FibrationProperties}.

\part*{Toward effective constructions}

\section{The explicit equations of Hopf operad actions}\label{section:ExplicitOperadActions}

\subsection{Introduction}
In this section we address the issue of constructing effectively
Hopf operad actions on the bar construction. In fact, our results
can supply explicit recursive constructions for the operations
$\nabla_\theta(q): B(A)^{\otimes r}\rightarrow B(A)$ associated to
operad elements $q\in\Q(r)$. Our purpose is to give this elementary
construction either as an illustration of our techniques or for a
direct application.

\medskip
To be explicit, we consider again the universal Hopf endomorphism
operad $\HopfEnd_B^{\P}$, defined
in~\ref{subsection:BarHopfEndomorphismOperad}, the operad of
universal bar operations $\HopfOp_B^{\P}$, defined
in~\ref{subsection:BarHopfOperations}, and the canonical morphism
$\nabla_\Theta: \HopfOp_B^{\P}\rightarrow\HopfEnd_B^{\P}$. According
to theorem~\ref{thm:BarHopfEndomorphismOperadRecall}, the action of
a unital Hopf operad~$\Q$ on the bar complex is equivalent to a
morphism $\nabla_\theta: \Q\rightarrow\HopfOp_B^{\P}$.
In~\ref{subsection:BarHopfOperations} we observe that the natural
morphism $\nabla_\Theta$ is an isomorphism. Consequently:

\begin{thmfact}\label{fact:BarOperationMorphismReduction}
For a $\Sigma_*$-projective operad $\P$, any morphism
$\nabla_\theta: \Q\rightarrow\HopfEnd_B^{\P}$ is equivalent to a
composite
\begin{equation*}
\xymatrix@C=3mm{ \Q\ar[rr]^(0.4){\nabla_\theta}\ar@{-->}[]!DR;[dr]_(0.35){\nabla_\rho} && \HopfEnd_B^{\P} \\
& \HopfOp_B^{\P}\ar[ur]!DL^(0.6){\simeq}_(0.7){\nabla_\Theta} & }.
\end{equation*}
As a byproduct, the action of a unital Hopf operad~$\Q$ on the bar
complex is actually determined by a morphism to the Hopf operad of
universal bar operations $\nabla_\rho: \Q\rightarrow\HopfOp_B^{\P}$
such that $\nabla_\theta = \nabla_\Theta\nabla_\rho$.
\end{thmfact}

Then, as an application of the structure results
of~\ref{section:OperadActionConstruction}, we obtain that the
operations $\nabla_\theta(q): B(A)^{\otimes r}\rightarrow B(A)$
satisfy the following characteristic properties:

\begin{thmfact}\label{fact:OperationStructure}
The components
\begin{equation*}
\xymatrix@C=24mm{ \Sigma A^{\otimes
m_1+\dots+m_r}\ar[r]^(0.6){\theta_{m_1,\dots,m_r}(q)}\ar@{^{(}->}[]!D-<0pt,4pt>;[d]
& \Sigma A \\
B(A)^{\otimes r}\ar[r]^(0.6){\nabla_\theta(q)} & B(A)\ar[u] }
\end{equation*}
of an operation~$\nabla_\theta(q)$ associated to an element
$q\in\Q(r)$ are given by the evaluation of operations
\begin{equation*}
\rho_{m_1,\dots,m_r}(q)\in\Lambda\P(m_1+\dots+m_r)
\end{equation*}
associated to~$q$.
\end{thmfact}

Recall that $\Lambda\P$ denotes the operadic suspension of~$\P$: the
operad whose algebras are suspensions $\Sigma A$ of $\P$-algebras
$A$.

\begin{thmfact}\label{fact:OperationExpansion}
The operation~$\nabla_\theta(q): B(A)^{\otimes r}\rightarrow B(A)$
associated to an element $q\in\Q(r)$ has an expansion of the form
\begin{equation*}
\nabla_\theta(q)(\alpha_1,\dots,\alpha_r) = \sum_{n=1}^{\infty}
\Bigl\{\sum\Bigl[
\underbrace{\theta(q^1)(\alpha_1^1,\dots,\alpha_r^1)\otimes\dots
\otimes\theta(q^n)(\alpha_1^n,\dots,\alpha_r^n)}_{\in\Sigma
A^{\otimes n}}\Bigr]\Bigr\},
\end{equation*}
for elements in the bar complex $\alpha_1,\dots,\alpha_r\in B(A)$,
where we consider the $n$-fold diagonals $\sum
\alpha_i^1\otimes\dots\otimes\alpha_i^n\in B(A)^{\otimes n}$ of the
tensors $\alpha_i\in B(A)$, the $n$-fold diagonals $\sum
q^1\otimes\dots\otimes q^n\in\Q(r)^{\otimes n}$ of the operation
$q\in\Q(r)$ and the homogeneous transformations
\begin{equation*}
\theta_{m^j_*}(q^j): \Sigma A^{\otimes
m^j_1}\otimes\dots\otimes\Sigma A^{\otimes m^j_r}\rightarrow\Sigma A
\end{equation*}
defined by the
operations~$\rho_{m^j_*}(q^j)\in\Lambda\P(m^j_1+\dots+m^j_r)$
associated to the elements $q^j\in\Q(r)$.
\end{thmfact}

The arguments are surveyed
in~\ref{subsection:BarHopfOperationExpansion}. According to these
statements, our problem is reduced to the construction of
appropriate maps $\rho_{m_*}:
\Q(r)\rightarrow\Lambda\P(m_1+\dots+m_r)$.
In~\ref{subsection:BarHopfOperationExpansion} we recall briefly the
structure of the Hopf operad of natural bar operations
$\HopfOp_B^{\P}$ and we make explicit the equations satisfied by the
maps $\rho_{m_*}$ associated to an operad morphism $\nabla_\rho:
\Q\rightarrow\HopfOp_B^{\P}$. Our results are recorded in
theorem~\ref{thm:OperationExplicitDefinition} in this subsection.

One can observe that these equations give rise to a recursive
definition for the operations $\nabla_\theta(q): B(A)^{\otimes
r}\rightarrow B(A)$ associated to elements of a cellular cofibrant
operad $\Q$. For the sake of precision, we state this recursive
definition in a theorem, namely
theorem~\ref{thm:OperationRecursiveDefinition}
in~\ref{subsection:BarHopfOperationConstruction}. Then we survey the
abstract lifting arguments set in
sections~\ref{section:CocellularCoalgebras}-\ref{section:OperadActionConstruction}
in order to prove this theorem and to give thorough justifications
for this recursive construction.

\medskip
To recapitulate: in~\ref{subsection:BarHopfOperationExpansion}, we
give an elementary interpretation, in term of operations, of our
abstract structure results;
in~\ref{subsection:BarHopfOperationConstruction}, we survey our
lifting arguments and we give a recursive construction of the action
of an operad on the bar complex.

\subsection{The expansion of operations on the bar complex}\label{subsection:BarHopfOperationExpansion}
As stated in the section introduction, the goal of this subsection
is to make explicit the equations satisfied by maps $\rho_{m_*}:
\Q(r)\rightarrow\Lambda\P(m_1+\dots+m_r)$ in order to obtain an
elementary and effective characterization of the action of a unital
Hopf operad on the bar construction. This result is obtained as a
consequence of constructions of the previous section. Namely such a
collection of maps is associated to coalgebra morphisms
$\nabla_\rho: \Q(r)\rightarrow\HopfOp_B^{\P}(r)$ that specify an
operad action on the bar complex if and only if they define a
morphism of unital Hopf operad. Therefore we relate simply the
equations of an operad morphism to equivalent properties for the
maps $\rho_{m_*}$. For this aim we recall briefly the definition of
the Hopf operad of universal bar operations $\HopfOp_B^{\P}$.

\subsubsection{Recalls: the coalgebra structure of the universal Hopf
operads} Precisely, recall that the Hopf endomorphism operad
$\HopfEnd_B^{\P}$, respectively the Hopf operad of universal bar
operations $\HopfOp_B^{\P}$, is defined by quasi-cofree coalgebras
such that
\begin{align*}
& \HopfEnd_B^{\P}(r) = (\Gamma(\PrimEnd_B^{\P}(r)),\partial),\\
&\qquad\text{where}\quad\PrimEnd_B^{\P}(r) = \prod_{m_1+\dots+m_r>0}
\DGHom_{A\in\P\Alg}(\Sigma A^{\otimes m_1+\dots+m_r},\Sigma A),\\
\intertext{respectively}
& \HopfOp_B^{\P}(r) = (\Gamma(\PrimOp_B^{\P}(r)),\partial),\\
&\qquad\text{where}\quad\PrimOp_B^{\P}(r) = \prod_{m_1+\dots+m_r>0}
\Lambda\P(m_1+\dots+m_r).
\end{align*}
The isomorphism
\begin{equation*}
\underbrace{(\Gamma(\PrimOp_B^{\P}),\partial)}_{\HopfOp_B^{\P}}
\xrightarrow[\simeq]{\nabla_\Theta}\underbrace{(\Gamma(\PrimEnd_B^{\P}),\partial)}_{\HopfEnd_B^{\P}}
\end{equation*}
is induced by the canonical morphisms
\begin{equation*}
\underbrace{\prod_{m_1+\dots+m_r>0}
\Lambda\P(m_1+\dots+m_r)}_{\PrimOp_B^{\P}(r)}
\xrightarrow{\Theta}\underbrace{\DGHom_{A\in\P\Alg}(\Sigma
A^{\otimes m_1+\dots+m_r},\Sigma A)}_{\HopfEnd_B^{\P}(r)}
\end{equation*}
which map an operad element $p\in\Lambda\P(m_1+\dots+m_r)$ to the
associated natural operation $\Theta(p): \Sigma A^{\otimes
m_1+\dots+m_r}\rightarrow\Sigma A$.

According to lemma~\ref{lemm:QuasiCofreeCoalgebraStructure}, a
coalgebra morphism $\nabla_\theta:
\Q(r)\rightarrow\HopfEnd_B^{\P}(r)$ is determined by a homogeneous
map of degree $0$
\begin{equation*}
\theta: \Q(r)\rightarrow\PrimEnd_B^{\P}(r)
\end{equation*}
such that $\theta = \pi\nabla_\theta$, where we consider the natural
projection
\begin{equation*}
\HopfEnd_B^{\P}(r) =
(\Gamma(\PrimEnd_B^{\P}(r)),\partial)\xrightarrow{\pi}\PrimEnd_B^{\P}(r).
\end{equation*}
Similarly, a coalgebra morphism $\nabla_\rho:
\Q(r)\rightarrow\HopfOp_B^{\P}(r)$ is determined by a homogeneous
map $\rho: \Q(r)\rightarrow\PrimOp_B^{\P}(r)$ such that $\theta =
\pi\nabla_\rho$. If we assume $\nabla_\theta =
\nabla_\Theta\nabla_\rho$, then we have also the relation $\theta =
\Theta\rho$.

The precise purpose of this subsection is to make explicit the
equations satisfied by a collection of maps $\rho:
\Q(r)\rightarrow\HopfOp_B^{\P}(r)$ so that the associated coalgebra
morphisms $\nabla_\rho: \Q(r)\rightarrow\HopfOp_B^{\P}(r)$ define a
morphism of unital Hopf operads.

\subsubsection{On the expansion of operations}
Before we explain briefly that, for the actual operation
$\nabla_\theta(q): B(A)^{\otimes r}\rightarrow B(A)$ associated to
an element $q\in\Q(r)$, the maps $\theta_{m_*}(q): \Sigma A^{\otimes
m_1+\dots+m_r}\rightarrow\Sigma A$ considered in
fact~\ref{fact:OperationStructure} are determined by the components
\begin{equation*}
\theta_{m_*}: \Q(r)\rightarrow\DGHom_{A\in\P\Alg}(\Sigma A^{\otimes
m_1+\dots+m_r},\Sigma A)
\end{equation*}
of the maps $\theta: \Q(r)\rightarrow\PrimEnd_B^{\P}(r)$. In fact,
this relationship is obtained
in~\ref{item:BarHopfEndomorphismEvaluation} where we make explicit
the evaluation product of a Hopf endomorphism operad
\begin{equation*}
\HopfEnd_{B(A)}(r)\otimes B(A)^{\otimes r}\rightarrow B(A).
\end{equation*}
The expansion given in fact~\ref{fact:OperationExpansion} for an
operation $\nabla_\theta(q): B(A)^{\otimes r}\rightarrow B(A)$ comes
also from this paragraph.

Recall simply that the Hopf endomorphism operad $\HopfEnd_{B(A)}$
consists of quasi-cofree coalgebras $\HopfEnd_{B(A)}(r) =
(\Gamma(\PrimEnd_{B(A)}(A)),\partial)$ like the universal Hopf
endomorphism operad $\HopfEnd_B^{\P}$ except that we replace the
modules of natural transformations
\begin{align*}
\PrimEnd_B^{\P}(r) & = \prod_{m_1+\dots+m_r>0}
\DGHom_{A\in\P\Alg}(\Sigma A^{\otimes m_1+\dots+m_r},\Sigma A)\\
\intertext{by the modules of homogeneous morphisms}
\PrimEnd_{B(A)}(r) & = \prod_{m_1+\dots+m_r>0} \DGHom(\Sigma
A^{\otimes m_1+\dots+m_r},\Sigma A)
\end{align*}
associated to the given $\P$-algebra $A$. The canonical morphism
$\HopfEnd_B^{\P}\rightarrow\linebreak\HopfEnd_{B(A)}$ is induced by
the obvious maps
\begin{equation*}
\DGHom_{A\in\P\Alg}(\Sigma A^{\otimes m_1+\dots+m_r},\Sigma
A)\rightarrow\DGHom(\Sigma A^{\otimes m_1+\dots+m_r},\Sigma A)
\end{equation*}
given by the specialization of a natural transformation to the given
algebra $A$.

Therefore the assertions of~\ref{item:BarHopfEndomorphismEvaluation}
together with fact~\ref{fact:BarOperationMorphismReduction} give
exactly the claim of fact~\ref{fact:OperationStructure} and the
expansion of fact~\ref{fact:OperationExpansion} for an operation
$\nabla_\theta(q): B(A)^{\otimes r}\rightarrow B(A)$.

\medskip
The operations $\rho_{m_*}(q)\in\Lambda\P(m_1+\dots+m_r)$ that occur
in facts~\ref{fact:OperationStructure}-\ref{fact:OperationExpansion}
are also yielded by the components $\rho_{m_*}:
\Q(r)\rightarrow\Lambda\P(m_1+\dots+m_r)$ of the maps $\rho:
\Q(r)\rightarrow\PrimOp_B^{\P}(r)$. Finally, we have the following
theorem that fulfils the objective of this subsection:

\begin{thm}\label{thm:OperationExplicitDefinition}
A collection of maps
\begin{equation*}
\rho_{m_*}: \Q(r)\rightarrow\Lambda\P(m_1+\dots+m_r)
\end{equation*}
determines a morphism of unital operads $\nabla_\rho:
\Q\rightarrow\HopfOp_B^{\P}$ if and only if the following properties
are satisfied:
\begin{enumerate}
\item
for the unit element $1\in\Q(1)$, we have
\begin{equation*}
\rho_{m_1}(1) = \begin{cases} 1\in\Lambda\P(1) & \text{if $m_1 =
1$},
\\ 0\in\Lambda\P(m_1) & \text{otherwise}; \end{cases}
\end{equation*}
\item
for any element $q\in\Q(r)$, we have the permutation relation
\begin{equation*}
\rho_{m_1,\dots,m_r}(w\cdot q) =
w(m_1,\dots,m_r)\cdot\rho_{m_{w(1)},\dots,m_{w(r)}}(q),
\end{equation*}
for all permutations $w\in\Sigma_r$, and the $\Lambda_*$-module
relation
\begin{equation*}
\rho_{m_1,\dots,m_r}(q) =
\rho_{m_1,\dots,\widehat{m_i},\dots,m_r}(q\circ_i *),
\end{equation*}
for all collections $m_*$ such that $m_i = 0$;
\item
for a composite operation $p\circ_i q\in\Q(s+t-1)$, where
$p\in\Q(s)$ and $q\in\Q(t)$, we have the composition relation
\begin{align*}
&\rho_{l_*}(p\circ_i q) = \\
&\quad\sum_{m}\Bigl\{\sum\shuffle_i(n^k_j)\cdot
\rho_{m_*}(p)\bigl(1,\dots,\rho_{n^1_*}(q^1),\dots,\rho_{n^m_*}(q^m),\dots,1\bigr)\Bigr\},
\end{align*}
for indices $(m_*),(n_*)$ such that
\begin{align*}
& m_k = l_k,\quad\text{for $k = 1,\dots,i-1$},\\
& m_i = m\quad\text{ranges over positive integers},\\
& m_k = l_{k+t-1}\quad\text{for $k = i+1,\dots,s$},\\
& n_k = l_{k+i-1}\quad\text{for $k = 1,\dots,t$}
\end{align*}
and where we consider the partitions $n^1_k+\dots+n^m_k = n_k$ and
the operations $\rho_{n^j_*}(q^j)$ associated to the $m$-fold
diagonals $\sum q^1\otimes\dots\otimes q^n\in\Q(s)^{\otimes n}$ of
the element $q\in\Q(s)$; on the right-hand side, the operations
$\rho_{n^j_*}(q^j)$ are substituted to the entries $t =
m_1+\dots+m_{i-1}+1,\dots, m_1+\dots+m_{i-1}+m$ of the operations
$\rho_{m_*}(p)$;
\item
for any element $q\in\Q(r)$, we have the differential relation
\begin{align*}
\delta(\rho_{m_*}(q)) & =
\rho_{m_*}(\delta(q)) + \rho_{m_*}(\partial(q))\\
& - \sum_n\Bigl\{\sum
\shuffle(m^j_i)\cdot\mu_n(\rho_{m^1_*}(q^1),\dots,\rho_{m^n_*}(q^n))\Bigr\}\\
&
\pm\sum_i\Bigl\{\sum_{n,t}\rho_{m'_*}(q)\circ_t\mu_n\Bigr\}\Bigr\};
\end{align*}
in the first summation we consider the $n$-fold diagonals $\sum
q^1\otimes\dots\otimes q^n\in\Q(r)^{\otimes n}$ of~$q$ and the
partitions $m^1_i+\dots+m^n_i = m_i$; in the second summation we
consider the collections $m'_*$ such that $m'_* = m_*$ for $*\not=i$
and $m'_i = m_i+n-1$ and $t$ ranges over the interval $t =
m_1+\dots+m_{i-1}+1,\dots, m_1+\dots+m_{i-1}+m_i$.
\end{enumerate}
\end{thm}

\begin{proof}
This theorem follows from direct applications of
lemma~\ref{lemm:QuasiCofreeCoalgebraStructure} and from the explicit
definition of the Hopf operad~$\HopfOp_B^{\P}$: in general, as a
Hopf operad structure is determined by coalgebra morphisms and since
the Hopf operad~$\HopfOp_B^{\P}$ consists of quasi-cofree
coalgebras, it is sufficient to check relations onto
$\PrimOp_B^{\P}$.

The properties (a-c) reflect the equations satisfied by an operad
morphism $\nabla_\rho: \Q(r)\rightarrow\HopfOp_B^{\P}(r)$.
Explicitly, the unit relation $\nabla_\rho(1) = 1$ is equivalent to
the commutativity of the diagram
\begin{equation*}
\xymatrix{ & \F\ar[dl]!UR_{\eta}\ar[dr]!UL^{\eta} & \\
\Q(1)\ar[rr]_(0.45){\nabla_\rho} && \HopfOp_B^{\P}(1) },
\end{equation*}
where we consider the coalgebra morphism $\eta: \F\rightarrow\Q(1)$,
respectively $\eta: \F\rightarrow\HopfOp_B^{\P}(1)$, specified by
the unit element $1\in\Q(1)$, respectively $1\in\HopfOp_B^{\P}(1)$.
As explained above, the identity of the coalgebra morphisms
$\nabla_\rho\eta = \eta$ is satisfied if and only if the composite
of these morphisms with the projection $\pi:
\HopfOp_B^{\P}(1)\rightarrow\PrimOp_B^{\P}(1)$ agree. Thus we obtain
the equation
\begin{equation*}
\pi\nabla_\rho(1) = \pi(1)\Leftrightarrow\rho(1) = \pi(1).
\end{equation*}
If we go back to the definition of the unit element
of~$\HopfOp_B^{\P}$ given in the proof of
claims~\ref{claim:BarOperationsUnitaryComposites}-\ref{claim:BarHopfOperationComposites},
then we obtain exactly the relation (a) of the theorem.

Similarly, the permutation relation $\nabla_\rho(w q) =
w\nabla_\rho(q)$ and the $\Lambda_*$-module relation
$\nabla_\rho(q\circ_i *) = \nabla_\rho(q)\circ_i *$ are equivalent
to the commutativity of diagrams of coalgebra morphisms. Namely:
\begin{equation*}
\xymatrix{ \Q(r)\ar[r]^(0.4){\nabla_\rho}\ar[d]^{w} &
\HopfOp_B^{\P}(r)\ar[d]^w \\
\Q(r)\ar[r]^(0.4){\nabla_\rho} & \HopfOp_B^{\P}(r)
},\quad\text{respectively}\quad\xymatrix{
\Q(r)\ar[r]^(0.4){\nabla_\rho}\ar[d]^{\partial_i} &
\HopfOp_B^{\P}(r)\ar[d]^{\partial_i} \\
\Q(r)\ar[r]^(0.4){\nabla_\rho} & \HopfOp_B^{\P}(r) }.
\end{equation*}
Then the commutativity of these diagrams is equivalent to the
relations
\begin{gather*}
\pi(\nabla_\rho(w q)) = \pi(w\nabla_\rho(q))\Leftrightarrow\rho(w q)
=
w\rho(q),\\
\text{respectively}\quad\pi(\nabla_\rho(q\circ_i *)) =
\pi(\nabla_\rho(q)\circ_i
*)\Leftrightarrow\rho(q\circ_i *) =
\partial_i(\rho(q)).
\end{gather*}
If we go back to the definition of the $\Lambda_*$-module structure
of $\PrimOp_B^{\P}$ given in the proof of
claim~\ref{claim:BarOperationsUnitaryComposites}, then we obtain
immediately the relation (b) of the theorem.

Notice that the statements (a-b) assert exactly that the maps $\rho:
\Q(r)\rightarrow\PrimOp_B^{\P}(r)$ define a morphism of unitary
$\Lambda_*$-modules.

The composition relation $\nabla_\rho(p\circ_i q) =
\nabla_\rho(p)\circ_i\nabla_\rho(q)$ is given by the diagram
\begin{equation*}
\xymatrix@C=24mm{
\Q(s)\otimes\Q(t)\ar[r]^(0.4){\nabla_\theta\otimes\nabla_\theta}\ar[d]_{\circ_i}
& \HopfOp_B^{\P}(s)\otimes\HopfOp_B^{\P}(t)\ar[d]_{\circ_i} \\
\Q(s+t-1)\ar[r]^(0.4){\nabla_\theta} & \HopfOp_B^{\P}(s+t-1) }.
\end{equation*}
Again this diagram commutes if and only if the involved morphisms
agree on the quotient object $\PrimOp_B^{\P}(s+t-1)$ of
$\HopfOp_B^{\P}(s+t-1)$. The quotient of the evaluation product
\begin{equation*}
\HopfOp_B^{\P}(s)\otimes\HopfOp_B^{\P}(t)\xrightarrow{\circ_i}\HopfOp_B^{\P}(s+t-1)
\end{equation*}
is made explicit in the proof of
claim~\ref{claim:BarHopfOperationComposites}. If we go back to this
definition, then we obtain exactly the relation (c) of the theorem.

Finally, recall that a morphism of quasi-cofree coalgebras, like
\begin{equation*}
\Q(r)\xrightarrow{\nabla_\rho}\underbrace{(\Gamma(\PrimOp_B^{\P}(r)),\partial)}_{\HopfOp_B^{\P}(r)},
\end{equation*}
commutes with differentials if and only if the associated map $\rho$
satisfies the equation
\begin{equation*}
\delta(\rho) + \beta\nabla_\rho = 0,
\end{equation*}
where $\beta: \Gamma(\PrimOp_B^{\P}(r))\rightarrow\PrimOp_B^{\P}(r)$
denotes the homogeneous map that determine the coderivation
$\partial =
\partial_\beta$ of the quasi-cofree coalgebra (see
lemma~\ref{lemm:QuasiCofreeCoalgebraStructure}). In the case of the
Hopf operad of bar operations $\HopfOp_B^{\P}(r) =
(\Gamma(\PrimOp_B^{\P}(r)),\partial)$ this coderivation $\partial$
is composed of components $\partial^h_i,\partial^v$ associated to
maps $\beta^h_i,\beta^v$ (see
claim~\ref{claim:BarHopfOperationDifferentials}). If we apply the
formula given in the proof of
claim~\ref{claim:BarHopfOperationDifferentials}, then we obtain
exactly the relation (d) of the theorem. The terms
\begin{equation*}
\sum_{n,t}\rho_{m'_*}(q)\circ_t\mu_n,\quad\text{respectively}\quad
\sum_n\Bigl\{\sum
\shuffle(m^j_i)\cdot\mu_n(\rho_{m^1_*}(q^1),\dots,\rho_{m^n_*}(q^n))\Bigr\},
\end{equation*}
represent precisely components of
$\beta^i_h(\rho(q))\in\PrimOp_B^{\P}(r)$, respectively
$\beta^v(\rho(q))\in\PrimOp_B^{\P}(r)$.
\end{proof}

\subsection{Lifting process and effective constructions}\label{subsection:BarHopfOperationConstruction}
According to results of~\ref{subsection:FibrationProperties}, for an
$E_\infty$-operad $\E$, a morphism $\nabla_\rho:
\Q\rightarrow\HopfOp_B^{\E}$, that supplies effectively an operad
action on the bar complex, can be obtained by lifting the morphism
$\nabla_\gamma: \C\rightarrow\HopfOp_B^{\C}$ associated to the
classical shuffle product for the bar complex of commutative
algebras. Explicitly, a morphism $\nabla_\rho:
\Q\rightarrow\HopfOp_B^{\E}$ is defined by a solution of the lifting
problem
\begin{equation*}
\xymatrix{  && \HopfOp_B^{\E}\ar@{->>}[d]^{\sim} \\
\Q\ar[r]\ar@{-->}[urr]^{\nabla_\rho} & \C\ar[r]_(0.4){\nabla_\gamma}
& \HopfOp_B^{\C} }.
\end{equation*}
The purpose of this subsection is to survey our arguments in order
to make this lifting process effective. As a byproduct, we obtain a
recursive definition of the operation $\nabla_\theta(q):
B(A)^{\otimes r}\rightarrow B(A)$ associated to an element
$q\in\Q(r)$. We state this result as a theorem in order to motivate
the study of this subsection.

In order to obtain effective results we need to have an effective
model of an $E_\infty$-operad $\E$. Explicitly, we suppose given a
computable strong deformation retract
\begin{equation*}
\xymatrix{ *+<2mm>{\C}\ar@<-2pt>[r]_{\eta} &
*+<2mm>{\E}\ar@<-2pt>[l]_{\epsilon}\ar@(ur,dr) "1,2"!UR;"1,2"!DR^{\nu}
},
\end{equation*}
where $\E$ is an $E_\infty$-operad of which $\epsilon:
\E\rightarrow\C$ is the augmentation. Similarly, we have to consider
a cofibrant Hopf operad $\Q$ provided with a manageable cellular
structure. To be explicit, we assume that $\Q$ is a connected unital
Hopf operad that arises as the colimit of a sequence of Hopf operads
\begin{equation*}
* = \Q^{-1}\rightarrow\Q^{0}\rightarrow\dots
\rightarrow\Q^d\rightarrow\dots\xrightarrow{}\colim_d\Q^d = \Q
\end{equation*}
obtained by pushouts
\begin{equation*}
\xymatrix{ \FOp_*(C^d)\ar[r]^{f^d}\ar[d]^{i^d} & \Q^{d-1}\ar@{-->}[d] \\
\FOp_*(D^d)\ar@{-->}[r] & \Q^{d} },
\end{equation*}
where $i^d: \FOp_*(C^d)\rightarrow\FOp_*(D^d)$ is a morphism of free
operads associated to a Reedy cofibration of Hopf
$\Lambda_*$-modules $i^d: C^d\rightarrow D^d$. In addition we shall
assume that $i^d$ splits $\Sigma_*$-equivariantly degreewise so that
$D^d$ is identified as a dg-$\Sigma_*$-module with a direct sum $D^d
= (C^d\oplus E^d,\partial)$ for a projective $\Sigma_*$-module $E^d$
defined effectively as a direct summand of a finitely generated free
$\Sigma_*$-module. For simplicity we can assume that $E^d$ is a free
$\Sigma_*$-module equipped with a finite set of generators
$\xi_\kappa$. Equivalently, the operad $\Q$ is associated to a
quasi-free reduced operad $\overline{\Q} =
(\FOp(\bigoplus_{d=0}^{\infty} E^d),\partial)$ for modules of
generators $E^d$ such that $E^d(r) = \bigoplus_{\kappa} \F[\Sigma_r]
\xi_\kappa$. The suboperads $\Q^d$ are determined by
$\overline{\Q}{}^d = (\FOp(\bigoplus_{i=0}^{d} E^i),\partial)$. The
differential of a basis element $\xi_\alpha\in E^d(r)$ is given by a
sum $(\delta+\partial)(\xi_\kappa)\in\Q^d(r)$, where
$\delta(\xi_\kappa)\in E^d(r)$ and
$\partial(\xi_\kappa)\in\Q^{d-1}(r)$. According to results of the
previous section, the Boardman-Vogt construction can supply explicit
operads $\Q = W^{\square}(\E)$ that satisfy these requirements.

In the abstract context, the morphism $\nabla_\theta:
\Q\rightarrow\HopfEnd_B^{\E}$ can be obtained effectively by
applications of cellular lifting properties. In an equivalent
elementary fashion, the maps $\rho_{m_*}:
\Q(r)\rightarrow\Lambda\E(m_1+\dots+m_r)$ that determine the action
of~$\Q$ on the bar construction are obtained inductively as
solutions of differential equations. In the effective context, we
can use the chain contraction $\nu: \E(m)\rightarrow\E(m)$ in order
to specify solutions of these equations. Then our result takes the
following form:

\begin{thm}\label{thm:OperationRecursiveDefinition}
In the context set in the previous paragraphs, an appropriate
collection of maps
\begin{equation*}
\rho_{m_*}: \Q(r)\rightarrow\Lambda\E(m_1+\dots+m_r)
\end{equation*}
that give rise to an action of~$\Q$ on the bar construction can be
defined recursively by the following requirements:
\begin{enumerate}
\item
for the unit element $1\in\Q(1)$, we have
\begin{equation*}
\rho_{m_1}(1) = \begin{cases} 1\in\Lambda\E(1) & \text{if $m_1 =
1$},
\\ 0\in\Lambda\E(m_1) & \text{otherwise}; \end{cases}
\end{equation*}
\item
for a composite operation $p\circ_i q\in\Q(s+t-1)$, where
$p\in\Q(s)$ and $q\in\Q(t)$, we have
\begin{align*}
&\rho_{l_*}(p\circ_i q) = \\
&\quad\sum_{m}\Bigl\{\sum\shuffle_i(n^k_j)\cdot
\rho_{m_*}(p)\bigl(1,\dots,\rho_{n^1_*}(q^1),\dots,\rho_{n^m_*}(q^m),\dots,1\bigr)\Bigr\},
\end{align*}
for indices $(m_*),(n_*)$ such that
\begin{align*}
& m_k = l_k\quad\text{for $k = 1,\dots,i-1$},\\
& m_i = m\quad\text{ranges over positive integers},\\
& m_k = l_{k+t-1}\quad\text{for $k = i+1,\dots,s$},\\
& n_k = l_{k+i-1}\quad\text{for $k = 1,\dots,t$}
\end{align*}
and where we consider the partitions $n^1_k+\dots+n^m_k = n_k$ and
the operations $\rho_{n^j_*}(q^j)$ associated to the $m$-fold
diagonals $\sum q^1\otimes\dots\otimes q^n\in\Q(s)^{\otimes n}$ of
the element $q\in\Q(s)$; on the right-hand side, the operations
$\rho_{n^j_*}(q^j)$ are substituted to the entries $t =
m_1+\dots+m_{i-1}+1,\dots, m_1+\dots+m_{i-1}+m$ of the operations
$\rho_{m_*}(p)$;
\item
for a permutation $w\in\Sigma_r$, we have the relation
\begin{equation*}
\rho_{m_1,\dots,m_r}(w\cdot q) =
w(m_1,\dots,m_r)\cdot\rho_{m_{w(1)},\dots,m_{w(r)}}(q),
\end{equation*}
for any operation $q\in\Q(r)$;
\item
for a generator $\xi_\kappa\in\Q(r)$, we set
\begin{align*}
\rho_{m_1,\dots,m_r}(\xi_\kappa) & =
\rho_{m_1,\dots,\widehat{m_i},\dots,m_r}(\xi_\kappa\circ_i
*)\quad\text{if $m_i = 0$ for some $i$},\\
\rho_{m_1,\dots,m_r}(\xi_\kappa) & =
\nu\Bigl\{\rho_{m_*}(\delta(\xi_\kappa)) + \rho_{m_*}(\partial(\xi_\kappa))\\
& - \sum_n\Bigl\{\sum
\shuffle(m^j_i)\cdot\mu_n(\rho_{m^1_*}(\xi_\kappa^1),\dots,\rho_{m^n_*}(\xi_\kappa^n))\Bigr\}\\
&
\pm\sum_i\Bigl\{\sum_{n,t}\rho_{m'_*}(\xi_\kappa)\circ_t\mu_n\Bigr\}\Bigr\}\quad\text{otherwise};
\end{align*}
in the first summation we consider the $n$-fold diagonals $\sum
\xi_\kappa^1\otimes\dots\otimes\xi_\kappa^n\in\Q(r)^{\otimes n}$ of
the operation $\xi_\kappa$ and the partitions $m^1_i+\dots+m^n_i =
m_i$; in the second summation we consider the collections $m'_*$
such that $m'_* = m_*$ for $*\not=i$ and $m'_i = m_i+n-1$ and $t$
ranges over the interval $t = m_1+\dots+m_{i-1}+1,\dots,
m_1+\dots+m_{i-1}+m_i$.
\end{enumerate}
\end{thm}

Thorough justifications of this recursive constructions can be
obtained directly from the assertions of
theorem~\ref{thm:OperationExplicitDefinition}. One checks
essentially that the terms on the right-hand side of the equation
of~$\rho_{m_1,\dots,m_r}(\xi_\kappa)$ are determined by
operations~$\rho_{n_1,\dots,n_r}(q)$ which are either associated to
elements $q\in\Q(r)$ of lower degree than $\xi_\kappa$ or such that
$n_1+\dots+n_r<m_1+\dots+m_r$.

In this subsection we give another proof of this theorem in order to
illustrate the abstract lifting arguments of
sections~\ref{section:CocellularCoalgebras}-\ref{section:OperadActionConstruction}.

\subsubsection{Connectedness assumptions}
Throughout this subsection we assume that $\Q$ is connected and
non-negatively graded. Consequently, our lifting problem can be
simplified according to observations of the introduction
of~\ref{section:OperadActionConstruction} (see
also~\ref{subsection:CommutativeHopfBarAction}).

Explicitly, recall that we have truncation functors
\begin{equation*}
\dg_{\Z}\Op^1_*\xrightarrow{\str^{\dg}_+}\dg_{\N}\Op^1_*\xrightarrow{\str^1_*}\dg_{\N}\Op^*_*
\end{equation*}
that target to the category formed by the non-negatively graded and
connected unital Hopf operads. Furthermore, by adjunction, our
lifting problem is equivalent to
\begin{equation*}
\xymatrix{ & \str^1_*\str^{\dg}_+(\HopfOp_B^{\E})\ar@{->>}[d]^{\sim}
\ar@{^{(}->}[]!R+<4pt,0pt>;[r] & \HopfOp_B^{\E}\ar[d] \\
\Q\ar@{-->}[ur]^{\exists\nabla_\rho}\ar[r] & \C\ar[r]^{=} & \C }
\end{equation*}
since we have $\str^1_*\str^{\dg}_+(\HopfOp_B^{\C}) = \C$. In the
construction of the next paragraphs we do not need to introduce
truncation functors explicitly. Thus we consider only a reduced
lifting problem:
\begin{equation*}
\xymatrix{ & \HopfOp_B^{\E}\ar[d] \\
\Q\ar@{-->}[ur]^{\exists\nabla_\rho}\ar[r] & \C }.
\end{equation*}
By adjunction, we are ensured that any lifting morphism in this
diagram factorizes through $\str^1_*\str^{\dg}_+(\HopfOp_B^{\E})$
and all the lifting problems considered in this paragraph are
equivalent.

\subsubsection{Lifting construction for a cellular cofibrant Hopf operad}
The idea is to define our lifting by induction by using cellular
structures. For this aim we assume that $\Q$ arises as the colimit
of a sequence of Hopf operads
\begin{equation*}
* = \Q^{-1}\rightarrow\Q^{0}\rightarrow\dots
\rightarrow\Q^d\rightarrow\dots\rightarrow\colim_d\Q^d = \Q
\end{equation*}
obtained by pushouts
\begin{equation*}
\xymatrix{ \FOp_*(C^d)\ar[r]^{f^d}\ar[d]^{i^d} & \Q^{d-1}\ar@{-->}[d] \\
\FOp_*(D^d)\ar@{-->}[r] & \Q^{d} },
\end{equation*}
where $i^d: \FOp_*(C^d)\rightarrow\FOp_*(D^d)$ is a morphism of free
operads induced by a Reedy cofibration of unitary
$\Lambda_*$-modules $i^d: C^d\rightarrow D^d$.

Then the lifting $\nabla_\rho: \Q\rightarrow\HopfOp_B^{\E}$ can be
obtained as the colimit of morphisms $\nabla_\rho = \nabla_\rho^d:
\Q^d\rightarrow\HopfOp_B^{\E}$ constructed by induction on $d$.
Accordingly, we are reduced to specify inductively a lifting in the
diagram of unitary Hopf $\Lambda_*$-modules
\begin{equation*}
\xymatrix{ C^d\ar[r]\ar[d]^{i^d} & \FOp_*(C^d)\ar[r] & \Q^{d-1}\ar[r] & \HopfOp_B^{\E}\ar[d] \\
D^d\ar@{-->}[urrr]^(0.4){\nabla_\rho}\ar[rrr] &&& \C },
\end{equation*}
for all $d\geq 0$.

\subsubsection{Cocellular structures and lifting constructions}
In order to obtain these Hopf $\Lambda_*$-module liftings, we
consider the cocellular decomposition of~$\HopfOp_B^{\E}$ introduced
in~\ref{subsection:FibrationProperties}. Explicitly, recall that
$\HopfOp_B^{\E}$ is the limit of a tower of quasi-cofree Hopf
$\Lambda_*$-modules $\ckcell_m\HopfOp_B^{\E}$ such that
$\ckcell_m\HopfOp_B^{\E} =
(\Gamma(\ckcell_m\PrimOp_B^{\E}),\partial)$. Furthermore, the
projection morphism
\begin{equation*}
\pk_m: \ckcell_m\HopfOp_B^{\E}\rightarrow\ckcell_{m-1}\HopfOp_B^{\E}
\end{equation*}
fits in a coextension diagram
\begin{equation*}
\xymatrix{ (\Gamma(K),\partial)\ar[r]\ar[d]
& \Gamma(\Delta^1\wedge K)\ar[d] \\
(\Gamma(L),\partial)\ar[r] & \Gamma(S^1\wedge K\times_{S^1\wedge
L}\Delta^1\wedge L) },
\end{equation*}
where we let $K = \ckcell_m\PrimOp_B^{\E}$ and $L =
\ckcell_{m-1}\PrimOp_B^{\E}$.

The lifting $\nabla_\rho: D^d\rightarrow\HopfOp_B^{\E}$ can be
obtained as the limit of morphisms $\nabla_\rho = \nabla_{\rho_m}:
D^d\rightarrow\ckcell_m\HopfOp_B^{\E}$ constructed by induction on
$m$. Accordingly, we are reduced to specify inductively a lifting in
the diagram of unitary $\Lambda_*$-modules
\begin{equation*}
\xymatrix{ C^d\ar[r]\ar[d]^{i^d} & (\Gamma(K),\partial)\ar[r] & \Delta^1\wedge K\ar[d] \\
D^d\ar[r]\ar@{-->}[urr]^{\tilde{\rho}} & (\Gamma(L),\partial)\ar[r]
& S^1\wedge K\times_{S^1\wedge L}\Delta^1\wedge L },
\end{equation*}
for all $m\geq 0$.

Recall that $\Delta^1\wedge K = (\x{01}\otimes K\oplus\x{1}\otimes
K,\partial)$, for homogeneous elements $\x{01}$ and $\x{1}$ such
that $\partial(\x{01}) = \x{1}$. Furthermore, we observe
in~\ref{subsection:CocellularCoalgebras} (see more especially
fact~\ref{fact:ConeMorphism}) that a morphism $\tilde{\rho}:
D^d\rightarrow\Delta^1\wedge K$ has the form
\begin{equation*}
\tilde{\rho}(\xi) = - \x{01}\otimes\delta(\rho)(\xi) +
\x{1}\otimes\rho(\xi),
\end{equation*}
for a homogeneous map $\rho: D^d\rightarrow K$ of degree $0$. In
fact, the morphisms $\nabla_\rho = \nabla_{\rho_m}:
D^d(r)\rightarrow\ckcell_m\HopfOp_B^{\E}(r)$ can be identified with
the coalgebra morphisms induced by these maps $\rho = \rho_m:
D^d(r)\rightarrow\ckcell_m\PrimOp_B^{\E}(r)$
(see~\ref{item:CofreeCoextensions} and
claim~\ref{claim:CofreeCoextension}). Accordingly, in an elementary
fashion, the lifting $\tilde{\rho}$ determines the components
$\rho_m(\xi) = \{\rho_{m_*}(\xi)\}\in\prod_{m\geq
m_1+\dots+m_r>0}\Lambda\E(m_1+\dots+m_r)$ of the collection of
operations associated to an element $\xi\in E^d(r)$.

In the next paragraph we survey the definition of the coextension
diagram above in order to give an elementary interpretation of this
lifting construction.

\subsubsection{The elementary interpretation of the lifting process}
Recall that $\Delta^1\wedge K = (\x{01}\otimes K\oplus\x{1}\otimes
K,\partial)$. We have similarly
\begin{equation*}
S^1\wedge K\times_{S^1\wedge L}\Delta^1\wedge L = (\x{01}\otimes
K\oplus\x{1}\otimes L,\partial)
\end{equation*}
and the canonical map
\begin{equation*}
\Delta^1\wedge K\xrightarrow{(\sigma\wedge
K,\Delta^1\wedge\pk_m)}S^1\wedge K\times_{S^1\wedge L}\Delta^1\wedge
L
\end{equation*}
can be identified with the obvious morphism
\begin{equation*}
(\x{01}\otimes K\oplus\x{1}\otimes
K,\partial)\xrightarrow{(\x{01}\otimes\Id,\x{01}\otimes\pk_m)}(\x{01}\otimes
K\oplus\x{1}\otimes L,\partial)
\end{equation*}
induced on one summand by the identity of~$K$ and on the other
summand by the projection morphism $\pk_m: K\rightarrow L$
(see~\ref{subsection:CocellularCoalgebras} and more especially
observation~\ref{obsv:ConeSuspensionFiberedProduct}).

The precise construction of the commutative square
\begin{equation*}
\xymatrix{ (\Gamma(K),\partial)\ar[r]\ar[d] & \Delta^1\wedge K\ar[d] \\
(\Gamma(L),\partial)\ar[r] & S^1\wedge K\times_{S^1\wedge
L}\Delta^1\wedge L }
\end{equation*}
is given in~\ref{item:CofreeCoextensions}. Recall that the morphism
\begin{equation*}
(\Gamma(K),\partial)\rightarrow\x{01}\otimes K\oplus\x{1}\otimes K
\end{equation*}
is defined on the summand $\x{01}\otimes K\subset\x{01}\otimes
K\oplus\x{1}\otimes K$ by the homogeneous maps $\beta =
\beta^h_i,\beta^v$ that determine the coderivations $\partial =
\partial^h_i,\partial^v$ of the quotient Hopf
operad~$\ckcell_m\HopfOp_B^{\E}$. The other component of this
morphism is given simply by the canonical projection of the cofree
coalgebra $\Gamma(K)$. The morphism
\begin{equation*}
(\Gamma(L),\partial)\rightarrow\x{01}\otimes K\oplus\x{1}\otimes L
\end{equation*}
is defined similarly once we observe that the homogeneous maps
$\beta = \beta^h_i,\beta^v$ that determine the coderivations of the
Hopf operad~$\HopfOp_B^{\E}$ admit factorizations
\begin{equation*}
\xymatrix{ \Gamma(\PrimOp_B^{\E})\ar[d]\ar[r]^{\beta} &
\PrimOp_B^{\E}\ar[d] \\
\Gamma(K)\ar[r]^{\beta}\ar[d] & K\ar[d] \\
\Gamma(L)\ar[r]^{\beta}\ar@{-->}[ur] & L }.
\end{equation*}
Accordingly, in the lifting problem considered in the previous
paragraph
\begin{equation*}
\xymatrix{ C^d\ar[r]\ar[d] & \x{01}\otimes K\oplus\x{1}\otimes
K\ar[d] \\
D^d\ar[r]_(0.3){\tilde{\rho}}\ar@{-->}[ur]^{\tilde{\rho}} &
\x{01}\otimes K\oplus\x{1}\otimes L },
\end{equation*}
the lower horizontal morphism is given by $\tilde{\rho}(\xi) =
\x{01}\otimes\beta\nabla_{\rho_{m-1}}(\xi) +
\x{1}\otimes\rho_{m-1}(\xi)$, where we consider the map $\rho_{m-1}:
D^d\rightarrow L = \ckcell_{m-1}\PrimOp_B^{\E}$ specified by
induction, the induced coalgebra morphism $\nabla_{\rho_{m-1}}:
D^d\rightarrow(\Gamma(L),\partial)$ and the map $\beta:
\Gamma(L)\rightarrow K$ defined by a quotient of the total
coderivation of~$\HopfOp_B^{\E}$.

Therefore the lifting morphism $\tilde{\rho}:
D^d\rightarrow\Delta^1\wedge K$ is associated to a map $\rho_{m-1}:
D^d\rightarrow K = \ckcell_m\PrimOp_B^{\E}$ that can be
characterized by the relation $\pk_m\rho_m = \rho_{m-1}$ and by the
equation
\begin{equation*}
\delta(\rho_m(\xi)) - \rho_m(\delta(\xi)) - \rho_m(\partial(\xi)) =
- \beta\nabla_{\rho_{m-1}}(\xi).
\end{equation*}
In fact, for a component $\rho_{m_*}(\xi)$ of~$\rho_m(\xi)$, we
recover the equation of
theorem~\ref{thm:OperationExplicitDefinition}:
\begin{multline*}
\delta(\rho_{m_*}(\xi)) - \rho_{m_*}(\delta(\xi)) -
\rho_{m_*}(\partial(\xi)) = \\- \sum_n\Bigl\{\sum
\shuffle(m^j_i)\cdot\mu_n(\rho_{m^1_*}(\xi^1),\dots,\rho_{m^n_*}(\xi^n))\Bigr\}\\
\pm \sum_i\Bigl\{\sum_{n,t}\rho_{m'_*}(\xi)\circ_i\mu_n\Bigr\}.
\end{multline*}
The important point addressed in this lifting process is that the
right-hand side of this equation consists of components
$\rho_{n_*}(q)$ such that $n_1+\dots+n_r\leq m-1$ and that are
already specified by the inductive construction. Similarly, as we
assume $\partial(\xi)\in C^d$, the term $\rho_{m_*}(\delta(\xi))$ is
determined by the map
\begin{equation*}
C^d\rightarrow\Q^{d-1}\xrightarrow{\nabla_\rho^{d-1}}\HopfOp_B^{\E}
\end{equation*}
specified at a previous stage of this inductive construction.

Finally, the lifting $\tilde{\rho}: D^d\rightarrow\Delta^1\wedge K$
has to be defined as a morphism of unitary $\Lambda_*$-modules.
Therefore, in the next paragraph, we survey the construction of
liftings in the Reedy model category of $\Lambda_*$-modules in order
to achieve our construction.

\subsubsection{On $\Lambda_*$-module liftings and the Reedy model structure}
In fact, in view of the definition of the Reedy model structure (see
more especially the proof of claim~\ref{claim:LambdaModulesLifting}
for properties M4.i-ii), the components of the lifting morphism
$\widetilde{\rho}: D^d(r)\rightarrow\Delta^1\wedge K(r)$ can be
obtained by induction on $r$ as lifting morphisms in the diagrams of
$\Sigma_r$-modules
\begin{equation*}
\xymatrix{ C^d(r)\ar[r]\ar[d] & M(r)\ar[d] \\
D^d(r)\ar[r]\ar@{-->}[ur] & \Match{M}(r)\times_{\Match{N}(r)} N(r)
},
\end{equation*}
where $M = K\wedge\Delta^1$ and $N = S^1\wedge K\times_{S^1\wedge
L}\Delta^1\wedge L$.

We make the cartesian product $\Match{M}(r)\times_{\Match{N}(r)}
N(r)$ explicit in order to make this lifting process effective. In
fact, we have clearly
\begin{multline*}
\Match{M}(r) =
\x{01}\otimes\Match{K}(r)\oplus\x{1}\otimes\Match{K}(r) \\
\text{and}\quad\Match{M}(r)\times_{\Match{N}(r)} N(r) =
\x{01}\otimes
K(r)\oplus\x{1}\otimes(\Match{K}(r)\times_{\Match{L}(r)} L(r)).
\end{multline*}
By definition of the $\Lambda_*$-module structure of
$\PrimOp_B^{\E}$, we obtain
\begin{equation*}
\Match{K}(r)\times_{\Match{L}(r)} L(r) =
{\prod_{m_*}}'\Lambda\E(m_1+\dots+m_r),
\end{equation*}
where $\prod'_{m_*}$ ranges over all collections $m\geq
m_1+\dots+m_r>0$ such that $m_i = 0$ for some $i$ or $m-1\geq
m_1+\dots+m_r>0$ (compare with
observation~\ref{obsv:LambdaModuleFibrations}). Furthermore, the
matching morphism $K(r)\rightarrow\Match{K}(r)\times_{\Match{L}(r)}
L(r)$ is given by the obvious projection of~$K(r) =
\prod_{m_*}\Lambda\E(m_1+\dots+m_r)$ onto the components
of~$\Match{K}(r)\times_{\Match{L}(r)} L(r)$.

As a conclusion the maps $\rho_{m_*}:
D^d(r)\rightarrow\Lambda\E(m_1+\dots+m_r)$ that determine our
lifting $\tilde{\rho}: D^d(r)\rightarrow\Delta^1\wedge K(r)$ can be
characterized by the following properties:\newline -- for the unit
element $1\in D^d(1)$, we have
\begin{equation*}
\rho_{m_1}(1) = \begin{cases} 1\in\Lambda\E(1) & \text{if $m_1 =
1$},
\\ 0\in\Lambda\E(m_1) & \text{otherwise} \end{cases}
\end{equation*}
so that $\tilde{\rho}$ defines a morphism of unitary
$\Lambda_*$-modules;\newline -- for a permutation $w\in\Sigma_r$, we
have the equivariance relations
\begin{equation*}
\rho_{m_1,\dots,m_r}(w\cdot\xi) =
w(m_1,\dots,m_r)\cdot\rho_{m_{w(1)},\dots,m_{w(r)}}(\xi)
\end{equation*}
which imply that $\tilde{\rho}$ commutes the action of~$w$; if $m_i
= 0$ for some $i$, then we have
\begin{equation*}
\rho_{m_1,\dots,m_r}(\xi) =
\rho_{m_1,\dots,\widehat{m_i},\dots,m_r}(\partial_i\xi)
\end{equation*}
so that the composite of $\tilde{\rho}$ with the matching morphism
$K(r)\rightarrow\Match{K}(r)\times_{\Match{L}(r)} L(r)$ matches the
previously defined components of $\tilde{\rho}$;\newline --
otherwise we have the differential equation
\begin{multline*}
\delta(\rho_{m_*}(\xi)) - \rho_{m_*}(\delta(\xi)) -
\rho_{m_*}(\partial(\xi)) = \\- \sum_n\Bigl\{\sum
\shuffle(m^j_i)\cdot\mu_n(\rho_{m^1_*}(\xi^1),\dots,\rho_{m^n_*}(\xi^n))\Bigr\}\\
\pm \sum_i\Bigl\{\sum_n\rho_{m'_*}(\xi)\circ_i\mu_{n,t}\Bigr\},
\end{multline*}
where in the first summation we consider the $n$-fold diagonals
$\sum \xi^1\otimes\dots\otimes\xi^n\in D^d(r)^{\otimes n}$ of the
element $\xi$ and the partitions $m^1_i+\dots+m^n_i = m_i$, in the
second summation we consider the collections $m'_*$ such that $m'_*
= m_*$ for $*\not=i$ and $m'_i = m_i+n-1$ and $t$ ranges over the
interval $t = m_1+\dots+m_{i-1}+1,\dots, m_1+\dots+m_{i-1}+m_i$.

Clearly, if we assume that the morphism $i^d: C^d(r)\rightarrow
D^d(r)$ splits equivariantly so that $D^d(r) = (C^d(r)\oplus
E^d(r),\partial)$ for a finitely generated free $\Sigma_r$-module
$E^d(r) = \bigoplus_\kappa \F[\Sigma_r]\xi_\kappa$, then one can fix
easily a map $\rho_m: D^d(r)\rightarrow K(r)$ that fulfils our
requirements. In fact, it is sufficient to specify solutions of our
differential equation for generators $\xi = \xi_\kappa$ of $E^d(r)$.
If the operad $\E(r)$ is equipped with a contracting homotopy $\nu:
\E(r)\rightarrow\E(r)$ as stated in the introduction, then a
solution can be specified by the formula of
theorem~\ref{thm:OperationRecursiveDefinition}:
\begin{multline*}
\rho_{m_1,\dots,m_r}(\xi_\kappa) =
\nu\Bigl\{\rho_{m_*}(\delta(\xi_\kappa)) + \rho_{m_*}(\partial(\xi_\kappa))\\
- \sum_n\Bigl\{\sum
\shuffle(m^j_i)\cdot\mu_n(\rho_{m^1_*}(\xi_\kappa^1),\dots,\rho_{m^n_*}(\xi_\kappa^n))\Bigr\}\\
\pm\sum_i\Bigl\{\sum_{n,t}\rho_{m'_*}(\xi_\kappa)\circ_t\mu_n\Bigr\}\Bigr\}.
\end{multline*}
Hence this comprehensive survey of our lifting constructions gives
thorough justifications for the recursive construction of
theorem~\ref{thm:OperationRecursiveDefinition}.\qed

\clearpage

\newcommand{\bibpreamble}{\textit{We put in a special reference for the
article~\cite{Bar1} which provides a background and motivations for
this research project. Otherwise numerical references to the
bibliography are given according to the author alphabetical order,
as usual.}\medskip}

\clearpage
\part*{Index and notation glossary}

\begin{center}{\scshape Index}\end{center}

\begin{trivlist}

  \indexspace\item $\Lambda_*$-module, \hyperpage{11}
    \subitem Hopf, \hyperpage{25}
    \subitem matching object of a, \hyperpage{19}
    \subitem non-unital, \hyperpage{11}
    \subitem Reedy model category of $\Lambda_*$-modules,
        \hyperpage{18}
    \subitem unital, \hyperpage{14}
    \subitem unitary, \hyperpage{11}

  \indexspace\item $\Sigma_*$-module, \hyperpage{11}
    \subitem non-unital, \hyperpage{11}
    \subitem unit element of a, \hyperpage{11}
    \subitem unitary, \hyperpage{11}

  \indexspace\item $r$-tree, \hyperpage{32}

  \indexspace\item Boardman-Vogt' construction, \hyperpage{29}

  \indexspace\item chain interval, \hyperpage{31}

  \indexspace\item coalgebra
    \subitem augmented, \hyperpage{76}
    \subitem augmented unitary, \hyperpage{76}
    \subitem cofree, \hyperpage{56}
    \subitem non-augmented, \hyperpage{76}
    \subitem quasi-cofree, \hyperpage{60}
    \subitem unit of a, \hyperpage{24}
    \subitem unitary, \hyperpage{25}

  \indexspace\item cofree
    \subitem coalgebra, \hyperpage{56}
    \subitem Hopf $\Lambda_*$-module, \hyperpage{67}

  \indexspace\item connected coalgebra, \hyperpage{83}
    \subitem morphism coalgebra of a, \hyperpage{83}

  \indexspace\item edge contraction, \hyperpage{33}

  \indexspace\item fibration
    \subitem of $\Z$-graded Hopf operads, \hyperpage{73}
    \subitem of a $\Z$-graded coalgebra, \hyperpage{54}
    \subitem of a $\Z$-graded Hopf object, \hyperpage{54}
    \subitem Reedy fibration of $\Lambda_*$-modules, \hyperpage{19}
    \subitem Reedy fibration of Hopf $\Lambda_*$-modules,
        \hyperpage{26}

  \indexspace\item free
    \subitem non-unital operad, \hyperpage{13}
    \subitem unital Hopf operad, \hyperpage{25}
    \subitem unital operad, \hyperpage{14}

  \indexspace\item Hopf
    \subitem algebra over a Hopf operad, \hyperpage{24}

  \indexspace\item Hopf $\Lambda_*$-module, \hyperpage{25}
    \subitem cofree, \hyperpage{67}
    \subitem matching object of a, \hyperpage{26}
    \subitem quasi-cofree, \hyperpage{68}
    \subitem Reedy model category of Hopf $\Lambda_*$-modules,
        \hyperpage{25}

  \indexspace\item Hopf endomorphism operad
    \subitem of a coalgebra, \hyperpage{86}
    \subitem of the bar construction, \hyperpage{87}, \hyperpage{92}

  \indexspace\item Hopf operad, \hyperpage{24}
    \subitem Boardman-Vogt' construction of a, \hyperpage{29}
    \subitem connected, \hyperpage{27}
    \subitem free unital, \hyperpage{25}
    \subitem of universal bar operations, \hyperpage{92}
    \subitem Reedy model category of unital Hopf operads,
        \hyperpage{27}
    \subitem unital, \hyperpage{24}

  \indexspace\item internal edge of a tree, \hyperpage{32}

  \indexspace\item leaves of a tree, \hyperpage{32}

  \indexspace\item length tensor, \hyperpage{33}
    \subitem module of length tensors, \hyperpage{33}

  \indexspace\item matching object
    \subitem of a $\Lambda_*$-module, \hyperpage{19}
    \subitem of a Hopf $\Lambda_*$-module, \hyperpage{26}
    \subitem of a quasi-cofree Hopf $\Lambda_*$-module, \hyperpage{68}

  \indexspace\item morphism coalgebra, \hyperpage{76}
    \subitem composites in a, \hyperpage{79}, \hyperpage{85}
    \subitem of a cofree coalgebra, \hyperpage{77}
    \subitem of a connected coalgebra, \hyperpage{83}
    \subitem of a quasi-cofree coalgebra, \hyperpage{81}
    \subitem of a quasi-cofree connected coalgebra, \hyperpage{84}
    \subitem of a tensor coalgebra, \hyperpage{84}

  \indexspace\item morphism of $r$-trees, \hyperpage{32}

  \indexspace\item non-unital
    \subitem $\Lambda_*$-module, \hyperpage{11}
    \subitem $\Sigma_*$-module, \hyperpage{11}
    \subitem operad, \hyperpage{11}

  \indexspace\item operad
    \subitem connected, \hyperpage{10}, \hyperpage{15}, \hyperpage{27}
    \subitem free non-unital, \hyperpage{13}
    \subitem free unital, \hyperpage{14}
    \subitem Hopf, \hyperpage{24}
      \subsubitem Hopf algebra over a, \hyperpage{24}
    \subitem Hopf endomorphism operad
      \subsubitem of a coalgebra, \hyperpage{86}
      \subsubitem of the bar construction, \hyperpage{87},
        \hyperpage{92}
    \subitem Hopf operad of universal bar operations, \hyperpage{92}
    \subitem non-unital, \hyperpage{11}
      \subsubitem free, \hyperpage{13}
    \subitem unit operation of an, \hyperpage{10}
    \subitem unital, \hyperpage{10}
      \subsubitem free, \hyperpage{14}
    \subitem unital Hopf, \hyperpage{24}
    \subitem unitary, \hyperpage{10}
      \subsubitem connected, \hyperpage{15}

  \indexspace\item quasi-cofree
    \subitem coalgebra, \hyperpage{60}
    \subitem connected coalgebra
      \subsubitem morphism coalgebra of a, \hyperpage{84}
    \subitem Hopf $\Lambda_*$-module, \hyperpage{68}

  \indexspace\item Reedy
    \subitem fibration
      \subsubitem of $\Lambda_*$-modules, \hyperpage{19}
      \subsubitem of $\Z$-graded Hopf operads, \hyperpage{73}
      \subsubitem of Hopf $\Lambda_*$-modules, \hyperpage{26}
    \subitem model category
      \subsubitem of $\Lambda_*$-modules, \hyperpage{18}
      \subsubitem of Hopf $\Lambda_*$-modules, \hyperpage{25}
      \subsubitem of unital Hopf operads, \hyperpage{27}
      \subsubitem of unital operads, \hyperpage{23}

  \indexspace\item root of a tree, \hyperpage{32}

  \indexspace\item tensor coalgebra
    \subitem morphism coalgebra of a, \hyperpage{84}

  \indexspace\item tree
    \subitem $n$-reduced, \hyperpage{33}
    \subitem $r$-tree, \hyperpage{32}
    \subitem cell metric, \hyperpage{33}
    \subitem edge contraction in a, \hyperpage{33}
    \subitem internal edge of a, \hyperpage{32}
    \subitem isomorphism of $r$-trees, \hyperpage{33}
    \subitem labeling of an $r$-tree, \hyperpage{38}
    \subitem leaves of a, \hyperpage{32}
    \subitem morphism of $r$-trees, \hyperpage{32}
    \subitem root of a, \hyperpage{32}

  \indexspace\item treewise tensor product, \hyperpage{38}

  \indexspace\item unit
    \subitem element of a $\Sigma_*$-module, \hyperpage{11}
    \subitem of a coalgebra, \hyperpage{24}
    \subitem operation of an operad, \hyperpage{10}

  \indexspace\item unital
    \subitem $\Lambda_*$-module, \hyperpage{14}
    \subitem Hopf operad, \hyperpage{24}
    \subitem operad, \hyperpage{10}
    \subitem operation, \hyperpage{10}

  \indexspace\item unital operad
    \subitem reduced operad of a, \hyperpage{10}

  \indexspace\item unitary
    \subitem $\Lambda_*$-module, \hyperpage{11}
    \subitem $\Sigma_*$-module, \hyperpage{11}
    \subitem coalgebra, \hyperpage{25}
    \subitem operad, \hyperpage{10}

\end{trivlist}

\clearpage
\begin{center}{\scshape Notation glossary}\end{center}

\begin{description}

\item[$\ctr^1_*$] the truncation functor, left adjoint to the category embedding $\itr^1_*$,
\hyperpage{15}, see also~\hyperpage{27} for Hopf objects

\item[$\CoAlg^a$] the category of non-augmented coassociative
coalgebras, \hyperpage{76}

\item[$\CoAlg^a_+$] the category of augmented coassociative
coalgebras, \hyperpage{24}, \hyperpage{76}

\item[$\CoAlg^a_*$] the category of augmented unitary coassociative
coalgebras, \hyperpage{25}, \hyperpage{76}

\item[$\CoAlg^a_0$] the category of connected augmented unitary coassociative
coalgebras, \hyperpage{76}, \hyperpage{83}

\item[$\Cube{\tau}$] the cubical chain complex build on the internal
edges of a tree $\tau$, \hyperpage{33}

\item[$\Delta^1\wedge E$] the cone of a dg-module $E$,
\hyperpage{61}

\item[$\partial_i$] the operadic composite with a unital operation in an operad or the corresponding operation in a $\Lambda_*$-module,
\hyperpage{11}

\item[$\partial_\alpha$] the coalgebra coderivation induced by a map $\alpha$, \hyperpage{60},
\hyperpage{68}

\item[$\partial^h_i$] the coalgebra coderivation induced by the bar coderivations on the source in
$\HopfEnd_{B(A)}$, $\HopfEnd_B^{\P}$ and $\HopfOp_B^{\P}$,
\hyperpage{88}, \hyperpage{89}, \hyperpage{93}

\item[$\partial^v$] the coalgebra coderivation induced by the bar coderivation on the target in
$\HopfEnd_{B(A)}$, $\HopfEnd_B^{\P}$ and $\HopfOp_B^{\P}$,
\hyperpage{88}, \hyperpage{89}, \hyperpage{93}

\item[$\epsilon$] the augmentation of an augmented object in a category,
\hyperpage{11}, of a unital operad, \hyperpage{12}, of a coalgebra
\hyperpage{24}, \hyperpage{76}, of the Boardman-Vogt construction,
\hyperpage{29}, of the chain interbal, \hyperpage{31}

\item[$E(\tau)$] the edge set of a tree, \hyperpage{32}

\item[$E'(\tau)$] the set of internal edges in a tree,
\hyperpage{32}

\item[$\eta_0$] the initial morphism of $\Lambda_*$, \hyperpage{14}

\item[$\eta^0$] the $0$-face of the chain interval, \hyperpage{31}

\item[$\eta^1$] the $1$-face of the chain interval, \hyperpage{31}

\item[$\ev$] the evaluation morphism for internal hom objects in general,
\hyperpage{78}

\item[$\ev_\Gamma$] the internal evaluation morphism for cofree coalgebras, \hyperpage{77},
\hyperpage{86}

\item[$\ev_T$] the internal evaluation morphism for tensor coalgebras, \hyperpage{84},
\hyperpage{86}

\item[$\F$] the ground field, \hyperpage{10}

\item[$\FOp(M)$] the free non-unital operad, \hyperpage{13}

\item[$\FOp_*(M)$] the free unital operad, \hyperpage{14}

\item[$\FOp_*(\Gamma)$] the free unital operad generated by a Hopf $\Lambda_*$-module,
\hyperpage{25}

\item[$\Gamma(V)$] the cofree coalgebra, \hyperpage{56}

\item[$\overline{\Gamma}(V)$] the unit cokernel of the cofree coalgebra,
\hyperpage{76}

\item[$\Gamma(M)$] the cofree Hopf $\Lambda_*$-module,
\hyperpage{67}

\item[$\HopfHom(L,M)$] the morphism coalgebra of unitary coalgebras $L$ and $M$,
\hyperpage{76}

\item[$\HopfEnd_\Gamma$] the Hopf endomorphism operad of a coalgebra,
\hyperpage{86}

\item[$\HopfEnd_B^{\P}$] the Hopf endomorphism operad of the bar construction,
\hyperpage{86}

\item[$\HopfOp_B^{\P}$] the Hopf operad of unversal bar operations,
\hyperpage{92}

\item[$\HopfOp^1_*$] the category of unitary Hopf operads,
\hyperpage{25}

\item[$\HopfOp^*_*$] the category of connected unital unitary Hopf operads,
\hyperpage{27}

\item[$\I$] the chain interval, \hyperpage{31}

\item[$I_v$] the entry set of a vertex in a tree, \hyperpage{32}

\item[$\itr^1_*$] the embedding for the subcategory of connected
objects, in the category of unital unitary operads, respectively
$\Lambda_*$-modules, \hyperpage{15}, see also~\hyperpage{27} for
Hopf objects

\item[$i^{\dg}_+$] the embedding of a category formed by $\N$-graded objects into the $\Z$-graded ones,
\hyperpage{60}

\item[$\iota_\infty$] the natural embedding of the cofree coalgebras into the cartesian product of tensor modules,
\hyperpage{56}

\item[$\K$] Stasheff's chain operad of associahedra, \hyperpage{71}

\item[$\overline{L}$] the unit cokernel of a unitary coalgebra,
\hyperpage{76}

\item[$\Lambda_*$] the category of injective maps, \hyperpage{11},
\hyperpage{18}

\item[$\Lambda_{r}^{s}$] the hom set in $\Lambda_*$, \hyperpage{11}

\item[$\Lambda^{\op}_*\HopfMod_0$] the category of non-unital
Hopf $\Lambda_*$-modules, \hyperpage{26}

\item[$\Lambda^{\op}_*\Mod_0$] the category of non-unital $\Lambda_*$-modules,
\hyperpage{11}

\item[$\Lambda^{\op}_*\Mod^1_0$] the category of non-unital unitary $\Lambda_*$-modules,
\hyperpage{11}

\item[$\Lambda^{\op}_*\Mod^1_0/\overline{\C}$] the category of non-unital unitary $\Lambda_*$-modules augmented over the underlying $\Lambda_*$-module of the reduced commutative operad,
\hyperpage{11}

\item[$\Lambda^{\op}_*\HopfMod^1_0$] the category of non-unital
unitary Hopf $\Lambda_*$-modules, \hyperpage{25}

\item[$\Lambda\P$] the suspension of an operad, \hyperpage{92}

\item[$\Match{M}$] the matching object of a $\Lambda_*$-module,
\hyperpage{19}

\item[$\Match{C}$] the matching object of a Hopf $\Lambda_*$-module,
\hyperpage{26}

\item[$\nabla_f$] the morphism of cofree coalgebras, \hyperpage{56},
respectively of quasi-cofree coalgebra, \hyperpage{61}, induced by a
map $f$, see also \hyperpage{68} for Hopf $\Lambda_*$-modules

\item[$\Op^1_0$] the category of non-unital unitary operads,
\hyperpage{11}

\item[$\Op^1_*$] the category of unital unitary operads,
\hyperpage{11}

\item[$\Op^*_*$] the category of connected unital unitary operads,
\hyperpage{15}

\item[$\Op^1_0/\overline{\C}$] the category of non-unital unitary operads augmented over the reduced commutative operad,
\hyperpage{11}

\item[$\overline{\P}$] the reduced operad of a unital operad,
\hyperpage{10}, \hyperpage{12}

\item[$\PrimEnd_{B(A)}$] the module of cogenerators of the quasi-cofree Hopf $\Lambda_*$-module $\HopfEnd_{B(A)}$,
\hyperpage{88}

\item[$\PrimEnd_B^{\P}$] the module of cogenerators of the quasi-cofree Hopf $\Lambda_*$-module $\HopfEnd_B^{\P}$,
\hyperpage{92}

\item[$\PrimOp_B^{\P}$] the module of cogenerators of the quasi-cofree Hopf $\Lambda_*$-module $\HopfOp_B^{\P}$,
\hyperpage{92}

\item[$\P(I)$] the numbering free version of an operad component,
\hyperpage{38}

\item[$\pi$] the universal morphism of the cofree coalgebra,
\hyperpage{56}

\item[$\str^1_*$] the truncation functor, right adjoint to the category embedding $\itr^1_*$,
\hyperpage{15}, see also~\hyperpage{27} for Hopf objects

\item[$\str^{\dg}_+$] the truncation functor, right adjoint to the category embedding $\itr^{\dg}_+$,
\hyperpage{60}, see also~\hyperpage{73} for the case of Hopf operads

\item[$\shuffle(m^j_i)$] a bloc shuffle permutation, \hyperpage{90}

\item[$S^1\wedge E$] the suspension of a dg-module $E$, \hyperpage{61}, see also~\hyperpage{69} for Hopf $\Lambda_*$-modules

\item[$\sigma\wedge E$] the cone projection, \hyperpage{62}, see also~\hyperpage{69} for Hopf $\Lambda_*$-modules

\item[$\Sigma_*\Mod_0$] the category of non-unital $\Sigma_*$-modules,
\hyperpage{11}

\item[$\Sigma_*\Mod^1_0$] the category of non-unital unitary $\Sigma_*$-modules,
\hyperpage{11}

\item[$\Sigma_*\Mod^1_0/\overline{\C}$] the category of non-unital unitary $\Sigma_*$-modules augmented over the underlying $\Sigma_*$-module of the reduced commutative operad,
\hyperpage{11}

\item[$T^c(V)$] the tensor coalgebra, \hyperpage{83}

\item[$\overline{T}^c(V)$] the unit cokernel of the tensor coalgebra,
\hyperpage{83}

\item[$\theta(r)$] the category of $r$-trees, \hyperpage{32}

\item[$\theta'(r)$] a poset formed by isomorphism classes of $r$-trees,
\hyperpage{33}

\item[$\theta'_d(r)$] a poset formed by isomorphism classes of $r$-trees with no more than $d$ vertices,
\hyperpage{33}

\item[$\theta''(r)$] a poset formed by isomorphism classes of reduced $r$-trees,
\hyperpage{33}

\item[$\tau/e_0$] an edge contraction in a tree, \hyperpage{33}

\item[$\tau_r$] the terminal $r$-tree, \hyperpage{33}

\item[$\tau(\P)$] the module spanned by the operadic labelings of a
tree $\tau$, \hyperpage{38}

\item[$V(\tau)$] the vertex set of a tree, \hyperpage{32}

\item[$W(\P)$] the Boardman-Vogt construction of an operad,
\hyperpage{29}

\end{description}


\begin{thebibliography}{99}

\bibitem[B1]{Bar1}
B. Fresse, \emph{The bar construction of an E-infinity algebra},
preprint
\href{http://arxiv.org/math.AT/0601085}{\texttt{arXiv:math.AT/0601085}}
(2006).

\bibitem{Berger}
C. Berger, \emph{Op\'erades cellulaires et espaces de lacets
it\'er\'es}, Ann. Inst. Fourier \textbf{46} (1996), 1125--1157.

\bibitem{BergerEnglish}
C. Berger, \emph{Combinatorial models for real configuration spaces
and $E_n$-operads}, \emph{in} ``Operads: proceedings of renaissance
conferences'', Contemp. Math. \textbf{202}, Amer. Math. Soc. (1997),
37--51.

\bibitem{BergerFresse}
C. Berger, B. Fresse, \emph{Combinatorial operad actions on
cochains}, Math. Proc. Camb. Philos. Soc. \textbf{137} (2004),
135--174.

\bibitem{BergerMoerdijk}
C. Berger, I. Moerdijk, \emph{Axiomatic homotopy theory for
operads}, Comment. Math. Helv. \textbf{78} (2003), 805--831.

\bibitem{BergerMoerdijkW}
\bysame, \emph{The Boardman-Vogt resolution of operads in monoidal
model categories}, preprint
\href{http://arxiv.org/math.AT/0502155}{\texttt{arXiv:math.AT/0502155}}
(2005).

\bibitem{BoardmanVogt}
J. Boardman, R. Vogt, \emph{Homotopy invariant algebraic structures
on topological spa\-ces}, Lecture Notes in Mathematics \textbf{347},
Springer-Verlag, 1973.

\bibitem{BousfieldFriedlander}
A.K. Bousfield, E.M. Friedlander, \emph{Homotopy theory of
$\Gamma$-spaces, spectra and bisimplicial sets}, \emph{in}
``Geometric applications of homotopy theory, II'', Lecture Notes in
Math. \textbf{658}, Springer-Verlag (1978), 80--130.

\bibitem{Ching}
M. Ching, \emph{Bar constructions for topological operads and the
Goodwillie derivatives of the identity},
\href{http://www.maths.warwick.ac.uk/gt/GTVol9/paper20.abs.html}{Geom.
Topol. \textbf{9} (2005), 833--933}.

\bibitem{Crans}
S. Crans, \emph{Quillen closed model structures for sheaves}, J.
Pure Appl. Alg. \textbf{101} (1995), 35--57.

\bibitem{DwyerSpalinski}
W. Dwyer, J. Spalinski, \emph{Homotopy theories and model
categories}, in ``Handbook of algebraic topology'', North-Holland
(1995), 73-126.

\bibitem{BarEinfinityNote}
B. Fresse, \emph{La construction bar d'une alg\`ebre comme alg\`ebre
de Hopf E-infini}, C.~R.~Acad. Sci. Paris S\'er. I \textbf{337}
(2003), 403--408.

\bibitem{OperadTextbook}
\bysame, \emph{Koszul duality of operads and homology of partition
posets}, \emph{in} ``Homotopy theory: relations with algebraic
geometry, group cohomology, and algebraic $K$-theory'', Contemp.
Math. \textbf{346}, Amer. Math. Soc. (2004), 115--215.

\bibitem{GetzlerGoerss}
E. Getzler, P. Goerss, \emph{A model category structure for
differential graded coalgebras}, preprint (1999).

\bibitem{GetzlerJones}
E. Getzler, J. Jones, \emph{Operads, homotopy algebra and iterated
integrals for double loop spaces}, preprint
\href{http://arxiv.org/hep-th/9403055}{\texttt{arXiv:hep-th/9403055}}
(1994).

\bibitem{GinzburgKapranov}
V. Ginzburg, M.M. Kapranov, \emph{Koszul duality for operads}, Duke
Math. J. \textbf{76} (1995), 203--272.

\bibitem{Hinich}
V. Hinich, \emph{Homological algebra of homotopy algebras}, Comm.
Algebra \textbf{25} (1997), 3291--3323.

\bibitem{Hirschhorn}
P. Hirschhorn, \emph{Model categories and their localizations},
Mathematical Surveys and Monographs \textbf{99}, American
Mathematical Society, 2003.

\bibitem{HoveyTextBook}
M. Hovey, \emph{Model categories}, Mathematical Surveys and
Monographs \textbf{63}, American Mathematical Society, 1999.

\bibitem{McClureSmith}
J. McClure, J.H. Smith, \emph{Multivariable cochain operations and
little $n$-cubes}, J. Amer. Math. Soc. \textbf{16} (2003), 681--704.

\bibitem{McClureSmithCosimplicial}
J. McClure, J.H. Smith, \emph{Cosimplicial objects and little
$n$-cubes, I}, Amer. J. Math. \textbf{126} (2004), 1109--1153.

\bibitem{QuillenBook}
D. Quillen, \emph{Homotopical algebra}, Lecture Notes in Mathematics
\textbf{43}, Sprin\-ger-Verlag, 1967.

\bibitem{Reedy}
C. Reedy, \emph{Homotopy theory of model categories}, unpublished
manuscript (1974).

\bibitem{Rezk}
C. Rezk, \emph{Spaces of algebra structures and cohomology of
operads}, PhD Thesis, Massachusetts Institute of Technology, 1996.

\bibitem{SmithJustin}
J.R. Smith, \emph{Cofree coalgebras over operads}, Topology Appl.
\textbf{133} (2003), 105--138.

\bibitem{Sweedler}
M. Sweedler, \emph{Hopf algebras}, Benjamin, 1969.

\end{thebibliography}
\end{document}